\documentclass{amsart}
\usepackage{amssymb}
\usepackage{psfig}
\newtheorem{theorem}{Theorem}[section]
\newtheorem{corollary}[theorem]{Corollary}
\newtheorem{lemma}[theorem]{Lemma}
\newtheorem{proposition}[theorem]{Proposition}
\theoremstyle{definition}
\newtheorem{definition}[theorem]{Definition}
\newtheorem{remark}[theorem]{Remark}
\newtheorem{scholium}[theorem]{Scholium}
\newtheorem{example}[theorem]{Example}
\theoremstyle{remark}
\newtheorem{observation}{Observation}[section]
\renewcommand{\theobservation}{\arabic{observation}}

\numberwithin{equation}{section}

\newlength{\customskipamount}
\def\func#1{\mathop{\rm #1}}%
\def\limfunc#1{\mathop{\rm #1}}%
\long\def\TeXButton#1#2{#2}%
\def\openone{\hbox{\upshape \small1\kern-3.3pt\normalsize1}}
\newbox\crossbox
\def\bigcross{\mathop{\mathchoice
{\setbox\crossbox=\hbox{\Huge \bfseries $\times$}
\raise\fontdimen22\textfont2%
\hbox{\raise0.5\dp\crossbox%
\hbox{\lower0.5\ht\crossbox%
\box\crossbox}}}
{\setbox\crossbox=\hbox{\huge \bfseries $\times$}
\raise\fontdimen22\textfont2%
\hbox{\raise0.5\dp\crossbox%
\hbox{\lower0.5\ht\crossbox%
\box\crossbox}}}
{\setbox\crossbox=\hbox{\LARGE \bfseries $\times$}
\raise\fontdimen22\scriptfont2%
\hbox{\raise0.5\dp\crossbox%
\hbox{\lower0.5\ht\crossbox%
\box\crossbox}}}
{\setbox\crossbox=\hbox{\Large \bfseries $\times$}
\raise\fontdimen22\scriptscriptfont2%
\hbox{\raise0.5\dp\crossbox%
\hbox{\lower0.5\ht\crossbox%
\box\crossbox}}}}}

\begin{document}
\title[Permutation representations of the Cuntz algebra]{Iterated function 
systems and permutation representations of the Cuntz algebra}
\author{Ola Bratteli}
\address{Mathematics Institute\\
University of Oslo\\
PB 1053 -- Blindern\\
N-0316 Oslo\\
Norway}
\email{bratteli@math.uio.no}
\thanks{Work supported in part by the U.S. National Science Foundation and the
Norwegian Research Council.}
\author{Palle E.T. Jorgensen}
\address{Department of Mathematics\\
The University of Iowa\\
14 MacLean Hall\\
Iowa City, IA 52242-1419\\
U.S.A.}
\email{jorgen@math.uiowa.edu}

\begin{abstract}
We study a class of representations of the Cuntz algebras $\mathcal{O}_N$, $%
N=2,3,\dots $, acting on $L^2\left( \mathbb{T}\right) $ where $\mathbb{T}=%
\mathbb{R}\diagup 2\pi \mathbb{Z}$. The representations arise in wavelet
theory, but are of independent interest. We find and describe the
decomposition into irreducibles, and show how the $\mathcal{O}_N$%
-irreducibles decompose when restricted to the subalgebra $\limfunc{UHF}%
\nolimits_N\subset \mathcal{O}_N$ of gauge-invariant elements; and we show
that the whole structure is accounted for by arithmetic and combinatorial
properties of the integers $\mathbb{Z}$. We have general results on a class
of representations of $\mathcal{O}_N$ on Hilbert space $\mathcal{H}$ such
that the generators $S_i$ as operators permute the elements in some
orthonormal basis for $\mathcal{H}$. We then use this to extend our results
from $L^2\left( \mathbb{T}\right) $ to $L^2\left( \mathbb{T}^d\right) $, $%
d>1 $; even to $L^2\left( \mathbf{T}\right) $ where $\mathbf{T}$ is some 
fractal version of the torus which carries more of the algebraic information 
encoded in our representations.
\end{abstract}
\maketitle
\tableofcontents
\listoffigures

\setcounter{section}{0}

\section{\label{Sec1}Introduction\label{Introduction}}

Let $\mathcal{H}$ be the Hilbert space $L^2\left( \mathbb{T}\right) $, $%
\mathbb{T}=\mathbb{R}\diagup 2\pi \mathbb{Z}$ the torus, and let $\left\{
z^n\right\} _{n=-\infty }^\infty $ be the usual orthonormal basis of Fourier
analysis: the convention is $z=e^{it}$ and $z^n=e^{itn}$, $t\in \mathbb{R}$, 
$n\in \mathbb{Z}$; and the Haar measure on $\mathbb{T}$ will be normalized.
The following representations of $\mathcal{O}_N$ (the Cuntz algebra of index 
$N$, \cite{Cun77}) arise in the study of filter-banks in wavelet theory; 
see, e.g., \cite
{Dau}, \cite{CoRy}, \cite{DDL95}, \cite{Jor95}, \cite{DaLa}, \cite{JoPe96},
and \cite{JoPe94}. Let functions $\left\{ m_i\right\} _{i=0}^{N-1}\subset
L^\infty \left( \mathbb{T}\right) $ be given such that the corresponding $%
N\times N$ matrix 
\begin{equation}
\left( m_k\left( e^{i\frac{2\pi l}N}z\right) \right) _{k,l=0}^{N-1}
\label{Eq1.1}
\end{equation}
is unitary. If, for example, $N=2$, the condition is unitarity of $\smash{%
\left( 
\begin{smallmatrix}
m_0\left( z\right) & m_0\left( -z\right) \\ 
m_1\left( z\right) & m_1\left( -z\right)
\end{smallmatrix}
\right) }$ $(\mathrm{a.a.}\,z\in \mathbb{T})$. If $m_0\in L^\infty \left( %
\mathbb{T}\right) $ is given, subject to the condition 
\begin{equation}
\left| m_0\left( z\right) \right| ^2+\left| m_0\left( -z\right) \right| ^2=1%
\text{\quad }\left( \mathrm{a.a.}\,z\in \mathbb{T}\right) \text{,}
\label{Eq1.2}
\end{equation}
then setting $m_1\left( z\right) =z\overline{m_0\left( -z\right) }$, the
unitarity property follows. Similarly, if $f\in L^\infty \left( \mathbb{T}%
\right) $ is unimodular, i.e., $\left| f\left( z\right) \right| =1$ $\mathrm{%
a.a.}\,z\in \mathbb{T}$, then 
\begin{equation}
m_1^f\left( z\right) :=zf\left( z^2\right) \overline{m_0\left( -z\right) }
\label{Eq1.3}
\end{equation}
is an admissible choice, and, in fact, conversely, any $m_1$ has this form;
see \cite{Dau}.

The announced representation of $\mathcal{O}_N$ is defined from the
functions $m_i$ as follows: 
\begin{equation*}
S_i\xi \left( z\right) =N^{\frac 12}m_i\left( z\right) \xi \left( z^N\right)
\end{equation*}
for $i\in \mathbb{Z}_N=\mathbb{Z}\diagup N\mathbb{Z}$, $\xi \in L^2\left( %
\mathbb{T}\right) $, and $z\in \mathbb{T}$. These representations were
introduced in \cite{Jor95}. The corresponding adjoint operators $S_i^{*}$
occur in more general contexts (although not as representations) in Ruelle's
theory \cite{Rue94} under the name ``transfer operators''.

Our analysis of the trigonometric basis $\left\{ z^n\right\} _{n\in %
\mathbb{Z}}$ (and the corresponding multi-variable case) is based
on another viewpoint which is more special in some respects and
more general in others. It is based on iterated function
systems (i.f.s.); and we refer to \cite{Mau95}, \cite{JoPe95}, and \cite
{Str95} for details on previous work. The theory of i.f.s. divides itself
into: (i) the fractals (iteration of ``small'' scales), and (ii) the
``large'' discrete systems. Our emphasis here will be (ii), and Strichartz 
\cite{Str95} has suggested the name ``fractals in the large'', or reverse
iterated function systems (r.i.f.s.) for the latter viewpoint. But let us call
attention to the duality between (i) and (ii)
(studied for example in \cite{Str95} and \cite{JoPe95}), in the
sense that (i) and (ii) serve as the two sides of Fourier duality. This
viewpoint is inspired by classical Pontrjagin duality (see, e.g., \cite{BrRo}
and \cite{HRW92}), but the new feature here is that neither side of the
(i)--(ii) setting in its fine structure alone admits group duality. Instead we
base our analysis on an associated discrete dynamical
system. This system in turn is tailored to the Cuntz algebra $\mathcal{O}_N$.

The simplest examples, for $N=2$, arise as follows: let $s$ be an odd
integer, and let the matrix be $2^{-\frac 12} 
\begin{pmatrix}
1 & 1 \\ 
z^s & -z^s
\end{pmatrix}
$. Then 
\begin{equation}
\begin{aligned} \left( S_0\xi \right) \left( z\right) &=\xi \left(
z^2\right) \text{,} \\ \left( S_1\xi \right) \left( z\right) &=z^s\xi \left(
z^2\right) \text{,\quad }\forall \xi \in L^2\left( \mathbb{T}\right)
\text{, }\mathrm{a.a.}\,z\in \mathbb{T}\text{,} \end{aligned}  \label{Eq1.4}
\end{equation}
and we note that 
\begin{equation}
\begin{split}
&S_i^{*}S_j^{}=\delta _{ij}\openone \text{\quad and} \\
&\sum_{i=0}^1S_i^{}S_i^{*}=\openone \text{,}
\end{split}
\label{Eq1.5}
\end{equation}
where $\openone$ is the identity operator in $\mathcal{H}=L^2\left( %
\mathbb{T}\right) $.

Similarly, if $N\in \left\{ 2,3,\dots \right\} $, and $D=\left\{ s_i\right\}
_{i=0}^{N-1}\subset \mathbb{Z}$ is a set of integers such that any two
distinct members are mutually incongruent modulo $N$, then the operators 
\begin{equation}
\left( S_i\xi \right) \left( z\right) =z^{s_i}\xi \left( z^N\right)
\label{Eq1.6}
\end{equation}
(for $\xi \in L^2\left( \mathbb{T}\right) $), satisfy the Cuntz relations 
\begin{equation}
\begin{split}
&S_i^{*}S_j^{}=\delta _{ij}\openone \text{\quad and} \\
&\sum_{i=0}^{\hbox to0pt{\hss $\scriptstyle N-1$ \hss} }S_i^{}S_i^{*}=%
\openone \text{.}
\end{split}
\label{Eq1.7}
\end{equation}

Recall that the Cuntz algebra $\mathcal{O}_N$ is the $C^{*}$-algebra
generated by $S_0,\dots ,S_{N-1}$ satisfying (\ref{Eq1.7}). This algebra is
simple \cite{Cun77}. Every system of operators $\left\{ S_i\right\} $ on a
Hilbert space $\mathcal{H}$, subject to these same relations, then
determines a representation $\pi $ of $\mathcal{O}_N$ on $\mathcal{H}$,
i.e., $\pi \in \limfunc{Rep}\left( \mathcal{O}_N,\mathcal{H}\right) $.
(See, e.g., \cite[Section 2.8]{EvKa}
and \cite{PoSt} for basic facts on representation theory.)

In Section \ref{Sec2} we will embed the representation defined in (\ref
{Eq1.6}) into a much more general class of representations
coming from certain iterated function systems, called
permutative, mul\-ti\-plic\-i\-ty-free representations, and one of the main
results of this paper is the embedding of all these into a universal
permutative, mul\-ti\-plic\-i\-ty-free representation; see Corollary \ref
{CorNew7.1}. This embedding is in terms of an embedding of our iterated
function systems into a certain ``coding space'' which
is a variant of a dynamical systems tool used also in \cite{Jor95}, \cite
{Pri87}, \cite{PoPr94}, \cite{BreJor}, \cite{EHK}
and more generally in \cite{Rue94}.
Another special way of obtaining such representations is to
replace $L^2\left( \mathbb{T}\right) $ with $L^2\left( \mathbb{T}^\nu
\right) $, and fix a matrix $\mathbf{N}$ with integer coefficients such that 
$\left| \det \left( \mathbf{N}\right) \right| =N$. Then choose a set $%
D=\left\{ d_1,\dots ,d_N\right\} $ of points in $\mathbb{Z}^\nu $ which are
incongruent modulo $\mathbf{N}\mathbb{Z}^\nu $, i.e., $d_i-d_j\notin \mathbf{%
N}\mathbb{Z}^\nu $ for $i\neq j$. As in elementary arithmetic, one proves
that the quotient mapping $\mathbb{Z}^\nu \rightarrow \mathbb{Z}^\nu \diagup 
\mathbf{N}\mathbb{Z}^\nu $ is $1$--$1$ when restricted to $D$, and if $d_i$ are
incongruent modulo $\mathbf{N}\mathbb{Z}^\nu $ then $d_i+\mathbf{N}l_i$ are
incongruent modulo $\mathbf{N}\mathbb{Z}^\nu $ for any choice of $l_i\in %
\mathbb{Z}^\nu $, and these vectors constitute the most general choice for $D
$. The associated representation of $\mathcal{O}_N$ is given by 
\begin{equation}
\left( S_i\xi \right) \left( z\right) =z^{d_i}\xi \left( z^{\mathbf{N}%
}\right) 
\label{Composition}
\end{equation}
where we define 
\begin{align*}
z& =\left( z_1,\dots ,z_\nu \right) \in \mathbb{T}^\nu  \\
z^{d_i}& =z_1^{d_{i1}}z_2^{d_{i2}}\cdots z_\nu ^{d_{i\nu }} \\
z^{\mathbf{N}}& =\left( z_1^{n_{11}}\cdots z_\nu ^{n_{1\nu
}},z_1^{n_{21}}\cdots z_\nu ^{n_{2\nu }},\dots ,z_1^{n_{\nu 1}}\cdots z_\nu
^{n_{\nu \nu }}\right) 
\end{align*}
if $\mathbf{N}=\left( n_{ij}\right) _{i,j=1}^\nu $. These types of
representations have been considered in \cite{JoPe96}, \cite{Ban91}, \cite
{Str94}. To make these representations mul\-ti\-plic\-i\-ty-free, one needs
the extra assumption 
\begin{equation*}
\bigcap_k\mathbf{N}^k\mathbb{Z}^\nu =0
\end{equation*}
on $\mathbf{N}$; see (\ref{Eq3.7}). In analyzing the properties of these
representations, the dynamical system defined by a certain $N$ to $1$ map $R:%
\mathbb{Z}^\nu \rightarrow \mathbb{Z}^\nu $ will play a dominant role. In
these cases, $R$ is defined by the requirement 
\begin{equation}
R\left( d_i+\mathbf{N}x\right) =x\text{,}
\label{Trans}
\end{equation}
and since (as we show) $R$ behaves like $\mathbf{N}^{-1}$ on a large scale,
both $\mathbf{N}$ and $D$ (up to permutation of the indices) can be
recovered from $R$, as a consequence of the fact that $R$ behaves like $%
\mathbf{N}^{-1}$ at large distances. See Scholium \ref{Sch2.4},
the discussion around (\ref{Eq3.6bis})--(\ref{Eq3.7}),
Proposition \ref{Pro3ins1},
and Subsection \ref{SubsecNew9.2}.

There are some similarities and some differences between the analysis
of the monomial representations (\ref{Composition}) in one dimension
(i.e., (\ref{Eq1.6})), and in higher dimensions. We study both cases
below: Section \ref{Examples} ($\nu =1$), and Section
\ref{ExamplesMatrixCase} ($\nu >1$). Our analysis in Section
\ref{Examples} in the case $N=2$ is based partly on factorization of
the Mersenne numbers, whereas the higher-dimensional cases (Section
\ref{ExamplesMatrixCase}) involve the geometry of a set $\mathbf{T}$
of generalized fractions constructed from a given pair $\left(
\mathbf{N},D\right) $ in $\nu $ dimensions where the matrix
$\mathbf{N}$ is as described, and $D$ represents generalized
``digits'', i.e. the vector version of such. See
(\ref{Eq3ins1})--(\ref{Eq3ins2}) for the precise definitions. One of
the conclusions coming out of our analysis in Sections
\ref{Examples}--\ref{ExamplesMatrixCase} is that there are infinite
families of non-isomorphic examples. See \cite{BaGe94} and
\cite{Gel96} on the isomorphism problem for the ``reptiles'' based on
$\left( \mathbf{N},D\right) $ systems. Some of the properties of these
reptiles can be understood from our study of periodic points of the
transformation $R$ from (\ref{Trans}), and also from the question of
which rank-$\nu $ lattices $\mathbb{L\kern0.5pt}$ make $\mathbf{T}$
tile $\mathbb{R}^\nu $. We also establish in Section
\ref{ExamplesMatrixCase} a connection between the last two, i.e.,
periodic points and lattices.  Another variation of the construction
in the previous paragraph is to replace $\ell ^2\left( \mathbb{Z}^\nu
\right) $ with the Hilbert space spanned by the monomials $z^n$ where
$n$ is specified to lie in a certain subset $S$ of $\mathbb{Z}^\nu $,
and we define the operators $S_i$ by the same formula as before, i.e.,
(\ref{Composition}), but with the function $\xi $ restricted to lie in
the corresponding closed subspace $\left( {}\subset L^2\left(
\mathbb{T}^\nu \right) \right) $.  Instead of assuming the elements of
$D$ be mutually incongruent $\bmod{\,\mathbf{N}\mathbb{Z}^\nu }$, we
now merely need to assume that $S=\bigcup _{i=1}^N\left(
d_i+\mathbf{N}S\right) $, and $\left( d_i+\mathbf{N}S\right) \cap
\left( d_j+\mathbf{N}S\right) =\emptyset$ for $i\neq j$.  This
structure occurs in some dynamical systems considered by \cite{Ke},
\cite{Od}, and \cite{LaWa96}, and is treated in detail in Remark
\ref{Rem3ins3}.

Let $\mathcal{B}\left( \mathcal{H}\right) $ denote the $C^{*}$-algebra of
all bounded operators on the Hilbert space $\mathcal{H}$. In \cite{Arv89}, 
\cite{Lac93a}, and \cite{BJP} it was noted that the study of the
endomorphisms of $\mathcal{B}\left( \mathcal{H}\right) $, i.e., $\limfunc{End%
}\left( \mathcal{B}\left( \mathcal{H}\right) \right) $, is equivalent to
that of $\limfunc{Rep}\left( \mathcal{O}_N,\mathcal{H}\right) $. If $S_i$
satisfy (\ref{Eq1.7}), then an endomorphism $\alpha $ is defined by 
\begin{equation}
\alpha \left( A\right) =\sum_{i=0}^{N-1}S_i^{}AS_i^{*}\text{,\quad }A\in 
\mathcal{B}\left( \mathcal{H}\right) \text{;}  \label{Eq1.8}
\end{equation}
and, conversely, every $\alpha \in \limfunc{End}\left( \mathcal{B}\left( 
\mathcal{H}\right) \right) $ is of this form.

Consider the automorphism group (``gauge group'') of $\mathcal{O}_N $,
denoted by $\left( \gamma _z\right) _{z\in \mathbb{T}}$, which is determined
by $\gamma _z\left( S_i\right) =zS_i$, $z\in \mathbb{T}$. The subalgebra $%
\limfunc{UHF}\nolimits_N$ of gauge-invariant elements is 
\begin{equation}
\limfunc{UHF}\nolimits_N:=\left\{ a\in \mathcal{O}_N:\gamma _z\left( a
\right) =a,\;\forall z\in \mathbb{T}\right\} \text{.}  \label{Eq1.9}
\end{equation}
It is a $\limfunc{UHF}$ algebra of Glimm type $\mathcal{M}_{N^\infty }$
(see, e.g., \cite{Cun77}, \cite{BrRo}, \cite{BJP}, and \cite{BJ}). For general
representations $\pi \in \limfunc{Rep}\left( \mathcal{O}_N,\mathcal{H}%
\right) $, the corresponding restrictions $\pi |_{\limfunc{UHF}\nolimits_N}$
have been studied in, e.g., \cite{BEEK} and \cite{Pow88}.

In those studies, the case when $\pi \left( \limfunc{UHF}\nolimits_N\right) $
is weakly*-dense in $\pi \left( \mathcal{O}_N\right) $ is predominant.
However, in the present examples this weak*-density does not hold. We will
show that, in the $D$-examples, the $\mathcal{O}_N$-irreducibles, when
restricted to $\limfunc{UHF}\nolimits_N$, break up as finite orthogonal sums
of mutually inequivalent $\limfunc{UHF}\nolimits_N$-irreducibles. Hence, we
get the $C^{*}$-algebra version of what in group theory is called ``Gelfand
pairs''; see, e.g., \cite{BeRa} and \cite{BJR}.

In group theory, the occurrence of Gelfand pairs appears to be rare, and
when it happens, the representation theory appears to account for the
underlying structure of the groups in question. Let $G$ be a group with a
subgroup $K$, and let $\pi $ denote some unitary representation of $G$ on a
Hilbert space. The restriction of $\pi $ to $K$ will be denoted $\pi |_K$.
We say that $\left( G,K\right) $ is a Gelfand pair if the irreducibles $\pi $
of $G$ are mul\-ti\-plic\-i\-ty-free when restricted to $K$. Typically, we
may not know \emph{all} the (equivalence classes of) irreducible
representations of $G$, and the definition then applies instead to a
suitably restricted family of irreducibles of $G$.

The analogy to the present setting is clear from this: we have $\limfunc{UHF}%
\nolimits_N\subset \mathcal{O}_N$ as a subalgebra, and any given permutative
mul\-ti\-plic\-i\-ty-free representation $\pi $ splits into a direct sum of
irreducible mutually inequivalent representations of $\mathcal{O}_N$, 
\begin{equation}
\pi =\sideset{}{^{\oplus}}{\sum}_i\pi _i\text{.}  \label{Eq1.10}
\end{equation}
Furthermore each $\pi _i|_{\limfunc{UHF}\nolimits_N}$ has a decomposition 
\begin{equation}
\pi _i|_{\limfunc{UHF}\nolimits_N}=\sideset{}{^{\oplus}}{\sum}_j\rho _j
\label{Eq1.11}
\end{equation}
where each $\rho _j$ is irreducible and the $\rho _j$ are mutually
inequivalent, also for different $i$; see Theorem \ref{Thm2.5}. Moreover, if
the representation $\pi $ comes from a subset $D=\left\{ s_1,\dots
,s_N\right\} $ of points in $\mathbb{Z}$ which are incongruent modulo $N$ as
described in (\ref{Eq1.6}), then all the decompositions (\ref{Eq1.10}) and (%
\ref{Eq1.11}) are finite; see Corollary \ref{Cor3.5} and Corollary \ref
{Cor3.10}. Both decompositions are also finite when $\pi $ comes from a
subset $D=\left\{ d_1,\dots ,d_N\right\} $ of points in $\mathbb{Z}^\nu $
which are incongruent modulo $\mathbf{N}$ and all the (complex) eigenvalues
of $\mathbf{N}$ have modulus greater than one; see Corollary \ref{Cor3.9}.

We obtain this double decomposition for $\pi $ in a very explicit form from
two equivalence relations on $\mathbb{N}$, $\sim $ and $\approx $ (to be
described below), such that the $\sim $ equivalence classes correspond to
the $\pi _i$'s, and (for fixed $i$) the $\approx $ (finer) equivalence
classes contained in each $\sim $ class correspond to the representations $%
\rho _j$ which occur in the decomposition (\ref{Eq1.11}) of $\pi _i|_{%
\limfunc{UHF}\nolimits_N}$.

For the given representation $\pi $ we will use the terminology that the
subrepresentations $\pi _i$ are the cycles in $\pi $, and the
subrepresentations $\rho _j$ of $\pi _i|_{\limfunc{UHF}\nolimits_N}$ are the
atoms in the cycles. (We use this terminology in a sense a little different
from that of Harish-Chandra (representation theory); and in our present use
it is for convenience, to catch some intuition behind the grouping
of irreducible representations into cyclic substrata.) If $%
\alpha $ is the endomorphism of $\mathcal{B}\left( L^2\left( \mathbb{T}%
\right) \right) $ corresponding to $\pi $, recall from \cite{BJP}, \cite
{Lac93a} that $\pi \left( \mathcal{O}_N\right) ^{\prime }$ is the fixed
point subalgebra under $\alpha $, and $\pi \left( \limfunc{UHF}%
\nolimits_N\right) ^{\prime }$ is the algebra at infinity $\bigcap_n\alpha
^n\left( \mathcal{B}\left( L^2\left( \mathbb{T}\right) \right) \right) $.
Thus, for our representations $\pi $, the fixed point algebra is an abelian
algebra with atomic spectrum, and for each minimal projection $E$ in this
algebra, the ergodic restriction of $\alpha $ to $E\mathcal{B}(L^2(\mathbb{T}%
))E$\/ has an algebra at infinity which is an abelian algebra with atomic
spectrum and with dimension equal to the number of atoms in the cycle
corresponding to $E$.

In Theorem \ref{Thm4.1} we will define an action of $\mathbb{Z}$ on the
atoms, such that the orbits under this action correspond to the cycles. This
picture gives our results the flavor of the familiar Kirillov-orbit picture
for the irreducible representations of nilpotent real Lie groups; see, e.g., 
\cite{BJR}.

Let us describe this action and the equivalence relations $\sim $ and $%
\approx $ on $\mathbb{N}$ further,
and their connection with the $N$ to $1$ map $R:%
\mathbb{N}\rightarrow \mathbb{N}$ defined in Scholium \ref{Sch2.4}. The most
instructive way of defining $\sim $ and $\approx $ is in Scholium \ref
{Sch3.8}. If $x,y\in \mathbb{N}$, then $x\approx y$ if and only if the tails
of the two orbit sequences $\left\{ x,Rx,R^2x,\dots \right\} $ and $\left\{
y,Ry,R^2y,\dots \right\} $ are equal, and $x\sim y$ if the tails of the two
sequences are equal up to translation of the sequences. This translation
thus defines the action of $\mathbb{Z}$ on the atoms alluded to in the
previous paragraph. Now, for a given $x\in \mathbb{N}$ two things may
happen: (i) all elements in the sequence $\left\{ x,Rx,R^2x,\dots \right\} $
may be distinct, in which case the cycle corresponding to $x$ has infinitely
many atoms, or (ii) there may be two elements in the sequence which are
identical, in which case the sequence asymptotically becomes periodic with
(minimal) period $\limfunc{Per}\left( x\right) $. In the latter case, there
are just $\limfunc{Per}\left( x\right) $ atoms in the cycle containing $x$, and
the corresponding representations are induced by so-called sub-Cuntz states
which will be described in Section \ref{Sec5}. Now, if $\mathbb{N}$ is replaced
by $\mathbb{Z}^\nu $, and $R$ is induced by a matrix $\mathbf{N}$, and a set 
$D$ as described above, we already mentioned that the large-scale behavior
of $R$ is the same as that of $\mathbf{N}^{-1}$. Thus if $\mathbf{N}^{-1}$
is contractive in some norm, all points will ultimately be mapped into a
periodic orbit near the origin, and the number of atoms can be explicitly
bounded; see Corollaries \ref{Cor3.9} and \ref{Cor3.10}. If more generally $%
\mathbf{N}^{-1}$ is hyperbolic, there will be a finite number of finite
orbits, the number being estimable, and in addition an infinite number of
orbits growing exponentially along the unstable directions of $\mathbf{N}%
^{-1}$; see Proposition \ref{Pro9.2}. If $\mathbf{N}^{-1}$ is not
hyperbolic, it is probably very hard in general to decide whether the
number of finite orbits is finite or not; see, e.g., \cite{Lag85} or
\cite[pp. 110--114]{Sen95}.

All the general definitions and results pertaining to permutative
mul\-ti\-plic\-i\-ty-free representations in Sections \ref{Sec2}, \ref{Sec4}%
, \ref{Sec5}, \ref{SecNew6}, and \ref{SecNew7} extend to the Cuntz algebra $%
\mathcal{O}_\infty $ with the obvious modifications, such as that the
combinatorial function $N\left( n\right) $ in Corollary \ref{CorNew6.3}, and
Corollary \ref{CorNew7.1}, is infinite for $n\in \mathbb{N}$ in the case $%
N=\infty $.

The present setting is chosen because it appears basic to Fourier analysis
and representation theory. The significance of the number $N$ for wavelet
theory is that it is the \emph{scaling,} the dyadic wavelets corresponding
to $N=2$. The main problem in wavelet theory is the basis problem for $%
L^2\left( \mathbb{R}\right) $, and not $L^2\left( \mathbb{T}\right) $, but $%
L^2\left( \mathbb{T}\right) $ plays a crucial intermediate role in the
construction of orthogonal wavelets, i.e., functions $\psi _i\in L^2\left( %
\mathbb{R}\right) $, $i=1,\dots ,N-1$, such that the triple-indexed family 
\begin{equation*}
\left\{ N^{\frac j2}\psi _i\left( N^jx+k\right) \right\} _{i,j,k}\subset
L^2\left( \mathbb{R}\right) \text{,\quad }j,k\in \mathbb{Z}\text{,}
\end{equation*}
is an orthonormal basis. (See \cite{Dau}, \cite{CoRy}, \cite{Jor95}, and 
\cite{Ho96} for details.)

In the papers \cite{JoPe96} and \cite{JoPe94}
(see also \cite{Ban91}, \cite{Ban96}, \cite{LaWa96}, \cite
{Ke}, \cite{Han92} and \cite{Str94}),
the higher dimensional version of this problem is
studied, i.e., $\mathbb{R}^\nu $\ $\left( \nu >1\right) $ in place of $%
\mathbb{R}$, and then the finite residues $\mathbb{Z\diagup }N\mathbb{Z}$
will be replaced with a set of representatives from $\mathbb{Z}^\nu \diagup 
\mathbf{N}\left( \mathbb{Z}^\nu \right) $, where now $\mathbf{N}=\left(
n_{ij}\right) _{i,j=1}^\nu $ is a $\nu $ by $\nu $ matrix with integral
entries, and $\left| \det \left( \mathbf{N}\right) \right| >1$. In this
setting, the digit set $D$ will be a subset of $\mathbb{Z}^\nu $ consisting
of mutually incongruent points modulo $\mathbf{N}\left( \mathbb{Z}^\nu
\right) $. The cardinality of $D$ must then equal $\left| \det \left( 
\mathbf{N}\right) \right| $. (See also \cite{LaWa96} and \cite{Ke}.) But the
case when the cardinality of $D$ is less than $\left| \det \left( \mathbf{N}%
\right) \right| $ is also interesting and studied in \cite{JoPe96} and \cite
{JoPe94} as well. In that case, the $\nu $-torus $\mathbb{T}^\nu :=\mathbb{R}%
^\nu \diagup \mathbb{Z}^\nu $ will be replaced with an associated fractal,
and the Haar measure on $\mathbb{T}^\nu $ with a corresponding fractal
measure $\mu $. The Hilbert space will be $L^2\left( \mu \right) $ in this
generalized setting; but the big difference, from the present setting to the
more general one, is that there is \emph{not} an analogue of the Fourier
basis $z^n=e^{int}$. In fact, we show in \cite{JoPe94} that, for the generic
fractal case, there will \emph{not} be an orthonormal basis for $L^2\left(
\mu \right) $ of pure frequencies $\left\{ e^{i\lambda t}\right\} $. We will
study these representations in a later paper.

\section{\label{Sec2}Permutative representations of $\mathcal{O}_N$\label
{Permutative}}

A representation of $\mathcal{O}_N$ on a Hilbert space $\mathcal{H}$ is said
to be \emph{permutative} if there is an orthonormal basis $\left\{
e_n\right\} _{n=1}^\infty $ for $\mathcal{H}$ such that 
\begin{equation}
S_ke_n\in \left\{ e_m\mid m\in \mathbb{N}\right\}  \label{Eq2.1}
\end{equation}
for $k\in \mathbb{Z}_N$, $n\in \mathbb{N}$, where the $S_k$'s are the
operators in $\mathcal{H}$ which define the above mentioned $\mathcal{O}_N$
representation. (Here of course $\mathbb{N}$ serves merely as a notation for
some generic countable index set.)
Then there exist maps $\sigma _k:\mathbb{N}\rightarrow %
\mathbb{N}$ for $k\in \mathbb{Z}_N$ such that 
\begin{equation}
S_ke_n=e_{\sigma _k\left( n\right) }  \label{Eq2.2}
\end{equation}
The Cuntz relations 
\begin{equation}
S_i^{*}S_j^{}=\delta _{ij}\openone \text{,\quad }\sum_{i=1}^NS_i^{}S_i^{*}=
\openone
\label{Eq2.3}
\end{equation}
immediately imply that 
\begin{subequations}
\label{Eq2.4}
\begin{align}
& \sigma _k:\mathbb{N}\rightarrow \mathbb{N}\text{ is injective, }k\in %
\mathbb{Z}_N\text{,} \label{Eq2.4a} \\
& \sigma _k\left( \mathbb{N}\right) \cap \sigma _l\left( \mathbb{N}\right)
=\emptyset \text{ for }k\neq l\text{,} \label{Eq2.4b} \\
& \bigcup_{\hbox to0pt{\hss $\scriptstyle k\in \mathbb{Z}_N$\hss}}\sigma
_k\left( \mathbb{N}\right) =\mathbb{N}\text{.} \label{Eq2.4c}
\end{align}
\end{subequations}
Conversely, if the maps $\sigma _k$ satisfy conditions (\ref{Eq2.4}) one
verifies that the operators $S_k$ defined by (\ref{Eq2.2}) satisfy the Cuntz
relations.

Next define a map $\sigma $ from the index set $\mathbb{N}$ into $\mathbb{Z}%
_N^\infty $ as follows. If $m\in \mathbb{N}$, $m$ corresponds to the
sequence $\left( j_1,j_2,\dots \right) $ defined inductively as follows: $%
j_1 $ is the unique $j$ such that 
\begin{equation}
S_j^{*}e_m\neq 0\text{,}  \label{Eq2.5}
\end{equation}
and then $S_j^{*}e_m\in \left\{ e_n\mid n\in \mathbb{N}\right\} $. When $%
j_1,\dots ,j_{k-1}$ are defined, let $j_k$ be the unique $j\in \mathbb{Z}_N$
such that 
\begin{equation}
S_{j_k}^{*}S_{j_{k-1}}^{*}\cdots S_{j_1}^{*}e_m\neq 0\text{.}  \label{Eq2.6}
\end{equation}
This definition has an obvious translation in terms of the maps $\sigma _k$,
since the sequence $\sigma \left( m\right) =\left( j_1,j_2,\dots \right) $
can be defined in terms of the $\sigma _i$: Define $j_1$ as the unique $j\in %
\mathbb{Z}_N$ such that there is an $m_1\in \mathbb{N}$ with $\sigma
_{j_1}\left( m_1\right) =m$, then $j_2$ as the unique $j$ such that there is
an $m_2\in \mathbb{N}$ with $\sigma _{j_2}\left( m_2\right) =m_1$, etc.

\begin{definition}
\label{Def2.1}We will say that a collection of $N$ maps $\sigma _1,\dots
,\sigma _N$ satisfying (\ref{Eq2.4}) is a \emph{branching function system,}
or simply a \emph{function system} of \emph{order} $N$. The associated map $%
\sigma :\mathbb{N}\rightarrow \mathbb{Z}_N^\infty $ is called the \emph{%
coding map} of the function system. We say that a function system is \emph{%
mul\-ti\-plic\-i\-ty-free} if the coding map is injective. We say that the
coding map $\sigma $ is \emph{partially injective} if it satisfies the
condition that if $n\in \mathbb{N}$ and $i_1,\dots ,i_k\in \mathbb{Z}_N$ and 
$\sigma \left( n\right) =\sigma \left( \sigma _{i_1}\cdots \sigma
_{i_k}\left( n\right) \right) $, then $n=\sigma _{i_1}\cdots \sigma
_{i_k}\left( n\right) $, and the function system is then said to be \emph{%
regular.}
\end{definition}

We will see later, in Theorem \ref{Thm2.5}, that mul\-ti\-plic\-i\-ty-free
function systems define representations of $\mathcal{O}_N$ and $\limfunc{UHF}%
\nolimits_N$ which are mul\-ti\-plic\-i\-ty-free in the sense of
representation theory. We will therefore already now say that a permutative
representation is \emph{mul\-ti\-plic\-i\-ty-free} if its function system is
so. Similarly the term \emph{regular} will be used for permutative
representations with regular function systems.

Following \cite{Str95}, the terminology ``branching reverse iterated
function system'' may be more appropriate than ``branching function
system'', but we will keep to the shorter term.

Clearly a mul\-ti\-plic\-i\-ty-free function system is regular.

Note that the notion of function system of order $N$ is closed under taking
disjoint union, while the notion of mul\-ti\-plic\-i\-ty-free function
system is clearly not; so this notion is a proper restriction.

We will show in Remark \ref{PermRepNotReg} that there are branching function
systems which are not regular.

We will now define two equivalence relations $\sim $ and $\approx $ on $%
\mathbb{N}$; and, by transporting these with $\sigma $, we get two
corresponding equivalence relations on $\sigma \left( \mathbb{N}\right)
\subset \mathbb{Z}_N^\infty $ which will also be denoted by $\sim $ and $%
\approx $ respectively. We say that $n\sim m$ if there are $I=\left(
i_1,\dots ,i_k\right) $ and $J=\left( j_1,\dots ,j_l\right) $ such that 
\begin{equation}
e_n=S_I^{}S_J^{*}e_m\text{.}  \label{Eq2.7}
\end{equation}
Then 
\begin{equation*}
e_m=S_J^{}S_I^{*}e_n
\end{equation*}
so $\sim $ is symmetric, and since the product of two monomials of the form $%
S_I^{}S_J^{*}$ is a monomial of the same form, $\sim $ is transitive. Thus $%
\sim $ is an equivalence relation.

We say that $n\approx m$ if furthermore $I,J$ may be chosen with $k=\left|
I\right| =\left| J\right| =l$. Then $\approx $ is clearly a stronger
equivalence relation.

We now characterize $\sim $ and $\approx $ on $\sigma \left( \mathbb{N}%
\right) \subset \mathbb{Z}_N^\infty $. Roughly, $n\approx m$ if the tails of
the corresponding sequences in $\mathbb{Z}_N^\infty $ are identical, and $%
n\sim m$ if the tails are identical up to translation.

\begin{proposition}
\label{Pro2.2}Assume that the coding map $\sigma :\mathbb{N}\rightarrow %
\mathbb{Z}_N^\infty $ is injective. The following two conditions are
equivalent for $n,m\in \mathbb{N}$.

\begin{enumerate}
\item  $\vphantom{\Bigl(}n\sim m$.\label{Pro2.2a(i)}

\item  There is a $\vphantom{\Bigl(}k\in \mathbb{Z}$ and an $n_0\in %
\mathbb{N}$ with $n_0>\left| k\right| $ such that \label{Pro2.2a(ii)} 
\begin{equation}
\sigma \left( n\right) _i=\sigma \left( m\right) _{i+k}  \label{Eq2.8}
\end{equation}
$\vphantom{\smash[t]{\Bigl(}}$for $i>n_0$.
\end{enumerate}

$\vphantom{\Bigl(}$Also, the following two conditions are equivalent for $%
n,m\in \mathbb{N}$.

\begin{enumerate}
\item  $\vphantom{\Bigl(}n\approx m$.\label{Pro2.2b(i)}

\item  There is an $\vphantom{\Bigl(}n_0\in \mathbb{N}$ such that \label
{Pro2.2b(ii)} 
\begin{equation}
\sigma \left( n\right) _i=\sigma \left( m\right) _i  \label{Eq2.9}
\end{equation}
$\vphantom{\smash[t]{\Bigl(}}$for $i>n_0$.
\end{enumerate}

$\vphantom{\smash[b]{\Bigl(}}$If the coding map is not assumed to be
injective, only the implications
\textup{(\ref{Pro2.2b(i)})}$\Rightarrow $\textup{(\ref{Pro2.2b(ii)})}
are valid.
\end{proposition}

\TeXButton{Begin Proof}{\begin{proof}}We first prove the first statement.
Assume that $n\sim m$, i.e., there are finite strings $I,J$ in $\mathbb{Z}_N$
with 
\begin{equation*}
e_n=S_I^{}S_J^{*}e_m\text{.}
\end{equation*}
But then 
\begin{equation*}
S_I^{*}e_n=S_J^{*}e_m\neq 0\text{.}
\end{equation*}
But this means 
\begin{align*}
\sigma \left( n\right) _i& =i_i\text{\quad for }i=1,\dots ,k\text{,} \\
\sigma \left( m\right) _j& =j_i\text{\quad for }i=1,\dots ,l\text{,}
\end{align*}
and thereafter the remaining parts of the strings $\sigma \left( n\right) $
and $\sigma \left( m\right) $ are identical,which is (\ref{Pro2.2a(ii)}).
This argument did not require the function system to be
mul\-ti\-plic\-i\-ty-free. Reverting this argument, if $\sigma \left(
n\right) =I\mathcal{K}$, $\sigma \left( m\right) =J\mathcal{K}$ where $%
\mathcal{K}$ is an infinite tail-string, it follows that if 
\begin{equation*}
e_{n^{\prime }}=S_I^{*}e_n\text{,\quad }e_{m^{\prime }}=S_J^{*}e_m\text{,}
\end{equation*}
then $\sigma \left( n^{\prime }\right) 
=\mathcal{K}=\sigma \left( m^{\prime }\right) $,
and hence if the function system is mul\-ti\-plic\-i\-ty-free, $%
n^{\prime }=m^{\prime }$ and 
\begin{equation*}
S_I^{*}e_n=S_J^{*}e_m\text{.}
\end{equation*}
Hence 
\begin{equation*}
e_n=S_I^{}S_J^{*}e_m\text{,}
\end{equation*}
or $n\sim m$.

The second statement in the proposition is proved in the same way, only with
the proviso of choosing $I,J$ with $\left| I\right| =\left| J\right| $
everywhere.\TeXButton{End Proof}{\end{proof}}

Another simple characterization of the equivalence relations $\sim $ and $%
\approx $, which does not depend on the function system being
mul\-ti\-plic\-i\-ty-free, is the following

\begin{proposition}
\label{Pro2.3}Assume that $\sigma _1,\dots ,\sigma _N$ is a branching
function system. The equivalence relation $\sim $ on $\mathbb{N}$ is the
equivalence relation generated by the relations 
\begin{equation}
\sigma _i\left( n\right) \sim n  \label{Eq2.10}
\end{equation}
for $i\in \mathbb{Z}_N$, $n\in \mathbb{N}$.

The equivalence relation $\approx $ on $\mathbb{N}$ is the equivalence
relation defined by the relations 
\begin{equation}
\sigma _{i_1}\sigma _{i_2}\cdots \sigma _{i_k}\left( n\right) \approx \sigma
_{j_1}\sigma _{j_2}\cdots \sigma _{j_k}\left( n\right)  \label{Eq2.11}
\end{equation}
for $i_l,j_l\in \mathbb{Z}_N$, $n\in \mathbb{N}$, $k=0,1,2,\dots $.
\end{proposition}

\TeXButton{Begin Proof}{\begin{proof}}Let us prove the first statement. One
clearly has $\sigma _i\left( n\right) \sim n$ by the defining relation (\ref
{Eq2.2}) 
\begin{equation*}
S_ie_n=e_{\sigma _i\left( n\right) }\text{,\quad i.e. }S_i^{*}e_{\sigma
_i\left( n\right) }=e_n\text{,}
\end{equation*}
for $\sigma _i$. Conversely if $n\sim m$, there are $I=\left( i_1,\dots
,i_k\right) $, $J=\left( j_1,\dots ,j_l\right) $ with 
\begin{equation*}
e_n=S_I^{}S_J^{*}e_m\text{,}
\end{equation*}
but then $m$ and $n$ are connected by $k+l$ elements in $\mathbb{N}$ such
that any two successive elements are related by applying some $\sigma _i$ to
either one or the other. This proves the first statement. For the second
statement, note that one has $\sigma _i\left( n\right) \approx \sigma
_j\left( n\right) $, since 
\begin{equation*}
e_{\sigma _i\left( n\right) }=S_ie_n\text{,\quad }e_{\sigma _j\left(
n\right) }=S_je_n
\end{equation*}
implies 
\begin{equation*}
e_{\sigma _i\left( n\right) }=S_i^{}S_j^{*}e_{\sigma _j\left( n\right) }%
\text{.}
\end{equation*}
Iterating this, one shows $\sigma _{i_1}\cdots \sigma _{i_k}\left( n\right)
\approx \sigma _{j_1}\cdots \sigma _{j_k}\left( n\right) $. Conversely, if $%
n\approx m$, there are $I=\left( i_1,\dots ,i_k\right) $, $J=\left(
j_1,\dots ,j_k\right) $ with 
\begin{equation*}
e_n=S_I^{}S_J^{*}e_m
\end{equation*}
and so 
\begin{equation*}
S_I^{*}e_n=S_J^{*}e_m\neq 0\text{.}
\end{equation*}
But if $l\in \mathbb{N}$ is such that 
\begin{equation*}
S_I^{*}e_n=S_J^{*}e_m=e_l\text{,}
\end{equation*}
then 
\begin{equation*}
e_n=S_Ie_l\text{,\quad }e_m=S_Je_l\text{,}
\end{equation*}
i.e., 
\begin{align*}
n& =\sigma _{i_1}\sigma _{i_2}\cdots \sigma _{i_k}\left( l\right) \text{,} \\
m& =\sigma _{j_1}\sigma _{j_2}\cdots \sigma _{j_k}\left( l\right) \text{,}
\end{align*}
and hence $n,m$ are connected applying the proposed relation $k$ times.%
\TeXButton{End Proof}{\end{proof}}

\begin{remark}
\label{MinimalSets}It follows from (\ref{Eq2.10}) and (\ref{Eq2.4}) that the 
$\sim $-equiv\-a\-lence classes can be characterized as the minimal nonempty
sets $S$ in $\mathbb{N}$ with the property that $\bigcup_{i=1}^N\sigma
_i\left( S\right) =S$. This is proved as follows: if $S$ has this property,
it follows immediately from (\ref{Eq2.4}) that 
\begin{equation*}
x\in S\Longleftrightarrow \sigma _i\left( x\right) \in S\text{,}
\end{equation*}
so $S$ is a union of $\sim $-equivalence classes by (\ref{Eq2.10}).
Conversely, any $\sim $-class $T$ (and thus any union of such classes) has
the property $\bigcup_{i=1}^N\sigma _i\left( T\right) =T$, and the claim
follows. This property was also considered in \cite{Str95}, \cite{Ban96},
and \cite{LaWa96}.
\end{remark}

\begin{scholium}
\label{Sch2.4}The equivalence classes of $\sim $ and $\approx $ can also be
described as the orbits of certain actions of two semigroups $G$ and $G_0$
on $\mathbb{N}$. Let $R:\mathbb{N}\rightarrow \mathbb{N}$ be defined
(uniquely) as the joint left inverse of the maps $\sigma _1,\dots ,\sigma _N$%
: 
\begin{equation*}
R\sigma _k\left( m\right) =m
\end{equation*}
for $k\in \mathbb{Z}_N$, $m\in \mathbb{N}$. Then 
\begin{equation*}
\left( \sigma _{i_1}\cdots \sigma _{i_k}R^n\right) \left( \sigma
_{j_1}\cdots \sigma _{j_l}R^m\right) =\begin{cases} \sigma _{i_1}\cdots
\sigma _{i_k}\sigma _{j_{n+1}}\cdots \sigma _{j_l}R^m & \text{if
}l>n\text{,} \\ \sigma _{i_1}\cdots \sigma _{i_k}R^{m+n-l} & \text{if }l\leq
n\text{,} \end{cases}
\end{equation*}
so $\left\{ \sigma _{i_1}\cdots \sigma _{i_k}R^n\mid k,n\in \left\{
0\right\} \cup \mathbb{N}\text{, }i_j\in \mathbb{Z}_N\right\} $ form a
semigroup $G$ with left inverses. Let $G_0$ be the sub-semigroup formed by
the elements $\sigma _{i_1}\cdots \sigma _{i_k}R^k$. Using Proposition \ref
{Pro2.3} and its proof, one verifies that $m\sim n$ if and only if there is
a $g\in G$ with $m=gn$, and $m\approx n$ if and only if there is a $g_0\in
G_0$ with $m=g_0n$. In particular, the $\sim $-equiv\-a\-lence classes are
just the orbits in $\mathbb{N}$ under the $G$ action, and the $\approx $%
-equiv\-a\-lence classes are the orbits under the $G_0$ action. This may
also be formulated as follows.
\end{scholium}

\begin{corollary}
\label{Cor2.4}Assume that $\sigma _1,\dots ,\sigma _N$ form a branching
function system. Let $R$ be as introduced in Scholium \ref{Sch2.4}. If $%
n,m\in \mathbb{N}$, then $n\sim m$ if and only if there are $k,l\in \left\{
0,1,\dots \right\} $ such that 
\begin{equation}
R^k\left( n\right) =R^l\left( m\right) \text{,}  \label{Eq2.11a}
\end{equation}
and $n\approx m$ if and only if there is a $k\in \left\{ 0,1,\dots \right\} $
such that 
\begin{equation}
R^k\left( n\right) =R^k\left( m\right) \text{.}  \label{Eq2.11b}
\end{equation}
\end{corollary}

It is clear from the definition of $\sim $ and $\approx $ that the closed
subspaces of $\mathcal{H}$ spanned by the $e_n$'s where $n$ runs through an
equivalence class in $\mathbb{N}$ are invariant under the action of $%
\mathcal{O}_N$ and $\limfunc{UHF}\nolimits_N$, respectively, and any vector $%
e_n$ is a cyclic vector for the corresponding subrepresentation. This in
itself is not enough to ensure irreducibility of these subrepresentations
(for example, the two matrices $\left( 
\begin{smallmatrix}
1 & 0 \\ 
0 & 1
\end{smallmatrix}
\right) $, $\left( \begin{smallmatrix}
0 & 1 \\ 
1 & 0
\end{smallmatrix}
\right) $ have the two vectors \smash[t]{$\left( 
\begin{smallmatrix}
1 \\ 
0
\end{smallmatrix}
\right) $, $\left( \begin{smallmatrix}
0 \\ 
1
\end{smallmatrix}
\right) $} as cyclic vectors, but do not generate an irreducible subalgebra
of $\mathcal{M}_2$). However, we can prove irreducibility of the
subrepresentations for permutative representations, and, moreover, if the
function system is mul\-ti\-plic\-i\-ty-free in the sense of Definition \ref
{Def2.1}, the representation is mul\-ti\-plic\-i\-ty-free in the usual
sense, i.e., the different subrepresentations are mutually non-equivalent.

\begin{theorem}
\label{Thm2.5}Consider a permutative representation of $\mathcal{O}_N$ on a
Hilbert space $\mathcal{H}$. Then the closure of any subspace of $\mathcal{H}
$ spanned by vectors $e_m$, where $m$ runs through a $\approx $-equivalence
class, is an irreducible $\limfunc{UHF}\nolimits_N$-module. If the function
system is regular, the closure of any subspace of $\mathcal{H}$ spanned by
vectors $e_m$, where $m$ runs through a $\sim $-equivalence class, is an
irreducible $\mathcal{O}_N$-module.

Furthermore, if the function system is mul\-ti\-plic\-i\-ty-free in the
sense of Definition \ref{Def2.1}, all the modules corresponding to different
equivalence classes are unitarily inequivalent in both cases.
\end{theorem}

\TeXButton{Begin Proof}{\begin{proof}}Assume first outright that the
function system is mul\-ti\-plic\-i\-ty-free. All the mentioned subspaces
are invariant, and hence to prove the statement it suffices to show that if $%
m_1,\dots ,m_n\in \mathbb{N}$ and $\psi _k\in \mathcal{H}\left( m_k\right)
\equiv $ subspace spanned by $e_l$ with $l\sim m_k$ (resp. $l\approx m_k$),
there is a polynomial $X$ in $S_I^{}S_J^{*}$'s (with $\left| I\right|
=\left| J\right| $ in the $\limfunc{UHF}\nolimits_N$ case) such that 
\begin{equation*}
Xe_{m_k}=\psi _k\text{.}
\end{equation*}
But for this it suffices, by linearity and re-indexing, to prove that for
any $m\in \mathbb{N}$ such that $m\sim m_1$ (resp. $m\approx m_1$) there is
a monomial $S_I^{}S_J^{*}$ (with $\left| I\right| =\left| J\right| $ in the $%
\limfunc{UHF}\nolimits_N$-case) such that 
\begin{equation}
S_I^{}S_J^{*}e_{m_i}=\begin{cases} e_m & \text{if }i=1\text{,} \\ 0 &
\text{otherwise.} \end{cases}  \label{Eq2.12}
\end{equation}
First note that as $m_1\sim m$ (or $m_1\approx m$) there is a monomial $%
S_I^{}S_J^{*}$ such that 
\begin{equation*}
S_I^{}S_J^{*}e_{m_1}=e_m
\end{equation*}
(with $\left| I\right| =\left| J\right| $ if $m_1\approx m$). Now, $J$ must
have the form 
\begin{equation*}
J=\left( j_1,\dots ,j_k\right)
\end{equation*}
where 
\begin{equation*}
\left( j_1,\dots ,j_k,j_{k+1},\dots \right) =\sigma \left( m_1\right)
\end{equation*}
is the sequence corresponding to $m_1$ in $\mathbb{Z}_N^\infty $. Now, by
injectivity of $\sigma $, we may choose $l\geq k$ so large that $\left(
j_1,\dots ,j_l\right) $ is different from the initial parts of the sequences 
$\sigma \left( m_i\right) $ for $i=2,\dots ,n$. Replacing $J$ by 
\begin{equation*}
\left( J,j_{k+1},\dots ,j_l\right) =\left( j_1,\dots ,j_l\right)
\end{equation*}
and $I$ by $\left( I,j_{k+1},\dots ,j_l\right) $ we still have 
\begin{equation*}
S_I^{}S_J^{*}e_{m_1}=e_m\text{,}
\end{equation*}
but now 
\begin{equation*}
S_I^{}S_J^{*}e_{m_i}=0
\end{equation*}
for $i=2,\dots ,n$. This ends the proof of (\ref{Eq2.12}) in the
mul\-ti\-plic\-i\-ty-free case.

Now, consider general branching function systems. Although the map $\sigma
:m\rightarrow \left( j_1,j_2,\dots \right) $ is then not necessarily
injective, it follows from Lemma \ref{Lem2.7}, below, that $\sigma $
restricted to each $\approx $-equivalence class is injective. Hence, if the $%
m_1,\dots ,m_n$ above are chosen from the same $\approx $-class, one argues
just as above that each $\approx $-class corresponds to an irreducible
representation of $\limfunc{UHF}\nolimits_N$ also in this case.

If $\sigma $ is regular, one likewise applies Lemma \ref{Lem2.7} on the $%
\sim $-equivalence classes.\TeXButton{End Proof}{\end{proof}}

\begin{lemma}
\label{Lem2.7} Assume that $\sigma :\mathbb{N}\rightarrow \mathbb{Z}%
_N^\infty $ is the coding map of a branching function system. Then $\sigma $ is
injective on each $\approx $-class, and if $\sigma $ is partially injective, 
$\sigma $ is injective on each $\sim $-class.
\end{lemma}

\TeXButton{Begin Proof}{\begin{proof}} If $n\sim m$, there are $I=\left(
i_1,\dots ,i_k\right) $, $J=\left( j_1,\dots ,j_l\right) $ (with $k=l$ if $%
n\approx m$) and a $p\in \mathbb{N}$ such that 
\begin{equation*}
n=\sigma _{i_1}\cdots \sigma _{i_k}\left( p\right) \text{,\quad }m=\sigma
_{j_1}\cdots \sigma _{j_l}\left( p\right)
\end{equation*}
but this means 
\begin{equation*}
\sigma \left( n\right) =IP\text{,\quad }\sigma \left( m\right) =JP
\end{equation*}
where $P=\sigma \left( p\right) $. Thus if $\sigma \left( n\right) =\sigma
\left( m\right) $, and $k=l$, then $I=J$, but then $n=\sigma _{i_1}\cdots
\sigma _{i_k}\left( p\right) =m$, and we have proved that $\sigma $ is
injective on $\approx $-classes. If, say, $k<l$, it follows from $IP=JP$
that $\sigma \left( p\right) =\sigma \left( \sigma _{j_{k+1}}\cdots \sigma
_{j_l}\left( p\right) \right) $. But partial injectivity then implies 
\begin{equation*}
p=\sigma _{j_{k+1}}\cdots \sigma _{j_l}\left( p\right)
\end{equation*}
and as $i_q=j_q$ for $q=1,\dots ,k$ it follows that 
\begin{equation*}
n=\sigma _I\left( p\right) =\sigma _I\sigma _{j_{k+1}}\cdots \sigma
_{j_l}\left( p\right) =\sigma _J\left( p\right) =m\text{,}
\end{equation*}
so $\sigma $ is injective on $\sim $-classes.\TeXButton{End Proof}
{\end{proof}}

\begin{remark}
\label{PermRepNotReg} Let us exhibit an example showing that not all
branching function systems are regular. Choose any $\sigma _1,\dots ,\sigma
_N$ satisfying (\ref{Eq2.4}), and assume that there is an $x\in \mathbb{N}$
with 
\begin{equation*}
\sigma _{i_1}\cdots \sigma _{i_k}\left( x\right) =x
\end{equation*}
where $k$ is odd (such examples abound by the initial remarks of Section \ref
{SecNew6}). Let $\Omega $ be the disjoint union of $\mathbb{N}$ with itself,
and let $\psi $ be the involution on $\Omega $ interchanging the two copies
of $\mathbb{N}$. Define $\bar{\sigma }_i$ on $\Omega $ by 
\begin{equation*}
\bar{\sigma }_i\left( y\right) =\sigma _i\left( \psi \left( y\right)
\right) =\psi \left( \sigma _i\left( y\right) \right)
\end{equation*}
for $y\in \Omega $. One checks readily that $\bar{\sigma }_i$ satisfy (%
\ref{Eq2.4}); and, if $y\in \Omega $, then $\sigma \left( y\right) $ is the
same regardless of whether $y$ is viewed as an element in $\Omega $ or as an
element in one of the two copies of $\mathbb{N}$ in $\Omega $. In
particular, if $x$ is the element mentioned above, sitting in one of the
copies of $\mathbb{N}$, then 
\begin{equation*}
\psi \left( x\right) =\bar{\sigma }_{i_1}\cdots \bar{\sigma }%
_{i_k}\left( x\right) \neq x\text{,}
\end{equation*}
but 
\begin{equation*}
\sigma \left( \psi \left( x\right) \right) =\sigma \left( x\right) \text{.}
\end{equation*}
Thus the function system $\bar{\sigma }_1,\dots ,\bar{\sigma }_N$
is not regular.

We have more remarks on nonregular function systems in Section
\ref{Concluding}, and we will see that the extension of Theorem
\ref{Thm2.5} to this case is nontrivial.
\end{remark}

\section{\label{Sec3}Monomial representations and integers modulo $N$\label
{Monomial}}

We now specialize to the case where we replace the index set of the basis $%
\left\{ e_n\right\} $ by $\mathbb{Z}$ instead of $\mathbb{N}$ (for
convenience) and define the maps $\sigma _k$ in (\ref{Eq2.4}) by 
\begin{equation}
\sigma _k\left( n\right) =nN+s_k\text{,}  \label{Eq3.1}
\end{equation}
where $s_1,\dots ,s_N$ are $N$ integers which are pairwise incongruent over $%
N$. The verification of the relations (\ref{Eq2.4}) (with $\mathbb{N}$
replaced by $\mathbb{Z}$ of course) is trivial. The associated
representation of $\mathcal{O}_N$ can then be realized on $L^2\left( %
\mathbb{T}\right) $ by the formula 
\begin{equation}
\left( S_k\xi \right) \left( z\right) =z^{s_k}\xi \left( z^N\right)
\label{Eq3.2}
\end{equation}
for $\xi \in L^2\left( \mathbb{T}\right) $, $z\in \mathbb{T}$, $k\in %
\mathbb{Z}_N$, and in this realization 
\begin{equation}
e_n\left( z\right) =z^n  \label{Eq3.3}
\end{equation}
for $n\in \mathbb{Z}$. Thus, this is one of the representations considered
in \cite{Jor95}. In this case the sequence $\sigma \left( m\right) =\left(
j_1,j_2,\dots \right) $ is determined by the requirement 
\begin{equation}
m=s_{j_1}+Ns_{j_2}+\dots +N^{k-1}s_{j_k} \mod{N^k}  \label{Eq3.4}
\end{equation}
for $k=1,2,\dots $. This follows from the computation 
\begin{equation}
S_k^{*}z^m= \begin{cases} z^n & \text{if }m=nN+s_k\text{ for some }n\in
\mathbb{Z}\text{,} \\ 0 & \text{otherwise,} \end{cases}  \label{Eq3.5}
\end{equation}
and thus 
\begin{setlength}{\multlinegap}{0pt}
\begin{multline}
S_{i_k}^{*}\cdots S_{i_1}^{*}z^m \\
=
\begin{cases}
0 & \text{if }m-s_{i_1}-Ns_{i_2}-\dots -N^{k-1}s_{i_k}\neq 0 \mod{N^k}%
\text{,} \\ 
z^{\frac{m-s_{i_1}-Ns_{i_2}-\dots -N^{k-1}s_{i_k}}{N^k}} & \text{if }%
m-s_{i_1}-Ns_{i_2}-\dots -N^{k-1}s_{i_k}=0 \mod{N^k}\text{.}
\end{cases}
  \label{Eq3.6}
\end{multline}
\end{setlength}

\begin{proposition}
\label{Pro3.1}The representation defined by \textup{(\ref{Eq3.1})} is
mul\-ti\-plic\-i\-ty-free.
\end{proposition}

\TeXButton{Begin Proof}{\begin{proof}}If $\sigma \left( m\right) =\sigma
\left( n\right) $ it follows from (\ref{Eq3.4}) that 
\begin{equation*}
m-n=0 \mod{N^k}
\end{equation*}
for $k=1,2,\dots $. Thus $m=n$ and $\sigma $ is injective.%
\TeXButton{End Proof}{\end{proof}}

The equivalence relations $\sim $ and $\approx $ on $\mathbb{Z}$ can in this
case be characterized by the remainders after $k$ divisions by $N$. More
precisely, if $m\in \mathbb{Z}$ and $\sigma \left( m\right) =\left(
j_1,j_2,\dots \right) $, define 
\begin{equation*}
R_k\left( m\right) =\left( m-s_{j_1}-Ns_{j_2}-\dots -N^{k-1}s_{j_k}\right)
N^{-k}\text{.}
\end{equation*}
Thus $R_k\left( m\right) $ is the remainder obtained after applying Euclid's
division algorithm $k$ times: 
\begin{equation*}
m=s_{j_1}+Ns_{j_2}+\dots +N^{k-1}s_{j_k}+N^kR_k\left( m\right) \text{.}
\end{equation*}
Because of the iterative nature of Euclid's algorithm, we have the semigroup
property 
\begin{equation*}
R_k\left( R_j\left( m\right) \right) =R_{k+j}\left( m\right)
\end{equation*}
and hence 
\begin{equation*}
R_k=\left( R_1\right) ^k=R^k
\end{equation*}
where $R=R_1$ is defined by 
\begin{equation*}
R\left( m\right) =\frac{m-s_k}N\text{\quad if }m=s_k\mod N\text{.}
\end{equation*}
Thus $R\sigma _k\left( m\right) =m$ for all $k\in \mathbb{Z}_N$, $m\in %
\mathbb{Z}$, so $R$ is a left inverse of each of the injective maps $\sigma
_1,\dots ,\sigma _N$. Thus $R$ is nothing but the map $R$ defined in general
in Scholium \ref{Sch2.4}. It follows immediately from Corollary \ref{Cor2.4}
that:

\begin{proposition}
\label{Pro3.2}Assuming \textup{(\ref{Eq3.1})}, the following two conditions are
equivalent.

\begin{enumerate}
\item  $\vphantom{\Bigl(}n\sim m$.\label{Pro3.2a(i)}

\item  $\vphantom{\Bigl(}R_k\left( m\right) =R_l\left( n\right) $ for some $%
k,l\in \mathbb{N}$\label{Pro3.2a(ii)}
\textup{(}and then $R_{k+i}\left( m\right)
=R_{l+i}\left( n\right) $
for $i=0,1,2,\dots \vphantom{\smash[t]{\Bigl(}}$\textup{)}.
\end{enumerate}

$\vphantom{\Bigl(}$Also the following two conditions are equivalent.

\begin{enumerate}
\item  $\vphantom{Bigl(}n\approx m$.\label{Pro3.2b(i)}

\item  $\vphantom{\Bigl(}R_k\left( m\right) =R_k\left( n\right) $ for some $%
k\in \mathbb{N}$. \label{Pro3.2b(ii)}
\end{enumerate}
\end{proposition}

Also, Proposition \ref{Pro2.3} has the following corollary:

\begin{corollary}
\label{Cor3.3}Assume \textup{(\ref{Eq3.1})}. Then the relation $\sim $ on $%
\mathbb{Z}$ is the equivalence relation generated by the relations 
\begin{equation*}
n\sim Nn+s_k
\end{equation*}
for $k\in \mathbb{Z}_N$, $n\in \mathbb{Z}$.
\end{corollary}

If $x$ is a nonnegative real number, let $\left\lfloor x\right\rfloor $
denote the largest nonnegative integer $\leq x$.

\begin{lemma}[Helge Tverberg]
\label{Lem3.4}Let $s=\max_k\left| s_k\right| $. Then the number of
equivalence classes in $\mathbb{Z}$ under $\sim $ is at most 
\begin{equation*}
1+2\left\lfloor \frac s{N-1}\right\rfloor \text{.}
\end{equation*}
\end{lemma}

\TeXButton{Begin Proof (due to Helge Tverberg)}
{\begin{proof}[Proof \textup{(}due to Helge Tverberg\/\textup{)}]}For
a given equivalence class,
let $a$ be the number in the class with smallest absolute value. If $s_k=a%
\mod N$, then $a$ is equivalent to $\frac{a-s_k}N$, and it follows that 
\begin{equation*}
\left| a\right| \leq \left| \frac{a-s_k}N\right| \text{,}
\end{equation*}
so 
\begin{equation*}
N\left| a\right| \leq \left| a-s_k\right| \text{.}
\end{equation*}
From this one deduces 
\begin{equation*}
\left| a\right| \leq \frac{\left| s_k\right| }{N-1}\text{,}
\end{equation*}
and hence 
\begin{equation*}
\left| a\right| \leq \left\lfloor \frac s{N-1}\right\rfloor
\end{equation*}
and the lemma follows.\TeXButton{End Proof}{\end{proof}}

\begin{corollary}
\label{Cor3.5}If $\pi $ is the permutative representation of $\mathcal{O}_N$
defined by \textup{(\ref{Eq3.1})}, then $\pi $ decomposes into at most $%
1+2\left\lfloor \frac s{N-1}\right\rfloor $ mutually inequivalent
irreducible representations.
Thus if $s=\max_k\left| s_k\right| <N-1$, then $\pi $ itself
is irreducible.
\end{corollary}

\TeXButton{Begin Proof}{\begin{proof}}This follows from Theorem \ref{Thm2.5}
in conjunction with Proposition \ref{Pro3.1} and Lemma \ref{Lem3.4}.%
\TeXButton{End Proof}{\end{proof}}

The description of the equivalence classes under $\approx $ is slightly more
complicated:

\begin{proposition}
\label{Pro3.6}Assume \textup{(\ref{Eq3.1}).} Then the following conditions are
equivalent for $n,m\in \mathbb{Z}$.

\begin{enumerate}
\item  $\vphantom{\Bigl(}n\approx m$.\label{Pro3.6(i)}

\item  There is a $\vphantom{\Bigl(}k\in \mathbb{N}$ and $i_1,\dots ,i_k\in %
\mathbb{Z}_N$ such that \label{Pro3.6(ii)} 
\begin{equation*}
n=s_{i_1}+Ns_{i_2}+\dots +N^{k-1}s_{i_k}+m-s_{j_1}-Ns_{j_2}-\dots
-N^{k-1}s_{j_k}\text{,}
\end{equation*}
$\vphantom{\Bigl(}$where $\left( j_1,j_2,\dots \right) =\sigma \left(
m\right) $.
\end{enumerate}
\end{proposition}

\TeXButton{Begin Proof}{\begin{proof}}This is a transcription of (\ref
{Eq2.11}) in Proposition \ref{Pro2.3} and Scholium \ref{Sch2.4} to the
present situation. In fact the relation in (\ref{Pro3.6(ii)}) just says
\newlength{\qedskip}
\settowidth{\qedskip}%
{$\displaystyle n=\sigma _{i_1}\sigma _{i_2}\cdots \sigma _{i_k}R^km\text{.}$}
\addtolength{\qedskip}{-\textwidth}
\begin{equation*}
n=\sigma _{i_1}\sigma _{i_2}\cdots \sigma _{i_k}R^km\text{.}%
\rlap{\hbox to-0.5\qedskip{\hfil \qed}}
\end{equation*}
\TeXButton{No QED Symbol}{\renewcommand{\qed}{}}\TeXButton{End Proof}
{\end{proof}}

\begin{remark}
\label{Rem3.7}Let $s=\max_k\left| s_k\right| $. Then Proposition \ref{Pro3.6}
can be used to show that the number of equivalence classes in $\mathbb{Z}$
under $\approx $ is at most $\left\lfloor \frac{4s}{N-1}\right\rfloor $, by
the following density argument:

For fixed $k\in \mathbb{N}$, put 
\begin{equation*}
F_k\left( m\right) =\left\{ \sigma _{i_1}\cdots \sigma _{i_k}R^km\mid
i_1,\dots ,i_k\in \mathbb{Z}_N\right\} \text{.}
\end{equation*}
Then all elements in $F_k\left( m\right) are$ $\approx $-equivalent to $m$,
and the $\approx $-equiv\-a\-lence class containing $m$ is $%
\bigcup_{k=0}^\infty F_k\left( m\right) $ by Proposition \ref{Pro3.6}. Note
that 
\begin{equation*}
\#\left( F_k\left( m\right) \right) =N^k\text{.}
\end{equation*}
Now 
\begin{equation*}
s+Ns+\dots +N^{k-1}s=s\frac{N^k-1}{N-1}<\frac s{N-1}N^k\text{,}
\end{equation*}
and it follows from Proposition \ref{Pro3.6}(\ref{Pro3.6(ii)}) that 
\begin{equation*}
F_k\left( m\right) \subseteq \left[ m-\frac{2s}{N-1}N^k,m+\frac{2s}{N-1}%
N^k\right] \text{.}
\end{equation*}
As $\#F_k\left( m\right) =N^k$, it follows that the mean density of $%
F_k\left( m\right) $ within 
\begin{equation*}
\left[ m-\frac{2s}{N-1}N^k,m+\frac{2s}{N-1}N^k\right]
\end{equation*}
is at most $\frac{N-1}{4s}$. Now if $m_1,\dots ,m_l$ are mutually non-$%
\approx $-equivalent elements in $\mathbb{Z}$, choose $k$ so large that $N^k$
is very large compared to the distances between the $m_i$'s. Then the
intervals above, for $m=m_i$, $i=1,\dots ,l$ are overlapping except for a
number of points which is very small compared to $N^k$. But as the $\approx $%
-equiv\-a\-lence classes corresponding to $m_1,\dots ,m_l$ are disjoint and
have mean density $\frac{N-1}{4s}$, it follows that 
\begin{equation*}
l\frac{N-1}{4s}\leq 1\text{,}
\end{equation*}
so 
\begin{equation*}
l\leq \frac{4s}{N-1}\text{.}
\end{equation*}
This ends the proof of the claim. As a corollary, %\label{Cor3.8} 
if $\pi $ is the restriction to $\limfunc{UHF}\nolimits_N$ of the
permutative representation of
$\mathcal{O}_N$ defined by \textup{(\ref{Eq3.1}),}
then $\pi $ decomposes into at most $\left\lfloor \frac{4s}{N-1}%
\right\rfloor $ mutually inequivalent irreducible representations. We will,
however, improve this estimate later by a completely different method; see
Corollary \ref{Cor3.10}.
\end{remark}

The method referred to in the end of the preceding remark might as well be
worked out in a setting where $\mathbb{Z}$ is replaced by $\mathbb{Z}^\nu $, 
$\mathbf{N}$ is a $\nu \times \nu $ matrix with integer entries and $\left|
\det \mathbf{N}\right| =N$, $D=\left\{ d_1,\dots ,d_N\right\} $ is a set of $%
N$ points in $\mathbb{Z}^\nu $ which are incongruent modulo $\mathbf{N}%
\mathbb{Z}^\nu $, and 
\begin{equation}
\sigma _i\left( n\right) =d_i+\mathbf{N}n  \label{Eq3.6bis}
\end{equation}
for $n\in \mathbb{Z}^\nu $. We assume that 
\begin{equation}
\bigcap_k\mathbf{N}^k\mathbb{Z}^\nu =0  \label{Eq3.7}
\end{equation}
and then the representation is mul\-ti\-plic\-i\-ty-free by the proof of
Proposition \ref{Pro3.1}; in fact (\ref{Eq3.7}) is both sufficient and
necessary for this. Note that (\ref{Eq3.7}) is not automatic: take for
example $\mathbf{N}=\left( n_{ij}\right) _{i,j=1}^\nu $ with $\nu \geq 2$, $%
n_{11}=N$, and $n_{ij}=\delta _{ij}$ otherwise. The condition (\ref{Eq3.7})
is fulfilled if $\left| \lambda \right| >1$ for all the eigenvalues $\lambda 
$ of $\mathbf{N}$, by using the Jordan canonical form for $\mathbf{N}$.
However, the condition (\ref{Eq3.7}) may be fulfilled even in cases when $%
\mathbf{N}$ has eigenvalues of modulus strictly less than $1$; see Section 
\ref{Appendix}.

It is convenient now to use the norm 
\begin{equation}
\left| n\right| =\sup_i\left| n_i\right|  \label{Eq3.8}
\end{equation}
on $\mathbb{Z}^\nu $. The corresponding norm on $\mathrm{GL}\left( \nu ,%
\mathbb{Z}\right) $ is 
\begin{equation}
\left\| \mathbf{N}\right\| =\sup_i\sum_j\left| n_{ij}\right| \text{.}
\label{Eq3.9}
\end{equation}

We first need

\begin{lemma}
\label{Lem3.7}If $\mathbf{N}$ is a $\nu \times \nu $ matrix with integer
entries and $\left| \det \mathbf{N}\right| =N$, and all the eigenvalues of $%
\mathbf{N}$ have modulus strictly greater than $1$, then there exists a
constant $C$ depending on $\mathbf{N}$ such that if $\sigma _i$ is defined by
\textup{(\ref{Eq3.6bis})}
and $d=\sup_i\left| d_i\right| $, then for any $x\in \mathbb{R}%
^\nu $ there is an $n_x\in \mathbb{N}$ such that 
\begin{equation*}
\left\| R^n\left( x\right) \right\| \leq Cd
\end{equation*}
for all $n\geq n_x$. Here $R:\mathbb{Z}^\nu \rightarrow \mathbb{Z}^\nu $ is
the map defined in Scholium \ref{Sch2.4}.
\end{lemma}

\TeXButton{Begin Proof}{\begin{proof}}In a suitable basis for $\mathbb{C}%
^\nu $, $\mathbf{N}$ may be written as a sum of blocks of the form 
\begin{equation*}
\begin{pmatrix}
\lambda & \varepsilon & 0 & 0 & \cdots & 0 & 0 \\ 
0 & \lambda & \varepsilon & 0 & \cdots & 0 & 0 \\ 
0 & 0 & \lambda & \varepsilon & \cdots & 0 & 0 \\ 
0 & 0 & 0 & \lambda & \cdots & 0 & 0 \\ 
\vdots & \vdots & \vdots & \vdots & \ddots & \vdots & \vdots \\ 
0 & 0 & 0 & 0 & \cdots & \lambda & \varepsilon \\ 
0 & 0 & 0 & 0 & \cdots & 0 & \lambda
\end{pmatrix}
\end{equation*}
where $\lambda $ is an eigenvalue and $\varepsilon $ is arbitrarily small.
Since the minimum of the modulus of the eigenvalues is strictly larger than
one, it follows that there is an
$\varepsilon >0$ and a norm on $\mathbb{C}^\nu 
$ such that 
\begin{equation*}
\left\| \mathbf{N}x\right\| \geq \left( 1+\varepsilon \right) \left\| x\right\|
\end{equation*}
for all $x\in \mathbb{C}^\nu $, in this norm. But this means 
\begin{equation*}
\left\| \mathbf{N}^{-1}\right\| <1
\end{equation*}
in the associated norm on $\mathcal{B}\left( \mathbb{C}^\nu \right) $.

We have 
\begin{equation*}
R\left( d_i+\mathbf{N}x\right) =R\left( \phi _i\left( x\right) \right) =x
\end{equation*}
for all $x\in \mathbb{Z}^\nu $, $i=1,\dots ,N$ so 
\begin{equation*}
R\left( d_i+\mathbf{N}x\right) -\mathbf{N}^{-1}\left( d_i+\mathbf{N}x\right)
=-\mathbf{N}^{-1}d_i
\end{equation*}
and hence 
\begin{equation*}
\left\| R\left( d_i+\mathbf{N}x\right) -\mathbf{N}^{-1}\left( d_i+\mathbf{N}%
x\right) \right\| \leq \left\| \mathbf{N}^{-1}\right\| \left\| d_i\right\|
\leq C^{\prime }d
\end{equation*}
for a constant $C^{\prime }$, where the last step use that all norms on $%
\mathbb{C}^\nu $ are equivalent. Thus 
\begin{equation*}
\left\| R\left( x\right) -\mathbf{N}^{-1}x\right\| \leq C^{\prime }d
\end{equation*}
for all $x\in \mathbb{Z}^\nu $. Hence 
\begin{align*}
\left\| R^n\left( x\right) -\mathbf{N}^{-n}\left( x\right) \right\| & \leq
\sum_{k=0}^{n-1}\left\| \mathbf{N}^{-k}\left( R^{n-k}\left( x\right) \right)
-\mathbf{N}^{-k-1}\left( R^{n-k-1}\left( x\right) \right) \right\| \\
& \leq \sum_{k=0}^{n-1}\left\| \mathbf{N}^{-1}\right\| ^k\left\| R\left(
R^{n-k-1}\left( x\right) \right) -\mathbf{N}^{-1}\left( R^{n-k-1}\left(
x\right) \right) \right\| \\
& \leq \sum_{k=0}^{n-1}\left\| \mathbf{N}^{-1}\right\| ^kC^{\prime }d \\
& \leq \frac 1{1-\left\| \mathbf{N}^{-1}\right\| }C^{\prime }d\text{.}
\end{align*}
But as $\lim_{n\rightarrow \infty }\left\| \mathbf{N}^{-n}\left( x\right)
\right\| =0$, it follows that 
\begin{equation*}
\limsup_{n\rightarrow \infty }\left\| R^n\left( x\right) \right\| \leq \frac{%
C^{\prime }}{1-\left\| \mathbf{N}^{-1}\right\| }d
\end{equation*}
and Lemma \ref{Lem3.7} is proved by noting that all norms on $\mathbb{C}^\nu 
$ are equivalent.\TeXButton{End Proof}{\end{proof}}

\begin{scholium}
\label{Sch3.8}Let 
\begin{equation*}
B=\left\{ x\in \mathbb{Z}^\nu \mid \left\| x\right\| \leq Cd\right\} \text{.}
\end{equation*}
The conclusion in Lemma \ref{Lem3.7} states that $R^n\left( x\right) $
ultimately is contained in $B$ as $n\rightarrow \infty $. But $B$ is a
finite subset of $\mathbb{R}^\nu $, and we have the estimate 
\begin{equation*}
\#\left( B\right) \leq \left( 2C+1\right) ^\nu d^\nu \text{.}
\end{equation*}
Since $B$ is finite, it follows that the sequence $R^n\left( x\right) $
ultimately becomes periodic, with period at most $\#\left( B\right) $. Call
this period $\limfunc{Per}\left( x\right) $. Now let $B_\infty $ be the set
of points $y\in B$ (equivalently $y\in \mathbb{Z}^\nu $) such that $%
n\rightarrow R^n\left( y\right) $ already is periodic. It follows that for
any $x\in \mathbb{Z}^\nu $, $R^n\left( x\right) \in B_\infty $ for
sufficiently large $n$. We have $\#\left( B_\infty \right) \leq \#\left(
B\right) \leq \left( 2C+1\right) ^\nu d^\nu $, and $B_\infty $ splits into a
finite number of finite orbits under $R|_{B_\infty }$. $R|_{B_\infty }$ is
actually a bijection.

Now consider points $x,y\in \mathbb{Z}^\nu $ such that $x\approx y$. It
follows from Proposition \ref{Pro2.3} and Scholium \ref{Sch2.4} that 
\begin{equation*}
y=\sigma _{i_1}\cdots \sigma _{i_n}R^nx
\end{equation*}
for suitable $n,i_1,\dots ,i_n$. By replacing $R^nx$ by $\sigma _{j_1}\cdots
\sigma _{j_m}R^{n+m}x$ where $\left( j_1,j_2,\dots \right) =\sigma \left(
R^nx\right) $ and $m$ is sufficiently large, we may assume that $n$ is so
large that $R^nx\in B_\infty $. Thus $x\approx y$ if and only if there is a $%
z\in B_\infty $ and $i_1,\dots ,i_n,j_1,\dots ,j_n\in \mathbb{Z}_N$ with 
\begin{equation*}
x=\sigma _{i_1}\cdots \sigma _{i_n}z\text{,\quad }y=\sigma _{j_1}\cdots
\sigma _{j_n}z\text{.}
\end{equation*}
But as $\sigma \left( z\right) $ is periodic with period $\limfunc{Per}%
\left( z\right) $, it follows also for example that 
\begin{equation*}
\sigma _{i_1}\cdots \sigma _{i_n}\sigma _{i_n+1}\cdots \sigma _{i_n+\limfunc{%
Per}\left( z\right) }\left( z\right) \approx \sigma _{j_1}\cdots \sigma
_{j_n}\left( z\right) \text{.}
\end{equation*}
It follows that the map $z\rightarrow z^{\approx }$ is a bijection between $%
B_\infty $ and the $\approx $-classes in $\mathbb{Z}^\nu $. In particular
there are exactly $\#\left( B_\infty \right) $ $\approx $-classes.

The set $B_\infty $ will be referred to as the \emph{set of periodic orbits}
in $\mathbb{Z}^\nu $. Applying $R$ gives a natural action of $\mathbb{Z}$ on 
$B_\infty $, and the cycles correspond to the orbits under this action. This
picture will be identified in the general setting of permutative
mul\-ti\-plic\-i\-ty-free representations in Section \ref{Sec4}, where also
nonperiodic orbits will occur.
\end{scholium}

\begin{corollary}
\label{Cor3.9}If $\mathbf{N}$ is a $\nu \times \nu $ matrix with integer
entries and $\left| \det \mathbf{N}\right| =N$, and all the eigenvalues of $%
\mathbf{N}$ have modulus greater than $1$, then the permutative
representation of $\limfunc{UHF}\nolimits_N$ defined by \textup{(\ref{Eq3.6})}
decomposes into a finite number of atoms. This number is dominated by $%
Cd^\nu $, where $C$ is a constant depending on $N$ and $d=\max_i\left\|
d_i\right\| $.
\end{corollary}

\TeXButton{Begin Proof}{\begin{proof}}By Theorem \ref{Thm2.5}, the number of
atoms is equal to the number of $\approx $-equiv\-a\-lence classes. But by
Scholium \ref{Sch3.8}, the number of $\approx $-equiv\-a\-lence classes is $%
\#\left( B_\infty \right) $, and 
\begin{equation*}
\#\left( B_\infty \right) \leq \left( 1+2C\right) ^\nu d^\nu
\end{equation*}
for the $C$ in Lemma \ref{Lem3.7}.\TeXButton{End Proof}{\end{proof}}

We will now prove the estimate announced in Remark
\ref{Rem3.7}. Note that this estimate on the number of atoms is better
than the estimate on the number of cycles in Lemma \ref{Lem3.4}.

\begin{corollary}
\label{Cor3.10}Consider the case $\nu =1$ and the permutative representation 
$\pi $ defined by \textup{(\ref{Eq3.1}).}
Put $D=\left\{ s_1,\dots ,s_N\right\} $ and
let $\limfunc{diam}\left( D\right) =\max D-\min D$ be the diameter of $D$.
Then the
number of equivalence classes in $\mathbb{Z}$ under $\approx $ is at most
equal to the number of integers in the interval
$\left[ -\frac {\max D}{N-1},-\frac {\min D}{N-1}\right]$,
and hence this number is dominated by $%
1+\left\lfloor \frac {\limfunc{diam}D}{N-1}\right\rfloor $.
Hence the restriction of $\pi $
to $\limfunc{UHF}\nolimits_N$ decomposes into at most
$1+\left\lfloor \frac {\limfunc{diam}\left( D\right) }{N-1}\right\rfloor $
irreducible representations, which are all mutually
non-equivalent.
\end{corollary}

\TeXButton{Begin Proof}{\begin{proof}}Put $s_-:=\min s_i$, $s_+:=\max s_i$.
A step for step repetition
of the proof of Lemma \ref{Lem3.7} in the case $\nu =1$, $\mathbf{N}=N$,
gives:
\begin{equation*}
-N^{-1}s_+\leq R\left( x\right) -N^{-1}x\leq -N^{-1}s_-\text{.}
\end{equation*}
Using
\begin{equation*}
R^n\left( x\right) -N^{-n}x=\sum_{k=0}^{n-1}N^{-k}%
\left( R^{n-k}\left( x\right) -N^{-1}R^{n-k-1}\left( x\right)\right) 
\end{equation*}
and
letting $n\rightarrow \infty $ one deduces that if $y$ is any limit point of 
$R^n\left( x\right) $ as $n\rightarrow \infty $ then
\begin{equation*}
-N^{-1}s_+\sum_{k=0}^\infty N^{-k}\leq y\leq -N^{-1}s_-\sum_{k=0}^\infty N^{-k}
\text{,}
\end{equation*}
so 
\begin{equation*}
-\frac{s_+}{N-1}\leq y\leq -\frac{s_-}{N-1}\text{.}
\end{equation*}
Thus the set $B_\infty $ of periodic points is contained in the interval
$\left[ -\frac{s_+}{N-1},-\frac{s_-}{N-1}\right] $.
This interval has length $\frac {\limfunc{diam}\left( D\right) }{N-1}$,
and thus 
\begin{equation*}
\#\left( B_\infty \right) \leq
1+\left\lfloor \frac {\limfunc{diam}\left( D\right) }{N-1}\right\rfloor
\text{.}
\end{equation*}
The rest follows from Scholium \ref{Sch3.8}.\TeXButton{End Proof}{\end{proof}}

Let the assumptions on the pair $\mathbf{N},D$ be as in Corollary \ref{Cor3.9}.
Then by a theorem of Bandt
(cf.~\cite{Ban91}, \cite{Ban96}, \cite{DDL95}, and \cite{Str94}) there is
a unique compact subset $\mathbf{T}=T\left( \mathbf{N},D\right)
\subset \mathbb{R}^\nu $ such that 
\begin{equation}
\mathbf{N}\left( \mathbf{T}\right) =\bigcup_{d\in D}\left( d+\mathbf{T}%
\right) \text{.}  \label{Eq3ins1}
\end{equation}
In fact 
\begin{equation}
\mathbf{T}=\left\{ \sum_{i=1}^\infty \mathbf{N}^{-i}d_i\biggm| d_i\in D\right\}
\text{.}  \label{Eq3ins2}
\end{equation}
Bandt showed that $\mathbf{T}$ is compact with non-empty interior in $%
\mathbb{R}^\nu $. If it is assumed that $\mathbf{T}$ has $\nu $-dimensional
Lebesgue measure equal to one, then it follows further that $\mathbf{T}$
must be a $\mathbb{Z}^\nu $-periodic tile. (See, e.g., \cite{Sen95},
\cite{Ke}, and \cite{GrMa92}.) 

Our next observation is that the asymptotics of $R^nx-\mathbf{N}^{-n}x$ for $%
x\in \mathbb{Z}^\nu $ derived in the proof of Lemma \ref{Lem3.7} above implies
that the lattice points in $-\mathbf{T}$ are exactly the periodic points under
$R$:

\begin{proposition}
\label{Pro3ins1}With the assumptions in Corollary \ref{Cor3.9}, we have 
\begin{equation}
B_\infty =\mathbb{Z}^\nu \cap \left( -\mathbf{T}\right) \text{.}
\label{Eq3ins3}
\end{equation}
\end{proposition}

\TeXButton{Begin Proof}{\begin{proof}}Assume first that $x\in B_\infty $.
We must show that $x\in -\mathbf{T}$.
But if $n=\limfunc{Per}\left( x\right) $, we have $R^nx=x$ and hence
$x=d_{i_1}+\mathbf{N}d_{i_2}+\dots +\mathbf{N}%
^{n-1}d_{i_n}+\mathbf{N}^nR^nx=d_{i_1}+\mathbf{N}%
d_{i_2}+\dots +\mathbf{N}%
^{n-1}d_{i_n}+\mathbf{N}^nx$ by the remark after Proposition \ref{Pro3.1}.
Thus 
\begin{equation}
-x=\mathbf{N}^{-1}d_{i_n}+\mathbf{N}^{-2}d_{i_{n-1}}+\dots +\mathbf{N}%
^{-n}d_{i_1}+\mathbf{N}^{-n}\left( -x\right)
\text{.}
\label{Eq3ins4}
\end{equation}
Iterating this expansion, and using (\ref{Eq3ins2}), we see that
$-x\in \mathbf{T}$,
and we have proved 
\begin{equation*}
B_\infty \subseteq \mathbb{Z}^\nu \cap \left( -\mathbf{T}\right) \text{.}
\end{equation*}

For the other inclusion of (\ref{Eq3ins3}), let $t\in \mathbf{T}\cap %
\mathbb{Z}^\nu $. Iterating (\ref{Eq3ins1}) in the form
$\mathbf{T}=\bigcup _{d\in D}\mathbf{N}^{-1}\left( d+\mathbf{T}\right) $,
we get, for all $n\in \mathbb{N}
$, points $d_1,\dots ,d_n\in D$, $t_n\in \mathbf{T}$, such that 
\begin{equation*}
t=\mathbf{N}^{-1}d_1+\mathbf{N}^{-2}d_2+\dots +\mathbf{N}^{-n}d_n+\mathbf{N}%
^{-n}t_n\text{,}
\end{equation*}
and therefore 
\begin{equation*}
-t_n=d_n+\mathbf{N}d_{n-1}+\dots +\mathbf{N}^{n-1}d_1+\mathbf{N}^n\left(
-t\right) \text{.}
\end{equation*}
Since the right-hand side is in $\mathbb{Z}^\nu $, we conclude that $t_n\in %
\mathbb{Z}^\nu $ for all $n\in \mathbb{N}$, and also 
\begin{equation}
R^n\left( -t_n\right) =-t\text{.}  \label{Eq3ins5}
\end{equation}
A recursion yields 
\begin{equation*}
R^n\left( -t_{n+m}\right) =-t_m\text{,}
\end{equation*}
$n,m\in \mathbb{N}$. Using that $t_n\in \mathbf{T}\cap \mathbb{Z}^\nu $ for
all $n\in \mathbb{N}$, we conclude that the pointset $\left\{ t_n\right\} $
must be finite. As a result, one of these points $s$, say, must satisfy $%
t_m=s=t_{n+m}$ for some $n\geq 1$. Therefore $R^n\left( -s\right) =R^n\left(
-t_{n+m}\right) =-t_m=-s$. Applying $R^m$ to both sides, we get 
\begin{equation*}
R^m\left( R^n\left( -s\right) \right) =R^m\left( -s\right) =-t\text{,}
\end{equation*}
and 
\begin{equation*}
R^m\left( R^n\left( -s\right) \right) =R^n\left( R^m\left( -s\right) \right)
=R^n\left( -t\right) \text{,}
\end{equation*}
proving $R^n\left( -t\right) =-t$; which is to say $-t\in B_\infty $. This
concludes the proof of the other inclusion in (\ref{Eq3ins3}).%
\TeXButton{End Proof}{\end{proof}}

We also have the following general observation on the effect of an arbitrary
integral translation of the set $D\subset \mathbb{Z}^\nu $ when the $\nu
\times \nu $ matrix $\mathbf{N}$ is expansive.

\begin{corollary}
\label{Cor3ins2}Let $\mathbf{N},D$ be given as in Corollary \ref{Cor3.9},
and let $\mathbf{T}%
=T\left( \mathbf{N},D\right) $ and $B_\infty =\left( -\mathbf{T}\right) \cap %
\mathbb{Z}^\nu $ be the corresponding attractors. Let $p\in \mathbb{Z}^\nu $%
, and set $D^{\prime }:=D+p$. Then $\mathbf{T}^{\prime }=T\left( \mathbf{N}%
,D^{\prime }\right) =\mathbf{T}+\left( \mathbf{N}-\openone \right) ^{-1}p$,
and it follows that
$B_\infty ^{\prime }=B_\infty -\left( \mathbf{N}-\openone
\right) ^{-1}p$ iff
$p\in \left( \mathbf{N}-\openone \right) \left( \mathbb{Z}%
^\nu \right) $.
\end{corollary}

\begin{remark}
\label{Rem3ins3}{\bfseries Radix representations.}
Let an expansive $\nu \times \nu $ matrix over $\mathbb{Z}$
be given as before, and set $N=\left| \det \mathbf{N}\right| $. Let $%
D=\left\{ d_i\right\} _{i=1}^N\subset \mathbb{Z}^\nu $ be given such that $%
0\in D$, and let $D^{*}$ be the set of finite sequences of elements of $D$,
and consider the radix-representation $\rho =\rho _D:D^{*}\rightarrow %
\mathbb{Z}^\nu $ which sends $\delta =\left( d_{i_0},d_{i_1},\dots
,d_{i_k}\right) $ to the point 
\begin{equation*}
\rho \left( \delta \right) =d_{i_0}+\mathbf{N}d_{i_1}+\dots +\mathbf{N}%
^kd_{i_k}\text{.}
\end{equation*}
Finally let $L:=\rho \left( D^{*}\right) $.
Kenyon \cite{Ke} assumes that $\rho $ is injective on $D^{*}$, and he
uses this for constructing more general attractors $\mathbf{T}=T\left( 
\mathbf{N},D\right) $ with non-periodic tiling properties. Note that Kenyon's
assumption is strictly weaker than our present assumption that 
\begin{equation}
d_i-d_j\notin \mathbf{N}\mathbb{Z}^\nu
\text{\quad for\quad }i\neq j\text{.}  \label{Eq3ins6}
\end{equation}
If we now define $\sigma _k:L\rightarrow L$ by
\begin{equation*}
\sigma _k\left( x\right) =d_k+\mathbf{N}x\text{,\quad }x\in L\text{,}
\end{equation*}
the conditions (\ref{Eq2.4}) follow immediately, with $\mathbb{N}$
replaced by $L$, and all the general results in Section \ref{Sec2} on
branching function systems apply. In particular the map $%
R:L\rightarrow L$ defined in Scholium \ref{Sch2.4} is now defined
by the requirement $R\left( d_i+\mathbf{N}x\right) =x$ for $x\in L$.
The results of Lemma \ref{Lem3.7} and Corollaries \ref{Cor3.9} and
\ref{Cor3.10}
go over to this new situation with the obvious modifications. Our
result from Proposition \ref{Pro3ins1} above then takes the form $B_\infty
=\left( -\mathbf{T}\right) \cap L$ where $B_\infty $ denotes the points in $L
$ which are of finite period under $R$. If $\nu =1$ and $N$ is a positive
integer, $N\geq 2$, then Kenyon \cite[Theorem 15]{Ke} gives a necessary and
sufficient condition on a subset $D=\left\{ d_i\right\} _{i=1}^N\subset %
\mathbb{Z}$ of non-negative integers that $\rho _D$ be injective. It is a
condition on the divisors of the polynomial $p_D\left( x\right)
=\sum_{i=1}^Nx^{d_i}$ and examples are given where $\rho _D$ is one-to-one
on $D^{*}$, but (\ref{Eq3ins6}) is not satisfied. Take for example $%
D=\left\{ 0,1,8,9\right\} $ and $N=4$. This $D$ has only two distinct
residues $\bmod{\,4}$. (This example was also studied in \cite{JoPe94},
\cite{Od}, \cite{LaWa96}, \cite{GrMa92}, and
earlier papers by Jorgensen and Pedersen.)

One could go even one step further and forget the way the subset
$L\subseteq \mathbb{Z}^\nu $ was constructed. That is, replace $L$ by
a general subset $S$ of $\mathbb{R}^\nu $ and define $\sigma
_k:S\rightarrow S$ from a set $D=\left\{ d_k\right\} $ as before, but
only assume the invariance $S=\bigcup _k\sigma _k\left( S\right) $ of
Remark \ref{MinimalSets} and property (\ref{Eq2.4b}), so that all the
properties (\ref{Eq2.4}) are valid. Again the different elements in
$D$ do not necessarily correspond to different $\bmod{\,\mathbf{N}}$
classes.  These situations were considered in \cite{Od} and
\cite{LaWa96}.  For $\nu =1$ \cite{Od} studies the selfsimilar compact
solution $\mathbf{T}$ to $\mathbf{T}=\bigcup_{d\in D}N^{-1}\left(
d+\mathbf{T}\right) $, and establishes the representation
$\mathbf{T}=\left[ 0,1\right] +E=\bigcup_{e\in E}\left[ 0,1\right]
+e$, disjoint union, where $\left[ 0,1\right] $ is the unit interval
in $\mathbb{R}$ and $E$ is some finite set determined from the digits
$D$.  Odlyzko restricts his analysis to the case where $D$ consists of
nonnegative integers and $0\in D$.  For $\nu >1$, \cite{LaWa96}
generalizes this project to higher dimensional ``octants'', but the
point is that $E$ can be related to the integral points in
$\mathbf{T}$.  Our approach to the integral points in $\mathbf{T}$ is
direct, general, and independent, and we hope to follow up on
implications for the Odlyzko-problem from our method.  We are indebted
to Yang Wang for explaining \cite{LaWa96} to us. The method which we
propose here for analyzing cycles is also related to (but different
from) the ``remainder-set algorithm'' used for radix-representations
of quadratic fields; see, e.g., \cite{Gil96}. Gilbert studies
expansions of the form $z=\sum _{j=0}^sd_{i_j}b^j$ for $z$ a Gaussian
integer, i.e., $z\in \mathbb{Z}\left[ i\right] $, where $b$ is a fixed
Gaussian integer, and $D=\left\{ d_i\mid i=1,2,\dots ,b\bar{b}\right\}
$ a digit set. Gilbert requires that the pair $\left\{ b,D\right\} $
be chosen such that every $z\in \mathbb{Z}\left[ i\right] $ has a
unique expansion. It is known that every $z\in \mathbb{C}$ will then
have an expansion $z=\sum _{j=-\infty }^sd_{i_j}b^j$, but the latter
are not unique. See \cite{Gil96} for definitions and prior
literature. Our analysis of the affine case in two real dimensions
below, Subsections \ref{TilePlane} and \ref{BoxSpline}, has
applications in the complex basis problem. The condition Gilbert
places on the pair $\left\{ b,D\right\} $ is strictly stronger than
$D$ being a full set of residues, as his fraction sets $\mathbf{T}$
will always have Lebesgue measure equal to one, and so his
$\mathbf{T}$'s will tile the plane with $\mathbb{Z}\left[ i\right]
$. Let $b=\alpha +i\beta $. Then complex multiplication $z\mapsto bz$
turns into matrix multiplication in $\mathbb{R}^2$ with
$\mathbf{N}_b=\left(
\begin{smallmatrix}
\alpha & -\beta \\
\beta & \alpha
\end{smallmatrix}
\right) $. Our residue sets $D$ in $\mathbb{Z}^2$ may also be viewed
as subsets of $\mathbb{Z}\left[ i\right] $. With our equivalence
relation $\sim $, it follows that Gilbert's condition on the pair
$\left\{ b,D\right\} $, if $0\in D$, implies that $0^\sim
=\mathbb{Z}^2$ where $0^\sim $ denotes the $\sim $-equivalence class
of the origin $0$. Hence, we are in the single-atom case discussed
more generally in Subsection \ref{SingleAtom} below.

There are other operator-theoretic approaches to fractional expansions
in the literature, e.g., the recent ones \cite{BoCo95} and
\cite{KMW96}. But they are different from ours both in ideas and
scope.
\end{remark}

\begin{remark}
If we restrict the attention to cycles, the estimates in the proof of Lemma 
\ref{Lem3.4} become, if $a=d_k\mod{\mathbf{N}\mathbb{Z}^\nu}$, 
\begin{equation*}
\left\| a\right\| \leq \left\| \mathbf{N}^{-1}\left( a-d_k\right) \right\|
\leq \left\| \mathbf{N}^{-1}\right\| \left\| a-d_k\right\| \text{.}
\end{equation*}
If $\left\| \mathbf{N}^{-1}\right\| <1$ we thus obtain 
\begin{equation*}
\left\| a\right\| \leq \frac{\left\| d_k\right\| }{\left\| \mathbf{N}%
^{-1}\right\| ^{-1}-1}\leq \frac d{\left\| \mathbf{N}^{-1}\right\| ^{-1}-1}
\end{equation*}
and hence the number of $\sim $-equiv\-a\-lence classes in $\mathbb{Z}^\nu $
under $\sim $ is at most 
\begin{equation*}
\left( 1+2\left\lfloor \frac d{\left\| \mathbf{N}^{-1}\right\|
^{-1}-1}\right\rfloor \right) ^\nu \text{.}
\end{equation*}
But the trouble with this argument is that we do not necessarily have $%
\left\| \mathbf{N}^{-1}\right\| <1$, even when all eigenvalues of $\mathbf{N}
$ have modulus greater than $1$, and we will have to operate with equivalent
norms in the case $\left\| \mathbf{N}^{-1}\right\| \geq 1$.
\end{remark}

\begin{remark}
\label{MeanDensity}The proof indicated in Remark \ref{Rem3.7},
becomes even more problematical in the case $\nu >1$. The obvious
modification is to use $\nu $-dimensional boxes in place of intervals, so if 
$m=\left( m_1,\dots ,m_\nu \right) \in \mathbb{Z}^\nu $, the subset $%
F_k\left( m\right) $ will be contained in the box 
\begin{equation*}
\left\{ x=\left( x_1,\dots ,x_\nu \right) \in \mathbb{Z}^\nu \bigm| m_i-%
\tfrac{2d}{\left\| \mathbf{N}^{-1}\right\| ^{-1}-1}\left\| \mathbf{N}%
\right\| ^k\leq x_i\leq m_i+\tfrac{2d}{\left\| \mathbf{N}^{-1}\right\|
^{-1}-1}\left\| \mathbf{N}\right\| ^k\right\}
\end{equation*}
centered at $m$. Again $\#\left( F_k\left( m\right) \right) =N^k$ so the
mean density of $F_k\left( m\right) $ within the box is now 
\begin{equation*}
N^k\left( \frac{4d}{\left\| \mathbf{N}^{-1}\right\| ^{-1}-1}\left\| \mathbf{N%
}\right\| ^k\right) ^{-\nu }M.
\end{equation*}
But here the problem turns up: the mean density contains the factor 
\begin{equation*}
\left( \frac N{\left\| \mathbf{N}\right\| ^\nu }\right) ^k
\end{equation*}
and 
\begin{equation}
\frac N{\left\| \mathbf{N}\right\| ^\nu }=\frac{\left| \det \mathbf{N}%
\right| }{\left\| \mathbf{N}\right\| ^\nu }\leq 1  \label{Eq3.11}
\end{equation}
with equality only in ultraspecial circumstances implying that all the
eigenvalues of $\mathbf{N}$ have modulus $N^{\frac 1\nu }$. (These cases are
however interesting for other reasons; see the examples in
Subsections \ref{SubsecNew9.2},
\ref{TilePlane}, \ref{SingleAtom}, \ref{Multiple},
and \ref{BoxSpline}.) 
Thus the mean density tends generically to zero as $k\rightarrow
\infty $, and the argument in the Remark is invalid as stated. However, if
the matrix $\mathbf{N}$ is diagonalizable as a matrix over the field $%
\mathbb{C}$, the situation can be remedied as follows: first note that if $%
m\in \mathbb{Z}^\nu $ and $\sigma \left( m\right) =\left( j_1,j_2,j_3,\dots
\right) $, then the set $F_k\left( m\right) $ consists of all points of the
form 
\begin{equation}
\sigma _{i_1}\cdots \sigma _{i_k}R^km=d_{i_1}+\mathbf{N}d_{i_2}+\dots +%
\mathbf{N}^{k-1}d_{i_k}+m-d_{j_1}-\mathbf{N}d_{j_2}-\dots -\mathbf{N}%
^{k-1}d_{j_k}\text{.}  \label{Eq3.13}
\end{equation}
Assume that $\mathbf{N}$ diagonalizes over $\mathbb{C}$, i.e., there is a
basis $\phi _1,\dots ,\phi _\nu $ for $\mathbb{C}^\nu $ such that 
\begin{equation}
\mathbf{N}\phi _\mu =\lambda _\mu \phi _\mu \text{,\quad }\mu =1,\dots ,\nu 
\text{.}  \label{Eq3.14}
\end{equation}
Then 
\begin{equation}
\prod_{\mu =1}^\nu \lambda _\mu =\det \left( \mathbf{N}\right) =N\text{.}
\label{Eq3.15}
\end{equation}
Assume furthermore that 
\begin{equation}
\left| \lambda _\mu \right| >1  \label{Eq3.16}
\end{equation}
for all $\mu $. This condition ensures the condition (\ref{Eq3.7}). Expand
the $d_i$ in the $\phi _\mu $: 
\begin{equation*}
d_i=\sum_\mu \delta _{i\mu }\phi _\mu \text{.}
\end{equation*}
Then 
\begin{equation*}
\mathbf{N}^kd_i=\sum_\mu \delta _{i\mu }\lambda _\mu ^k\phi _\mu \text{,}
\end{equation*}
so 
\begin{equation}
\left| \left( \mathbf{N}^kd_i\right) _\mu \right| \leq \limfunc{const.}\cdot\;
\left| d_i\right| \left| \lambda _\mu \right| ^k\leq \limfunc{const.}\cdot\;
d\left| \lambda _\mu \right| ^k\text{,}  \label{Eq3.17}
\end{equation}
where the constant only depends on $\mathbf{N}$. It follows that 
\begin{align}
\left| \left( \sigma _{i_1}\cdots \sigma _{i_k}R^k\left( m\right) -m\right)
_\mu \right| &\leq \limfunc{const.}\cdot\; 
d\left( 1+\left| \lambda _\mu \right|
+\left| \lambda _\mu \right| ^2+\dots +\left| \lambda _\mu \right|
^{k-1}\right)  \label{Eq3.18} \\
&=\limfunc{const.}\cdot\; d\frac{\left| \lambda _\mu \right| ^k-1}{\left|
\lambda _\mu \right| -1}\leq \limfunc{const.}\cdot\; d\left| \lambda _\mu
\right| ^k\text{,}  \notag
\end{align}
where the constants (which are different) only depend on $\mathbf{N}$. It
follows that $F_k\left( m\right) $ is contained in a ``complex
parallelepiped'' in $\mathbb{C}^\nu \supset \mathbb{Z}^\nu $ centered on $m$
whose linear dimension in the $\mu $'th coordinate direction is dominated by
a constant multiple of $d\left| \lambda _\mu \right| ^k$. Let $\eta _1,\dots
,\eta _\nu $ be the dual basis of $\phi _1,\dots ,\phi _\nu $, i.e., $\eta
_\mu $ is the linear functional defined by 
\begin{equation*}
\eta _\mu \left( \smash{\sum_{\mu ^{\prime }}\delta _{\mu ^{\prime }}\phi _{\mu
^{\prime }}}\vphantom{\sum \delta _{\mu ^{\prime }}\phi _{\mu
^{\prime }}}\right) =\delta _\mu \text{.}
\end{equation*}
Then $\eta _1,\dots ,\eta _\nu $ are linearly independent in $\mathbb{C}%
^{\nu *}\cong \mathbb{C}^\nu $. But if $\eta _1,\dots ,\eta _\nu $ is any
linearly independent set of vectors in $\mathbb{C}^\nu $, then at least one
of the $2^\nu $ combinations 
\begin{equation*}
\TeXButton{ReIm}{\genfrac{}{}{0pt}{}{\func{Re}}{\func{Im}}}\eta _1,%
\TeXButton{ReIm}{\genfrac{}{}{0pt}{}{\func{Re}}{\func{Im}}}\eta _2,\dots ,%
\TeXButton{ReIm}{\genfrac{}{}{0pt}{}{\func{Re}}{\func{Im}}}\eta _\nu
\end{equation*}
of vectors in $\mathbb{R}^\nu $ is linearly independent. This can be seen as
follows: viewing $\eta _1,\dots ,\eta _\nu $ as column vectors, we have 
\begin{equation*}
\det \left( \eta _1\eta _2\cdots \eta _\nu \right) \neq 0\text{.}
\end{equation*}
Decomposing $\eta _\mu =\limfunc{Re}\eta _\mu +i\func{Im}\eta _\mu $, this
determinant decomposes into a linear combination of the $2^\nu $
determinants 
\begin{equation*}
\det \left( \TeXButton{ReIm}{\genfrac{}{}{0pt}{}{\func{Re}}{\func{Im}}}\eta
_1\cdots \TeXButton{ReIm}{\genfrac{}{}{0pt}{}{\func{Re}}{\func{Im}}}\eta
_\nu \right) \text{,}
\end{equation*}
and hence one of the latter must be nonzero. Now, replacing each of our
original $\eta _\mu $'s by $\limfunc{Re}\eta _\mu $ or $\func{Im}\eta _\mu $%
, we may thus find $\nu $ $\mathbb{R}$-linearly independent functionals on $%
\mathbb{R}^\nu $, from now on called $\eta _1,\dots ,\eta _\nu $, such that 
\begin{equation*}
n\in F_k\left( m\right) \Rightarrow \left| \eta _\mu \left( n-m\right)
\right| \leq \limfunc{const.}\cdot\; d\left| \lambda _\mu \right| ^k\text{.}
\end{equation*}
It follows that $F_k\left( m\right) $ is contained in a parallelepiped in $%
\mathbb{R}^\nu $ whose volume grows like $d^\nu \left| \lambda _1\lambda
_2\cdots \lambda _\nu \right| ^k=d^\nu N^k$, where the last equality follows
from (\ref{Eq3.16}). But since $\#\left( F_k\left( m\right) \right) =N^k$ it
follows that the mean density of $F_k\left( m\right) $ in this
parallelepiped is dominated by $\limfunc{const.}\cdot\; d^{-\nu }$, where the
constant is independent of $k$. As $\left| \lambda _\mu \right| >1$ for all $%
\mu $, we conclude as in the Remark following Proposition \ref{Pro3.6} that
if $\pi $ is the restriction to $\limfunc{UHF}\nolimits_N$ of the
permutative representation of $\mathcal{O}_N$ defined by (\ref{Eq3.6}), then 
$\pi $ decomposes into at most $\limfunc{const.}\cdot\; d^\nu $ mutually
inequivalent irreducible representations, where $d=\max_i\left| d_i\right| $
and the constant only depends on $\mathbf{N}$.
\end{remark}

We finally remark on the difficulty of extending this density estimate
method to matrices $\mathbf{N}$ which are not diagonalizable: $\mathbf{N}$
is then equivalent to a matrix in Jordan canonical form, i.e., a direct sum
of matrices of the form 
\begin{equation*}
\begin{pmatrix}
\lambda & 1 & 0 & 0 & \cdots & 0 & 0 \\ 
0 & \lambda & 1 & 0 & \cdots & 0 & 0 \\ 
0 & 0 & \lambda & 1 & \cdots & 0 & 0 \\ 
0 & 0 & 0 & \lambda & \cdots & 0 & 0 \\ 
\vdots & \vdots & \vdots & \vdots & \ddots & \vdots & \vdots \\ 
0 & 0 & 0 & 0 & \cdots & \lambda & 1 \\ 
0 & 0 & 0 & 0 & \cdots & 0 & \lambda
\end{pmatrix}
\end{equation*}
The $k$'th power of this matrix has the first row 
\begin{equation*}
\left( \lambda ^k,\binom k1\lambda ^{k-1},\binom k2\lambda ^{k-2},\dots
\right) \text{,}
\end{equation*}
and hence in the estimate (\ref{Eq3.18}) the right-hand side will contain a
polynomial in $k$ of order at least $1$ as an extra factor. Thus the
parallelepiped containing $F_k\left( m\right) $ will grow so fast that the
mean density tends to zero.

We have seen how the set of $\approx $-classes provides a finite partition
of $\mathbb{Z}$; and we gave general conditions for finiteness in the case
of $\mathbb{Z}^\nu $. It is known (see \cite[pp. 11--13]{Fur81}) in general
that, for \emph{any} finite partition $\mathbb{Z}^\nu =\bigcup_{i=1}^qB_i$,
at least one of the $B_i$'s must contain arbitrarily long arithmetic
progressions: in fact, if $C\subset \mathbb{Z}^\nu $ is any finite
configuration, then some $B_i$ must contain a set of the form $aC+b$ for
some $a\in \mathbb{N}$, $b\in \mathbb{Z}^\nu $. The case $\nu =1$ is van der
Waerden's theorem to the effect that, for every $l\in \mathbb{N}$ some $B_i$
must contain a progression $\left\{ b,b+a,b+2a,\dots ,b+la\right\} $ with
the further property that $a\in B_i$. Our present finite partitions from the 
$\approx $ relation have much more explicit properties which we will take up
in Sections \ref{SecNew8} and \ref{SecNew9},
where we calculate various examples where
the relevant arithmetic progressions will be displayed with explicit
formulas.

\section{\label{Sec4}Cycles of irreducible $\mathcal{O}_N$ representations
and their \\atoms of $\limfunc{UHF}\nolimits_N$ representations\label{Cycles}}

Let us for the moment return to the general setup in Section \ref{Sec2}. Let 
$\sim $ and $\approx $ be the equivalence relations on $\mathbb{N}$
described by Proposition \ref{Pro2.3}, Scholium \ref{Sch2.4}, and Corollary 
\ref{Cor2.4}, and if $x\in \mathbb{N}$, let $x^{\sim }$ and $x^{\approx }$
denote its corresponding equivalence classes. Since $\approx $ is a finer
relation than $\sim $ we will use the terminology that the $\sim $-classes
are cycles, the $\approx $-classes are atoms, and an atom $x^{\approx }$ is
included in a cycle $y^{\sim }$ if $x^{\approx }\subseteq y^{\sim }$. This
terminology can also be transposed to the subrepresentations of $%
\mathcal{O}_N$, $\limfunc{UHF}\nolimits_N$ corresponding to the cycles and
atoms, respectively. We will now introduce an action of $\mathbb{Z}$ on the
atoms, which transforms each atom to an atom inside the same cycle, and
which acts transitively on the atoms within each cycle.

\begin{theorem}
\label{Thm4.1}Assume that the maps $\sigma _1,\dots ,\sigma _N:\mathbb{N}%
\rightarrow \mathbb{N}$ form a branching function system. Then 
\begin{equation}
\tau ^j\left( x^{\approx }\right) =\begin{cases} \left( R_{}^jx\right)
^{\approx } & \text{if }j\geq 0\text{,} \\ \left( \sigma _1^{\left| j\right|
}x\right) ^{\approx } & \text{if }j<0\text{,} \end{cases}  \label{Eq4.1}
\end{equation}
for $j\in \mathbb{Z}$, $x\in \mathbb{N}$ defines an action of $\mathbb{Z}$
on the atoms, which leaves the atoms within each cycle globally invariant,
and which is transitive on the atoms within each cycle.
\end{theorem}

\begin{remark}
The subindex $1$ chosen in the definition of $\tau ^j$ is not significant.
In fact, 
\begin{equation}
\tau ^{-j}\left( x^{\approx }\right) =\left( \sigma _{i_1}\sigma
_{i_2}\cdots \sigma _{i_j}\left( x\right) \right) ^{\approx }  \label{Eq4.2}
\end{equation}
for each $j>0$, $i_k\in \mathbb{Z}_N$.
\end{remark}

\TeXButton{Begin Proof of Theorem 4.1}
{\textit{Proof of Theorem \ref{Thm4.1}.\enspace}}Recall from Proposition \ref
{Pro2.3} and Scholium \ref{Sch2.4} that $x\sim y$ if and only if there is a
monomial $M$ in $R$'s and $\sigma _i$'s such that $y=Mx$, and $x\approx y$
if and only if this monomial can be chosen to contain the same number of $R$%
-factors as the number of $\sigma $-factors. Thus if $x\sim y$, then $Rx\sim
Ry$ and $\sigma _ix\sim \sigma _jy$ for any $i,j\in \mathbb{Z}_N$. We now
prove the same for $\approx $.

\begin{observation}
\label{Sec4Obs1}If $x\approx y$, then $Rx\approx Ry$.
\end{observation}

\TeXButton{Begin Proof}{\begin{proof}}If $x=y$, there is nothing more to
prove. If $x\neq y$ there is a $k\in \mathbb{N}$ and $i_1,\dots ,i_k\in %
\mathbb{Z}_N$ such that 
\begin{equation*}
x=\sigma _{i_1}\cdots \sigma _{i_k}R^ky\text{.}
\end{equation*}
Using $R\sigma _{i_1}=\limfunc{id}$, we obtain 
\begin{equation*}
Rx=\sigma _{i_2}\cdots \sigma _{i_k}R^{k-1}Ry\text{,}
\end{equation*}
and hence $Rx\approx Ry$.\TeXButton{No QED Symbol}{\renewcommand{\qed}{}}%
\TeXButton{End Proof}{\end{proof}}

\begin{observation}
\label{Sec4Obs2}If $x\approx y$, then $\sigma _ix\approx \sigma _jy$ for any 
$i,j\in \mathbb{Z}_N$.
\end{observation}

\TeXButton{Begin Proof}{\begin{proof}}There is a $k\in \mathbb{N}\cup
\left\{ 0\right\} $ and $i_1,\dots ,i_k\in \mathbb{Z}_N$ such that 
\begin{equation*}
x=\sigma _{i_1}\cdots \sigma _{i_k}R^ky\text{;}
\end{equation*}
but using $R\sigma _j=\limfunc{id}$, we then have 
\begin{equation*}
\sigma _ix=\sigma _i\sigma _{i_1}\cdots \sigma _{i_k}R^kR\sigma _jy\text{,}
\end{equation*}
so $\sigma _ix\approx \sigma _jy$.\TeXButton{No QED Symbol}
{\renewcommand{\qed}{}}\TeXButton{End Proof}{\end{proof}}

Observations \ref{Sec4Obs1} and \ref{Sec4Obs2} imply that the maps $\tau ^j$
in (\ref{Eq4.1}) are well defined on $\approx $-equiv\-a\-lence classes,
and, furthermore, the remark subsequent to the theorem is valid. Also, since
the analogous $\tau ^j$ on $\sim $-equiv\-a\-lence classes is the identity
map, the $\tau ^j$ maps the atoms within each cycle into themselves. It
remains to verify that $\tau $ defines an action of $\mathbb{Z}$, i.e.,

\begin{observation}
\label{Sec4Obs3}$\tau ^{i+j}=\tau ^i\tau ^j$ for all $i,j\in \mathbb{Z}$,
and $\tau ^0=\limfunc{id}$.
\end{observation}

\TeXButton{Begin Proof}{\begin{proof}}The last statement is trivial, and the
first statement is true if $i,j$ have the same sign, by (\ref{Eq4.1}). If $%
i\geq -j\geq 0$, then 
\begin{equation*}
\tau ^i\tau ^j\left( x^{\approx }\right) =\tau ^i\left( \sigma
_1^{-j}x\right) ^{\approx }=\left( R^i\sigma _1^{-j}x\right) ^{\approx
}=\left( R^{i+j}x\right) ^{\approx }=\tau ^{i+j}\left( x^{\approx }\right) 
\text{.}
\end{equation*}
If $-j\geq i\geq 0$, then 
\begin{equation*}
\tau ^i\tau ^j\left( x^{\approx }\right) =\left( R^i\sigma _1^{-j}x\right)
^{\approx }=\left( \sigma _1^{-j-i}x\right) ^{\approx }=\tau ^{i+j}\left(
x^{\approx }\right) \text{.}
\end{equation*}
If $j\geq -i\geq 0$, then 
\begin{align*}
\tau ^i\tau ^j\left( x^{\approx }\right) &=\tau ^i\left( R^jx\right)
^{\approx }=\left( \sigma _1^{-i}R^jx\right) ^{\approx } \\
&=\left( \sigma _1^{-i}R^{-i}R^{j+i}x\right) ^{\approx }=\left(
R^{i+j}x\right) ^{\approx }=\tau ^{i+j}\left( x^{\approx }\right) \text{.}
\end{align*}
If $-i\geq j\geq 0$, then 
\begin{equation*}
\tau ^i\tau ^j\left( x^{\approx }\right) =\left( \sigma _1^{-i}R^jx\right)
^{\approx }=\left( \sigma _1^{-i-j}\sigma _1^jR^jx\right) ^{\approx }=\left(
\sigma _1^{-i-j}x\right) ^{\approx }=\tau ^{i+j}\left( x^{\approx }\right) 
\text{.}
\end{equation*}
\TeXButton{No QED Symbol}{\renewcommand{\qed}{}}\TeXButton{End Proof}
{\end{proof}}

Thus $\tau $ is an action of $\mathbb{Z}$, and it remains to establish:

\begin{observation}
\label{Sec4Obs4}The action $\tau $ acts transitively on the atoms within
each cycle.
\end{observation}

\TeXButton{Begin Proof}{\begin{proof}}If $x\sim y$, then 
\begin{equation*}
x=\sigma _{i_1}\cdots \sigma _{i_k}R^ly
\end{equation*}
for suitable $i_1,\dots ,i_k\in \mathbb{Z}_N$, $k,l\in \mathbb{N}\cup
\left\{ 0\right\} $. But then 
\begin{equation*}
x^{\approx }=\tau ^{l-k}\left( y^{\approx }\right)
\end{equation*}
by the definitions (\ref{Eq4.1}), (\ref{Eq4.2}), and the relations $R\sigma
_i=\limfunc{id}$.\TeXButton{End Proof}{\end{proof}}

\section{\label{Sec5}Relations of finite type and sub-Cuntz states\label
{Relations}}

Let $\sim $ and $\approx $ be the relations on $\mathbb{N}$ coming from a
permutative representation. If $n\in \mathbb{N}$, let $\limfunc{Per}\left(
n\right) \in \mathbb{N}\cup \left\{ \infty \right\} $ be the period of $%
n^{\approx }$ under the action $\tau $ of $\mathbb{Z}$ on the $\approx $%
-equiv\-a\-lence classes defined in Theorem \ref{Thm4.1}, i.e., $\limfunc{Per%
}\left( n\right) $ is the number of atoms in the cycle consisting of the $%
\approx $-classes inside $n^{\sim }$. It follows from Theorem \ref{Thm4.1}
that $\limfunc{Per}\left( n\right) $ only depends on $n^{\sim }$. In this
section we will show that if $\limfunc{Per}\left( n\right) $ is finite, the
corresponding representation of $\mathcal{O}_N$ comes from a so-called
sub-Cuntz state, which is moreover a vector state $\left( e_m,\,\cdot
\,e_m\right) $ where we may choose $m\approx n$, and there are $i_1,\dots
,i_k\in \mathbb{Z}_N$ with 
\begin{equation}
S_{i_1}\cdots S_{i_k}e_m=e_m\text{.}  \label{Eq5.1}
\end{equation}

Let us first describe the sub-Cuntz states. If $m\in \mathbb{N}$, there is a
canonical embedding of $\mathcal{O}_{N^m}$ into $\mathcal{O}_N$ given by 
\begin{equation}
s_{i_1+Ni_2+\dots +N^{m-1}i_m}^{\left( N^m\right) }\hookrightarrow
s_{i_1}^{\left( N\right) }s_{i_2}^{\left( N\right) }\cdots s_{i_m}^{\left(
N\right) }\text{,}  \label{Eq5.2}
\end{equation}
where the super-index on $s$ refers to the sub-index on $\mathcal{O}$.
Because of the relation $\sum_is_i^{}s_i^{*}=\openone $,
the corresponding embedding of $\limfunc{UHF}%
\nolimits_{N^m} $ into $\limfunc{UHF}\nolimits_N$ is actually an
isomorphism. In the tensor product decomposition, this isomorphism is given
by 
\begin{equation}
\left( x_1\otimes \dots \otimes x_N\right) \otimes \left( x_{N+1}\otimes
\dots \otimes x_{2N}\right) \otimes \cdots \hookrightarrow x_1\otimes
x_2\otimes \cdots  \label{Eq5.3}
\end{equation}
for all $x_i\in \mathcal{M}_N$. Now we say that a state $\omega $ on $%
\mathcal{O}_N$ is a sub-Cuntz state of order $m$ if its restriction to $%
\mathcal{O}_{N^m}$ is a Cuntz state. (For general background facts on the
Cuntz algebra, its subalgebras, and its states, we refer to \cite{BrRo}.)
Let us summarize some known facts about Cuntz states in this setting:

\begin{proposition}
\label{Pro5.1}Let $\omega $ be a state on $\mathcal{O}_N$, let $\left( 
\mathcal{H},\pi ,\Omega \right) $ be the corresponding cyclic
representation,
let $m\in \mathbb{N}$ and let $\lambda _{i_1,\dots ,i_m}$ \textup{(}%
for $i_k\in \mathbb{Z}_N$\textup{)} be complex numbers. 
Then the following conditions
are equivalent. \TeXButton{Proposition 5.1}
{\begin{align}
&\Omega =\sum_{\hbox to0pt{\hss $\scriptstyle i_1,\dots ,i_m$\hss}
}\lambda _{i_1,\dots ,i_m}S_{i_1}\cdots
S_{i_m}\Omega \text{.}  \label{Eq5.4} \\
&S_{i_m}^{*}\cdots S_{i_1}^{*}\Omega =\lambda _{i_1,\dots ,i_m}\Omega
\text{\quad
for all }i_1,\dots ,i_m\in \mathbb{Z}_N\text{.}  \label{Eq5.5} \\
\intertext{Furthermore, these conditions imply that }
&\sum_{\hbox to0pt{\hss $\scriptstyle i_1,\dots ,i_m$\hss}
}\left| \lambda _{i_1,\dots ,i_m}\right| ^2=1\text{,}
\label{Eq5.6}
\end{align}} and the restriction of $\omega $ to $\limfunc{UHF}%
\nolimits_{N^m}=\limfunc{UHF}\nolimits_N$ is the pure product state
determined by the unit vector $\sum_{i_1,\dots ,i_m}\lambda _{i_1,\dots
,i_m}e_{i_1}\otimes \dots \otimes e_{i_m}$ in $\left( \mathbb{C}^N\right)
^{\otimes m}$.
\end{proposition}

\TeXButton{Begin Proof}{\begin{proof}}See \cite{Cun77}, \cite{BJ}, \cite{BJP}.%
\TeXButton{End Proof}{\end{proof}}

We will say that the state on $\mathcal{O}_N$ defined by a vector $\Omega $
with the properties of Proposition \ref{Pro5.1} is a sub-Cuntz state of
order $m$. Thus sub-Cuntz states of order $1$ are ordinary Cuntz states. It
should be emphasized that the conditions (\ref{Eq5.4})--(\ref{Eq5.6}) only
determine the restriction of $\omega $ to $\mathcal{O}_{N^m}$, and hence a
sub-Cuntz state is just any extension to $\mathcal{O}_N$ of an ordinary
Cuntz state on $\mathcal{O}_{N^m}$. The sub-Cuntz states we will consider in
the following will, however, be defined on all of $\mathcal{O}_N$ by
specific requirements.

\begin{proposition}
\label{Pro5.2}Let $\sim $ and $\approx $ be the equivalence relations on $%
\mathbb{N}$ coming from a branching function system. Let $n\in \mathbb{N}$,
and assume that $k=\limfunc{Per}\left( n\right) $ is finite. Then there
exists an $m\approx n$ and $i_1,\dots ,i_k\in \mathbb{Z}_N$ such that 
\begin{equation}
\sigma _{i_1}\cdots \sigma _{i_k}\left( m\right) =m\text{.}  \label{Eq5.7}
\end{equation}
Furthermore, $m$ and $\left( i_1\dots ,i_k\right) $ are uniquely determined
by $n$.
\end{proposition}

\TeXButton{Begin Proof}{\begin{proof}}By (\ref{Eq4.1}) we have 
\begin{equation*}
R^kn\approx n\text{,}
\end{equation*}
and then by (\ref{Eq2.11b}), 
\begin{equation*}
R^{k+l}n=R^ln
\end{equation*}
for all sufficiently large $l$. But again by Theorem \ref{Thm4.1}, $%
R^{jk}n\approx n$ for $j=1,2,\dots $, and hence putting 
\begin{equation*}
m=R^{jk}n
\end{equation*}
for some large $j$, we have 
\begin{equation*}
n\approx m=R^km\text{.}
\end{equation*}
But as 
\begin{equation*}
\sigma \left( m\right) =\left( i_1,i_2,\dots \right) \text{,}
\end{equation*}
we then have 
\begin{equation*}
m=\sigma _{i_1}\sigma _{i_2}\cdots \sigma _{i_k}R^km=\sigma _{i_1}\sigma
_{i_2}\cdots \sigma _{i_k}m\text{,}
\end{equation*}
and the proposition is proved.\TeXButton{End Proof}{\end{proof}}

\begin{corollary}
\label{Cor5.3}If $\sim $ and $\approx $ are the equivalence relations coming
from a permutative, regular representation, and $n\in \mathbb{N}$ with $%
\limfunc{Per}\left( n\right) <+\infty $, then there is a vector in the
representation of $\mathcal{O}_N$ corresponding to $n^{\sim }$ which defines
a sub-Cuntz state of order $k=\limfunc{Per}\left( n\right) $. Furthermore,
this vector can be taken to be $e_m$ where $e_m$ is in the subspace
associated with the irreducible representation of $\limfunc{UHF}\nolimits_N$
corresponding to $n^{\approx }$, and there are $i_1,\dots ,i_k\in \mathbb{Z}%
_N$ such that 
\begin{equation}
S_{i_1}\cdots S_{i_k}e_m=e_m\text{.}  \label{Eq5.8}
\end{equation}
\end{corollary}

\TeXButton{Begin Proof}{\begin{proof}}By Theorem \ref{Thm2.5}, this follows
from Proposition \ref{Pro5.2}. The relation (\ref{Eq5.8}) is a transcription
of (\ref{Eq5.7}), using (\ref{Eq2.2}).\TeXButton{End Proof}{\end{proof}}

\begin{corollary}
\label{Cor5.4}If $\mathbb{N}$ is replaced by $\mathbb{Z}$, and the $\sigma
_k $ are defined by \textup{(\ref{Eq3.1}),} 
\begin{equation*}
\sigma _k\left( n\right) =nN+s_k\text{,}
\end{equation*}
where $s_1,\dots ,s_N$ are pairwise incongruent over $N$, then the
corresponding representation of $\mathcal{O}_N$ decomposes into a finite sum
of pairwise non-equivalent irreducible representations which are all induced
by sub-Cuntz states.
\end{corollary}

\TeXButton{Begin Proof}{\begin{proof}}The first part of the corollary is a
restatement of Theorem \ref{Thm2.5}, in conjunction with Corollary \ref
{Cor3.5}. But by Corollary \ref{Cor3.10}, there is only a finite number of $%
\approx $-equiv\-a\-lence classes, so any element in $\mathbb{Z}$ has finite
period under $\tau $. Thus the ulterior statement of the corollary follows
from Corollary \ref{Cor5.3}.\TeXButton{End Proof}{\end{proof}}

\begin{remark}
\label{SubCuntzMatrix}{\bfseries Periodic points.}
Corollary \ref{Cor5.4} is of course also true in the
setting defined by (\ref{Eq3.6bis}) and (\ref{Eq3.7}), i.e., $\mathbb{N}$ is
replaced by $\mathbb{Z}^\nu $, $\mathbf{N}$ is a $\nu \times \nu $ matrix
with integer entries and $\left| \det \mathbf{N}\right| =N$, $D=\left\{
d_1,\dots ,d_N\right\} $ is a set of $N$ points in $\mathbb{Z}^\nu $ which
are incongruent modulo $\mathbf{N}\mathbb{Z}^\nu $, 
\begin{equation*}
\sigma _i\left( n\right) =d_i+\mathbf{N}n
\end{equation*}
for $n\in \mathbb{Z}^\nu $, and $\bigcap_k\mathbf{N}^k\mathbb{Z}^\nu =0$.
The latter condition entails that $\left( \openone -\mathbf{N}^k\right) ^{-1}$
exists on $\mathbb{R}^\nu $ for $k=1,2,\dots $. In this case the equation (%
\ref{Eq5.7}) is 
\begin{equation}
d_{i_1}+\mathbf{N}d_{i_2}+\dots +\mathbf{N}^{k-1}d_{i_k}+\mathbf{N}^km=m%
\text{,}  \label{Eq5.9}
\end{equation}
i.e., 
\begin{equation}
m=\left( \openone -\mathbf{N}^k\right) 
^{-1}\left( d_{i_1}+\mathbf{N}d_{i_2}+\dots +%
\mathbf{N}^{k-1}d_{i_k}\right) \text{.}  \label{Eq5.10}
\end{equation}
But this gives a method to find all the atoms of finite period: just check
for all $k\in \mathbb{N}$, $i_1,i_2,\dots ,i_k\in \mathbb{Z}_N$ whether the
right-hand side of (\ref{Eq5.10}) is in $\mathbb{Z}^\nu $. Corollary \ref
{Cor3.9} and Proposition \ref{Pro9.2} give a bound on how high-order $k$ one
has to check this. This method will be discussed further
in Section \ref{SecNew9}. We also note that (\ref{Eq5.10})
identifies the solution $m$ to (\ref{Eq5.9}) as an element in the
tile $-T\left( \mathbf{N},D\right) $ from (\ref{Eq3ins1}) as follows:
recalling that arbitrary elements $t$ in
$\mathbf{T}=T\left( \mathbf{N},D\right) $
are given as radix-representations 
\begin{equation}
t=\sum_{j=1}^\infty\mathbf{N}^{-j}\delta_j\text{\quad
with\quad }\delta_j\in D \label{Eq5.11}
\end{equation}
cf. (\ref{Eq3ins2}). Substitution of
$\left( \openone -\mathbf{N}^k\right) ^{-1}=-\sum_{j=1}^\infty\mathbf{N}^{-jk}$
into (\ref{Eq5.10})
then yields the following periodic form for the radix-representation
of $-m\in \mathbf{T}$, 
\begin{equation}
-m=\left( \smash{\underbrace{d_{i_1},\dots ,d_{i_k}},\underbrace{d_{i_1},\dots
,d_{i_k}},\underbrace{\vphantom{d_{i_1},\dots ,d_{i_k}}\dots
\vphantom{d_{i_1},\dots ,d_{i_k}}},\underbrace{\vphantom{d_{i_1},\dots
,d_{i_k}}\dots \vphantom{d_{i_1},\dots ,d_{i_k}}},\dots
}\vphantom{d_{i_1},\dots ,d_{i_k}}\right) 
\vphantom{\underbrace{d_{i_1},\dots ,d_{i_k}}}\label{Periodic}
\end{equation}
relative to the expansion (\ref{Eq5.11}), specifically $\delta _{jk+r}=d_{i_r}$
for $j=0,1,2,\dots $, and $r=1,2,\dots ,k$, corresponding to $t=-m$.
\end{remark}

\begin{scholium}
\label{Sch5.5}The canonical embedding $\mathcal{O}_{N^m}$ into $\mathcal{O}%
_N $ defines the same gauge-invariant algebra $\limfunc{UHF}\nolimits_N$ for
the two algebras. If $\pi $ is a permutative, mul\-ti\-plic\-i\-ty-free
representation of $\mathcal{O}_N$, then the restriction $\pi _0$ of $\pi $
to $\mathcal{O}_{N^m}$ is the permutative representation defined by the maps 
\begin{equation*}
\sigma _{i_1}\sigma _{i_2}\cdots \sigma _{i_m}\text{\quad and\quad }R^m
\end{equation*}
instead of 
\begin{equation*}
\sigma _i\text{\quad and\quad }R\text{.}
\end{equation*}
Thus the two representations $\pi $ and $\pi _0$ have the same atoms, but if 
$\tau $ is the canonical action of $\mathbb{Z}$ on the atoms such that the
cycles of $\pi $ correspond to the $\tau $-orbits, then the corresponding
action associated to $\pi _0$ is $\tau ^m$. Thus we get a double-carousel
picture: each $\pi $-cycle with $n$ atoms splits into $\gcd \left(
m,n\right) $ $\pi _0$-cycles each containing $n\diagup \gcd \left(
m,n\right) $ atoms.
\end{scholium}

\section{\label{SecNew6}The shift representation\label{Shift}}

In the relations (\ref{Eq2.4}), the set $\mathbb{N}$ just plays the role of
a countable set. Further insight into the structure of the maps $\sigma _i$
can be obtained by equipping $\mathbb{N}$ with extra structure, and we will
do that by embedding $\mathbb{N}$ as a countable subset $\Omega $ of $%
\mathbb{Z}_N^\infty =\bigcross_{k=1}^\infty \mathbb{Z}_N$. For this we will,
through this section, assume that the coding map $\sigma $ is injective, and
the embedding of $\mathbb{N}$ in $\Omega $ is just the coding map $\sigma $
defined by (\ref{Eq2.5}) and (\ref{Eq2.6}). Let $R:\mathbb{N}\rightarrow %
\mathbb{N}$ be the map defined in Scholium \ref{Sch2.4}, and denote $\sigma
\sigma _i\sigma ^{-1}$ and $\sigma R\sigma ^{-1}$ again by $\sigma _i$ and $%
R $. Then it follows immediately from the definitions that 
\begin{equation}
\sigma _i\left( x_1,x_2,x_3,\dots \right) =\left( i,x_1,x_2,x_3,\dots \right)
\label{EqNew6.1}
\end{equation}
and 
\begin{equation}
R\left( x_1,x_2,x_3,\dots \right) =\left( x_2,x_3,x_4,\dots \right)
\label{EqNew6.2}
\end{equation}
for all $\left( x_1,x_2,x_3,\dots \right) \in \Omega $, $i\in \mathbb{Z}_N$.
Conversely, if $\Omega $ is any countable subset of $\bigcross_{k=1}^\infty %
\mathbb{Z}_N$ which is invariant under the two maps (\ref{EqNew6.1}) and (%
\ref{EqNew6.2}), then one immediately verifies that the relations (\ref
{Eq2.4}) are fulfilled. The representation associated with $\sigma _i$ can
then be realized on $\ell _2\left( \Omega \right) $ by mapping the basis
element $e_n$ into $\delta _{\sigma \left( n\right) }$, and hence 
\begin{equation}
S_k\delta _x=\delta _{\sigma _k\left( x\right) }  \label{EqNew6.3}
\end{equation}
for $k\in \mathbb{Z}_N$, $x\in \Omega $.

Transposing the characterizations in Propositions \ref{Pro2.2} and \ref
{Pro2.3} and Corollary \ref{Cor2.4} to $\Omega $, we obtain

\begin{corollary}
\label{CorNew6.1}In the shift representation, the following two conditions
are equivalent.

\begin{enumerate}
\item  $\vphantom{\Bigl(}x\sim y$.\label{CorNew6.1a(1)}

\item  There is a $\vphantom{\Bigl(}k\in \mathbb{Z}$ and an $n_0\in %
\mathbb{N}$ with $n_0>\left| k\right| $ such that\label{CorNew6.1a(2)} 
\begin{equation*}
x_i=y_{i+k}
\end{equation*}
$\vphantom{\smash[t]{\Bigl(}}$for $i>n_0$.
\end{enumerate}

$\vphantom{\Bigl(}$Also, the following two conditions are equivalent.

\begin{enumerate}
\item  $\vphantom{\Bigl(}x\approx y$.\label{CorNew6.1b(1)}

\item  There is an $\vphantom{\Bigl(}n_0\in \mathbb{N}$ such that \label
{CorNew6.1b(2)} 
\begin{equation*}
x_i=y_i
\end{equation*}
$\vphantom{\smash[t]{\Bigl(}}$for $i>n_0$.
\end{enumerate}
\end{corollary}

Thus the $\sim $-equiv\-a\-lence classes are characterized as classes of $x$
in $\Omega $ with the same tail up to translation, and the $\approx $%
-equiv\-a\-lence classes are characterized as classes of $x$ in $\Omega $
with the same tail.

The map $\tau ^j$ on the $\approx $-equiv\-a\-lence classes defined in
Theorem \ref{Thm4.1} is then defined by translation of the tail sequence by $%
j$ places towards the left if $j\geq 0$, and by $\left| j\right| $ places
towards the right if $j<0$. Thus the atoms in the cycle defined by a fixed
tail sequence is defined by the various translates of this tail, so in this
picture Theorem \ref{Thm4.1} becomes almost trivial. Furthermore,
Proposition \ref{Pro5.2} and Corollary \ref{Cor5.3} get the following form:

\begin{corollary}
\label{CorNew6.2}In the shift representation, the following conditions are
equivalent for an $x\in \Omega $.

\begin{enumerate}
\item  $\vphantom{\Bigl(}k=\func{Per}\left( x\right) $ is finite.\label
{CorNew6.2(1)}

\item  The sequence $\vphantom{\Bigl(}\left( x_1,x_2,x_3,\dots \right) $ has
a periodic tail with minimal period $k$. \label{CorNew6.2(2)}
\end{enumerate}

Furthermore, if $\vphantom{\smash[b]{\Bigl(}}l\in \mathbb{N}$ is so large
that the periodic tail is 
\begin{equation*}
\left( \smash{\dots ,\underbrace{x_{lk+1},x_{lk+2},\dots
,x_{lk+k}},\underbrace{x_{lk+1},\dots ,x_{lk+k}},\dots}\right)
\end{equation*}
then put 
\begin{equation*}
\left( i_1,\dots ,i_k\right) =\left( x_{lk+1},\dots ,x_{lk+k}\right)
\end{equation*}
and the element 
\begin{equation*}
y=\left( y_1,y_2,\dots \right) =\left( i_1,\dots ,i_k,i_1,\dots
,i_k,i_1\dots ,i_k,\dots \right)
\end{equation*}
has the properties $y\approx x$, and 
\begin{equation*}
\sigma _{i_1}\cdots \sigma _{i_k}\left( y\right) =y\text{.}
\end{equation*}
Hence 
\begin{equation*}
S_{i_1}\cdots S_{i_k}e_y=e_y
\end{equation*}
so $e_y$ defines a sub-Cuntz state of order $k$. The cycle has $k$ such
states, corresponding to the $k$ different translates of the tail sequence,
and each of these $k$ states represents an atom.
\end{corollary}

It is now clear which subsets $\Omega \subseteq \bigcross_{k=1}^\infty %
\mathbb{Z}_N$ can occur in this manner: There is a countable (or finite)
sequence $y_1,y_2,y_{3,}\dots $ in $\Omega $ such that $\Omega $ consists of
all $x\in \bigcross_{k=1}^\infty \mathbb{Z}_N$ with corresponding tail
sequence equal to a translate of the tail sequence of some $y_k$. For any
countable sequence $y_1,y_2,\dots $, we may define a countable $\Omega $ in
this way; and conversely, if $\Omega $ is given, we may simply take $%
y_1,y_2,\dots $ to be an enumeration of all the elements in $\Omega $.

Note that the number of periodic sequences in $\mathbb{Z}_N^\infty $ with
period $k$ is $N^k$. If $N^{\prime }\left( k\right) $ denotes the number of
periodic sequences with \emph{minimal} period $k$, we have 
\begin{multline*}
N^{\prime }\left( k\right) =N^k-\smash{\sum \{ N^{\prime }\left( n\right)
\mid n\text{ factor of }k\text{, }} \\
\text{inclusive of }n=1\text{ but exclusive of }n=k\}
\end{multline*}
Let $\sim $ denote the equivalence on the set of periodic sequences given by
translation. The number of $\sim $-classes of sequences with minimal period $%
k$ is then 
\begin{equation*}
N\left( k\right) =N^{\prime }\left( k\right) \diagup k
\end{equation*}
and hence 
\begin{multline*}
N\left( k\right) =\smash{\Big( N^k-\sum \{ nN\left( n\right) \mid n\text{
factor of }k\text{, }} \\
\text{inclusive of }n=1\text{ but exclusive of }n=k\}\smash{\Big)} \diagup k
\end{multline*}
In particular, 
\begin{alignat*}
{2} N\left( k\right) &=\left( N^k-N\right) \diagup k & \quad &\text{ if }k%
\text{ is prime;} \\
N\left( 1\right) &=N\text{;} & \quad & \\
N\left( k^2\right) &=\left( N^{k^2}-N^k\right) \diagup k^2 & \quad &\text{
if }k\text{ is prime; and} \\
N\left( kn\right) &=\left( N^{kn}-N^k-N^n+N\right) \diagup kn & \quad &\text{
if }k, n\text{ are distinct primes, etc.}
\end{alignat*}

Put $N\left( \infty \right) =\infty $. From the remarks of the last two
paragraphs, we deduce

\begin{corollary}
\label{CorNew6.3}For a given permutative mul\-ti\-plic\-i\-ty-free
representation of $\mathcal{O}_N$ there is a countable or finite number of
cycles, and the number of cycles with exactly $k$ atoms does not exceed $%
N\left( k\right) $ for $k=1,2,\dots ,\infty $. Conversely, any countable
family of sets each containing countably or finitely many points, can be
represented as the cycle-atom structure of a permutative
mul\-ti\-plic\-i\-ty-free representation of $\mathcal{O}_N$ if the number of
sets in the family containing $k$ points does not exceed $N\left( k\right) $
for $k=1,2,3,\dots ,\infty $.
\end{corollary}

We will remark on the results corresponding to those in this section for
general permutative representations in Section \ref{Concluding}.

\section{\label{SecNew7}The universal permutative mul\-ti\-plic\-i\-ty-free
representation\label{Universal}}

\noindent From now, put 
\begin{equation*}
\Omega =\bigcross_{k=1}^\infty \mathbb{Z}_N
\end{equation*}
and equip $\Omega $ with the discrete topology instead of the usual compact
Cantor one. Define maps $\bar{\sigma }_i$ on $\Omega $ by (\ref
{EqNew6.1}): 
\begin{equation*}
\bar{\sigma }_i\left( x_1,x_2,\dots \right) =\left( i,x_1,x_2,\dots
\right)
\end{equation*}
for $i\in \mathbb{Z}_N$, $\left( x_1,x_2,\dots \right) \in \Omega $. We
obtain a representation $\bar{\pi }$ of $\mathcal{O}_N$ on the
non-separable Hilbert space $\ell _2\left( \Omega \right) $ by putting 
\begin{equation*}
S_k\delta _x=\delta _{\bar{\sigma }_k\left( x\right) }
\end{equation*}
for $k\in \mathbb{Z}_N$, $x\in \Omega $. One verifies as in Section \ref
{SecNew6} that this representation is permutative (albeit not on a separable
Hilbert space) and mul\-ti\-plic\-i\-ty-free. Furthermore, if another
permutative mul\-ti\-plic\-i\-ty-free representation $\pi $, defined by the
maps $\sigma _i$ satisfying (\ref{Eq2.4}), is given, then the map $%
e_n\rightarrow \delta _{\sigma \left( n\right) }$ defines a unitary
equivalence between $\pi $ and a subrepresentation of $\bar{\pi }$.
Therefore we call $\bar{\pi }$ the universal permutative
mul\-ti\-plic\-i\-ty-free representations. Again the Hilbert spaces of the
irreducible subrepresentations of $\bar{\pi }$ are spanned by $\delta
_x $ with $x$ running through $\sim $-equiv\-a\-lence classes in $\Omega $,
and the Hilbert spaces of the irreducible subrepresentations of 
$\bar{\pi }|_{\limfunc{UHF}\nolimits_N}$ are spanned by the respective 
$\delta _x$%
's with $x$ running through $\approx $-equiv\-a\-lence classes in $\Omega $.
Here $x\sim y$ if $x$, $y$ have the same tail sequence up to translation,
and $x\approx y$ if $x$, $y$ have the same tail sequence. In summary,

\begin{corollary}
\label{CorNew7.1}There exists a universal permutative
mul\-ti\-plic\-i\-ty-free representation $\bar{\pi }$ of $\mathcal{O}_N$
on a Hilbert space of continuum dimension, with the property that 
$\bar{\pi }$ contains any permutative mul\-ti\-plic\-i\-ty-free representation.
For $k=1,2,\dots $, the representation $\bar{\pi }$ contains $N\left(
k\right) $ cycles with $k$ atoms, and the subspace corresponding to each of
these atoms contains a vector defining a sub-Cuntz state of order $k$,
namely the unique vector of the form 
\begin{equation*}
\left( \smash{\underbrace{x_1,x_2,\dots ,x_k},\underbrace{x_1,x_2,\dots
,x_k},\underbrace{x_1,\dots},\dots}\right) \vphantom{\underbrace{x_1,x_2,%
\dots ,x_k}}
\end{equation*}
in the associated $\approx $-class. The representation $\bar{\pi }$
contains a continuum of cycles with a countable infinity of atoms, one for
each of the translates of one element in the associated
$\sim $-class \textup{(}which
is not asymptotically periodic\/\textup{).}
\end{corollary}

We will comment on the extension of the results in Sections \ref
{SecNew6} and \ref{SecNew7} to permutative representations $\pi $ of $%
\mathcal{O}_N$ which are not necessarily mul\-ti\-plic\-i\-ty-free, and
their associated branching function system
$\sigma _1,\dots ,\sigma _N$ in Section \ref{Concluding}. In
this case the map $e_n\rightarrow \delta _{\sigma \left( n\right) }$ defines
merely an intertwiner between $\pi $ and the universal
mul\-ti\-plic\-i\-ty-free representation.

\section{\label{Sec6}Some specific examples of the cycle and atom structure:
\\The $\bmod{\,N}$ case\label{Examples}\label{SecNew8}}

We will now compute the cycle and atom structure for some of the monomial
representations defined by (\ref{Eq3.1}). First consider the unitary
operator $U_m$ defined on $\mathcal{H}$ by 
\begin{equation}
U_me_n=e_{n+m}  \label{Eq6.1}
\end{equation}
for $n\in \mathbb{Z}$. We have 
\begin{align*}
U_m^{*}S_k^{}U_m^{}e_n^{}& =U_m^{*}S_k^{}e_{n+m}^{}=U_m^{*}e_{N\left(
n+m\right) ^{}+s_k^{}} \\
& =e_{N\left( n+m\right) +s_k-m}^{}=e_{Nn+\left( N-1\right) m+s_k}^{}\text{,}
\end{align*}
so the monomial representation defined by (\ref{Eq3.1}) is unitarily
equivalent to that defined by 
\begin{equation}
\sigma _k^{\left( m\right) }\left( n\right) =nN+\left( s_k+\left( N-1\right)
m\right) \text{.}  \label{Eq6.2}
\end{equation}
This is particularly convenient if $N=2$, since we then may assume that one
of the $s_k$ is zero. In general we cannot do so (for example, $N=3$, $%
\left\{ s_i\right\} =\left\{ 1,3,5\right\} $). Another useful symmetry is
defined by the unitary self-adjoint operator $V$ given on $\mathcal{H}$ by 
\begin{equation}
Ve_n=e_{-n}  \label{Eq6.3}
\end{equation}
for $n\in \mathbb{Z}$. We have 
\begin{equation*}
VS_k^{}Ve_n^{}=e_{nN-s_k}^{}
\end{equation*}
so the monomial representation defined by (\ref{Eq3.1}) is unitarily
equivalent to that given by 
\begin{equation}
\sigma _k^{\prime }\left( n\right) =nN-s_k\text{.}  \label{Eq6.4}
\end{equation}
\medskip

\subsection{The case $N=2$}\label{SubsecNew8.1}

\ \bigskip

By the remark following (\ref{Eq6.2}), when $N=2$ we may assume $s_0=0$, and
then $s_1=p$ where $p$ is an odd integer. By (\ref{Eq6.4}) we may assume $%
p>0 $. Then 
\begin{equation}
\begin{aligned} \sigma _0\left( n\right) &=2n\text{,} \\ \sigma _1\left(
n\right) &=2n+p\text{,} \end{aligned}  \label{Eq6.5}
\end{equation}
and 
\begin{equation}
R\left( n\right) = \begin{cases} \frac n2 & \text{if }n\text{ is even,} \\
\frac{n-p}2 & \text{if }n\text{ is odd.} \end{cases}  \label{Eq6.6}
\end{equation}
For $m\in \mathbb{Z}$, $i_1,\dots ,i_k\in \mathbb{Z}_2=\left\{ 0,1\right\} $
we have 
\begin{align}
\sigma _{i_1}\cdots \sigma _{i_k}\left( m\right) & =s_{i_1}+2s_{i_2}+\dots
+2^{k-1}s_{i_k}+2^km  \label{Eq6.7} \\
& =p\left( i_1+2i_2+\dots +2^{k-1}i_k\right) +2^km  \notag \\
& =pi+2^km\text{,}  \notag
\end{align}
where $i$ can be any integer in $\left\{ 0,1,\dots ,2^k-1\right\} $. Thus,
by Proposition \ref{Pro2.3}, $n_1\approx n_2$ if and only if $n_1$ and $n_2$
both belong to a set of the form 
\begin{equation}
X\left( m,k\right) =\left\{ 2^km+pi\bigm|i=0,1,\dots ,2^k-1\right\}
\label{Eq6.8}
\end{equation}
for some $m\in \mathbb{Z}$, $k\in \left\{ 0\right\} \cup \mathbb{N}$. The
case $p=1 $ is special:

\begin{proposition}
\label{Pro6.1}\label{ProNew8.1}If $N=2$ and $p=1$, there are exactly two
atoms and two cycles, corresponding to the equivalence classes $\left\{
0,1,2,\dots \right\} $ and $\left\{ \dots ,-3,-2,-1\right\} $. Thus all
atoms have period $1$, and the corresponding eigen relations are 
\begin{equation}
\sigma _0\left( 0\right) =0  \label{Eq6.9}
\end{equation}
for the first class, and 
\begin{equation}
\sigma _1\left( -1\right) =-1  \label{Eq6.10}
\end{equation}
for the second.
\end{proposition}

\TeXButton{Begin Proof}{\begin{proof}}This is immediate by putting $m=0$ and
letting $k\rightarrow \infty $ in (\ref{Eq6.8}), and verifying that the two
indicated sets are invariant under both $\sigma _0$, $\sigma _1$ and $R$, so
they are both $\approx $- and $\sim $-equiv\-a\-lence classes.%
\TeXButton{End Proof}{\end{proof}}

\begin{proposition}
\label{Pro6.2}\label{ProNew8.2}If $N=2$ and $p\in \left\{ 3,5,7,9,\dots
\right\} $, there are exactly $p+1$ atoms, given by the $\approx $%
-equiv\-a\-lence classes 
\begin{equation}
X_m=\left\{ n\in \mathbb{Z}\mid n=m\mod p\right\}  \label{Eq6.11}
\end{equation}
for $m=1,2,\dots ,p-1$, in conjunction with the two classes 
\begin{equation}
\begin{aligned} X_{0^{+}} &=\left\{ np\mid n=0,1,2,\dots \right\} \text{,}
\\ X_{0^{-}} &=\left\{ -np\mid n=1,2,3,\dots \right\} \text{.} \end{aligned}
\label{Eq6.12}
\end{equation}
The latter two $\approx $-equiv\-a\-lence classes are also $\sim $%
-equiv\-a\-lence classes, while the cycle of $X_m$ for $m=1,2,\dots ,p-1$
consists of the atoms $X_{2^km\mod p}$, $k=0,1,2,\dots $. Thus the period of
the atom $X_m$ under the $\tau $ action is the order of $2$ modulo $p\diagup
\gcd \left( m,p\right) $, i.e., this period is the smallest positive integer 
$k$ such that $p$ divides $m\left( 2^k-1\right) $.
\end{proposition}

\begin{remark}
\label{RemNew8.2}In the present case, Corollary \ref{Cor3.10}
gives the upper bound
$1+\frac p{2-1}=1+p$ on the number of atoms, so by Proposition \ref{ProNew8.2},
this bound is optimal in this case. Also, scrutinizing the proof of
Corollary \ref{Cor3.10},
we note that the set $B_\infty $ of points in $\mathbb{Z}$ which
are periodic under $R$
must be exactly 
\begin{equation*}
B_\infty =\left\{ -p,-p+1,\dots ,-1,0\right\} 
\end{equation*}
in this case.
\end{remark}

\TeXButton{Begin Proof of Proposition 8.2}
{\begin{proof}[Proof of Proposition \ref{ProNew8.2}]}We will need
Euler's familiar theorem on
divisors (see, e.g., \cite{And94}): for $n\in \mathbb{N}$, define the Euler
function $\phi \left( n\right) $ as the number of integers $m\in \left\{
1,2,\dots ,n-1\right\} $ such that $\gcd \left( m,n\right) =1$. If $a\in %
\mathbb{N}$ and $\gcd \left( a,n\right) =1$, then $a^{\phi \left( n\right)
}=1\mod n$, hence $a^{k\phi \left( n\right) }=1\mod n$ for all $k\in %
\mathbb{N}$, and then 
\begin{equation*}
ma^{k\phi \left( n\right) }=m\mod n
\end{equation*}
for all $m,k\in \mathbb{N}$. If $a=2$ and $p\in \mathbb{N}$ is odd, we have
in particular 
\begin{equation*}
m2^{k\phi \left( p\right) }=m\mod p\text{.}
\end{equation*}
Thus, for each fixed $m\in \left\{ -\left( n-1\right) ,-\left( n-2\right)
,\dots ,-1\right\} $, the sets $X\left( m,k\phi \left( p\right) \right) $
form an increasing sequence in $k$, and 
\begin{equation*}
\bigcup_kX\left( m,k\phi \left( p\right) \right) =X_{p+m}\text{.}
\end{equation*}
On the other hand it follows immediately from (\ref{Eq6.8}) and its
preceding remark that 
\begin{equation*}
n_1\approx n_2\Rightarrow n_1=n_2\mod p\text{,}
\end{equation*}
and hence $X_m$ is a $\approx $-equiv\-a\-lence class for $m=1,2,\dots ,p-1$.

Next note that $k\mapsto X\left( 0,k\right) $ form an increasing sequence
with union $X_{0^{+}}$, and $k\mapsto X\left( -p,k\right) $ form an
increasing sequence with union $X_{0^{-}}$, so each of the two sets are
subsets of a $\approx $-equiv\-a\-lence class. But $X_{0^{-}}\cup
X_{0^{+}}=X_0=p\mathbb{Z}$. If one of the sets $X\left( m,k\right) $
intersects $X_0$, one has 
\begin{equation*}
2^km=0 \mod p\text{,}
\end{equation*}
i.e., $p$ is a factor of $2^km$, but as $p$ is odd, $p$ is then a factor of $%
m$. But if $m$ is a multiple of $p$, then either $X\left( m,k\right) \in
\left\{ \dots ,-3,-2,-1\right\} $ (if $m<0$) or $X\left( m,k\right) \in
\left\{ 0,1,2,\dots \right\} $. Thus, if $X\left( m,k\right) $ intersects $%
X_{0^{\pm }}$, then $X\left( m,k\right) \subseteq X_{0^{\pm }}$. It follows
that $X_{0^{+}}$ and $X_{0^{-}}$ really are $\approx $-equiv\-a\-lence
classes.

To compute the cycles, we must compute the action $\tau $ of $\mathbb{Z}$ on
the $\approx $-classes defined in Theorem \ref{Thm4.1}. Computing modulo $p$%
, the cycle containing $X_m$ for $m=1,2,\dots ,p-1$ thus consists of the
atoms $X_{2^km}$ for $k=0,1,2,\dots $, and the number of atoms in this cycle
is the smallest $k=1,2,\dots $ such that $2^km=m \mod p$, i.e., $\left(
2^k-1\right) m=0 \mod p$. Dividing out by $\gcd \left( m,p\right) $, this is
the smallest $k$ such that $p\diagup \gcd \left( m,p\right) $ divides $2^k-1$%
, i.e., $k$ is the order of $2$ modulo $p\diagup \gcd \left( m,p\right) $.
(We recall the Mersenne numbers $2^k-1$; see \cite{LeV96}, \cite{BLSTW},
\cite{Rib89}.)%
\TeXButton{End Proof}{\end{proof}}

\textit{Summary.} When $N=2$, one thus always has the two one-atom cycles $%
\left\{ X_{0^{+}}\right\} $ and $\left\{ X_{0^{-}}\right\} $. The remaining
cycles are: \newlength{\saveleftmargini} \setlength{\saveleftmargini}{%
\leftmargini} \settowidth{\leftmargini}{$p=8388607$} \addtolength{%
\leftmargini}{2.7778pt} \newlength{\saveleftmarginii} \setlength{%
\saveleftmarginii}{\leftmarginii} \settowidth{\leftmarginii}{$\gcd(m,p)=$} 
\addtolength{\leftmarginii}{2.7778pt} 
\begin{list}{}{\setlength{\leftmargin}{\leftmargini} 
\setlength{\itemsep}{\customskipamount}
\setlength{\topsep}{\customskipamount}
\setlength{\labelsep}{2.7778pt}}
\item[\hss\llap{$p=1$:}]  None.

\item[\hss\llap{$p=3$:}]  $\left\{ X_1,X_2\right\} $.

\item[\hss\llap{$p=5$:}]  $\left\{ X_1,X_2,X_3,X_4\right\} $.

\item[\hss\llap{$p=7$:}]  $\left\{ X_1,X_2,X_4\right\} $, 
$\left\{ X_3,X_5,X_6\right\} $.

\item[\hss\llap{$p=9$:}]  $\left\{ X_1,X_2,X_4,X_5,X_7,X_8\right\} $, $\left\{
X_3,X_6\right\} $.

\item[\hss\llap{$p=11$:}]  $\left\{ X_1,X_2,\dots ,X_{10}\right\} $.

\item[\hss\llap{$p=13$:}]  $\left\{ X_1,X_2,\dots ,X_{12}\right\} $.

\item[\hss\llap{$p=15$:}]  $\left\{ X_1,X_2,X_4,X_8\right\} $, $\left\{
X_3,X_6,X_9,X_{12}\right\} $,
$\left\{ X_5,X_{10}\right\} $, $\left\{ X_7,X_{11},X_{13},X_{14}\right\} $.

\item[\hss\llap{$p=21$:}]  $2$ cycles with $6$ atoms, $2$ with $3,$ 
and $1$ with $2$.

\item[\hss\llap{$p=25$:}]  $1$ cycle with $20$ atoms and $1$ with $4$.

\item[\hss\llap{$p=27$:}]  $1$ cycle with $18$ atoms, $1$ with $6$, 
and $1$ with $2$.

\item[\hss\llap{$p=41$:}]  $2$ cycles with $20$ atoms.

\item[\hss\llap{$p=49$:}]  $2$ cycles with $21$ atoms and $2$ with $3$.

\item[\hss\llap{$p=81$:}]  $1$ cycle with $54$ atoms, $1$ with $18$, 
$1$ with $6$,
and $1$ with $2$.

\item[\hss\llap{$p=105$}]  $=3\cdot 5\cdot 7$:
$6$ cycles with $12$ atoms, $2$ with $6$, $3$
with $4$, $2$ with $3$, and $1$ with $2$.

\item[\hss\llap{$p=451$}]  $=11\cdot 41$: 
Note that $\gcd \left( 2^{20}-1,451\right) =451$,
and $11$ is a factor of the Mersenne number $2^k-1$ if and only if $k$ is a
multiple of $10$, and $41$ is a factor of the Mersenne number $2^k-1$ if and
only if $k$ is a multiple of $20$. Thus all numbers in 
$\mathbb{Z}_{451}\diagdown \left\{ 0\right\} $ 
have period $20$ under multiplication by 
$2$, except for the integer multiples of $41$, which have period $10$.
Conclusion:
\begin{list}{}{\setlength{\leftmargin}{\leftmarginii}
\setlength{\itemsep}{0\customskipamount}
\setlength{\topsep}{\customskipamount}}
\item[ ]  $1$ cycle with $10$ atoms,

\item[ ]  $22$ cycles with $20$ atoms.
\end{list}

\item[\hss\llap{$p=1387$}]  $=19\cdot 73$: 
Note that the order of $2$ modulo $19$ is $18$, and the
order of $2$ modulo $73$ is $9$, so there are $73-1=72$ atoms of order $9$,
while all others have order $18$. Conclusion:
\begin{list}{}{\setlength{\leftmargin}{\leftmarginii}
\setlength{\itemsep}{0\customskipamount}
\setlength{\topsep}{\customskipamount}}
\item[ ]  $8$ cycles with $9$ atoms,

\item[ ]  $73$ cycles with $18$ atoms.
\end{list}

\item[\hss\llap{$p=4095$}]  $=2^{12}-1=3\cdot 3\cdot 5\cdot 7\cdot 13$: 
Thus all $m\in \mathbb{Z}_{4095}\diagdown
\left\{ 0\right\} $ have a period dividing $12$. Looking at a table of the
prime decomposition of the Mersenne numbers $2^k-1$ for $k=2,3,4,\dots ,12$,
one deduces:
\begin{list}{}{\setlength{\leftmargin}{\leftmarginii}
\setlength{\itemsep}{\customskipamount}
\setlength{\topsep}{\customskipamount}
\setlength{\labelsep}{2.7778pt}}
\item[\hss\llap{$\gcd (m,p)=$}]  $3\cdot 5\cdot 7\cdot 13$:\ 

$2$ atoms;

$p\diagup \gcd =3\Rightarrow $ period $=2$;

$1$ cycle with $2$ atoms.

\item[\hss\llap{$\gcd(m,p)=$}]  $3\cdot 3\cdot 7\cdot 13$:\ 

$4$ atoms;

$p\diagup \gcd =5\Rightarrow $ period $=4$;

$1$ cycle with $4$ atoms.

\item[\hss\llap{$\gcd(m,p)=$}]  $3\cdot 3\cdot 5\cdot 13$:\ 

$6$ atoms;

$p\diagup \gcd =7\Rightarrow $ period $=3$;

$2$ cycles with $3$ atoms.

\item[\hss\llap{$\gcd(m,p)=$}]  $5\cdot 7\cdot 13$:\ 

$3\cdot 3-1=8$ atoms;

$2$ of these accounted for before; $6$ remain;

$p\diagup \gcd =3\cdot 3\Rightarrow $ period $=6$;

$1$ cycle with $6$ atoms.

\item[\hss\llap{$\gcd(m,p)=$}]  $3\cdot 7\cdot 13$:\ 

$3\cdot 5-1=14$ atoms;

$2+4$ of these accounted for before; $8$ remain;

$p\diagup \gcd =3\cdot 5\Rightarrow $ period $=4$;

$2$ cycles with $4$ atoms.

\item[\hss\llap{$\gcd(m,p)=$}]  $3\cdot 5\cdot 13$:\ 

$3\cdot 7-1=20$ atoms;

$2+6$ of these accounted for before; $12$ remain;

$p\diagup \gcd =3\cdot 7\Rightarrow $ period $6$;

$2$ cycles with $6$ atoms.

\item[\hss\llap{$\gcd(m,p)=$}]  $5\cdot 13$:\ 

$3\cdot 3\cdot 7-1=62$ atoms;

$2+6+6+12=26$ accounted for before; $36$ remain;

$p\diagup \gcd =3\cdot 3\cdot 7\Rightarrow $ period $6$;

$6$ cycles with $6$ atoms each.

\item[\hss\llap{\hbox to\leftmargin{\hskip\labelsep All remaining
$m$:\hss}}]  \ 

period $12$.
\end{list}

Thus we have
\begin{list}{}{\setlength{\leftmargin}{\leftmarginii}
\setlength{\itemsep}{0\customskipamount}
\setlength{\topsep}{\customskipamount}}
\item[ ]  $1$ cycle with $2$ atoms;

\item[ ]  $2$ cycles with $3$ atoms;

\item[ ]  $3$ cycles with $4$ atoms;

\item[ ]  $9$ cycles with $6$ atoms;

\item[ ]  $335$ cycles with $12$ atoms.
\end{list}
\end{list}
\setlength{\leftmargini}{\saveleftmargini} \setlength{\leftmarginii}{%
\saveleftmarginii}

This last example has an interesting optimal property: if $N\left( k\right) $
is the maximal number of cycles with $k$ atoms, in the general situation of
a permutative mul\-ti\-plic\-i\-ty-free representation, one computes (by the
recipe prior to Corollary \ref{CorNew6.3}) that (for $N=2$): 
{\allowdisplaybreaks
\begin{align*}
N\left( 1\right) &=2 \\
N\left( 2\right) &=\left( 2^2-2\right) \diagup 2=1 \\
N\left( 3\right) &=\left( 2^3-2\right) \diagup 3=2 \\
N\left( 4\right) &=\left( 2^4-2^2\right) \diagup 4=3 \\
N\left( 5\right) &=\left( 2^5-2\right) \diagup 5=6 \\
N\left( 6\right) &=\left( 2^6-2^3-2^2+2\right) \diagup 6=9 \\
N\left( 12\right) &=\left( 2^{12}-6N\left( 6\right) -4N\left( 4\right)
-3N\left( 3\right) -2N\left( 2\right) -N\left( 1\right) \right) \diagup
12=335
\end{align*}%
}%

Thus in this example, $p=4095=2^{12}-1$, for each $k$ that divides $12$ the
number of cycles with $k$ atoms is as large as it can be by the general
theory of Corollary \ref{CorNew6.3}. It follows that the image of $\mathbb{Z}
$ by the map $\sigma $ (defined by (\ref{Eq2.5})--(\ref{Eq2.8})) consists of 
\emph{all} sequences in $\bigcross_{k=1}^\infty \mathbb{Z}_2$ with a tail
which is periodic with period $12$ (i.e., the minimal period divides $12$).
In particular, this means that, for any sequence $i_1,i_2,\dots ,i_{12}$ of
length $12$, there is a vector $e_m$ in the canonical basis for the
representation of $\mathcal{O}_2$ such that 
\begin{equation*}
S_{i_1}S_{i_2}\cdots S_{i_{12}}e_m=e_m\text{.}
\end{equation*}
(This vector is unique if the periodic sequence of period $12$ defined by $%
\left( i_1,\dots ,i_{12}\right) $ has \emph{minimal} period $12$; if the
minimal period $d$ is just a factor of $12$, then the vector is only unique
up to translation by $d$, so there are $12\diagup d$ such vectors.)

This special property comes from the fact that $4095=2^{12}-1$ is a Mersenne
number, and looking back at the examples $p=3=2^2-1$, $p=7=2^3-1$, $%
p=15=2^4-1$ we see that this phenomenon occurs also there. Thus if $2\left(
k\right) $ is the function defined prior to Corollary \ref{CorNew6.3} for $%
N=2$, we have

\begin{proposition}
\label{ProNew8.3}If $N=2$, and $p=2^k-1$ $\left( k=1,2,3,\dots \right) $,
then the range of the map $\sigma :\mathbb{Z}\rightarrow \bigcross%
_{k=1}^\infty \mathbb{Z}_2$ consists of all sequences which have a periodic
tail with minimal period
$n$ dividing $k$ \textup{(}inclusive of $n=1$ and $n=k$\textup{)
(}i.e., the tails
are all periodic sequences with period $k$\textup{).} Thus, if $n$
divides $k$ in this situation, and $\left( i_1,\dots ,i_n\right) $ is a
sequence in $\mathbb{Z}_2^n$ such that the associated periodic sequence $%
\left( \smash{\dots ,\underbrace{i_1,\dots ,i_n},\underbrace{i_1,\dots ,i_n},%
\underbrace{\vphantom{i_1,\dots ,i_n}\dots 
\vphantom{i_1,\dots ,i_n}}}\right) $ has minimal period $n$%
, then the eigenvalue equation 
\begin{equation*}
S_{i_1}\cdots S_{i_n}e_m=e_m
\end{equation*}
has a unique solution $m$. Thus, for each $n$ dividing $k$, there are $%
2\left( n\right) $ cycles with $n$ atoms, and these are all the cycles
occurring.
\end{proposition}

\TeXButton{Begin Proof}{\begin{proof}}We merely note that if $n$ divides $k$%
, then $2^n-1$ divides $2^k-1$. But this is just because $2^n=1\mod q$, then 
$2^{mn}=\underbrace{2^n\cdot 2^n\cdots 2^n}_m=1\mod q$ for $m=1,2,\dots $;
i.e., if $k=nm$, then $2^k-1=\left( 2^n-1\right) \left( 1+2^n+\dots
+2^{n\left( m-1\right) }\right) $.\TeXButton{End Proof}{\end{proof}}

In general, if $p$ is a prime number and $k$ is the order of $2$ modulo $p$,
there will be: 
\setlength{\saveleftmargini}{\leftmargini}
\settowidth{\leftmargini}{$p=8388607$}
\addtolength{\leftmargini}{2.7778pt}
\setlength{\saveleftmarginii}{\leftmarginii}
\settowidth{\leftmarginii}{$\gcd(m,p)=$}
\addtolength{\leftmarginii}{2.7778pt}
\begin{list}{}{\setlength{\leftmargin}{\leftmargini}
\addtolength{\leftmargin}{\leftmarginii}
\setlength{\itemsep}{\customskipamount}
\setlength{\topsep}{\customskipamount}
\setlength{\labelsep}{2.7778pt}}
\item  $\left( p-1\right) \diagup k$ cycles with $k$ atoms each.
\end{list}
But it is not necessary for $p$ to be prime for us to have this structure. 
Example: \begin{list}{}{\setlength{\leftmargin}{\leftmargini} 
\setlength{\itemsep}{\customskipamount}
\setlength{\topsep}{\customskipamount}
\setlength{\labelsep}{2.7778pt}}
\item[\hss\llap{$p=8388607$}]  $=47\cdot 178481=2^{23}-1$. 
Since the order of $2$ modulo $47$, and the order
of $2$ modulo $178481$, are both $23$, all the atoms in this example (except
for the two fixed points) have period $23$. Conclusion:
\begin{list}{}{\setlength{\leftmargin}{\leftmarginii}
\setlength{\itemsep}{0\customskipamount}
\setlength{\topsep}{\customskipamount}}
\item[ ]  $364722$ cycles with $23$ atoms each.
\end{list}
\end{list}
A computer check done by Brian Treadway shows that the numbers
$p$ up to $10000$ with this property are:
\begin{list}{}{\setlength{\leftmargin}{\leftmargini} 
\setlength{\itemsep}{\customskipamount}
\setlength{\topsep}{\customskipamount}
\setlength{\labelsep}{2.7778pt}}
\item[\hss\llap{$p=2047$}]  $=23\cdot 89=2^{11}-1$: $186$ cycles with $11$ 
atoms each.

\item[\hss\llap{$p=3277$}]  $=29\cdot 113$: $117$ cycles with $28$ atoms each.

\item[\hss\llap{$p=4033$}]  $=37\cdot 109$: $112$ cycles with $36$ atoms each.

\item[\hss\llap{$p=8321$}]  $=53\cdot 157$: $160$ cycles with $52$ atoms each.
\end{list}
\setlength{\leftmargini}{\saveleftmargini}
\setlength{\leftmarginii}{\saveleftmarginii}
(The three last numbers are not Mersenne numbers.) A transcript of
this calculation is available by anonymous {\ttfamily ftp} from
{\ttfamily ftp.math.uiowa.edu}: look in the directory {\ttfamily
pub/jorgen/PermRepCuntzAlg} for the file {\ttfamily TwOrdP.log}.
\medskip

\subsection{The case $N$ general, $s_i=i$ for $i\in \mathbb{Z}_N=\left\{
0,1,\dots ,N-1\right\} $}\label{SubsecNew8.2}

\ \bigskip

One verifies exactly as in Proposition \ref{Pro6.1} that in this case there
are two atoms $\left\{ 0,1,2,\dots \right\} $ and $\left\{ \dots
,-3,-2,-1\right\} $, and these are also cycles.
\medskip

\subsection{The case $N=3$, $\left\{ s_i\right\}
=\left\{ 1,3,5\right\} $}\label{SubsecNew8.3}

\ \bigskip

This is one case where one cannot use (\ref{Eq6.2}) to reduce to the case
that one $s_i$ is $0$. We have $s_i=1+2i$, $i=0,1,2$, and hence 
\begin{align*}
\sigma _{i_1}\cdots \sigma _{i_k}\left( m\right) & =\left( 1+2i_1\right)
+3\left( 1+2i_2\right) +\dots +3^{k-1}\left( 1+2i_k\right) +3^km \\
& =\frac{3^k-1}{3-1}+2i+3^km=\frac{3^k\left( 1+2m\right) -1}2+2i\text{,}
\end{align*}
where $i=i_1+3i_2+\dots +3^{k-1}i_k$ can assume the values $0,1,2,\dots
,3^k-1$. Thus $n_1\approx n_2$ if and only if both $n_1$ and $n_2$ belong to
one of the arithmetic progressions 
\begin{equation*}
X\left( m,k\right) =\left\{ \frac{3^k\left( 1+2m\right) -1}2+2i\biggm|%
i=0,1,\dots ,3^k-1\right\} \text{.}
\end{equation*}
Putting in particular $m=-1$, $X\left( -1,k\right) $ is the step $2$
arithmetic progression between $-\frac{3^k+1}2$ and $\frac{3\cdot 3^k-5}2$.
But one checks that $\frac{3^k+1}2$ is even for $k=1,3,5,\dots $ and odd for 
$k=2,4,6,\dots $. Since $n_1\approx n_2$ implies that $n_1=n_2\mod 2$, one
deduces that there are two $\approx $-equiv\-a\-lence classes, namely the
set of even integers and the set of odd integers. These two sets are
interchanged by the maps 
\begin{equation*}
\sigma _i\left( m\right) =3m+1+2i\text{,}
\end{equation*}
and hence there is only one cycle, containing the two atoms. Thus the
corresponding representation of $\mathcal{O}_3$ is irreducible, and its
restriction to $\limfunc{UHF}\nolimits_3$ decomposes into two inequivalent
irreducible representations.

\section{\label{SecNew9}Some specific examples of the cycle and atom
structure: \\The $\bmod{\,\mathbf{N}}$ case\label{ExamplesMatrixCase}}

The main step from Section \ref{Examples} to the present one is the
increase in the dimension $\nu $ from one to higher; however, since a
general and comprehensive theory is not yet in the cards, we work out
the case $\nu =2$ in considerable detail. While we shall consider
several prototypical cases for the $2\times 2$ matrix $\mathbf{N}$
over $\mathbb{Z}$, our emphasis is the variety of dynamical systems
which result as the residue sets (or digit sets) $D$ range over an
infinite family.

But we also have two general developments in the generality where
$\left( \mathbf{N},D\right) $ is specified in $\nu $ dimensions,
assuming $0\in D$, and $D$ being a full set of residues for
$\mathbb{Z}^\nu \diagup \mathbf{N}\mathbb{Z}^\nu $. We consider the
setting when the iteration limit $\mathbf{T}$ from (\ref{Eq3ins1})
(i.e., the generalized fractions) is compact, and we compare our
representations of $\mathcal{O}_N$ ($N=\left| \det \mathbf{N}\right| $)
on the sequence spaces $\ell ^2$ with corresponding representations
on $L^2\left( \mathbf{T}\right) $, where this latter Hilbert space is
now with respect to the restricted Lebesgue measure from
$\mathbb{R}^\nu $. Our second general development concerns the set
$B_\infty $ of periodic points for $R$, and the corresponding
equivalence classes $p^\sim $ in $\mathbb{Z}^\nu $, defined for
$p\in B_\infty $. We shall study when there is a minimal subset $P$,
$0\in P\subset B_\infty $ such that $\bigcup _{p\in P}p^\sim $ is a
rank-$\nu $ lattice $\mathbb{L\kern0.5pt}$ which turns $\mathbf{T}$
into an $\mathbb{L\kern0.5pt}$-periodic tile for $\mathbb{R}^\nu $.
\medskip

\subsection{A matrix case}\label{SubsecNew9.1}

\ \bigskip

\begin{example}
\label{JordanMatrix}Let us compute the cycle and atom structure
for one of the monomial
representations defined by (\ref{Eq3.6bis}). We specify 
\begin{equation*}
\mathbf{N}=
\begin{pmatrix}
2 & 1 \\ 
0 & 2
\end{pmatrix}
\text{,}
\end{equation*}
acting on $\mathbb{Z}^2$, and then $N=\det \mathbf{N}=4$. As the set $%
\left\{ d_i,i=1,\dots ,4\right\} \subset \mathbb{Z}^2$ we take 
\begin{equation*}
d_{ij}=
\begin{pmatrix}
i \\ 
j
\end{pmatrix}
\text{,\quad }i,j=0,1\text{.}
\end{equation*}
Define 
\begin{equation*}
\sigma _{ij}\left( x\right) =d_{ij}+\mathbf{N}x
\end{equation*}
for $x\in \mathbb{Z}^2$. If $p$ is the parity function on $\mathbb{Z}$,
i.e., 
\begin{equation*}
p\left( m\right) =\begin{cases} 0 & \text{if }m\text{ is even,} \\ 1 &
\text{if }m\text{ is odd,} \end{cases}
\end{equation*}
a computation shows that 
\begin{align*}
R
\begin{pmatrix}
x_1 \\ 
x_2
\end{pmatrix}
& =\mathbf{N}^{-1}\left( 
\begin{pmatrix}
x_1 \\ 
x_2
\end{pmatrix}
-
\begin{pmatrix}
p\left( x_1-\frac{x_2-p\left( x_2\right) }2\right) \\ 
p\left( x_2\right)
\end{pmatrix}
\right) \\
& =\frac 14
\begin{pmatrix}
2x_1-2p\left( x_1-\frac{x_2-p\left( x_2\right) }2\right) -x_2+p\left(
x_2\right) \\ 
2x_2-2p\left( x_2\right)
\end{pmatrix}
\text{.}
\end{align*}
Now we have 
\begin{equation*}
\left\| R
\begin{pmatrix}
x_1 \\ 
x_2
\end{pmatrix}
-\mathbf{N}^{-1}
\begin{pmatrix}
x_1 \\ 
x_2
\end{pmatrix}
\right\| \leq \frac 12
\end{equation*}
for any $\left( \begin{smallmatrix}
x_1 \\ 
x_2
\end{smallmatrix}
\right) \in \mathbb{Z}^2$ because 
\begin{equation*}
R\left( d_{ij}+\mathbf{N}
\begin{pmatrix}
y_1 \\ 
y_2
\end{pmatrix}
\right) -\mathbf{N}^{-1}\left( d_{ij}+\mathbf{N}
\begin{pmatrix}
y_1 \\ 
y_2
\end{pmatrix}
\right) =\mathbf{N}^{-1}\left( d_{ij}\right)
\end{equation*}
and 
\begin{equation*}
\left\| \mathbf{N}^{-1}d_{ij}\right\| =\left\| \frac 14
\begin{pmatrix}
2i-j \\ 
2j
\end{pmatrix}
\right\| \leq \frac 12\text{.}
\end{equation*}
Thus 
\begin{align*}
& \left\| R^n
\begin{pmatrix}
x_1 \\ 
x_2
\end{pmatrix}
-\mathbf{N}^{-n}
\begin{pmatrix}
x_1 \\ 
x_2
\end{pmatrix}
\right\| \\
& \qquad \leq \sum_{k=0}^{n-1}\left\| \mathbf{N}^{-k}\left( R^{n-k}
\begin{pmatrix}
x_1 \\ 
x_2
\end{pmatrix}
\right) -\mathbf{N}^{-k-1}\left( R^{n-k-1}
\begin{pmatrix}
x_1 \\ 
x_2
\end{pmatrix}
\right) \right\| \\
& \qquad \leq \sum_{k=0}^{n-1}\left\| \mathbf{N}^{-k}\right\| \frac 12\text{.%
}
\end{align*}
\begin{figure}
\makebox[\textwidth]{\hfill
\begin{minipage}{0in}
\llap{$\mathbf{N}^2\left( \mathbf{T}\right) $ showing}
\llap{action of $\mathbf{N}^{-1}$}
\llap{\ }
\end{minipage}\quad
\begin{minipage}{3.6in}
\psfig{file=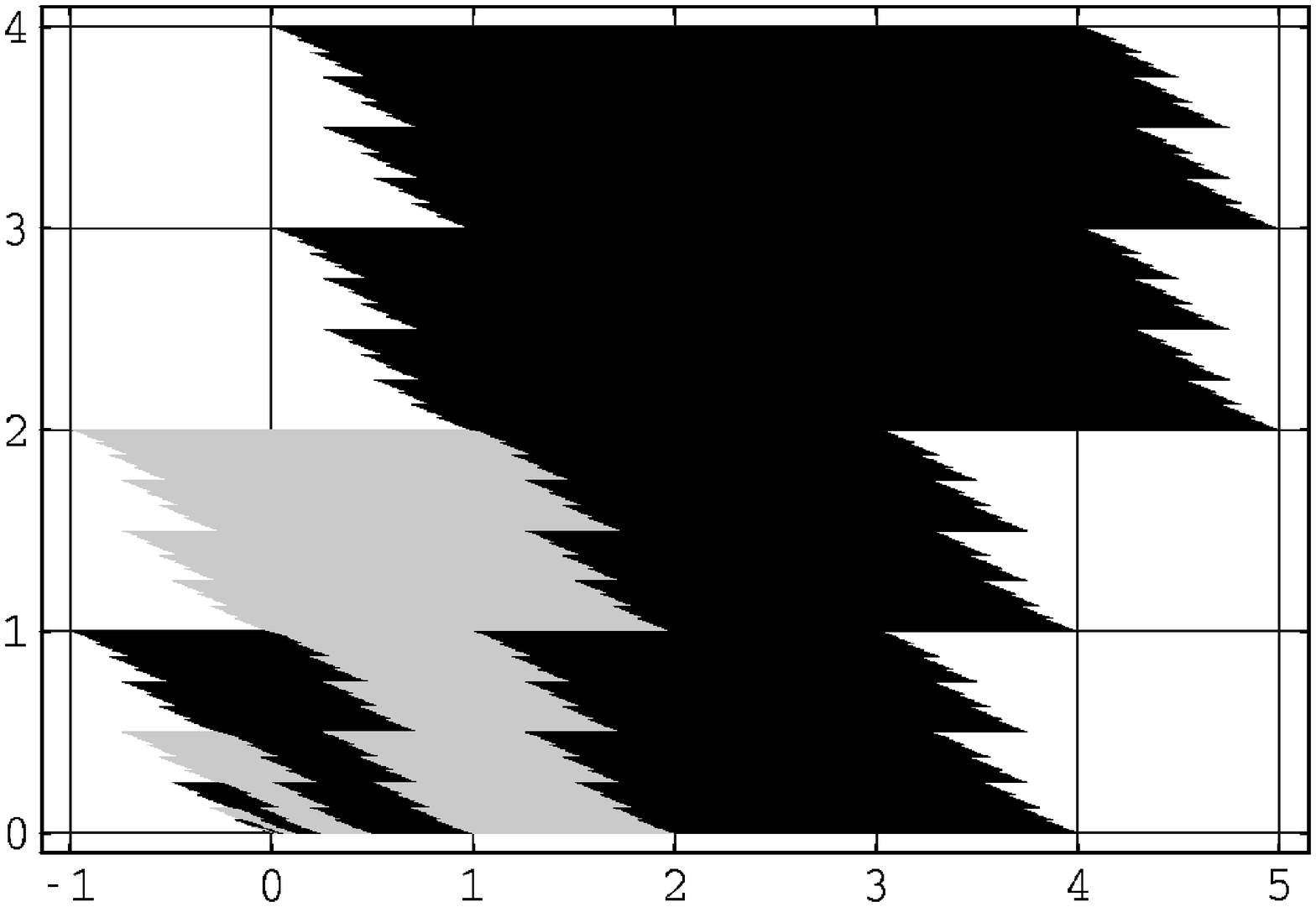,width=259bp}
\end{minipage}
}
\makebox[\textwidth]{
\begin{picture}(360,40)
\multiput(109.25,40)(-13.9,-6.95){6}{\line(-2,-1){10}}
\multiput(190.75,40)(13.9,-6.95){6}{\line(2,-1){10}}
\end{picture}
}
\makebox[\textwidth]{
\begin{minipage}{0in}
\rlap{$\mathbf{T}$}
\rlap{\ }
\end{minipage}\hfill
\begin{minipage}{3.6in}
\psfig{file=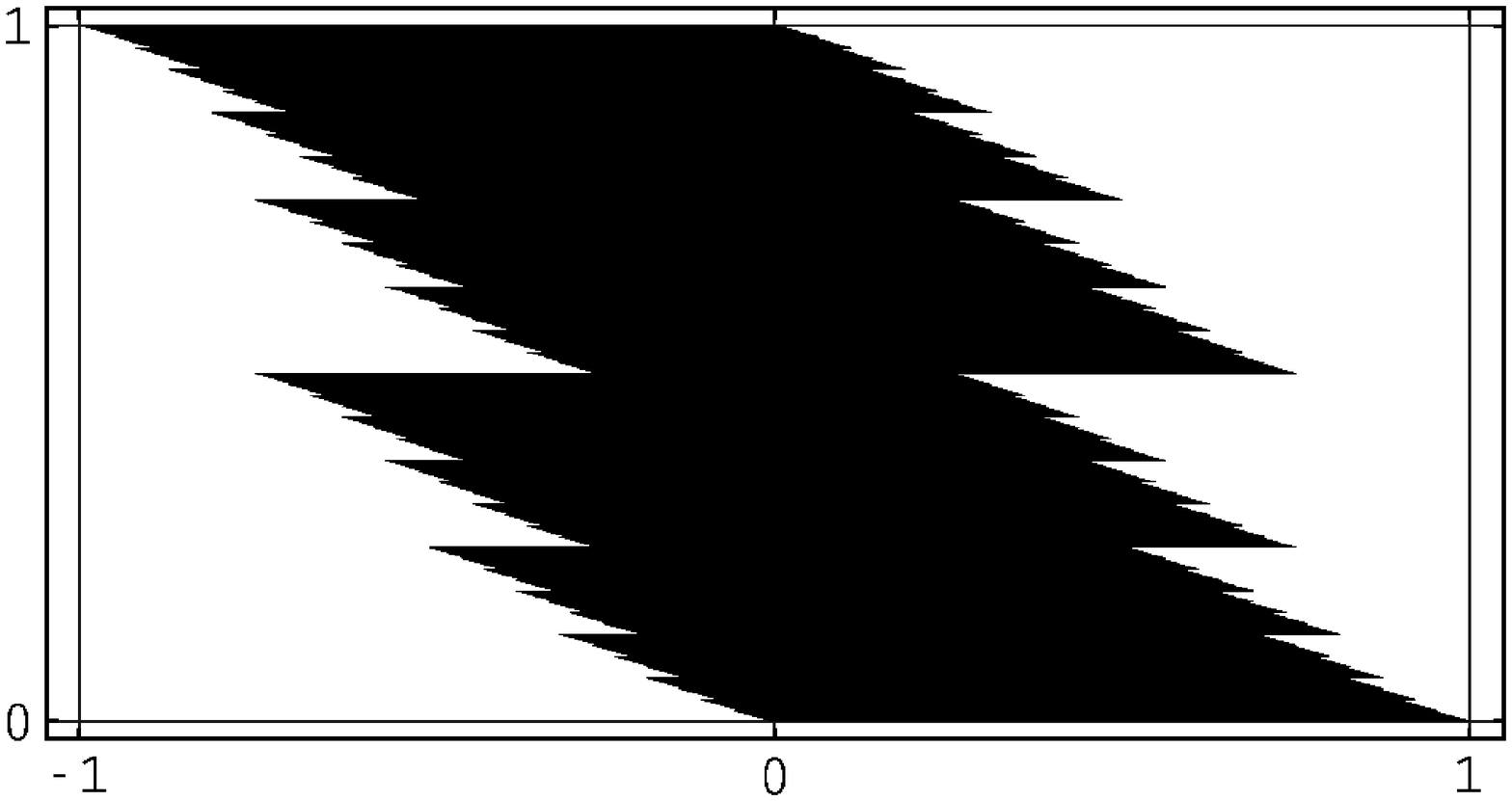,width=259bp}
\end{minipage}%
\makebox[83pt]{}%
}
\raisebox{-9pt}[12pt][0pt]{\makebox[\textwidth]{%
$\mathbf{N}=\left(
\begin{smallmatrix}
2 & 1 \\
0 & 2
\end{smallmatrix}
\right) $, $D=\left\{ \left(
\begin{smallmatrix}
0 \\
0
\end{smallmatrix}
\right), \left(
\begin{smallmatrix}
1 \\
0
\end{smallmatrix}
\right), \left(
\begin{smallmatrix}
1 \\
1
\end{smallmatrix}
\right), \left(
\begin{smallmatrix}
0 \\
1
\end{smallmatrix}
\right) \right\} $%
}}
\caption{Shark-Jawed Parallelogram (Example \protect \ref{JordanMatrix})}
\label{Reptile}
\end{figure}
But 
\begin{equation*}
\mathbf{N}^{-k}=\frac 1{4^k}
\begin{pmatrix}
2^k & -k2^{k-1} \\ 
0 & 2^k
\end{pmatrix}
\end{equation*}
so 
\begin{equation*}
\left\| \mathbf{N}^{-k}\right\| =\frac 1{4^k}\left( 2^k+k2^{k-1}\right)
=\frac 1{2^k}\left( 1+\frac k2\right) \text{.}
\end{equation*}
Hence 
\begin{align*}
\left\| R^n
\begin{pmatrix}
x_1 \\ 
x_2
\end{pmatrix}
-\mathbf{N}^{-n}
\begin{pmatrix}
x_1 \\ 
x_2
\end{pmatrix}
\right\| & \leq \frac 12\sum_{k=0}^{n-1}\frac 1{2^k}\left( 1+\frac
k2\right) \\
& =\frac 32+\left( n-3\right) \left( \frac 12\right) ^{n+1}\text{.}
\end{align*}
But $\lim_{n\rightarrow \infty }\mathbf{N}^{-n}\left( 
\begin{smallmatrix}
x_1 \\ 
x_2
\end{smallmatrix}
\right) =0$, and hence 
\begin{equation*}
\limsup_{n\rightarrow \infty }\left\| R^n
\begin{pmatrix}
x_1 \\ 
x_2
\end{pmatrix}
\right\| \leq \frac 32
\end{equation*}
for any $\left( \begin{smallmatrix}
x_1 \\ 
x_2
\end{smallmatrix}
\right) \in \mathbb{Z}^2$. Thus $R^n\left( 
\begin{smallmatrix}
x_1 \\ 
x_2
\end{smallmatrix}
\right) $ is contained in the $9$-point set 
\begin{equation*}
\left\{ 
\begin{pmatrix}
y_1 \\ 
y_2
\end{pmatrix}
\biggm| y_1,y_2=-1,0,1\right\}
\end{equation*}
for $n$ large enough (compare with Lemma \ref{Lem3.7} and Scholium \ref
{Sch3.8}). Let us compute how $R$ acts on this $9$-point set: 
{\allowdisplaybreaks
\begin{align*}
\begin{pmatrix}
-1 \\ 
-1
\end{pmatrix}
& \overset{\raise1pt\hbox{$\scriptstyle R$}}{\rightarrow }
\begin{pmatrix}
0 \\ 
-1
\end{pmatrix}
\overset{\raise1pt\hbox{$\scriptstyle R$}}{\rightarrow }
\begin{pmatrix}
0 \\ 
-1
\end{pmatrix}
\rightarrow \cdots \\
\begin{pmatrix}
-1 \\ 
0
\end{pmatrix}
& \overset{\raise1pt\hbox{$\scriptstyle R$}}{\rightarrow }
\begin{pmatrix}
-1 \\ 
0
\end{pmatrix}
\overset{\raise1pt\hbox{$\scriptstyle R$}}{\rightarrow }
\begin{pmatrix}
-1 \\ 
0
\end{pmatrix}
\rightarrow \cdots \\
\begin{pmatrix}
-1 \\ 
1
\end{pmatrix}
& \overset{\raise1pt\hbox{$\scriptstyle R$}}{\rightarrow }
\begin{pmatrix}
-1 \\ 
0
\end{pmatrix}
\overset{\raise1pt\hbox{$\scriptstyle R$}}{\rightarrow }
\begin{pmatrix}
-1 \\ 
0
\end{pmatrix}
\rightarrow \cdots \\
\begin{pmatrix}
0 \\ 
-1
\end{pmatrix}
& \overset{\raise1pt\hbox{$\scriptstyle R$}}{\rightarrow }
\begin{pmatrix}
0 \\ 
-1
\end{pmatrix}
\rightarrow \cdots \\
\begin{pmatrix}
0 \\ 
0
\end{pmatrix}
& \overset{\raise1pt\hbox{$\scriptstyle R$}}{\rightarrow }
\begin{pmatrix}
0 \\ 
0
\end{pmatrix}
\rightarrow \cdots \\
\begin{pmatrix}
0 \\ 
1
\end{pmatrix}
& \overset{\raise1pt\hbox{$\scriptstyle R$}}{\rightarrow }
\begin{pmatrix}
0 \\ 
0
\end{pmatrix}
\overset{\raise1pt\hbox{$\scriptstyle R$}}{\rightarrow }
\begin{pmatrix}
0 \\ 
0
\end{pmatrix}
\rightarrow \cdots \\
\begin{pmatrix}
1 \\ 
-1
\end{pmatrix}
& \overset{\raise1pt\hbox{$\scriptstyle R$}}{\rightarrow }
\begin{pmatrix}
1 \\ 
-1
\end{pmatrix}
\rightarrow \cdots \\
\begin{pmatrix}
1 \\ 
0
\end{pmatrix}
& \overset{\raise1pt\hbox{$\scriptstyle R$}}{\rightarrow }
\begin{pmatrix}
0 \\ 
0
\end{pmatrix}
\overset{\raise1pt\hbox{$\scriptstyle R$}}{\rightarrow }
\begin{pmatrix}
0 \\ 
0
\end{pmatrix}
\rightarrow \cdots \\
\begin{pmatrix}
1 \\ 
1
\end{pmatrix}
& \overset{\raise1pt\hbox{$\scriptstyle R$}}{\rightarrow }
\begin{pmatrix}
0 \\ 
0
\end{pmatrix}
\overset{\raise1pt\hbox{$\scriptstyle R$}}{\rightarrow }
\begin{pmatrix}
0 \\ 
0
\end{pmatrix}
\rightarrow \cdots
\end{align*}%
}%
Thus $R$ has $4$ fixed points, 
\begin{equation*}
\begin{pmatrix}
0 \\ 
0
\end{pmatrix}
,
\begin{pmatrix}
0 \\ 
-1
\end{pmatrix}
,
\begin{pmatrix}
-1 \\ 
0
\end{pmatrix}
,
\begin{pmatrix}
1 \\ 
-1
\end{pmatrix}
\text{,}
\end{equation*}
and $R^n\left( \begin{smallmatrix}
x_1 \\ 
x_2
\end{smallmatrix}
\right) $ is ultimately mapped into one of these fixed points as $%
n\rightarrow \infty $. Thus the set $B_\infty $ defined in Scholium \ref
{Sch3.8} just consists of these four fixed points, and there are $4$ atoms
in the representation, corresponding to ordinary Cuntz states on $\mathcal{O}%
_4$, by Proposition \ref{Pro5.2} and Corollary \ref{Cor5.3}. These states
can be computed explicitly using 
\begin{align*}
\sigma _{00}
\begin{pmatrix}
0 \\ 
0
\end{pmatrix}
& =
\begin{pmatrix}
0 \\ 
0
\end{pmatrix}
, \\
\sigma _{11}
\begin{pmatrix}
0 \\ 
-1
\end{pmatrix}
& =
\begin{pmatrix}
0 \\ 
-1
\end{pmatrix}
, \\
\sigma _{10}
\begin{pmatrix}
-1 \\ 
0
\end{pmatrix}
& =
\begin{pmatrix}
-1 \\ 
0
\end{pmatrix}
, \\
\sigma _{01}
\begin{pmatrix}
1 \\ 
-1
\end{pmatrix}
& =
\begin{pmatrix}
1 \\ 
-1
\end{pmatrix}
\text{.}
\end{align*}
Since these points are fixed by $R$, it follows that the corresponding atoms
have period $1$ under the action $\tau $ of $\mathbb{Z}$ on $\left( %
\mathbb{Z}^2\right) ^{\approx }$ defined by Theorem \ref{Thm4.1}, and hence
there are also $4$ cycles, each containing one of the atoms. The four points in
$B_\infty =\left( -\mathbf{T}\right) \cap \mathbb{Z}^2$ for this example
are illustrated in Figure \ref{Reptile}. (See Proposition \ref{Pro3ins1}.)

Keeping $\mathbf{N}=\left( 
\begin{smallmatrix}
2 & 1 \\
0 & 2
\end{smallmatrix}\right) $ the same, but varying $D$, we get a
different transformation $R$, see (\ref{Trans}), and therefore
different fixed points: if for example
$D^{\prime } =\left\{ \left(
\begin{smallmatrix}
0 \\
0
\end{smallmatrix}
\right) ,\left(
\begin{smallmatrix}
1 \\
0
\end{smallmatrix}
\right) ,\left(
\begin{smallmatrix}
1 \\
1
\end{smallmatrix}
\right) ,\left(
\begin{smallmatrix}
2 \\
1
\end{smallmatrix}
\right) \right\} $, then we will have the corresponding four periodic points
$B^{\prime } _\infty =\left\{ \left(
\begin{smallmatrix}
0 \\
0
\end{smallmatrix}
\right) ,\left(
\begin{smallmatrix}
-1 \\
0
\end{smallmatrix}
\right) ,\left(
\begin{smallmatrix}
0 \\
-1
\end{smallmatrix}
\right) ,\left(
\begin{smallmatrix}
-1 \\
-1
\end{smallmatrix}
\right) \right\} $. But this second version is an $\mathbf{N}^2$-scaled
version of the previous one. Note that $\mathbf{N}=\left( 
\begin{smallmatrix}
2 & 1 \\
0 & 2
\end{smallmatrix}\right) $ yields the transformation
\begin{equation*}
\mathbf{N}^{-1}+\mathbf{N}^{-2}+\dots 
=\left( \mathbf{N}-\openone \right) ^{-1}=
\begin{pmatrix}
1 & -1 \\
0 & 1
\end{pmatrix}
\end{equation*}
which, for both choices of $D$, maps the points in $D$ onto the
vertices in $\mathbf{T}$, and therefore it gets us the integral points
in $\mathbf{T}$. This is because all points in $B_\infty $ are fixed
points in this case; see (\ref{Eq3.8}).

Scholium \ref{Sch3.8} is a universal algorithm allowing an explicit
calculation of all of the $R$-periodic points $B_\infty =\left(
-\mathbf{T}\right) \cap \mathbb{Z}^\nu $ for general pairs
$\mathbf{N},D$ (as described) in $\nu $ dimensions, and
(\ref{Eq5.9})--(\ref{Periodic}) in Remark \ref{SubCuntzMatrix} is in
principle an abstract and general formula yielding all the integral
points in $\mathbf{T}$ in any dimension $\nu $. But the only general
way to solve (\ref{Eq5.9}) seems to be to use Scholium
\ref{Sch3.8}. In the examples of the Figures, $k$ from (\ref{Eq5.9})
is $1$. In Example \ref{Exa8ins2}, we also need $k=2$.

\end{example}

For other related matrix applications of this type, we refer to
\cite{Ban96}, \cite{Pri87}, \cite{BreJor}, \cite{Han92}, and \cite{PoPr94}.
\medskip

\subsection{Another matrix case: full spectrum}\label{SubsecNew9.2}

\ \bigskip

An instance where the estimate (\ref{Eq3.11}) is in fact an equality in an
equivalent norm is when the matrix $\mathbf{N}$ satisfies $\mathbf{N}^m=p%
\openone_\nu $ for some integer values $m$, $p$, i.e., $m\in \left\{
2,3,\dots \right\} $, $p\in \mathbb{Z}$. This is a special case of the
situation in Remark \ref{SubCuntzMatrix},
but with the added simplification that $\mathbf{N}^m$ is now a scalar.

\begin{remark}
\label{Rem8ins1}{\bfseries Semisimple matrices $\mathbf{N}$.}
Suppose there are integers $m$, $p$ as described such that $%
\mathbf{N}^m=p\openone_\nu $.

\begin{enumerate}
\item  \label{Rem8ins1(i)}Then the matrix $\mathbf{N}$ is diagonalizable
with spectrum $\left\{ \lambda _k\right\} $ of the form $\lambda
_k=N^{\frac 1\nu }e^{i\theta _k}=\left| p\right| ^{\frac 1m }e^{i\theta _k}$,
where $N=\left| \det \mathbf{N}\right| $,
and necessarily $m\theta _k\in 2\pi \mathbb{Z}$.

\item  \label{Rem8ins1(ii)}Let $D\subset \mathbb{Z}^\nu $ be a full set of
residues for $\mathbf{N}$ as described. If, for some $d\in D$, the solution
in $\mathbb{R}^\nu $ to $d+\mathbf{N}x=x$ is integral (i.e., $x\in \mathbb{Z}%
^\nu $) (and then $x=R\left( x\right) $), then the point 
\begin{equation*}
d^{\prime }:=d+\mathbf{N}d+\dots +\mathbf{N}^{m-1}d\in D^{\prime }:=D+%
\mathbf{N}D+\dots +\mathbf{N}^{m-1}D
\end{equation*}
satisfies $d^{\prime }+px=x$; in particular, $x\in \left( 1-p\right)
^{-1}D^{\prime }$.

\item  \label{Rem8ins1(iii)}If there are points $d_1,\dots ,d_k\in D$
such that the solution $x$ to 
\begin{equation*}
d_1+\mathbf{N}d_2+\dots +\mathbf{N}^{k-1}d_k+\mathbf{N}^kx=x
\end{equation*}
is integral, i.e., $x\in \mathbb{Z}%
^\nu $ (and then $x=R^k\left( x\right) $)
then there are also points $d_1^{\prime },\dots ,d_k^{\prime }\in
D^{\prime }$ such that 
\begin{equation*}
d_1^{\prime }+pd_2^{\prime }+\dots +p^{k-1}d_k^{\prime }+p^kx=x\text{.}
\end{equation*}
\end{enumerate}
\end{remark}

\TeXButton{Begin Proof}{\begin{proof}}Points $\lambda _k\in \limfunc{spec}%
\left( \mathbf{N}\right) $ must satisfy $\lambda _k^m=p$, so the
root-of-unity assertion in (\ref{Rem8ins1(i)}) follows. The matrix $\mathbf{N%
}$ is diagonalizable because its minimal polynomial must divide $x^m-p$.
Note that the eigenvalues of $\mathbf{N}$ need not have uniform 
multiplicity over
$\mathbb{Z}_m$ by Example \ref{Exa8ins1} below.

It follows that, if $m\in \mathbb{N}$ is chosen minimal $\mathrm{s.t.}$ $%
\mathbf{N}^m=p\openone_\nu $ for some $p\in \mathbb{Z}$, then the minimal
polynomial of $\mathbf{N}$ equals $x^m-p$ iff $\limfunc{spec}\left( \mathbf{N%
}\right) =\left\{ N^{\frac 1\nu }e^{i2\pi \frac km}\right\} _{k=0}^{m-1}$,
with the minimal polynomial defined as the monic minimal degree polynomial $%
\phi \left( x\right) \in \mathbb{C}\left[ x\right] $ satisfying $\phi \left( 
\mathbf{N}\right) =0$. Note that, in general, we have $\phi \left( x\right) $
dividing $x^m-p$; and $\phi \left( x\right) $ dividing $\det \left( x\openone%
_\nu -\mathbf{N}\right) $. If $p\neq 0$ is square-free and 
$p\neq \pm 1$, then $%
x^m-p$ is irreducible over $\mathbb{Q}$, by Eisenstein; and so $\phi \left(
x\right) =x^m-p$ also divides $\det \left( x\openone_\nu -\mathbf{N}\right) $%
.

For Examples \ref{Exa8ins1}--\ref{Exa8ins2} below, $\mathbf{N}=\left( 
\begin{smallmatrix}
1 & -1 \\
1 & 1
\end{smallmatrix}\right) $, resp., $\mathbf{N}=\left( 
\begin{smallmatrix}
0 & 1 \\
4 & 0
\end{smallmatrix}\right) $, we have the minimality condition satisfied for $%
m=4$, $p=-4$; resp., $m=2$, $p=4$. In the second example, all three
polynomials $\phi \left( x\right) $, $x^m-p$, and
$\det \left( x\openone_\nu -\mathbf{N}\right) $
coincide; whereas in the first example
$x^m-p=x^4+4$ while $\phi \left( x\right)
=\det \left( x\openone_\nu -\mathbf{N}\right) =x^2-2x+2$, and $\limfunc{%
spec}\left( \mathbf{N}\right) =\left\{ \sqrt{2}e^{\pm i\frac \pi 4}\right\} $%
.

(\ref{Rem8ins1(ii)}): Given the points $d$, $x$ as specified in (\ref
{Rem8ins1(ii)}), i.e., $d+\mathbf{N}x=x$, we have 
\begin{align*}
\left( 1-p\right) x &=\left( \openone -\mathbf{N}^m\right) x \\
&=\left( \openone +\mathbf{N}+\dots +\mathbf{N}^{m-1}\right)
\left( \openone -\mathbf{N}%
\right) x \\
&=\left( \openone +\mathbf{N}+\dots +\mathbf{N}^{m-1}\right) d \\
&=d+\mathbf{N}d+\dots +\mathbf{N}^{m-1}d=d^{\prime }\text{,}
\end{align*}
which is the assertion.

(\ref{Rem8ins1(iii)}): Setting $e:=d_1+\mathbf{N}d_2+\dots +\mathbf{N}%
^{k-1}d_k$, and iterating the given equation, we get 
\begin{align*}
x &=e+\mathbf{N}^kx=e+\mathbf{N}^ke+\mathbf{N}^{2k}x=\dots \\
&=e+\mathbf{N}^ke+\mathbf{N}^{2k}e+\dots +\mathbf{N}^{k\left( m-1\right) }e+%
\mathbf{N}^{km}x \\
&=\left( \openone +\mathbf{N}^k+\dots +\mathbf{N}^{k\left( m-1\right) }\right)
e+p^kx%
\text{.}
\end{align*}
The individual terms in the sum are 
\begin{equation*}
\sum_{i<m}\mathbf{N}^{ik}e=\sum_{i<m}\sum_{j<k}\mathbf{N}^{ik+j}d_j=%
\sum_{j<k}\sum_{i<m}\mathbf{N}^{ik+j}d_j=\sum_{l<k}\sum_{r<m}\mathbf{N}%
^{lm+r}d_{l\left( r\right) }=\sum_{l<k}p^ld_l^{\prime }\text{,}
\end{equation*}
where $d_l^{\prime }=\sum_{r<m}\mathbf{N}^rd_{l\left( r\right) }$.%
\TeXButton{End Proof}{\end{proof}}

The equation $\mathbf{N}^m=p\openone _\nu $ arises in the study 
of tilings as they
are used in wavelet theory; see, e.g., \cite{DDL95}, \cite{Ban91}, \cite
{Hou94}, \cite{HRW92}, and \cite{LaWa96}. Candidates
for dual objects for our discrete orbits in $\mathbb{Z}^\nu $, and
corresponding to the given data $\left( \mathbf{N},D\right) $, are
the sets introduced in (\ref{Eq3ins1}) and (\ref{Eq3ins2}):
\begin{equation*}
\mathbf{T}=T\left( \mathbf{N},D\right) 
:=\left\{ \sum_{i=1}^\infty \mathbf{N}^{-i}d_i%
\biggm| d_i\in D\right\} \text{.}
\end{equation*}
Recall that if the $\nu $-dimensional
Lebesgue measure of $\mathbf{T}$ is one, i.e., $\mu _\nu \left( T\left( 
\mathbf{N},D\right) \right) =1$, then $\mathbf{T}$ is a periodic tile, i.e., 
\begin{equation*}
\bigcup_{n\in \mathbb{Z}^\nu }\left( \mathbf{T}+n\right) =\mathbb{R}^\nu
\end{equation*}
and $\mu _\nu \left( \mathbf{T\cap }\left( \mathbf{T}+n\right) \right) =0$
for $n\in \mathbb{Z}^\nu $, $n\neq 0$. The $\mu _\nu $-measure-one
condition is true in our case since there are
$m$, $p$ such that $\mathbf{N}%
^m=p\openone $. Recall also from Proposition \ref{Pro3ins1} that we may
recover $B_\infty $ as
$B_\infty =\mathbb{Z}^\nu \cap \left\{ -\mathbf{T}\right\} $.
For Example \ref{JordanMatrix}, above, the tile $\mathbf{T}$ is the
parallelogram in Figure \ref{Reptile}, above, and we see from
this figure that $B_\infty =-\left\{ \left(
\begin{smallmatrix}
0 \\
0
\end{smallmatrix}
\right) ,\left(
\begin{smallmatrix}
0 \\
1
\end{smallmatrix}
\right) ,\left(
\begin{smallmatrix}
1 \\
0
\end{smallmatrix}
\right) ,\left(
\begin{smallmatrix}
-1 \\
1
\end{smallmatrix}
\right) \right\} $ in accordance with Proposition \ref{Pro3ins1}.
The appearance of this tile can be explained as follows: from
(\ref{Eq3ins2}) and the expression for $\mathbf{N}^{-k}$ in Example
\ref{JordanMatrix}, it follows that
\begin{equation*}
\sum _{k=1}^\infty \mathbf{N}^{-k}d_k=
\begin{pmatrix}
x-\phi \left( y\right) \\
y
\end{pmatrix}
\end{equation*}
if $d_k=
\begin{pmatrix}
i_k \\
j_k
\end{pmatrix}
$, $i_k,j_k=0,1$, where 
\begin{equation*}
x=\sum _k\frac 1{2^k}i_k\text{,\quad }y=\sum _k\frac 1{2^k}j_k
\end{equation*}
are arbitrary numbers in the interval $\left[ 0,1\right] $,
given in dyadic expansion. The function $\phi $, defined by 
\begin{equation*}
\phi \left( y\right) =\sum _k\frac k{2^{k+1}}j_k\text{,}
\end{equation*}
is strictly speaking not a function of $y\in \left[ 0,1\right] $
alone, but depends on the particular dyadic expansion used for $y$
(that is, $\phi $ is a function on the Cantor set obtained by cutting
and doubling points at each dyadic rational in $\left[ 0,1\right] $). 
Now, if $y$ is the dyadic rational $y=p2^{-k}$, where 
$1\leq p\leq 2^k-1$ and $p$ is odd, then $j_k=1$, and
\begin{equation*}
y=\frac {p-1}{2}2^{-k+1}+2^{-k}=\frac {p-1}{2}2^{-k+1}+
\sum _{l=k+1}^\infty 2^{-l}\text{.}
\end{equation*}
Thus
\begin{equation*}
\phi \left( y_+\right) =\phi \left( y\right) 
=\phi \left( \frac {p-1}{2}2^{-k+1}\right) +k2^{-k-1}
\end{equation*}
and
\begin{equation*}
\phi \left( y_-\right) =\phi \left( \frac {p-1}{2}2^{-k+1}\right) 
+\sum _{l=k+1}^\infty l2^{-l-1}\text{.}
\end{equation*}
Thus $\phi $ is continuous everywhere except at the dyadic rationals,
and $\phi $ is continuous to the right even at the dyadic
rationals. If $y=p2^{-k}$ with $p$ odd, we have
\begin{equation*}
\phi \left( y_+\right) -\phi \left( y_-\right) 
=k2^{-k-1}-\sum _{l=k+1}^\infty l2^{-l-1}=-2^{-k}\text{.}
\end{equation*}
Now, since
\begin{equation*}
\mathbf{T}=\left\{ 
\begin{pmatrix}
x-\phi \left( y\right) \\ 
y
\end{pmatrix}
\biggm|
\begin{pmatrix}
x \\
y
\end{pmatrix}
\in \left[ 0,1\right] \times \left[ 0,1\right] \right\} 
=\bigcup _{y\in \left[ 0,1\right] }
\begin{pmatrix}
\left[ -\phi \left( y\right) ,-\phi \left( y\right) +1\right] \\
y
\end{pmatrix}
\end{equation*}
this explains the appearance of the shark-jawed parallelogram in
Figure \ref{Reptile} above. In particular, the teeth constitute the
graph of the function $\phi $, and one sees that at $y=\frac 12$ the
graph makes a jump of $-\frac 12$, and at $y=\frac 14$ and $\frac 34$
the graph makes a jump of $-\frac 14$, etc., in accordance with the
description above of $\phi $.
\medskip

\subsection{Self-similarity and tiles}\label{SelfSimTile}

\ \bigskip

The more detailed tiling properties of $\mathbf{T}$ give us a precise
way of relating properties of the mapping $R$ from (\ref{Trans}) to
our Fourier representation (\ref{Composition}) for the
$\mathcal{O}_N$-representation, and then back to the $\ell ^2$-basis
viewpoint in Section \ref{Sec4}. If $\mathbf{T}$ is in fact a
$\mathbb{Z}^\nu $-tile for $\mathbb{R}^\nu $, then it follows from
\cite{Jor95} that the $L^2\left( \mathbb{T}^\nu \right) $
representation in (\ref{Composition}) is unitarily equivalent to a
representation on $L^2\left( \mathbf{T}\right) $, where the measure on
$\mathbf{T}$ is simply the restriction of the $\mathbb{R}^\nu
$-Lebesgue measure. The mapping $z\mapsto z^\mathbf{N}$ of
$\mathbb{T}^\nu $ into itself in (\ref{Composition}) is then
equivalent to the endomorphism $\tau $ on $\mathbf{T}$ defined as
follows: let $t\in \mathbf{T}$; pick the unique representation:
$\mathbf{N}t=s+n$, $s\in \mathbf{T}$, $n\in \mathbb{Z}^\nu $, and then
define $\tau \left( t\right) =s$.  This gives a concrete realization
of our symbolic coding mappings $\sigma _i$ and $R$ from
(\ref{EqNew6.1}) and (\ref{EqNew6.2}), respectively. For, if $D$ is
represented as $D=\left\{ d_i\right\} _{i=1}^N$, and $\tau _i\left(
x\right) :=\mathbf{N}^{-1}\left( x+d_i\right) $, then $\tau
_i:\mathbf{T}\rightarrow \mathbf{T}$, and $\tau \circ \tau
_i=\limfunc{id}_\mathbf{T}$ for all $i=1,\dots ,N$. So $\tau $ is a
fractal version of the integral mapping $R$ in (\ref{EqNew6.2}).

More generally, suppose $0\in D$, and $D$ is a full set of residues,
i.e., a selection of a point in $\mathbb{Z}^\nu $ for each of the $N$
cosets of $\mathbb{Z}^\nu \diagup \mathbf{N}\mathbb{Z}^\nu $. Then let
$\mathbb{L\kern0.5pt}:=\mathbb{Z}\left( \mathbf{N},D\right) $ be the
smallest $\mathbf{N}$-invariant lattice containing $D$, and suppose
$\mathbf{T}$ tiles by $\mathbb{L\kern0.5pt}$. Then, using that
$\mathbf{N}\mathbb{L\kern0.5pt}\subset \mathbb{L\kern0.5pt}$, note
that $\tau :\mathbf{T}\rightarrow \mathbf{T}$ may be defined as
before: $\mathbf{N}t=s+n$, $t,s\in \mathbf{T}$, $n\in
\mathbb{L\kern0.5pt}$, $\tau \left( t\right) =s$; and we will still
have the formula $\tau \circ \tau _i=\limfunc{id}_\mathbf{T}$ from
above, but now in the more general case.  In fact (\ref{Eq3ins1})
shows that $\tau _i$ and $\tau $ may be defined on $\mathbf{T}$
independently of $\mathbb{L\kern0.5pt}$. (See also the discussion in
Scholium \ref{LagerLattice} below, especially the conclusion
(\ref{Ortho}) there. We may then define a representation of
$\mathcal{O}_N$ on $L^2\left( \mathbf{T}\right) $ by the following
explicit formulas for the operators:
\begin{align}
S_i\psi &=N^\frac 12\chi _{\tau _i\left( 
\mathbf{T}\right) }\psi \circ \tau
\label{Up} \\
\intertext{and}
S^*_i\psi &=N^{-\frac 12}\psi \circ \tau _i\quad \left( 
\psi \in L^2\left( \mathbf{T}\right) \right) \text{.}
\label{Down}
\end{align}

The argument for why the $\mathcal{O}_N$-relations (\ref{Eq2.3}) are
indeed satisfied for the operators $S_i$ defined in (\ref{Up}) is from
\cite{JoPe94}. Indeed the orthogonality relation in (\ref{Eq2.3}) is
clear, and $\sum _iS^{}_iS^*_i=\openone $ follows from (\ref{Eq3ins1})
combined with the formula (\ref{Down}) for $S^*_i$. We also use the
basic properties of the Lebesgue measure, when restricted to
$\mathbf{T}$: we have, for $\psi \in L^2\left( \mathbf{T}\right) $,
\begin{align*}
\sum _i\left\langle \psi \mid S^{}_iS^*_i\psi \right\rangle 
& =\sum _i \left\| S^*_i\psi \right\| ^2 
=\sum _iN^{-1} \int _\mathbf{T}\left| \psi \circ \tau _i \right| ^2 \,dx \\
& =\sum _iN^{-1} \int _\mathbf{T}\left| \psi \left( 
\mathbf{N}^{-1}\left( x+d_i\right) \right) \right| ^2 \,dx \\
& =\int _\mathbf{T}\left| \psi \right| ^2 \,dx
=\left\langle \psi \mid \psi \right\rangle \text{,}
\end{align*}
where (\ref{Eq3ins1}), and the corresponding formula for Lebesgue
measure, was used in the last step.
If $\mathbf{T}$ is not a $\mathbb{Z}^\nu $-tile, it is still a finite
covering of the torus $\mathbb{T}^\nu $, as long as $D$ is assumed to
be a full set of residues for $\mathbb{Z}^\nu \diagup
\mathbf{N}\mathbb{Z}^\nu $, as follows from \cite[Theorem 6.2]{JoPe92}, 
which also gives an expression for the covering
number. We refer to \cite[Theorem 6.2]{JoPe92} for the definition of a
\emph{finite covering:} what is implied here is that
$\mathbf{T}=\bigcup _{i=1}^qE_i$ where the $E_i$'s are $\mathbb{Z}^\nu
$-periodic tiles, and where the intersections $E_i\cap E_j$, $i\neq
j$, have $\nu $-dimensional Lebesgue measure equal to zero. In
particular, we conclude that the Lebesgue measure of $\mathbf{T}$ must
be an integer. To apply Theorem 6.2 in \cite{JoPe92} we first show
that $\left\{ e_n\left( x\right) :=e^{i2\pi n\cdot x}\right\} $, $x\in
\mathbf{T}$, $n\in \mathbb{Z}^\nu $, is an orthogonal family in
$L^2\left( \mathbf{T}\right) $, relative to the restricted Lebesgue
measure. In general, $\left\{ e_n|_\mathbf{T}\right\} _{n\in
\mathbb{Z}^\nu }$ will not be a basis for $L^2\left( \mathbf{T}\right)
$. For that it is necessary and sufficient that $\mathbf{T}$ be a
$\mathbb{Z}^\nu $-tile. But the functions are always mutually
orthogonal in $L^2\left( \mathbf{T}\right) $. For the general theory,
see \cite{LiMa}.

We now turn to the specifics of the aforementioned application of
\cite[Theorem 6.2]{JoPe92}. It is based on the following general fact. 
\begin{proposition}
\label{FiniteCover}Let $\mathbf{N},D,\mathbf{T}$ be as described
above \textup{(}i.e., with the eigenvalues $\lambda _i$ of
$\mathbf{N}$ satisfying $\left| \lambda _i\right| >1$, $D$ a full set
of residues for $\mathbb{Z}^\nu \diagup \mathbf{N}\mathbb{Z}^\nu $,
and $\mathbf{T}$ the fractal determined by
\textup{(}$\ref{Eq3ins1}$\textup{)).}  We then have
\begin{equation*}
\int _\mathbf{T}e_n\left( x\right) \,dx=0
\end{equation*}
for all $n\in \mathbb{Z}^\nu \,\diagdown \left\{ 0\right\} $; 
and therefore $\mathbf{T}$ is a finite cover space for the torus
$\mathbb{T}^\nu =\mathbb{R}^\nu \diagup \mathbb{Z}^\nu $.
\end{proposition}

\begin{proof}
Let $\mathbf{N},D$ be given as specified, and let $\mathbf{M}$
denote the transposed matrix $\mathbf{N}^{\limfunc{tr}}$. 
For $x\in \mathbb{R}^\nu $, let
\begin{equation*}
F_D\left( x\right) 
:=\sum _{d\in D}e^{i2\pi x\cdot d}=\sum _{d\in D}e_d\left( x\right) \text{.}
\end{equation*}
An iteration of (\ref{Eq3ins1}) in the form $\mathbf{T}=\bigcup _{d\in
D}\mathbf{N}^{-1}\left( d+\mathbf{T}\right) $ then yields the
following product formula:
\begin{equation*}
\int _\mathbf{T}e_n\left( x\right) \,dx=F\left( \mathbf{M}^{-1}n\right) 
F\left( \mathbf{M}^{-2}n\right) \cdots F\left( \mathbf{M}^{-k}n\right) 
\int _{\mathbf{N}^{-k}\mathbf{T}}e_n\left( x\right) \,dx\text{,}
\end{equation*}
valid for all $n\in \mathbb{Z}^\nu $, $k=1,2,\dots $, and where $dx$
denotes the $\nu $-dimensional Lebesgue measure. Assuming $n\in
\mathbb{Z}^\nu \,\diagdown \left\{ 0\right\} $, use (\ref{Eq3.6bis})
for the matrix $\mathbf{M}$, and pick $k$ large enough to get a
non-trivial residue for $n$, taking into account (\ref{Eq3.7}). We
will then have
\begin{equation*}
F\left( \mathbf{M}^{-1}n\right) =\cdots =F\left(
\vphantom{\mathbf{M}^{-1}}\smash{\mathbf{M}^{-\left(
k-1\right) }n}\right) =N\text{,}
\end{equation*}
and $F\left( \mathbf{M}^{-k}n\right) =0$; and the result from the
Proposition follows.  Indeed, if $n\in \mathbb{Z}^\nu \,\diagdown
\left\{ 0\right\} $, then there is some $k\in \mathbb{N}$, and
$n^{\prime },r\in \mathbb{Z}^\nu $, such that $n=\mathbf{M}^kn^{\prime
}+\mathbf{M}^{k-1}r$ and $r\notin \mathbf{M}\left( \mathbb{Z}^\nu
\right) $. We then get
\begin{equation*}
F_D\left( \mathbf{M}^{-k}n\right) =F_D\left( 
n^{\prime }+\mathbf{M}^{-1}r\right) =\sum _{d\in D}e_d\left( 
n^{\prime }+\mathbf{M}^{-1}r\right) =\sum _{d\in D}e_d\left( 
\mathbf{M}^{-1}r\right) \text{.}
\end{equation*}
But now we can use $D$ as a set of representatives for the elements in
the finite group $\mathbb{Z}^\nu \diagup \mathbf{N}\left(
\mathbb{Z}^\nu \right) $. Since $\mathbf{M}^{-1}r\notin \mathbb{Z}^\nu
$, the last sum is an average of a non-trivial character on a finite
group, and so it is zero by a standard fact about group
characters. Indeed, replacing $D=\left\{ d_i\right\} _{i=1}^N$ with
$d_i^{\prime }=d_i+\mathbf{N}l_i$, for any $l_i\in \mathbb{Z}^\nu $,
will not affect the sum. It then follows from \cite[Theorem
6.2]{JoPe92} that $\mathbf{T}$ is a finite cover of the standard
torus.
\end{proof}

As a consequence, in the general case, even when $\mathbf{T}$ is not a
$\mathbb{Z}^\nu $-tile, our $\mathcal{O}_N$-rep\-re\-sen\-ta\-tions
are realized on $L^2\left( \mathbf{T}\right) $ in the same manner as
in the special case (\ref{Up})--(\ref{Down}). The representation on
$L^2\left( \mathbf{T}\right) $ which is defined by (\ref{Up}) is not
unitarily equivalent to our original permutative representation on
$\ell ^2\left( \mathbb{Z}^\nu \right) $ from (\ref{Composition}) (or,
equivalently, from (\ref{Eq3.6bis})). But when (\ref{Up}) is adjusted
by a cocycle as described in \cite[Section 2]{BJP}, then the old
representation on $\ell ^2\left( \mathbb{Z}^\nu \right) \simeq
L^2\left( \mathbb{T}^\nu \right) $ will be \emph{contained in} the
corresponding (adjusted) representation on $L^2\left(
\mathbf{T}\right) $ where the intertwining operator $W:\ell ^2\left(
\mathbb{Z}^\nu \right) \rightarrow L^2\left( \mathbf{T}\right) $ is
the one which takes the basis $\left\{ e_n:n\in \mathbb{Z}^\nu
\right\} $ for $\ell ^2\left( \mathbb{Z}^\nu \right) $ into the
functions $f_n\left( x\right) =e^{i2\pi n\cdot x}$, restricted to
$\mathbf{T}$, i.e., $x\in \mathbf{T}:We_n=f_n$, $n\in \mathbb{Z}^\nu
$. As an application of Proposition \ref{FiniteCover}, we note that
the two representations, the ``old'' one (\ref{Eq3.6bis}), and the new
one on $L^2\left( \mathbf{T}\right) $, will be unitarily equivalent
iff $\mathbf{T}$ is a $\mathbb{Z}^\nu $-tile. The present section has
examples of both $\mathbf{T}$ having this $\mathbb{Z}^\nu $-tiling
property, and not: Example $D_1$ has it (see Figure
\ref{RedCrossFigure}), while $D_3$, $D_5$, and $D_9$ do not (see
Figures \ref{CloudNineFigure}, \ref{CloudThreeFigure}, and
\ref{CloudFiveFigure}).

\begin{theorem}
\label{Intertwines}Let $\left( \mathbf{N},D\right) $ be as in
Proposition \ref{FiniteCover}. Let
$\mathbf{M}:=\mathbf{N}^{\limfunc{tr}}$ be the transposed matrix, and
let $G=\left\{ r_j\right\} _{j=1}^N\subset \mathbb{Z}^\nu $ be a full
set of residues for $\mathbf{M}$ such that $N^{-\frac 12}\left(
e^{i2\pi r_j\cdot \mathbf{N}^{-1}d_k}\right) _{jk}$ is unitary as an
$N\times N$ matrix. Given $D$, a corresponding $G$ can be found by
group duality. Let $\left\{ V_j\right\} $ be the representation of
$\mathcal{O}_N$ which is determined by $\left( \mathbf{M},G\right) $
and \textup{(\ref{Eq3.6bis})}, i.e., $V_je_n=e_{r_j+\mathbf{M}n}$,
$n\in \mathbb{Z}^\nu $, and let $W:\ell ^2\left( \mathbb{Z}^\nu
\right) \rightarrow L^2\left( \mathbf{T}\right) $ be the isometry
which is given by $e_n\mapsto f_n$ according to Proposition
\ref{FiniteCover}. Let $m_{ji}:=N^{-\frac 12}f_{r_i}\circ \tau _j$,
and $T_i:=\sum _{j=1}^NS_jm_{ji}$, where $\left\{ S_j\right\} $ is
given by \textup{(\ref{Up})}. Then $W$ intertwines, i.e.,
\begin{equation}
WV_i=T_iW\text{,\quad }i=1,\dots ,N\text{.}
\label{Inter}
\end{equation}
\end{theorem}

\begin{proof}
See the formulas in Section \ref{Monomial} for the $\mathbf{N}$
representations, and Section 2 in \cite{BJP}. The modified formulas
(\ref{Up})--(\ref{Down}) for the $L^2\left( \mathbf{T}\right)
$-representations may be written in the following explicit form: let
$\left( \mathbf{N},D\right) $ and $\left( \mathbf{M},G\right) $ be
chosen as before in duality, i.e., $D=\left\{ d_i\right\} _{i=1}^N$
and $G=\left\{ r_i\right\} _{i=1}^N$ residue sets in $\mathbb{Z}^\nu $
as described. The duality condition will be satisfied for the pair if,
e.g.,
\begin{equation}
r_j\cdot \mathbf{N}^{-1}\left( d_k\right) \in\frac {\left( j-1\right) 
\left( k-1\right) }N+\mathbb{Z}
\label{Polar}
\end{equation}
for all $j,k\in \left\{1,2,\dots ,N\right\} $, where $N=\left| \det
\mathbf{N}\right| $ as usual. This choice will work if $N$ is a
prime. Then we claim that the representation $\left\{ T_i\right\} $ on
$L^2\left( \mathbf{T}\right) $ \emph{contains} a copy of the original
$\mathcal{O}_N$-representation on $\ell ^2\left( \mathbb{Z}^\nu
\right) $. It is determined by:
\begin{align}
T_i\psi &=f_{r_i}\left. \psi \circ \tau \right. \text{,}
\label{UpPrime}\tag{$\ref{Up}^{\prime }$}\\
\intertext{and}
T^*_i\psi &=\frac 1N\sum _{j=1}^N\overline{f_{r_i}\circ \tau _j}\left. 
\psi \circ \tau _j \right. \quad\left( \psi \in L^2\left( \mathbf{T}\right) 
\right) \text{,}
\label{DownPrime}\tag{$\ref{Down}^{\prime }$}
\end{align}
where $\tau $ and $\tau _j$ are the mappings on $\mathbf{T}$ from
$\left( \mathbf{N},D\right) $ as in (\ref{Up})--(\ref{Down}). Notice
that (\ref{UpPrime}) should be compared to (\ref{Composition}) from
the Introduction, and the unitarity condition for the $N\times N$
matrix in Theorem \ref{Intertwines} is a generalization of the
condition on the matrix (\ref{Eq1.1}).The basic idea behind the
isometric intertwiner
\begin{equation}
We_n=f_n=e^{i2\pi n\cdot x}\text{,\quad }n\in \mathbb{Z}^\nu ,\;x\in 
\mathbf{T}\text{,}
\label{Isometry}
\end{equation}
is from \cite{JoPe92} and \cite{JoPe94} (which also contain more
details on this point).

Firstly, it is immediate from (\ref{UpPrime}) that the intertwining
relation (\ref{Inter}) holds on $\ell ^2\left( \mathbb{Z}^\nu \right)
$, where the $\mathcal{O}_N$-representation $\left\{ V_i\right\} $,
acting on $\ell ^2\left( \mathbb{Z}^\nu \right) $, is given by the
formula
\begin{equation}
V_ie_n=e_{r_i+\mathbf{M}n}\text{,\quad }n\in \mathbb{Z}^\nu 
\label{Dual}
\end{equation}
on the canonical basis vectors. Indeed, for $n\in \mathbb{Z}^\nu $, 
$x\in \mathbf{T}$,
\begin{equation*}
T_iWe_n\left( x\right) =f_{r_i}\left( x\right) f_n\left( \mathbf{N}x\right) 
=f_{r_i}\left( x\right) f_{\mathbf{M}n}\left( x\right) 
=f_{r_i+\mathbf{M}n}\left( x\right) =WV_ie_n\left( x\right) \text{,}
\end{equation*}
proving (\ref{Inter}).

We now show that $\left\{ T_i\right\} $ satisfies the
$\mathcal{O}_N$-relations (\ref{Eq2.3}), and that, in addition to
(\ref{Inter}), we also have 
\begin{equation*}
T^*_iW=WV^*_i\text{,\quad }i=1,\dots ,N\text{.}
\end{equation*}

To see that $\sum _iT^{}_iT^*_i=\openone $ on $L^2\left( \mathbf{T}\right) $,
we proceed as follows: let $\psi \in L^2\left( \mathbf{T}\right) $. Then 
\begin{align*}
\sum _i\left\| T^*_i\psi \right\| ^2 &= \sum _i\frac 1{N^2}\int 
_{\mathbf{T}}\left| \smash{\sum _j\overline{f_{r_i}\circ \tau _j}\left. 
\psi \circ \tau _j\right. }\vphantom{\sum f}\right| ^2 \,dx \\
&= \frac 1{N^2}\sum _j\sum _k\int _{\mathbf{T}}\sum _i\overline{f_{r_i}\circ 
\tau _j}\left. f_{r_i}\circ \tau _k\right. \left. \psi \circ \tau _j\right. 
\overline{\psi \circ \tau _k}\,dx \\
&= \frac 1N\sum _j\sum _k\delta _{jk}\int _{\mathbf{T}}\left| \psi \right| 
^2\circ \tau _j \,dx= \frac 1N\sum _j\int _{\mathbf{T}}\left| \psi \right| 
^2\circ \tau _j \,dx \\ 
&= \int _{\mathbf{T}}\left| \psi \right| ^2 \,dx
\end{align*}
which is the desired conclusion.

We claim that $\ell ^2\left( \mathbb{Z}^\nu \right) $ in fact defines
a subrepresentation of $\mathcal{O}_N$ in $L^2\left( \mathbf{T}\right)
$ in the strict sense: indeed, let $\ell ^2\left( \mathbb{Z}^\nu
\right) $ be viewed as a subspace of $L^2\left( \mathbf{T}\right) $
via the isometry $W$, and let $V_i$ define the original $\ell ^2\left(
\mathbb{Z}^\nu \right) $ representation by (\ref{Dual}). There are
then operators $U_i$, $B_i$, defined on $L^2\left( \mathbf{T}\right)
\ominus \ell ^2\left( \mathbb{Z}^\nu \right) $, such that
\begin{equation*}
T_i=
\begin{pmatrix}
V_i & B_i \\
0 & U_i
\end{pmatrix}
\text{,\quad }i=1,\dots ,N\text{.}
\end{equation*}
The exact form of the complementary operators $B_i$, $U_i$ can be read
off from (\ref{UpPrime})--(\ref{DownPrime}), and they are both zero
precisely when $\mathbf{T}$ is a $\mathbb{Z}^\nu $-tile. Since both
systems $\left\{ V_i\right\} $ and $\left\{ T_i\right\} $ satisfy the
Cuntz relations (\ref{Eq2.3}) in the respective Hilbert spaces, we get
\begin{equation*}
\sum _i\left(V^{}_iV^*_i+B^{}_iB^*_i\right) =\openone \text{,}
\end{equation*}
and therefore $B_i=0$ for all $i$. It follows that $T_i=\left( 
\begin{smallmatrix}
V_i & 0 \\
0 & U_i
\end{smallmatrix}
\right) $, and that the operators $U_i$ also define an
$\mathcal{O}_N$-representation, now on $L^2\left( \mathbf{T}\right)
\ominus \ell ^2\left( \mathbb{Z}^\nu \right) $.
\end{proof}

\begin{scholium}
\label{LagerLattice}In summary, we have shown that, for every pair $\left(
\mathbf{M},G\right) $ in $\mathbb{Z}^\nu $ (i.e., $G$ a full set of
residues for $\mathbb{Z}^\nu \diagup \mathbf{M}\mathbb{Z}^\nu $),
there is a dual pair $\left( \mathbf{N},D\right) $ (see (\ref{Polar}))
such that the $\mathcal{O}_N$-representation defined from $\left(
\mathbf{M},G\right) $ on $\ell ^2\left( \mathbb{Z}^\nu \right) $ by
(\ref{Dual}) is unitarily equivalent to a subrepresentation of the
$\mathcal{O}_N$-representation on $L^2\left( \mathbf{T}\right) $
defined from $\left( \mathbf{N},D\right) $ via formulas
(\ref{UpPrime})--(\ref{DownPrime}). The intertwining isometry $\ell
^2\left( \mathbb{Z}^\nu \right)
\overset{\smash[b]{\lower1pt\hbox{$\scriptstyle W$}}}{\hookrightarrow
}L^2\left( \mathbf{T}\right) $ is the one defined from Proposition
\ref{FiniteCover}; see (\ref{Isometry}).

We also note that the orthogonality condition in Proposition
\ref{FiniteCover} generally holds for a \emph{larger} lattice than
$\mathbb{Z}^\nu $. In fact, given $\left( \mathbf{N},D\right) $, $0\in
D$, as before, we introduced the smallest lattice
$\mathbb{L\kern0.5pt}$ containing $D$ and invariant under
$\mathbf{N}$, and the corresponding \emph{dual lattice}
\begin{equation*}
\mathbb{L\kern0.5pt}^*=\left\{ \xi \in \mathbb{R}^\nu \mid \xi \cdot \lambda 
\in \mathbb{Z}, \forall \lambda \in \mathbb{L\kern0.5pt}\right\} \text{.}
\end{equation*}
Using the same argument as from the proof of Proposition \ref{FiniteCover}
and (\ref{Polar}), we can show that, in fact, 
\begin{equation}
\int _{\mathbf{T}}f_\xi \left( x\right) \,dx=0\text{,\quad }\forall \xi \in 
\mathbb{L\kern0.5pt}^*\,\diagdown \left\{ 0\right\} \text{.}
\label{Ortho}
\end{equation}

For those special cases when $\mathbf{T}$ is already an
$\mathbb{L\kern0.5pt}$-periodic tile (for example, for some of the
``cloud'' examples in Examples \ref{CloudNine} and \ref{Clouds} below,
e.g. $D_3$ and $D_5$ in (\ref{CloudD})), we then conclude from
formulas (\ref{UpPrime}), (\ref{Inter}), and (\ref{Dual}) above that
the $\left\{ T_i\right\}$-representation of $\mathcal{O}_N$ on
$L^2\left( \mathbf{T}\right) $ is then permutative. The argument is
also based on the following simple implication: $\mathbf{N}\left(
\mathbb{L\kern0.5pt}\right) \subset \mathbb{L\kern0.5pt} \Rightarrow
\mathbf{M}\left( \mathbb{L\kern0.5pt}^*\right) \subset
\mathbb{L\kern0.5pt}^*$, when
$\mathbf{M}=\mathbf{N}^{\limfunc{tr}}$. Also recall $\left(
\mathbf{N}\mathbb{Z}^\nu \right) ^*=\mathbf{M}^{-1}\left(
\mathbb{Z}^\nu \right) $. Further recall that
$\mathbb{L\kern0.5pt}\subset \mathbb{Z}^\nu $, and therefore
$\mathbb{Z}^\nu \subset \mathbb{L\kern0.5pt}^*$. For the statement of
the result (\ref{Ortho}), we need to observe that both formulas
(\ref{Isometry}) and (\ref{Dual}) make sense if $n$ is taken to be an
element in the (generally) larger lattice $\mathbb{L\kern0.5pt}^*$,
rather than in $\mathbb{Z}^\nu $. In considering the general form of
(\ref{Isometry}), we think of $e_n$ as an element in $\ell ^2\left(
\mathbb{L\kern0.5pt}^*\right) $ when $n\in \mathbb{L\kern0.5pt}^*$;
and similarly, for the generalized (\ref{Dual}), we use that $n\mapsto
r_i+\mathbf{M}n$ is a mapping of $\mathbb{L\kern0.5pt}^*$ into itself
for all $r_i\in G$, and defines an iterated function system. We also
have properties (\ref{Eq2.4}) for these mappings, relative to
$\mathbb{L\kern0.5pt}^*$.

The reference above to reptiles $\mathbf{T}$ from residue sets $D_3$
and $D_5$ refers to a construction of a two-parameter family of
examples to be taken up in detail in Subsection \ref{TilePlane} below,
and we refer to the $\mathbf{T}$'s as ``Cloud'' examples because of
their appearance; see Figures \ref{CloudNineFigure},
\ref{CloudThreeFigure}, and \ref{CloudFiveFigure}. The explicit
numbers for the various residue sets $D$ for $\mathbf{N}=\left(
\begin{smallmatrix}
1 & 2 \\
-2 & 1
\end{smallmatrix}
\right) $ are listed in (\ref{CloudD}), and the subscript $5$ in $D_5$
and $\mathbf{T}_5$ (for example) is the number of points in the
corresponding $B_\infty $-set, in this case $B_\infty ^{(5)}$.
\end{scholium}

We noted above that $\mathbf{T}_3$ and $\mathbf{T}_5$ tile the plane
(i.e., $\mathbb{R}^2$) by the respective lattices
$\mathbb{L\kern0.5pt}_3$ and $\mathbb{L\kern0.5pt}_5$ where
\begin{equation*}
\mathbb{L\kern0.5pt}_3=\mathbb{Z}\left[ 
\begin{pmatrix}
2 \\
0
\end{pmatrix}%
,%
\begin{pmatrix}
0 \\
1
\end{pmatrix}%
\right] \text{\quad and\quad }
\mathbb{L\kern0.5pt}_5=\mathbb{Z}\left[ 
\begin{pmatrix}
1 \\
0
\end{pmatrix}%
,%
\begin{pmatrix}
0 \\
2
\end{pmatrix}%
\right] \text{.}
\end{equation*}
See Example \ref{Clouds} below for more details. The proof of this
(and a corresponding negative observation for $\mathbf{T}_9$) is based
on the following general
\begin{lemma}
\label{TilesByL}Let 
$\mathbf{N}$ be a $\nu \times \nu $ expansive matrix with integer
entries, and let $D$ be a full set of residues for $\mathbb{Z}^\nu
\diagup \mathbf{N}\mathbb{Z}^\nu $, and assume that $0\in D$. Let
$\sim $ denote the $\mathcal{O}_N$-equivalence relation on
$\mathbb{Z}^\nu $ from Scholium \ref{Sch2.4}, and let $\mathbf{T}$ be
the fractal determined by \textup{(\ref{Eq3ins1})}. If the $\sim
$-class of $0$ in $\mathbb{Z}^\nu $ is a rank-$\nu $ lattice
$\mathbb{L\kern0.5pt}$, then $\mathbf{T}$ tiles $\mathbb{R}^\nu $ by
$\mathbb{L\kern0.5pt}$.
\end{lemma}
\begin{proof}
We use an argument from Bandt \cite{Ban91}, specifically (i) in his
proof (p. 554) of his Theorem 1. This argument gets us a $c\in
\mathbb{R}_+$ and a finite subset $F\subset \mathbb{L\kern0.5pt}$ such
that $\mathbf{U}_c=:\left\{y\in \mathbb{R}^\nu \mid \left\| y\right\|
\leq c\right\} \subset \mathbf{T}+F$ ($=\bigcup_{\lambda \in
F}\mathbf{T}+\lambda $). If $x\in \mathbb{R}^\nu $, the condition on
$\mathbf{N}$ then gets us $k\in \mathbb{N}$ such that
$\mathbf{N}^{-k}x\in \mathbf{U}_c$, and therefore there are $t_0\in
\mathbf{T}$ and $\lambda _0\in F$ ($\subset \mathbb{L\kern0.5pt}$)
such that $\mathbf{N}^{-k}x=t_0+\lambda _0$. Using (\ref{Eq3ins1}),
there are indices $i_1,\dots ,i_k$ and $t_k\in \mathbf{T}$ such that
$t_0=\tau _{i_k}\left( \tau _{i_{k-1}}\left( \cdots \left( \tau
_{i_1}\left( t_k\right) \right) \cdots \right) \right) $. Hence
\begin{equation*}
x=\mathbf{N}^kt_0+\mathbf{N}^k\lambda _0=d_{i_1}+\mathbf{N}d_{i_2}+\dots 
+\mathbf{N}^{k-1}d_{i_k}+t_k+\mathbf{N}^k\lambda _0\text{.}
\end{equation*}
But
\begin{equation*}
d_{i_1}+\mathbf{N}d_{i_2}+\dots +\mathbf{N}^{k-1}d_{i_k}\in 0^\sim 
=\mathbb{L\kern0.5pt}\text{,}
\end{equation*}
and $\mathbf{N}^k\lambda _0\in \mathbb{L\kern0.5pt}$. By assumption,
$\mathbb{L\kern0.5pt}+\mathbb{L\kern0.5pt}\subset
\mathbb{L\kern0.5pt}$, so we conclude that $x\in
\mathbf{T}+\mathbb{L\kern0.5pt}$. Also, from the open-set property in
\cite{Ban91}, the intersections of the distinct parts $\tau
_{i_k}\circ \dots \circ \tau _{i_1}\left( \mathbf{T}\right) $ in
$\mathbf{T}$ have empty interior, so it follows, by the same scaling
argument, that the interior of $\left( \mathbf{T}+\lambda \right) \cap
\mathbf{T}$, for $\lambda \in \mathbb{L\kern0.5pt}\diagdown \left\{
0\right\}$, is empty. Here we used again that the identity
$\mathbb{L\kern0.5pt}=0^\sim $ is assumed. This concludes our proof
that $\mathbf{T}$ tiles $\mathbb{R}^\nu $ by the lattice
$\mathbb{L\kern0.5pt}$, under the stated assumptions.
\end{proof}

The application to the Cloud examples (Example \ref{Clouds}) is 
now based on the following
\renewcommand{\theobservation}{\unskip}
\begin{observation}
Consider the Cloud examples $\mathbf{T}_i$, $i=1,3,5,9$, and let $\sim
_i$ denote the corresponding equivalence relation on
$\mathbb{Z}^2$. Then the lattice condition from Lemma \ref{TilesByL}
is satisfied for $i=1,3,5$, but not for $i=9$. In fact
\begin{align*}
\begin{pmatrix}
0 \\
0
\end{pmatrix}%
^{\sim _1} &=\mathbb{Z}^2\text{,} \\
\begin{pmatrix}
0 \\
0
\end{pmatrix}%
^{\sim _3} &=\left\{ 
\begin{pmatrix}
2m \\
n
\end{pmatrix}
\biggm| m,n\in \mathbb{Z}\right\} =\mathbb{L\kern0.5pt}_3\text{,} \\
\begin{pmatrix}
0 \\
0
\end{pmatrix}%
^{\sim _5} &=\left\{ 
\begin{pmatrix}
m \\
2n
\end{pmatrix}
\biggm| m,n\in \mathbb{Z}\right\} =\mathbb{L\kern0.5pt}_5\text{,} \\
\intertext{while}
\begin{pmatrix}
0 \\
0
\end{pmatrix}%
^{\sim _9} &\text{ is a proper subset of }\mathbb{L\kern0.5pt}_9=\left\{ 
\begin{pmatrix}
m \\
2n
\end{pmatrix}
\biggm| m,n\in \mathbb{Z}\right\} \text{ that is not a lattice.}
\end{align*}
\end{observation}
\begin{proof}
The proof of each assertion is based on Corollary \ref{Cor3.3} above,
but in the $\mathbb{Z}^2$-variant, and an inductive argument. A more
detailed analysis of $\left(
\begin{smallmatrix}
0 \\
0
\end{smallmatrix}
\right) ^{\sim _9}$ may also be gained from an extension of the
flow-diagram in Figure \ref{CloudNineCycle} below. See also
Proposition \ref{CloudNineTiles}.
\end{proof}

Our present approach to tilings is somewhat different from that of
\cite{Rud88}--\cite{Rud89}. Our starting point is the permutative
representation, and that in turn encodes the action of arbitrary
symbols on basis vectors in some Hilbert space. In contrast, the
viewpoint in the two mentioned papers by Rudolph is to understand
$\mathbb{R}^n$-actions from tilings: in \cite{Rud88}, the orbits of an
$\mathbb{R}^n$-action are measurably tiled by rectangles of $2^n$
different shapes; and in \cite{Rud89} it is shown that, in some cases,
the tiles can be canonically labeled by $\mathbb{Z}^n$.
\medskip

\subsection{Tiles in the plane}\label{TilePlane}

\ \bigskip

The following two examples, for $\nu =2,$ also illustrate
the tiling property:

\begin{figure}
\begin{minipage}{126bp}
\psfig{file=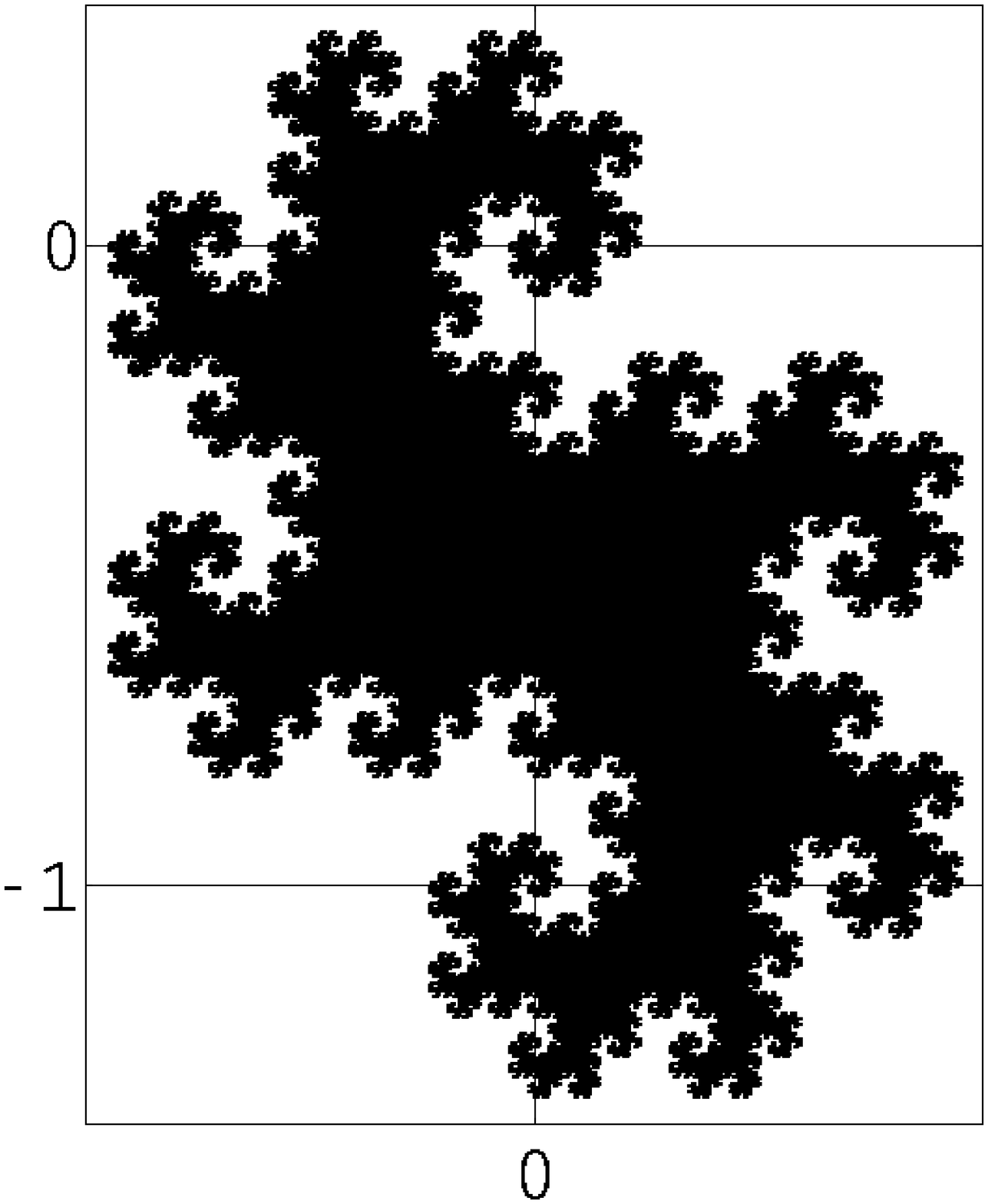,width=126bp}
\makebox[1.75in]{$\mathbf{T}$}
\end{minipage}\hfill%
\begin{minipage}{176bp}
\psfig{file=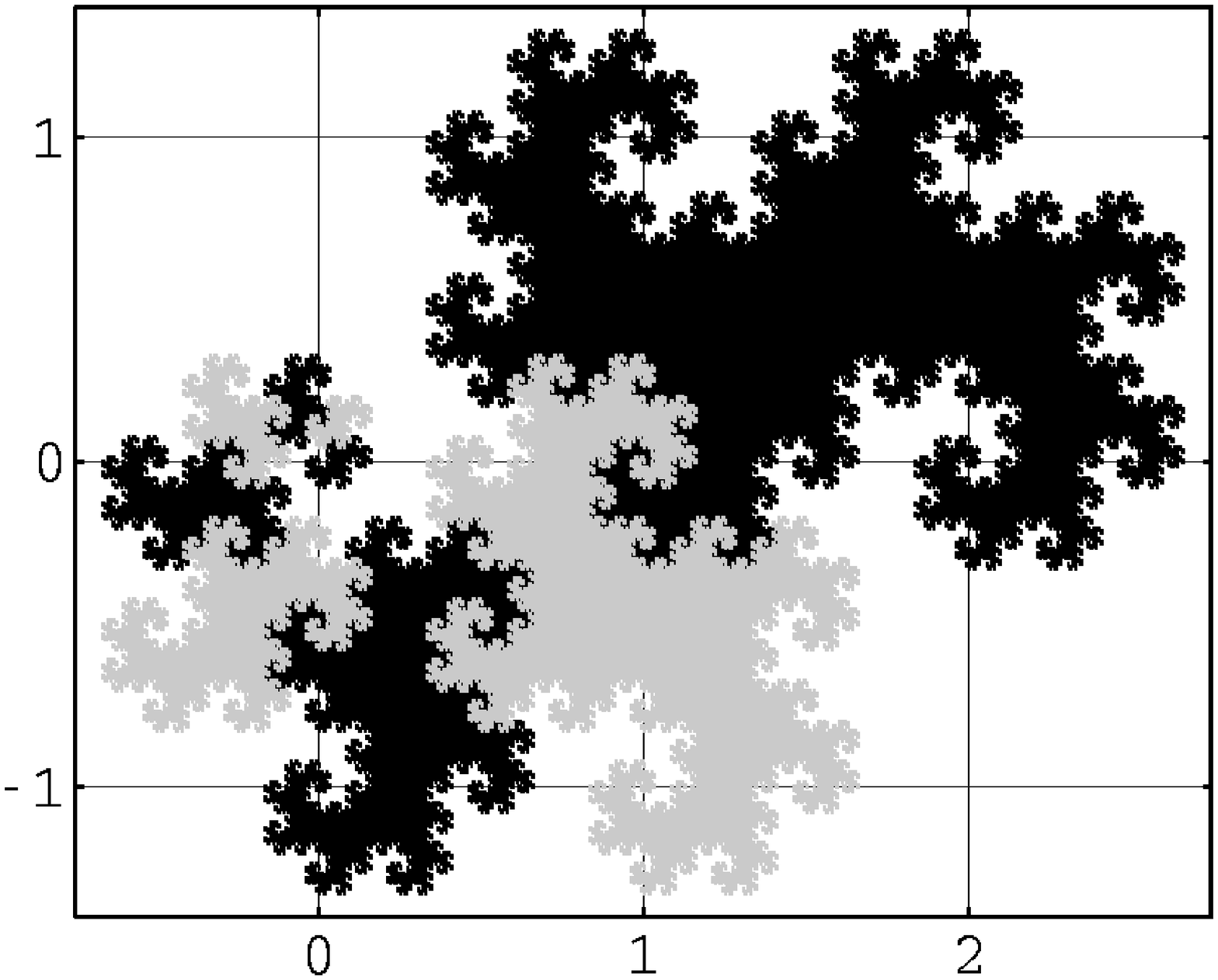,width=176bp}
\makebox[176bp]{$\mathbf{N}^2\left( \mathbf{T}\right) $ showing 
action of $\mathbf{N}^{-1}$}
\end{minipage}
\raisebox{-9pt}[12pt][0pt]{\makebox[\textwidth]{%
$\mathbf{N}=\left(
\begin{smallmatrix}
1 & -1 \\
1 & 1
\end{smallmatrix}
\right) $, $D=\left\{ \left(
\begin{smallmatrix}
0 \\
0
\end{smallmatrix}
\right), \left(
\begin{smallmatrix}
1 \\
0
\end{smallmatrix}
\right) \right\} $%
}}
\caption{Twin-Dragon (Example \protect \ref{Exa8ins1})}
\label{Twin-Dragon}
\end{figure}

\begin{example}
\label{Exa8ins1}Let $\mathbf{N}=\left( 
\begin{smallmatrix}
1 & -1 \\
1 & 1
\end{smallmatrix}\right) $, $D=\left\{ \left( 
\begin{smallmatrix}
0 \\
0
\end{smallmatrix}\right) ,\left( 
\begin{smallmatrix}
1 \\
0
\end{smallmatrix}\right) \right\} $. Then $\mathbf{N}^4=-4\openone _2=\left( 
\begin{smallmatrix}
-4 & 0 \\
0 & -4
\end{smallmatrix}\right) $, and $\mathbf{N}^8=2^4\openone _2=\left( 
\begin{smallmatrix}
16 & 0 \\
0 & 16
\end{smallmatrix}\right) $; and $T\left( \mathbf{N},D\right) $ is the
twin-dragon set (\cite{Ban91}, \cite{Ban96}, \cite{DDL95})
which has fractal boundary. (See Figure \ref
{Twin-Dragon} for a sketch.) Property (\ref{Rem8ins1(ii)}) is
illustrated in this example where $\left( 
\begin{smallmatrix}
1 \\
0
\end{smallmatrix}\right) +\mathbf{N}\left( 
\begin{smallmatrix}
0 \\
1
\end{smallmatrix}\right) =\left( 
\begin{smallmatrix}
0 \\
1
\end{smallmatrix}\right) $ and
\begin{equation*}
D^{\prime }=D+\mathbf{N}D+\mathbf{N}^2D+\mathbf{N}^3D=\left\{ 
\begin{pmatrix}
i_0+i_1-2i_3 \\ 
i_1+2i_2+2i_3
\end{pmatrix}
\biggm|i_0,i_1,i_2,i_3\in \left\{ 0,1\right\} \right\}\text{.}
\end{equation*}
So we have two order-one Cuntz states
given by the fixed points $\left(
\begin{smallmatrix}
0 \\
0
\end{smallmatrix}
\right) $ and $\left(
\begin{smallmatrix}
0 \\
1
\end{smallmatrix}
\right) $,
and using the method of Scholium \ref{Sch3.8}, one can
show that these are all the atoms of this representation.
These points can also be read off Figure \ref
{Twin-Dragon}.

A result of Knuth \cite[p. 190, Fig. 1]{Knu81} may be translated into
the observation that the related $\left( \mathbf{N},D\right) $-system
with $\mathbf{N}=\left(
\begin{smallmatrix}
-1 & -1 \\
1 & -1
\end{smallmatrix}
\right) $ and $D=\left\{ \left( 
\begin{smallmatrix}
0 \\
0
\end{smallmatrix}
\right) ,\left( 
\begin{smallmatrix}
1 \\
0
\end{smallmatrix}
\right) \right\} $ has just a single Cuntz state, corresponding to $\left( 
\begin{smallmatrix}
0 \\
0
\end{smallmatrix}
\right) ^\sim =\mathbb{Z}^2$. The effect on our present $\mathbf{T}$
in Figure \ref{Twin-Dragon} is an affine transformation wiping out
$\left(
\begin{smallmatrix}
0 \\
1
\end{smallmatrix}
\right) $ in $B_\infty $, i.e., this point won't be in the rotated
$-\mathbf{T}$ from Knuth's example.
\end{example}

\begin{figure}
\begin{minipage}[b]{64bp}
\psfig{file=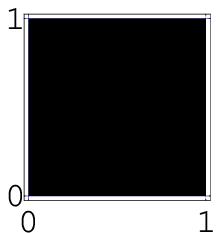}
\makebox[64bp]{$\mathbf{T}$}
\end{minipage}\hfill%
\begin{minipage}[b]{259bp}
\psfig{file=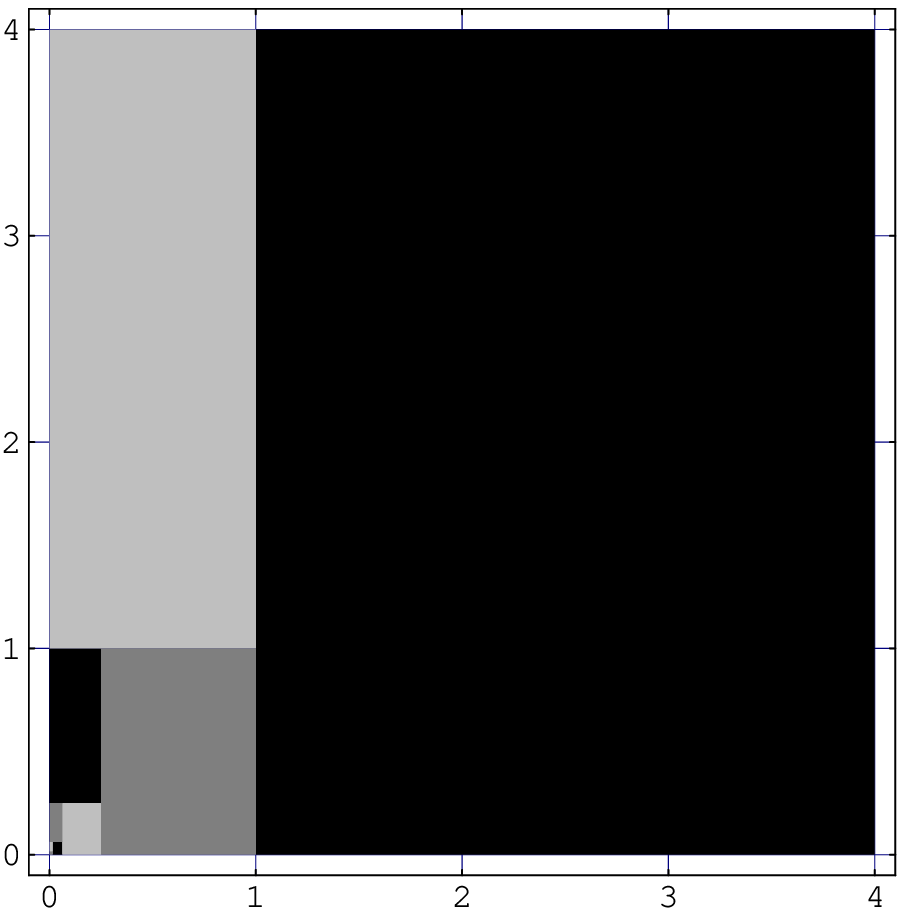}
\makebox[259bp]{$\mathbf{N}^2\left( \mathbf{T}\right) $ showing 
action of $\mathbf{N}^{-1}$}
\end{minipage}
\raisebox{-9pt}[12pt][0pt]{\makebox[\textwidth]{%
$\mathbf{N}=\left(
\begin{smallmatrix}
0 & 1 \\
4 & 0
\end{smallmatrix}
\right) $, $D=\left\{ \left(
\begin{smallmatrix}
0 \\
0
\end{smallmatrix}
\right), \left(
\begin{smallmatrix}
0 \\
1
\end{smallmatrix}
\right), \left(
\begin{smallmatrix}
0 \\
2
\end{smallmatrix}
\right), \left(
\begin{smallmatrix}
0 \\
3
\end{smallmatrix}
\right) \right\} $%
}}
\caption{Unit Square (Example \protect \ref{Exa8ins2})}
\label{UnitSquare}
\end{figure}

\begin{example}
\label{Exa8ins2}$\mathbf{N}=\left( 
\begin{smallmatrix}
0 & 1 \\
4 & 0
\end{smallmatrix}\right) $, $D=\left\{ \left( 
\begin{smallmatrix}
0 \\
k
\end{smallmatrix}\right) \mid k=0,1,2,3\right\} $, where $\mathbf{N}%
^2=2^2\openone _2=\left( \begin{smallmatrix}
4 & 0 \\
0 & 4
\end{smallmatrix}\right) $, and $D^{\prime }=D+\mathbf{N}D=\left\{ \left( 
\begin{smallmatrix}
i \\
j
\end{smallmatrix}\right) \mid i,j\in \left\{ 0,1,2,3\right\} \right\} $.
Here we have the non-trivial order-one Cuntz state corresponding to the
solution 
\begin{equation*}
\begin{pmatrix}
0 \\ 
3
\end{pmatrix}
+\mathbf{N}
\begin{pmatrix}
-1 \\ 
-1
\end{pmatrix}
=
\begin{pmatrix}
-1 \\ 
-1
\end{pmatrix}
\text{.}
\end{equation*}

The complete cycle-atom structure for this example is computed
prior to (\ref{EqNew8.29}) from which we see that $B_\infty =-\left\{ \left(
\begin{smallmatrix}
0 \\
0
\end{smallmatrix}
\right) ,\left(
\begin{smallmatrix}
1 \\
1
\end{smallmatrix}
\right) ,\left(
\begin{smallmatrix}
1 \\
0
\end{smallmatrix}
\right) ,\left(
\begin{smallmatrix}
0 \\
1
\end{smallmatrix}
\right) \right\} $.
This can also be read off the square in Figure \ref{UnitSquare}.

The iteration for $\mathbf{N}^2=\left( 
\begin{smallmatrix}
4 & 0 \\
0 & 4
\end{smallmatrix}\right) $ separates coordinates: 
\begin{equation*}
\begin{cases} i_0+4i_1+4^2i_2+4^3i_3+\cdots \\
j_0+4j_1+4^2j_2+4^3j_3+\cdots \text{,}
\end{cases}
\end{equation*}
and we get the following equation satisfied by the tile $\mathbf{T}=T\left( 
\mathbf{N},D\right) $ which is the unique compact solution to the equation: 
\begin{equation*}
4\mathbf{T}=\bigcup_{i,j}\left\{ \mathbf{T}+
\begin{pmatrix}
i \\ 
j
\end{pmatrix}
\biggm|i,j\in \left\{ 0,1,2,3\right\} \right\}
\end{equation*}
meaning that $\mathbf{T}$ is the unit square in this case, i.e., 
\begin{equation*}
\mathbf{T}=I\times I=\left\{ 
\begin{pmatrix}
x \\ 
y
\end{pmatrix}
\biggm|0\leq x\leq 1,0\leq y\leq 1\right\} \text{.}
\end{equation*}
So the boundary is not a fractal in Example \ref{Exa8ins2}.
\end{example}

Since, in Example \ref{Exa8ins1}, 
\begin{equation*}
\mathbf{N}\left( \mathbb{Z}^2\right) =\left\{ 
\begin{pmatrix}
n_1 \\ 
n_2
\end{pmatrix}
\in \mathbb{Z}^2\biggm|n_1+n_2\in 2\mathbb{Z}\text{, }n_1-n_2\in 2%
\mathbb{Z}\right\} \text{,}
\end{equation*}
and, in Example \ref{Exa8ins2}, 
\begin{equation*}
\mathbf{N}\left( \mathbb{Z}^2\right) =\left\{ 
\begin{pmatrix}
n_1 \\ 
n_2
\end{pmatrix}
\in \mathbb{Z}^2\biggm|n_2\in 4\mathbb{Z}\right\} \text{,}
\end{equation*}
the candidates we gave for the respective residue sets $D\subset \mathbb{Z}%
^2 $ are easily seen to have the desired property: $d-d^{\prime }\notin 
\mathbf{N}\left( \mathbb{Z}^2\right) $ whenever $d,d^{\prime }\in D$ and $%
d\neq d^{\prime }$. But, more importantly, we can use our method from above
in one dimension to find all the cycles and atoms for the other more
complicated choices of full residue sets $D$.

In particular, Remark \ref{Rem8ins1}(\ref{Rem8ins1(iii)}) above shows how
the calculations for $\nu =1$, following Proposition \ref{ProNew8.3}, may be
used in the analysis of cycle--atom structure for Examples \ref{Exa8ins1}
and \ref{Exa8ins2} above, and for related examples $\mathrm{s.t.\;}\mathbf{N}%
^m=p\openone_\nu $ holds for some $m$, $p$.
\medskip

\subsection{Single-atom cases}\label{SingleAtom}

\ \bigskip

In Example \ref{JordanMatrix}, the matrix $\mathbf{N}=\left( 
\begin{smallmatrix}
2 & 1 \\
0 & 2
\end{smallmatrix}
\right) $ is special in that $\left( \mathbf{N}-\openone \right)
^{-1}\mathbf{N}$ is also an integral matrix. An application of
Proposition \ref{Admissible}, below, then yields the result that, for
any other $D^{\prime } \subset \mathbb{Z}^2$ consisting of $4$
mutually $\mathbf{N}$-incongruent points, the corresponding $B^{\prime
} _\infty $ will have no fewer than $4$ points, so there are at least
$4$ equivalence classes for the corresponding $\approx ^{\prime } $
relation. The same argument, applied to Example \ref{Exa8ins1}, shows
that the corresponding number is $2$ there. Example \ref{Exa8ins2}
provides an illustration of a case where the points in $\mathbf{T}\cap
\mathbb{Z}^2$ are obtained from the second iteration of (\ref{Eq5.10})
in Remark \ref{SubCuntzMatrix}, needing both the cases $k=1$ and
$k=2$.

The assertions made here about the lower bound on the set of periodic
points in the respective examples follow from the following:

\begin{proposition}
\label{Admissible}Assume 
that $\mathbf{N}$ is a $\nu \times \nu $ matrix with integer
entries and $N=\left| \det \mathbf{N}\right| >0$. Assume that all the
\textup{(}complex\/\textup{)} eigenvalues of $\mathbf{N}$ have modulus
greater than one. If $D$ is a set of $N$ points in $\mathbb{Z}^\nu $
which are mutually incongruent modulo $\mathbf{N}\mathbb{Z}^\nu $, and
\begin{equation}
D_0=\left\{ d\in D\mid \smash{\left( \openone -\mathbf{N}\right)
^{-1}}d\in \mathbb{Z}^\nu \right\}
\label{Incl}
\end{equation}
then the map $R$ associated with the admissible pair $\left(
\mathbf{N},D\right) $ by Scholium \ref{Sch2.4} \textup{(}or
\textup{(\ref{Trans}))} has exactly $\#\left( D_0\right) $ fixed
points, so
\begin{equation*}
\#\left( B_\infty \right) \geq \#\left( D_0\right) \text{.}
\end{equation*}
In particular, if $\left( \openone -\mathbf{N}\right) ^{-1}$ has 
integral entries, then
\begin{equation}
\#\left( B_\infty \right) \geq N=\#\left( D\right) \text{.}
\label{Est}
\end{equation}
\end{proposition}

\TeXButton{Begin Proof}{\begin{proof}}
If $d\in D$ and $x=\left( \openone -\mathbf{N}\right) ^{-1}d\in
\mathbb{Z}^\nu $, then $R\left( x\right) =x$ and $x\in B_\infty $. As
$\left( \openone -\mathbf{N}\right) ^{-1}$ is invertible, distinct
$d$'s give rise to distinct $x$'s in this way.  
\TeXButton{End Proof}{\end{proof}}

This raises the question of finding geometric properties of admissible
pairs $\left( \mathbf{N},D\right) $ in $\nu $ dimensions $\nu \geq 2$;
for example, when $\mathbf{N}$ is given, finding conditions on
$D\subset \mathbb{Z}^\nu $ such that $\left( \mathbf{N},D\right) $ is
admissible, and the corresponding $B_\infty =\left( -\mathbf{T}\right)
\cap \mathbb{Z}^\nu $ is a singleton. This means that the associated
representation of $\limfunc{UHF}\nolimits_N$ is irreducible on
$L^2\left( \mathbb{T}^\nu \right) $.

\begin{figure}
\begin{minipage}{1.5in}
\psfig{file=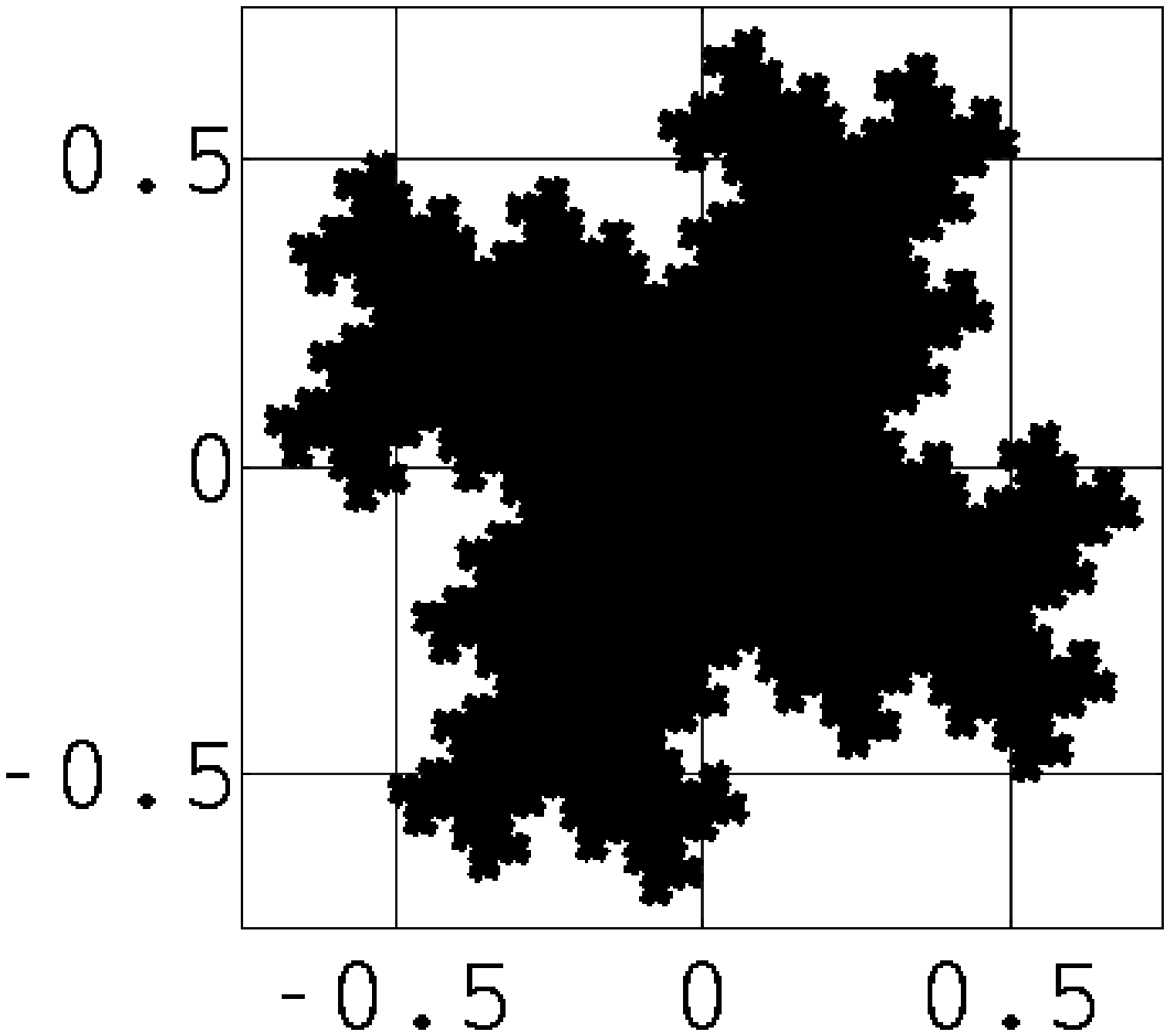,width=108bp}
\makebox[1.5in]{$\mathbf{T}$}
\end{minipage}\hfill%
\begin{minipage}{3in}
\psfig{file=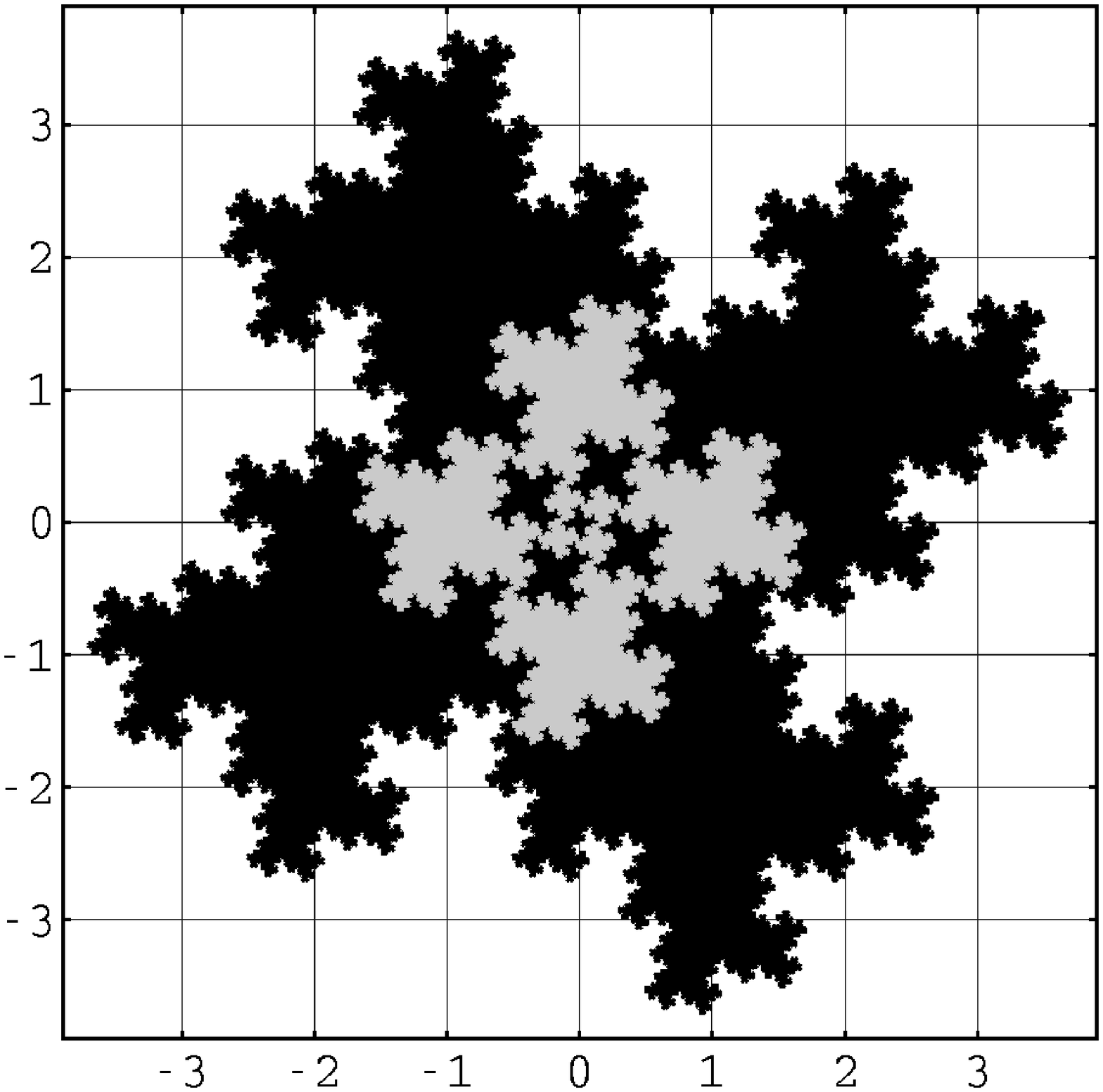,width=216bp}
\makebox[3in]{$\mathbf{N}^2\left( \mathbf{T}\right) $ showing action of
$\mathbf{N}^{-1}$}
\end{minipage}
\raisebox{-9pt}[12pt][0pt]{\makebox[\textwidth]{%
$\mathbf{N}=\left(
\begin{smallmatrix}
1 & 2 \\
-2 & 1
\end{smallmatrix}
\right) $, $D=\left\{ \left(
\begin{smallmatrix}
0 \\
0
\end{smallmatrix}
\right), \left(
\begin{smallmatrix}
\pm 1 \\
0
\end{smallmatrix}
\right), \left(
\begin{smallmatrix}
0 \\
\pm 1
\end{smallmatrix}
\right) \right\} $%
}}
\caption{Fractal Red Cross (Example \protect \ref{FractalRedCross})}
\label{RedCrossFigure}
\end{figure}

\begin{example}
\label{FractalRedCross}We note that Proposition \ref{Admissible}
does not apply to $\mathbf{N}=\left( 
\begin{smallmatrix}
1 & 2 \\
-2 & 1
\end{smallmatrix}
\right) $, as $\left( \mathbf{N}-\openone \right) ^{-1}\mathbf{N}=\left( 
\begin{smallmatrix}
1 & -\frac 12 \\
\frac 12 & 1
\end{smallmatrix}
\right) $ does not have integer entries.
Using estimates from Section \ref{Monomial}, we now show that the answer is
``yes'' for $\mathbf{N}=\left( 
\begin{smallmatrix}
1 & 2 \\
-2 & 1
\end{smallmatrix}
\right) $, i.e., it admits $D$ with $B_\infty $ a singleton. There are
also examples in one dimension if $N>2$, and in higher dimensions with
$B_\infty $ a singleton. Take for example $N=3$ and $D=\left\{
2k-1,2k,2k+1\right\} $ with $k\in \mathbb{Z}$; then it follows from
Corollary \ref{Cor3.10} that $B_\infty =\left\{ -k\right\} $. When
$N=2$ it follows from Subsection \ref{SubsecNew8.1} that there are no
examples with $B_\infty $ a one-point set. There are simple examples
in two dimensions when we get non-$\mathbb{Z}^2$-lattice points in
$\mathbf{T}$ upon applying $\left( \mathbf{N}-\openone \right) ^{-1}$
to the vectors in $D$. One such is the ``Four-Clover'', or the
``Fractal Red Cross''. This version of $\mathbf{T}$ arises from the
following admissible pair $\left( \mathbf{N},D\right) $:
\begin{equation*}
\mathbf{N}=
\begin{pmatrix}
1 & 2 \\
-2 & 1
\end{pmatrix}
\text{,\quad } D=\left\{ 
\begin{pmatrix}
0 \\
0
\end{pmatrix}
,
\begin{pmatrix}
\pm 1 \\
0
\end{pmatrix}
,
\begin{pmatrix}
0 \\
\pm 1
\end{pmatrix}
\right\} \text{.}
\end{equation*}
Here $\left( \mathbf{N}-\openone \right) ^{-1}=\left( 
\begin{smallmatrix}
0 & -\frac 12 \\
\frac 12 & 0
\end{smallmatrix}
\right) $, i.e., non-integral, and 
\begin{equation*}
\left( \mathbf{N}-\openone \right) ^{-1}D=\left\{ 
\begin{pmatrix}
0 \\
0
\end{pmatrix}
,
\begin{pmatrix}
0 \\
\pm \frac 12
\end{pmatrix}
,
\begin{pmatrix}
\pm \frac 12 \\
0
\end{pmatrix}
\right\} \text{.}
\end{equation*}
In fact, it can be checked that 
\begin{equation*}
\mathbf{T}\cap \mathbb{Z}^2=\left\{ 
\begin{pmatrix}
0 \\
0
\end{pmatrix}
\right\} \text{,}
\end{equation*}
meaning that here we get just a single atom. In the terminology of
Powers \cite{Pow88}, this means (as we recalled in Section \ref{Sec1})
that the corresponding endomorphism of $\mathcal{B}\left( L^2\left(
\mathbb{T}^2\right) \right) $ is a (Powers) shift.

To verify that $\mathbf{T}\cap \mathbb{Z}^2=\left\{ \left( 
\begin{smallmatrix}
0 \\
0
\end{smallmatrix}
\right) \right\} $, use Corollary \ref{Cor3.9}, formula (\ref{Eq5.9})
in Remark \ref{SubCuntzMatrix}, and the observation that $\mathbf{N}$
here is diagonalizable with eigenvalues $\lambda _\pm =1\pm
i2$. Specifically, since $\left| \lambda _+\right| =\left|
\lambda_-\right| =\sqrt 5$, the argument from the proof of Corollary
\ref{Cor3.10} also shows that every point $x=\left(
\begin{smallmatrix}
i \\
j
\end{smallmatrix}
\right) \in B_\infty =\left( -\mathbf{T}\right) \cap \mathbb{Z}^2$
must satisfy $\left\| x\right\| \leq \frac 1{\sqrt 5-1}<1$ where
$\left\| x\right\| =\sqrt {i^2+j^2}$. This norm is appropriate since
$\mathbf{N}=\openone + \left(
\begin{smallmatrix}
0 & 2 \\
-2 & 0
\end{smallmatrix}
\right) $ and $\left(
\begin{smallmatrix}
0 & 2 \\
-2 & 0
\end{smallmatrix}
\right) $ is skew-adjoint, so the eigenvectors of $\mathbf{N}$ in
$\mathbb{C}^2$ are orthogonal. 
It follows that $B_\infty =\left\{ \left( 
\begin{smallmatrix}
0 \\
0
\end{smallmatrix}
\right) \right\} $ as claimed. Note further that the discussion here
falls within the class of ultraspecial examples in Remark
\ref{MeanDensity} where the estimate (\ref{Eq3.11}) is in fact an
equality, since $\left\| \mathbf{N}^{-1}\right\| =5^{-\frac 12}$. (The
same argument, for Example \ref{Exa8ins1} (the Twin Dragon), shows
that points $x$ in $\mathbf{T}$ there must satisfy $\left\| x\right\|
\leq \frac 1{\sqrt 2-1}$.)
\end{example}

A more refined construction of single-atom examples (analogous to
$D_1$) in higher even dimensions may be based on the known Hadamard
matrices from combinatorics; see, e.g., \cite{JoPe96} and
\cite{Haa95}. If $\nu $ is even, a $\nu \times \nu $ matrix
$\mathbf{N}$ with entries $\pm 1$ is said to be Hadamard if
$\mathbf{N}^*\mathbf{N}=\nu \openone _\nu $. Then of course, the
spectrum of $\mathbf{N}$ must be in $\left\{ \lambda \in
\mathbb{C}\mid \left| \lambda \right| =\sqrt \nu \right\} $ and
$\left\| \mathbf{N}^{-1}\right\| =\nu ^{-\frac 12}$. For the
determinant, we have
\begin{equation*}
N_\nu =\left| \det \mathbf{N}_\nu \right| =\left( \sqrt \nu \right) 
^\nu \text{.}
\end{equation*}
We shall need $\nu =8$ below, in which case $N_8=8^4=2^{12}$, and
$\left\| \mathbf{N}_8^{-1}\right\| =\frac 1{2\sqrt 2}$. The case $\nu
=2$, and $\mathbf{N}_2=\left(
\begin{smallmatrix}
1 & 1 \\
1 & -1
\end{smallmatrix}
\right) $, is essentially our present Example \ref{Exa8ins1}, where
the cardinality of $B_\infty $ is $2$. For $\nu =8$, we can take
$\mathbf{N}=\mathbf{N}_2\otimes \mathbf{N}_2\otimes \mathbf{N}_2$. Now
choose $D=\left\{ d_i\in \mathbb{Z}^8\mid i=1,2,\dots ,2^{12}\right\}
$ such that the pair $\left( \mathbf{N},D\right) $ is admissible in
$\mathbb{Z}^\nu $, and
\begin{equation}
\max _i \left\| d_i\right\| <2\sqrt 2-1\text{.} 
\label{HadamardBound}
\end{equation}
Then the argument from before shows that every point $x\in
\mathbf{T}=-\mathbf{T}$ must satisfy $\left\| x\right\| <1$; and
therefore, by Proposition \ref{Pro3ins1}, $B_\infty =\left\{ 0\right\} $
for those examples.

It is known (see, e.g., \cite{Haa95}) that there are other
higher-order Hadamard matrices than those which arise from the simple
tensor construction.

We also note that the tensor construction works for matrices more
general than the Hadamard matrices: if $\mathbf{N}_i$ are given
integral matrices of size $\nu _i$, $i=1,2$, and if there are $p_i$'s
$>1$, such that $\mathbf{N}^*_i\mathbf{N}_i=p_i\openone _{\nu _i}$,
then set $\mathbf{M}=\mathbf{N}_1\otimes \mathbf{N}_2$, and we will
have $\mathbf{M}^*\mathbf{M}$ represented as a diagonal matrix which
yields us the estimate $\left\| \mathbf{M}^{-1}\right\| \leq \frac
1{\min \left( p_1,p_2\right) }<1$. Therefore, if $D\subset
\mathbb{Z}^{\nu _1\nu _2}$ is chosen such that $\left(
\mathbf{N},D\right) $ is admissible, and also assuming
\begin{equation}
\max _j\left\| d_j\right\| <\min \left( p_1,p_2\right) -1\text{,}
\label{TensorDBound}
\end{equation}
then we will have an example $\left( \mathbf{M},D\right) $ in $\nu
_1\nu _2$ dimensions with $B_\infty =\left\{ 0\right\} $. But the
selection of the set $D$ of residues in $\nu _1\nu _2$ dimensions is
easier if the individual matrices $\mathbf{N}_i$, $i=1,2$, are
Hadamard. In the more general cases, the explicit construction of
admissible residue sets in $\nu _1\nu _2$ dimensions may be
difficult. In particular, the conditions (\ref{HadamardBound}) and
(\ref{TensorDBound}) show that a more detailed classification of the
single-atom cases must depend also on the asymptotics of the counting
function
\begin{equation}
N_\nu \left( R,a\right) =\#\left\{ x\in \mathbb{Z}^\nu \mid \left| 
x-a\right| \leq R\right\}
\label{CountingFunction}
\end{equation}
where $\left| x-a\right| =\left( \sum _i\left( x_i-a_i\right)
^2\right) ^\frac 12$. Even for $\nu =2$, this is a problem of current
interest; see, e.g., \cite{BlDy94}.
\medskip

\subsection{Multiple coverings}\label{Multiple}

\ \bigskip

\begin{figure}[p]
\raisebox{-9pt}[0pt][0pt]{\makebox[\textwidth][l]{%
{\itshape Sometimes we see a cloud that's dragonish...}
(Shakespeare, {\itshape Antony and Cleopatra})%
}}
\makebox[\textwidth][r]{\psfig{file=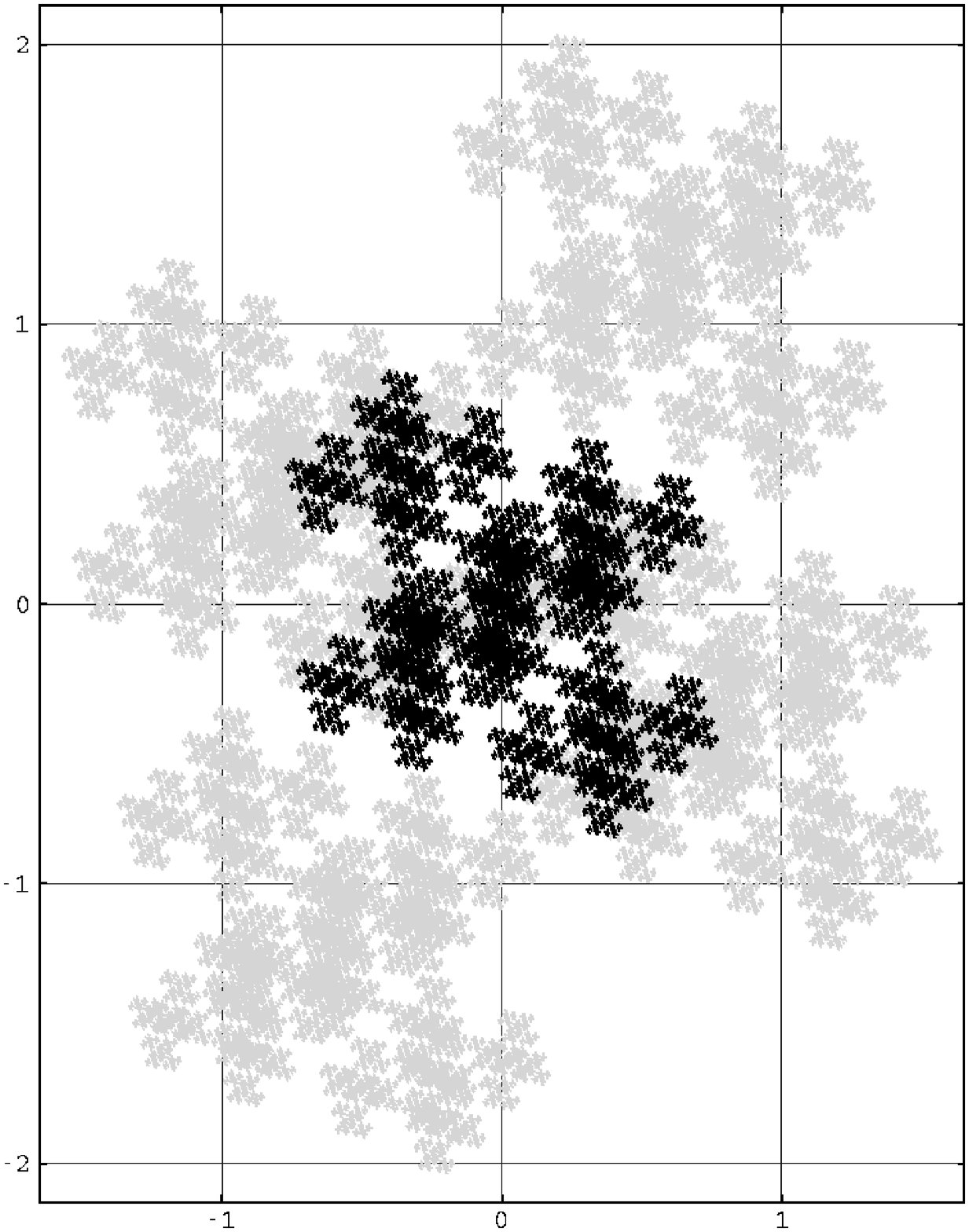,width=360bp}}
\raisebox{15pt}[0pt]{\parbox[t]{\textwidth}{Shading added to illustrate how
$\mathbf{T}$ is formed as a ``reptile'', that is, the union of the five
subtiles, the central one $\tau _{\left( 
\begin{smallmatrix}
0 \\
0
\end{smallmatrix}
\right) }\left( \mathbf{T}\right) =\mathbf{N}^{-1}\mathbf{T}$ (dark),
and the two pairs, $\tau _{\left( 
\begin{smallmatrix}
\pm 3 \\
0
\end{smallmatrix}
\right) }\left( \mathbf{T}\right) =\pm \tau _{\left( 
\begin{smallmatrix}
3 \\
0
\end{smallmatrix}
\right) }\left( \mathbf{T}\right) $ and $\tau _{\left( 
\begin{smallmatrix}
0 \\
\pm 2
\end{smallmatrix}
\right) }\left( \mathbf{T}\right) =\pm \tau _{\left( 
\begin{smallmatrix}
0 \\
2
\end{smallmatrix}
\right) }\left( \mathbf{T}\right) $ in symmetry (light), with respective
interiors of distinct subtiles non-overlapping, i.e., ``just-touching''
subtiles.}}
\raisebox{-9pt}[9pt][0pt]{\makebox[\textwidth]{%
$\mathbf{N}=\left(
\begin{smallmatrix}
1 & 2 \\
-2 & 1
\end{smallmatrix}
\right) $, $D=\left\{ \left(
\begin{smallmatrix}
0 \\
0
\end{smallmatrix}
\right), \left(
\begin{smallmatrix}
\pm 3 \\
0
\end{smallmatrix}
\right), \left(
\begin{smallmatrix}
0 \\
\pm 2
\end{smallmatrix}
\right) \right\} $%
}}
\caption{Cloud Nine (Example \protect \ref{CloudNine})}
\label{CloudNineFigure}
\end{figure}

The single-atom case is relevant to, among other things, the
complex-base problem \cite{Gil96}. See Remark \ref{Rem3ins3}, where we
also noted that the single-atom case makes the corresponding
$\mathbf{T}$ tile $\mathbb{R}^\nu $ by the unit lattice
$\mathbb{Z}^\nu $. We now turn to the case where $B_\infty $ is not a
singleton, but where we still retain the lattice geometry for
$\mathbf{T}$. Then we will expect $\mathbf{T}$ to tile $\mathbb{R}^\nu $
with a proper sublattice $\mathbb{L\kern0.5pt}\subset \mathbb{Z}^\nu $,
and, hence, $\mathbf{T}$ will be a finite covering.  When
$\mathbf{N}=\left(
\begin{smallmatrix}
1 & 2 \\
-2 & 1
\end{smallmatrix}
\right) $, the associated fractal $\mathbf{T}=\mathbf{T}\left(
D\right) $ can easily be computed for the following candidates for
$D$:
\begin{equation}
\begin{aligned}
D_1&=\left\{
\begin{pmatrix}
0 \\
0
\end{pmatrix}
,
\begin{pmatrix}
\pm 1 \\
0
\end{pmatrix}
,
\begin{pmatrix}
0 \\
\pm 1
\end{pmatrix}
\right\} & &\text{\quad (Fractal Red Cross, Fig. \ref{RedCrossFigure})} \\
D_9&=\left\{
\begin{pmatrix}
0 \\
0
\end{pmatrix}
,
\begin{pmatrix}
\pm 3 \\
0
\end{pmatrix}
,
\begin{pmatrix}
0 \\
\pm 2
\end{pmatrix}
\right\} & &\text{\quad (Cloud Nine, Fig. \ref{CloudNineFigure})} \\
D_3&=\left\{
\begin{pmatrix}
0 \\
0
\end{pmatrix}
,
\begin{pmatrix}
0 \\
\pm 1
\end{pmatrix}
,
\begin{pmatrix}
0 \\
\pm 2
\end{pmatrix}
\right\} & &\text{\quad (Cloud Three, Fig. \ref{CloudThreeFigure})} \\
D_5&=\left\{
\begin{pmatrix}
0 \\
0
\end{pmatrix}
,
\begin{pmatrix}
\pm 3 \\
0
\end{pmatrix}
,
\begin{pmatrix}
\pm 1 \\
0
\end{pmatrix}
\right\} & &\text{\quad (Cloud Five, Fig. \ref{CloudFiveFigure})}
\end{aligned}
\label{CloudD}
\end{equation}

\begin{example}
\label{CloudNine}Choosing, for example, $D_9 =\left\{ 
\left( 
\begin{smallmatrix}
0 \\
0
\end{smallmatrix}
\right) 
,
\left( 
\begin{smallmatrix}
\pm 3 \\
0
\end{smallmatrix}
\right) 
,
\left( 
\begin{smallmatrix}
0 \\
\pm 2
\end{smallmatrix}
\right) \right\} $, gets us two fixed points $\left( 
\begin{smallmatrix}
\pm 1 \\
0
\end{smallmatrix}
\right) $ in $B^{(9)} _\infty $ in addition to the one point in
$B_\infty =\left\{ \left( 
\begin{smallmatrix}
0 \\
0
\end{smallmatrix}
\right) \right\} $ for the other choice of $D=\left\{ \left( 
\begin{smallmatrix}
0 \\
0
\end{smallmatrix}
\right) , \left( 
\begin{smallmatrix}
\pm 1 \\
0
\end{smallmatrix}
\right) , \left( 
\begin{smallmatrix}
0 \\
\pm 1
\end{smallmatrix}
\right) \right\} $. More specifically, there are nine points in the set
$B^{(9)}_\infty =\left( -\mathbf{T}_9\right) \cap \mathbb{Z}^2=\left\{ \left( 
\begin{smallmatrix}
0 \\
0
\end{smallmatrix}
\right) ,\left( 
\begin{smallmatrix}
\pm 1 \\
0
\end{smallmatrix}
\right) ,\left( 
\begin{smallmatrix}
0 \\
\pm 1
\end{smallmatrix}
\right) ,\left( 
\begin{smallmatrix}
\pm 1 \\
\pm 1
\end{smallmatrix}
\right) \right\} $ (see Figure \ref{CloudNineFigure}). We still have
the symmetry $\mathbf{T}_9=-\mathbf{T}_9$ when $\mathbf{T}_9$ is
determined from (\ref{Eq3ins1}) corresponding to $D_9$, and we have
displayed both $\mathbf{T}_1$ and $\mathbf{T}_9$ in Figures
\ref{RedCrossFigure} and \ref{CloudNineFigure}. We now prove the
assertion about $B^{(9)}_\infty $ ``the hard way'' (as opposed to
``merely'' reading the picture). First the argument from Example
\ref{JordanMatrix} shows that $R^{\prime n}\left( x\right) $ is
eventually contained in the euclidean sphere of radius
\begin{equation*}
\frac 3{\sqrt 5}\left( \vphantom{\sum \left( \frac 1{\sqrt 5}\right) 
^n}\smash{\sum _{n=0}^\infty \left( \frac 1{\sqrt 5}\right) ^n}\right) 
=\frac 3{\sqrt 5-1}=2.427\dots 
\end{equation*}
regardless of which starting point $x$ in $\mathbb{Z}^2$ is
chosen. Here we use the notation $R^{\prime }$ for the mapping given
by (\ref{Trans}) with $D_9 =\left\{ \left(
\begin{smallmatrix}
0 \\
0
\end{smallmatrix}
\right) 
,
\left( 
\begin{smallmatrix}
\pm 3 \\
0
\end{smallmatrix}
\right) 
,
\left( 
\begin{smallmatrix}
0 \\
\pm 2
\end{smallmatrix}
\right) \right\} $, and it is immediate that $R^{\prime }
\left( 
\begin{smallmatrix}
0 \\
0
\end{smallmatrix}
\right) =
\left( 
\begin{smallmatrix}
0 \\
0
\end{smallmatrix}
\right) $, $R^{\prime }
\left( 
\begin{smallmatrix}
\pm 1 \\
0
\end{smallmatrix}
\right) =
\left( 
\begin{smallmatrix}
\pm 1 \\
0
\end{smallmatrix}
\right) $; i.e., three points of period one. A direct calculation, 
as in Example \ref{JordanMatrix} above, also yields
\begin{equation*}
\begin{pmatrix}
0 \\ 
1
\end{pmatrix}
\overset{\raise1pt\hbox{$\scriptstyle R^{\prime }$}}{\rightarrow }
\begin{pmatrix}
-1 \\ 
-1
\end{pmatrix}
\overset{\raise1pt\hbox{$\scriptstyle R^{\prime }$}}{\rightarrow }
\begin{pmatrix}
1 \\ 
-1
\end{pmatrix}
\overset{\raise1pt\hbox{$\scriptstyle R^{\prime }$}}{\rightarrow }
\begin{pmatrix}
0 \\ 
-1
\end{pmatrix}
\overset{\raise1pt\hbox{$\scriptstyle R^{\prime }$}}{\rightarrow }
\begin{pmatrix}
1 \\ 
1
\end{pmatrix}
\overset{\raise1pt\hbox{$\scriptstyle R^{\prime }$}}{\rightarrow }
\begin{pmatrix}
-1 \\ 
1
\end{pmatrix}
\overset{\raise1pt\hbox{$\scriptstyle R^{\prime }$}}{\rightarrow }
\begin{pmatrix}
0 \\ 
1
\end{pmatrix}
\end{equation*}
(see also Figure \ref{CloudNineCycle}), and thus,
\begin{equation*}
\left( R^{\prime }\right) ^6
\begin{pmatrix}
0 \\ 
1
\end{pmatrix}
=
\begin{pmatrix}
0 \\ 
1
\end{pmatrix}
\text{.}
\end{equation*}
Hence the point $\left( 
\begin{smallmatrix}
0 \\
1
\end{smallmatrix}
\right) $ has (minimal) period equal to $6$.

\begin{figure}
\psfig{file=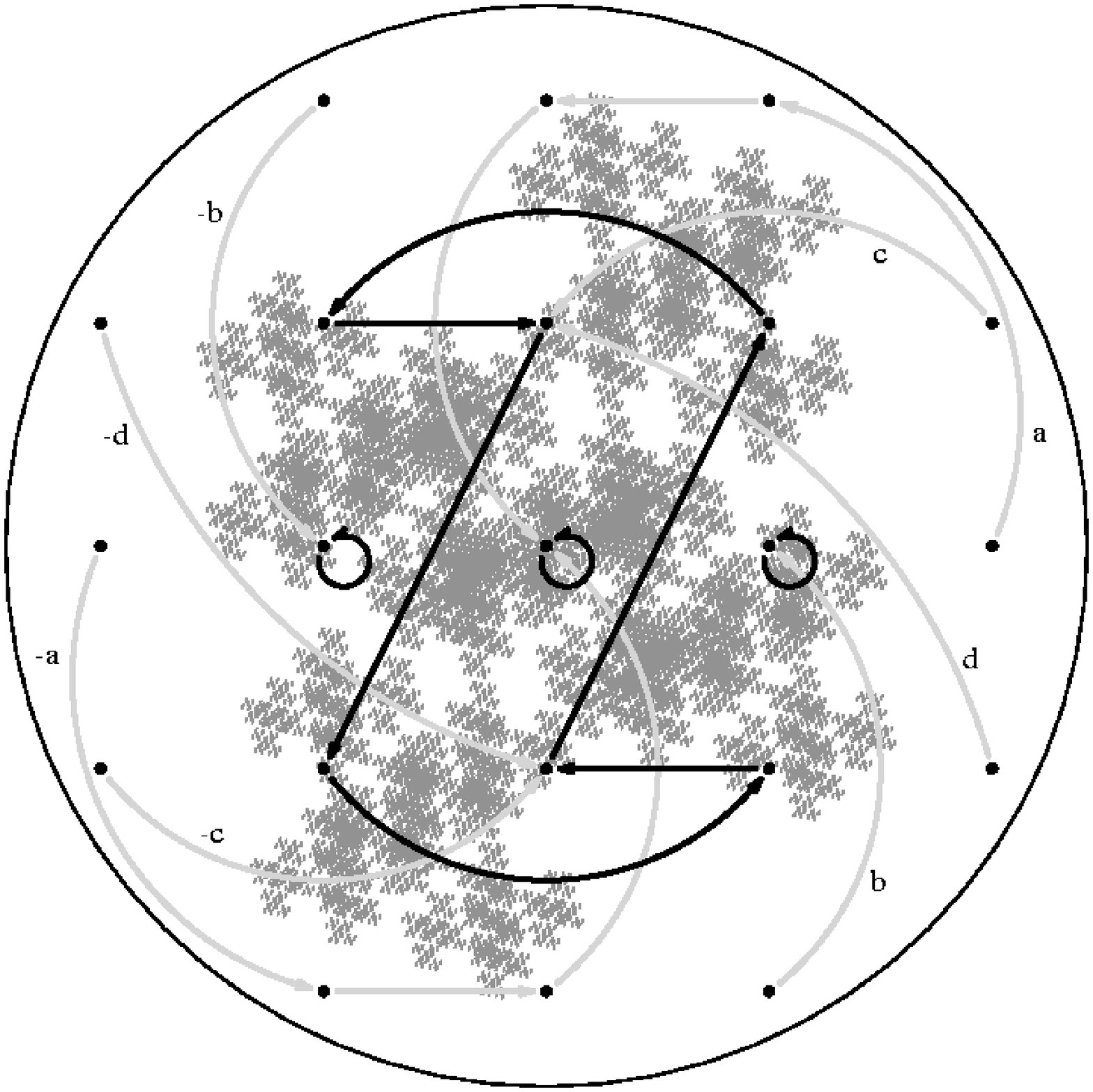,width=360bp}
\caption{Paths of points within spectral circle for Cloud Nine
(Example \protect \ref{CloudNine})}
\label{CloudNineCycle}
\end{figure}

Furthermore, we can check directly that the additional points within
the spectral circle have the following paths:
{\allowdisplaybreaks
\begin{gather}
\begin{pmatrix}
2 \\ 
0
\end{pmatrix}
\overset{\raise1pt\hbox{$\scriptstyle R^{\prime }$}}{\rightarrow }
\begin{pmatrix}
1 \\ 
2
\end{pmatrix}
\overset{\raise1pt\hbox{$\scriptstyle R^{\prime }$}}{\rightarrow }
\begin{pmatrix}
0 \\ 
2
\end{pmatrix}
\overset{\raise1pt\hbox{$\scriptstyle R^{\prime }$}}{\rightarrow }
\begin{pmatrix}
0 \\ 
0
\end{pmatrix}
\text{,} \tag{$\mathrm{a}$}\\
\begin{pmatrix}
-2 \\ 
0
\end{pmatrix}
\overset{\raise1pt\hbox{$\scriptstyle R^{\prime }$}}{\rightarrow }
\begin{pmatrix}
-1 \\ 
-2
\end{pmatrix}
\overset{\raise1pt\hbox{$\scriptstyle R^{\prime }$}}{\rightarrow }
\begin{pmatrix}
0 \\ 
-2
\end{pmatrix}
\overset{\raise1pt\hbox{$\scriptstyle R^{\prime }$}}{\rightarrow }
\begin{pmatrix}
0 \\ 
0
\end{pmatrix}
\text{,} \tag{$-\mathrm{a}$}\\
\begin{pmatrix}
1 \\ 
-2
\end{pmatrix}
\overset{\raise1pt\hbox{$\scriptstyle R^{\prime }$}}{\rightarrow }
\begin{pmatrix}
1 \\ 
0
\end{pmatrix}
\overset{\raise1pt\hbox{$\scriptstyle R^{\prime }$}}{\rightarrow }
\begin{pmatrix}
1 \\ 
0
\end{pmatrix}
\text{,} \tag{$\mathrm{b}$}\\
\begin{pmatrix}
-1 \\ 
2
\end{pmatrix}
\overset{\raise1pt\hbox{$\scriptstyle R^{\prime }$}}{\rightarrow }
\begin{pmatrix}
-1 \\ 
0
\end{pmatrix}
\overset{\raise1pt\hbox{$\scriptstyle R^{\prime }$}}{\rightarrow }
\begin{pmatrix}
-1 \\ 
0
\end{pmatrix}
\text{,} \tag{$-\mathrm{b}$}\\
\begin{pmatrix}
2 \\ 
1
\end{pmatrix}
\overset{\raise1pt\hbox{$\scriptstyle R^{\prime }$}}{\rightarrow }
\begin{pmatrix}
0 \\ 
1
\end{pmatrix}
\overset{\raise1pt\hbox{$\scriptstyle R^{\prime }$}}{\rightarrow }
6\text{-cycle,} \tag{$\mathrm{c}$}\\
\begin{pmatrix}
-2 \\ 
-1
\end{pmatrix}
\overset{\raise1pt\hbox{$\scriptstyle R^{\prime }$}}{\rightarrow }
\begin{pmatrix}
0 \\ 
-1
\end{pmatrix}
\overset{\raise1pt\hbox{$\scriptstyle R^{\prime }$}}{\rightarrow }
6\text{-cycle,} \tag{$-\mathrm{c}$}\\
\begin{pmatrix}
2 \\ 
-1
\end{pmatrix}
\overset{\raise1pt\hbox{$\scriptstyle R^{\prime }$}}{\rightarrow }
\begin{pmatrix}
0 \\ 
1
\end{pmatrix}
\overset{\raise1pt\hbox{$\scriptstyle R^{\prime }$}}{\rightarrow }
6\text{-cycle,} \tag{$\mathrm{d}$}\\
\begin{pmatrix}
-2 \\ 
1
\end{pmatrix}
\overset{\raise1pt\hbox{$\scriptstyle R^{\prime }$}}{\rightarrow }
\begin{pmatrix}
0 \\ 
-1
\end{pmatrix}
\overset{\raise1pt\hbox{$\scriptstyle R^{\prime }$}}{\rightarrow }
6\text{-cycle.} \tag{$-\mathrm{d}$}
\end{gather}%
}% 
The flow diagram illustrates the fact that some points in the
pattern are stable points of attraction in the sense of \cite{LiMa}.

Hence we have all $21$ points in the spectral circle covered: three
have period $1$, six period $6$, and the other twelve are not
periodic. (See Figure \ref{CloudNineCycle}.) Thus, of the $21$
$\mathbb{Z}^2$-points inside the spectral circle with radius $r=\frac
3{\sqrt 5-1}$, $9$ are in $B^{(9)}_\infty $: three single-atom cycles
on the $x$-axis, and one six-atom cycle (Z-shaped). All of the three
single-atom cycles have flow paths linking up to them and tracing a
portion of the points inside the circle (paths
$(\mathrm{a})$--$(\mathrm{b})$ in Figure \ref{CloudNineCycle}.) Two of
the six atoms in the $6$-cycle attract separate flow lines from
respective points inside the circle (paths
$(\mathrm{c})$--$(\mathrm{d})$ in Figure \ref{CloudNineCycle}), and an
inspection shows that all $21$ points are thus accounted for. We have
described all the interior flow paths, and we note the reflection
symmetry of paths under $\left(
\begin{smallmatrix}
x \\
y
\end{smallmatrix}
\right) \mapsto \left( 
\begin{smallmatrix}
-x \\
-y
\end{smallmatrix}
\right) $.
\end{example}

\begin{proposition}
\label{CloudNineTiles}Cloud Nine $(\mathbf{T}_9)$ tiles $\mathbb{R}^2$
by the lattice $\mathbb{L\kern0.5pt}=\mathbb{Z}\left[ \left( 
\begin{smallmatrix}
1 \\
0
\end{smallmatrix}
\right) ,\left( 
\begin{smallmatrix}
0 \\
2
\end{smallmatrix}
\right) \right] $.
\end{proposition}
\begin{proof}
In the Observation after Lemma \ref{TilesByL}, we consider the
equivalence class $\left( 
\begin{smallmatrix}
0 \\
0
\end{smallmatrix}
\right) ^\sim $ where $\sim $ is the $\mathcal{O}_5$-equivalence
relation determined on $\mathbb{Z}^2$ in Section \ref{Monomial} from
$\mathbf{N}$ and $D_9=\pm \left\{ \left(
\begin{smallmatrix}
0 \\
0
\end{smallmatrix}
\right) ,\left( 
\begin{smallmatrix}
3 \\
0
\end{smallmatrix}
\right) ,\left( 
\begin{smallmatrix}
0 \\
2
\end{smallmatrix}
\right) \right\} $. We noted after Lemma \ref{TilesByL} that $\left( 
\begin{smallmatrix}
0 \\
0
\end{smallmatrix}
\right) ^\sim $ is \emph{not} a lattice. For example, from path
($\mathrm{a}$) in Figure \ref{CloudNineCycle}, note that the two
points $\left(
\begin{smallmatrix}
1 \\
2
\end{smallmatrix}
\right) $ and $\left( 
\begin{smallmatrix}
2 \\
0
\end{smallmatrix}
\right) $ are in $\left( 
\begin{smallmatrix}
0 \\
0
\end{smallmatrix}
\right) ^\sim $, but their sum $\left( 
\begin{smallmatrix}
3 \\
2
\end{smallmatrix}
\right) $ is not. In fact, $\left( 
\begin{smallmatrix}
3 \\
2
\end{smallmatrix}
\right) \in \left( 
\begin{smallmatrix}
-1 \\
0
\end{smallmatrix}
\right) ^\sim $, as follows from the extended path ($-\mathrm{b}$), $\left( 
\begin{smallmatrix}
3 \\
2
\end{smallmatrix}
\right) \overset{\raise1pt\hbox{$\scriptstyle R^{\prime
}$}}{\rightarrow }\left(
\begin{smallmatrix}
-1 \\
2
\end{smallmatrix}
\right) \overset{\raise1pt\hbox{$\scriptstyle R^{\prime
}$}}{\rightarrow }\left(
\begin{smallmatrix}
-1 \\
0
\end{smallmatrix}
\right) $. We claim that 
\begin{equation}
\label{ClassUnion}
\begin{pmatrix}
0 \\
0
\end{pmatrix}%
^\sim \cup 
\begin{pmatrix}
1 \\
0
\end{pmatrix}%
^\sim \cup 
\begin{pmatrix}
-1 \\
0
\end{pmatrix}%
^\sim =\mathbb{L\kern0.5pt}
\end{equation}
where $\mathbb{L\kern0.5pt}$ is the lattice
$\mathbb{L\kern0.5pt}=\left\{ \left(
\begin{smallmatrix}
m \\
2n
\end{smallmatrix}
\right) \mid m,n\in \mathbb{Z}\right\} $. It is clear from an
inspection that the three equivalence classes are contained in
$\mathbb{L\kern0.5pt}$, and an induction argument shows that we get
all of $\mathbb{L\kern0.5pt}$ as a disjoint union of the three $\sim
$-classes, proving (\ref{ClassUnion}). Recall that the three points
$\left(
\begin{smallmatrix}
0 \\
0
\end{smallmatrix}
\right) $ and $\left( 
\begin{smallmatrix}
\pm 1 \\
0
\end{smallmatrix}
\right) $ are precisely the solutions $x$ in $\mathbb{Z}^2$ to
$R^{\prime }x=x$, i.e., the points in $B^{(9)}_\infty $ of period
one. Using now (\ref{ClassUnion}) in place of the corresponding
condition in Lemma \ref{TilesByL} above, we conclude that
$\mathbf{T}_9$ tiles $\mathbb{R}^2$ by the lattice
$\mathbb{L\kern0.5pt}$ as claimed. We omit details, but they are
essentially in the proof of Lemma \ref{TilesByL}.
\end{proof}
\begin{remark}
The proof also suggests a possible general scheme for getting lattices
$\mathbb{L\kern0.5pt}$ associated with given ``reptiles'' $\mathbf{T}$
in $\mathbb{R}^\nu $, i.e., the (\ref{Eq3ins1}) ones. Given $\left(
\mathbf{N},D\right) $ in $\nu $ dimensions with $0\in D$, and
specified as in Section \ref{Monomial}, then introduce $R$, the
equivalence relation $\sim $, and the corresponding set $B_\infty $ of
periodic points. Then look for minimal subsets $P$ of $B_\infty $
containing $0$ such that $\bigcup _{p\in P}p^\sim $ is a rank-$\nu $
lattice $\mathbb{L\kern0.5pt}\subset \mathbb{Z}^\nu $, and show that
$\mathbf{T}$ tiles $\mathbb{R}^\nu $ by $\mathbb{L\kern0.5pt}$.
\end{remark}

The Fractal Red Cross is graphed in Figure \ref{RedCrossFigure} and
discussed in \cite{Str94}, and its boundary is fractal in a sense
which is specified there; its scaling similarity dimension is
$d_s=\frac {2\ln 3}{\ln 5}$. We note that the boundaries of
$\mathbf{T}$ for our Examples \ref{JordanMatrix} and \ref{Exa8ins1}
are also fractal. In Figure \ref{Reptile}, for example, the boundary
consists of two intervals, and the rest is fractal. (See also
\cite{Ban96}.) While the extent of the overlap in $\mathbf{T}=\bigcup
_{d\in D}\mathbf{N}^{-1}\left( d+\mathbf{T}\right) $ from
(\ref{Eq3ins1}) isn't quite clear in Example \ref{CloudNine} (Figure
\ref{CloudNineFigure}), there is some overlap which can be entered
into the formula used for the similarity scaling dimension $d_s$ of
$\partial \mathbf{T}$ for Example \ref{FractalRedCross}. If the
relative overlap is specified by a fraction $0<\beta <1$, then the
scaling dimension is $\frac {2\ln \left( 3+\beta \right) }{\ln 5}$ by
the Strichartz argument used in Example \ref{FractalRedCross}.

\begin{figure}
\psfig{file=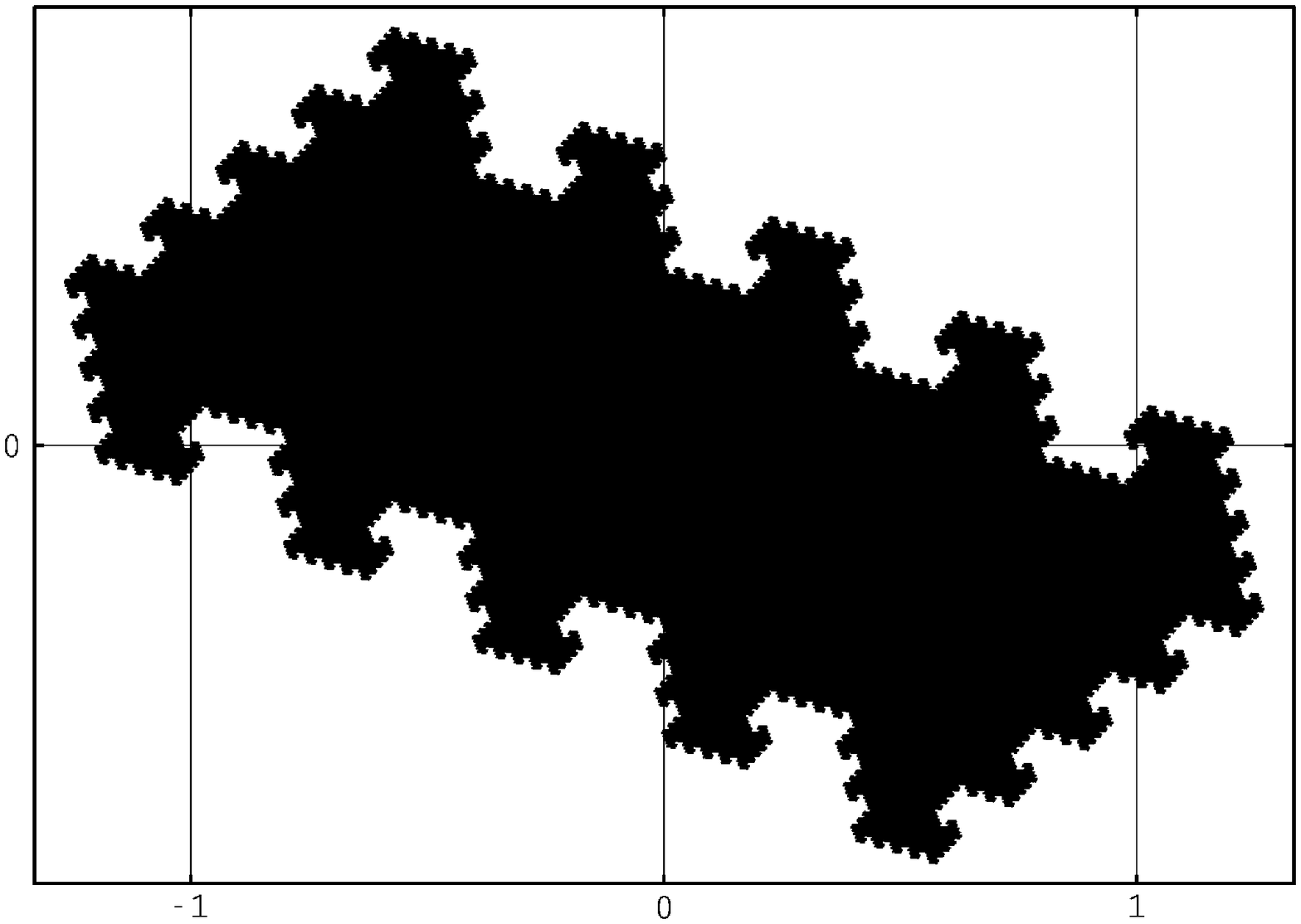,width=5in}
\raisebox{-9pt}[12pt][0pt]{\makebox[\textwidth]{%
$\mathbf{N}=\left(
\begin{smallmatrix}
1 & 2 \\
-2 & 1
\end{smallmatrix}
\right) $, $D=\left\{ \left(
\begin{smallmatrix}
0 \\
0
\end{smallmatrix}
\right), \left(
\begin{smallmatrix}
0 \\
\pm 1
\end{smallmatrix}
\right), \left(
\begin{smallmatrix}
0 \\
\pm 2
\end{smallmatrix}
\right) \right\} $%
}}
\caption{Cloud Three (Example \protect \ref{Clouds})}
\label{CloudThreeFigure}
\end{figure}

\begin{figure}
\psfig{file=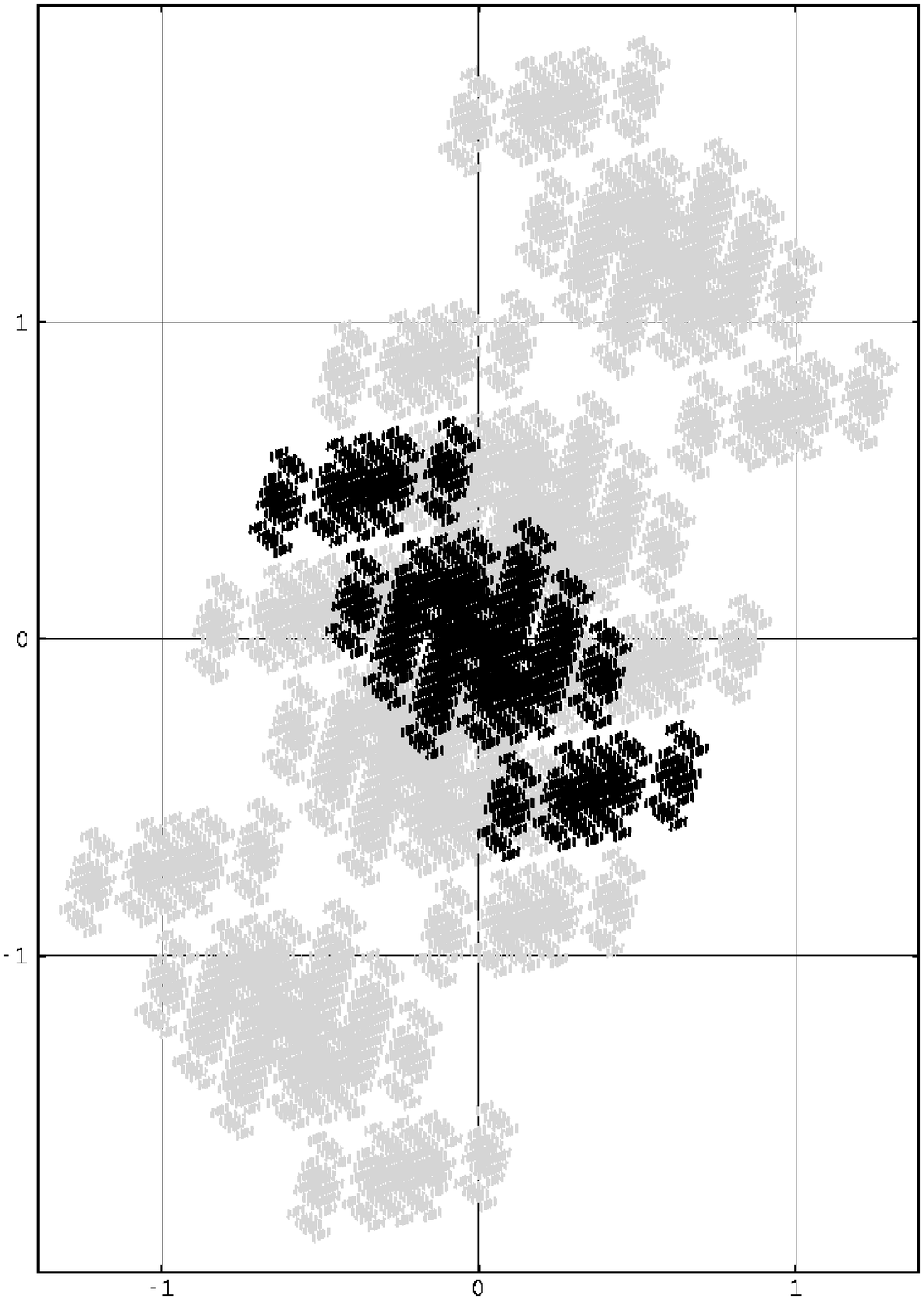,width=5in}
\makebox[\textwidth]{%
We use the same shading convention here as in Figure \ref{CloudNineFigure}.}
\raisebox{-9pt}[12pt][0pt]{\makebox[\textwidth]{%
$\mathbf{N}=\left(
\begin{smallmatrix}
1 & 2 \\
-2 & 1
\end{smallmatrix}
\right) $, $D=\left\{ \left(
\begin{smallmatrix}
0 \\
0
\end{smallmatrix}
\right), \left(
\begin{smallmatrix}
\pm 3 \\
0
\end{smallmatrix}
\right), \left(
\begin{smallmatrix}
\pm 1 \\
0
\end{smallmatrix}
\right) \right\} $%
}}
\caption{Cloud Five (Example \protect \ref{Clouds})}
\label{CloudFiveFigure}
\end{figure}

\begin{example}
\label{Clouds}\label{CloudThree}\label{CloudFive}A simpler version 
of the fractal $\mathbf{T}_9$ derives from taking $D_3=\left\{ \left(
\begin{smallmatrix}
0 \\
0
\end{smallmatrix}
\right) ,\left(
\begin{smallmatrix}
0 \\
\pm 1
\end{smallmatrix}
\right) ,\left(
\begin{smallmatrix}
0 \\
\pm 2
\end{smallmatrix}
\right) \right\} $, still with $\mathbf{N}=\left(
\begin{smallmatrix}
1 & 2 \\
-2 & 1
\end{smallmatrix}
\right) $. Then $B^{(3)}_\infty =\left\{ \left(
\begin{smallmatrix}
0 \\
0
\end{smallmatrix}
\right) ,\left(
\begin{smallmatrix}
\pm 1 \\
0
\end{smallmatrix}
\right) \right\} $, and there are three cycles each with one atom. 
$\mathbf{T}$ is then as in Figure \ref{CloudThreeFigure}.

It follows from a theorem of Lagarias and Wang \cite{LaWa96} that when 
$\mathbf{N}=\left(
\begin{smallmatrix}
1 & 2 \\
-2 & 1
\end{smallmatrix}
\right) $, and $D\subset \mathbb{Z}^2$ is chosen as a full set of
residues for $\mathbb{Z}^2\diagup \mathbf{N}\mathbb{Z}^2$, then the
corresponding fractal $\mathbf{T}=\mathbf{T}\left( D\right) $ is a
$\mathbb{Z}^2$-periodic tile iff $\mathbb{Z}\left( \mathbf{N},D\right)
=\mathbb{Z}^2$ where $\mathbb{Z}\left( \mathbf{N},D\right) $ is the
smallest $\mathbf{N}$-invariant lattice containing $\left\{
d-d^{\prime } \mid d,d^{\prime } \in D\right\} $.  The examples with
$D=D_i$ as given in (\ref{CloudD}) represent the first four cases in a
double-indexed family of examples: $\mathbf{N}=\left(
\begin{smallmatrix}
1 & 2 \\
-2 & 1
\end{smallmatrix}
\right) $, and $D\left( p,q\right) =\mathbf{M}\left( p,q\right) 
\left(D_1\right) $, a full set of residues, where the matrix 
$\mathbf{M}$ is 
\begin{equation}
\mathbf{M}\left( p,q\right) =
\begin{pmatrix}
1-p & 3q \\
2p & 1-q
\end{pmatrix}
\text{,}
\label{Mpq}
\end{equation}
and $\left( p,q\right) \in \mathbb{Z}^2$.  The numbers of elements in
the corresponding sets $B_\infty $ of periodic points are $1$, $9$,
$3$, and $5$, respectively. By checking $\mathbb{Z}\left(
\mathbf{N},D_i\right) $, and using the Lagarias-Wang theorem, we note
that, of the four fractals $\mathbf{T}_i$, $i=1,9,3$, and $5$, only
$\mathbf{T}_1$ (Fractal Red Cross) is a $\mathbb{Z}^2$-tile, i.e., it
is the only one of the four which is a periodic tile of $\mathbb{R}^2$
by $\mathbb{Z}^2$. Specifically, $\left(
\begin{smallmatrix}
1 \\
0
\end{smallmatrix}
\right) \notin \mathbb{Z}\left( \mathbf{N},D_3\right) $, and the point $\left(
\begin{smallmatrix}
0 \\
1
\end{smallmatrix}
\right) $ is not in either of the two lattices $\mathbb{Z}\left(
\mathbf{N},D_5\right) $ or $\mathbb{Z}\left( \mathbf{N},D_9\right) $.
The scaling self-similarity dimensions of the respective boundaries
$\partial \mathbf{T}_i$ also vary from one to the next.  (The graph of
$\mathbf{T}_5$ is shown in Figure \ref{CloudFiveFigure}.)

\begin{figure}
\psfig{file=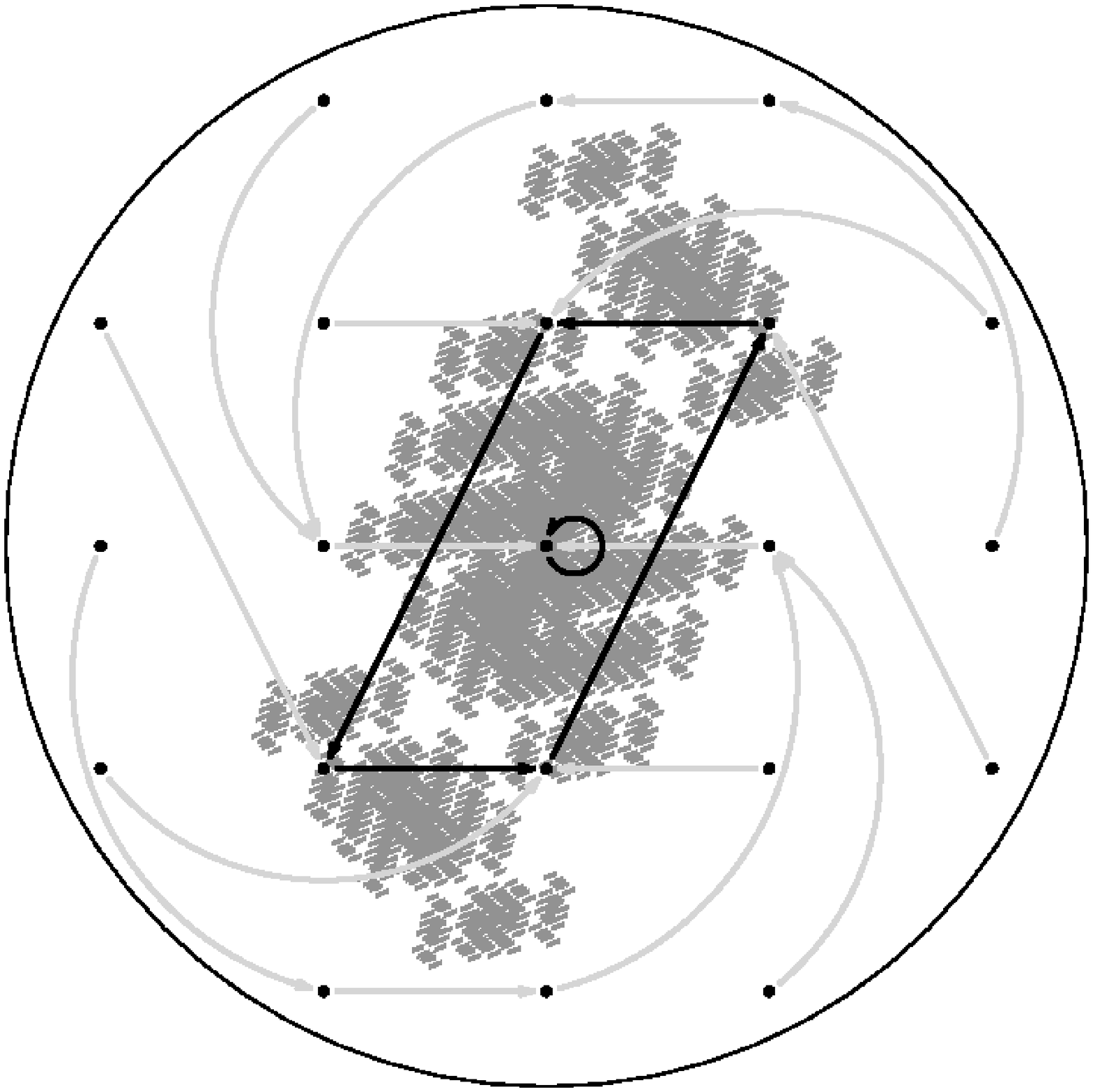,width=360bp}
\caption{Paths of points within spectral circle for Cloud Five
(Example \protect \ref{CloudFive})}
\label{CloudFiveCycle}
\end{figure}

All four $D_i$'s have been selected such that $0\in D_i$, and
$D_i=-D_i$. As a result, $0\in \mathbf{T}_i$, and
$\mathbf{T}_i=-\mathbf{T}_i$, with the point $0$ making a single-atom
cycle $\left\{ \left(
\begin{smallmatrix}
0 \\
0
\end{smallmatrix}
\right) \right\} $. For the two corresponding to $i=5$ and $9$, we
have higher-order cycles in $B^{(i)}_\infty =\mathbf{T}_i\cap
\mathbb{Z}^2$, and they surround the single-atom cycle $\left\{ \left(
\begin{smallmatrix}
0 \\
0
\end{smallmatrix}
\right) \right\} $ with winding number one. The higher-order closed
cycles have the respective orders $6$, for $\mathbf{T}_9$, as computed
in Example \ref{CloudNine} (Cloud Nine, cf. Figure
\ref{CloudNineCycle}), and $4$ for $\mathbf{T}_5$ (Cloud Five,
cf. Figure \ref{CloudFiveCycle}). Let us look closer at
$\mathbf{T}_5$, and repeat the analysis in Example
\ref{CloudNine}. There are $21$ $\mathbb{Z}^2$-points inside the
spectral circle with radius $r=\frac 3{\sqrt 5-1}$, and, of these, $5$
are in $B_\infty ^{(5)}$: one single-atom cycle at the origin $0$, and
the parallelogram-shaped four-atom cycle
\begin{equation*}
\begin{pmatrix}
1 \\ 
1
\end{pmatrix}
\rightarrow 
\begin{pmatrix}
0 \\ 
1
\end{pmatrix}
\rightarrow 
\begin{pmatrix}
-1 \\ 
-1
\end{pmatrix}
\rightarrow 
\begin{pmatrix}
0 \\ 
-1
\end{pmatrix}
\rightarrow 
\begin{pmatrix}
1 \\ 
1
\end{pmatrix}%
\text{.}
\end{equation*}
The single-atom cycle $0$ has multiple flow paths from points inside
the circle, and linking up at $0$. All four atoms in the $4$-cycle
similarly attract flow paths from respective points inside (and
outside) the circle, and an inspection shows that all $21$ points are
accounted for. We have described the following flow paths for Cloud
Five (inside the circle) (see Figure \ref{CloudFiveCycle}):
%{\allowdisplaybreaks
\begin{gather*}
\begin{pmatrix}
1 \\ 
-2
\end{pmatrix}
\rightarrow 
\begin{pmatrix}
1 \\ 
0
\end{pmatrix}
\rightarrow 
\begin{pmatrix}
0 \\ 
0
\end{pmatrix}
\text{,} \\
\begin{pmatrix}
2 \\ 
0
\end{pmatrix}
\rightarrow 
\begin{pmatrix}
1 \\ 
2
\end{pmatrix}
\rightarrow 
\begin{pmatrix}
0 \\ 
2
\end{pmatrix}
\rightarrow 
\begin{pmatrix}
-1 \\ 
0
\end{pmatrix}
\rightarrow 
\begin{pmatrix}
0 \\ 
0
\end{pmatrix}
\text{,} \\
\begin{pmatrix}
1 \\ 
-1
\end{pmatrix}
\rightarrow 
\begin{pmatrix}
0 \\ 
-1
\end{pmatrix}
\text{,} \\
\begin{pmatrix}
2 \\ 
-1
\end{pmatrix}
\rightarrow 
\begin{pmatrix}
1 \\ 
1
\end{pmatrix}
\text{,} \\
\begin{pmatrix}
2 \\ 
1
\end{pmatrix}
\rightarrow 
\begin{pmatrix}
0 \\ 
1
\end{pmatrix}
\text{,}
\end{gather*}%
%}%
and in addition we also have the reflected paths under $\left( 
\begin{smallmatrix}
x \\
y
\end{smallmatrix}
\right) \mapsto \left( 
\begin{smallmatrix}
-x \\
-y
\end{smallmatrix}
\right) $.

We have instances here in $2$ dimensions where $\mathbf{N}$ is fixed,
but different choices of full residue sets $D_i\subset \mathbb{Z}$,
$i=1,3,5,9$ result in dynamical systems $\left(
\mathbb{Z}^2,R_i\right) $ which are mutually non-conjugate. This
follows from the fact (see \cite{LiMa}) that the counting of periodic
points for a dynamical system is a conjugacy invariant, i.e., the
numbers $\#\left\{ x\in X\mid R^n\left( x\right) =x\right\} $ when
some system $R:X\rightarrow X$ is given. In one dimension, such
non-conjugate examples abound; see, e.g., Subsection
\ref{SubsecNew8.1} above.

Note that the Lagarias-Wang theorem (\cite{LaWa96},
\cite{LaWa96Number}) based on $\mathbb{Z}\left( \mathbf{N},D\right) $
does not apply to $\mathbf{N}=\left(
\begin{smallmatrix}
2 & 1 \\
0 & 2
\end{smallmatrix}
\right) $ in Example \ref{JordanMatrix} (Figure \ref{Reptile}) or to
$\mathbf{N}=\left(
\begin{smallmatrix}
0 & 1 \\
4 & 0
\end{smallmatrix}
\right) $ in Example \ref{Exa8ins2} (Figure \ref{UnitSquare}) because
the determinant there is not a prime. On the other hand, it is known
\cite{Str94} that the $\mathbf{T}$'s for those examples are in fact
$\mathbb{Z}^2$-tiles of the plane. But the story for
$\mathbf{N}=\left(
\begin{smallmatrix}
2 & 1 \\
0 & 2
\end{smallmatrix}
\right) $ is completely different if $D$ from Example
\ref{JordanMatrix} is changed to $D^{\prime } :=\left\{ \left(
\begin{smallmatrix}
0 \\
0
\end{smallmatrix}
\right) , \left(
\begin{smallmatrix}
3 \\
0
\end{smallmatrix}
\right) , \left(
\begin{smallmatrix}
0 \\
1
\end{smallmatrix}
\right) , \left(
\begin{smallmatrix}
3 \\
1
\end{smallmatrix}
\right) \right\} $. Then the corresponding $\mathbf{T}\left( D^{\prime
} \right) $ does not tile with $\mathbb{Z}^2$, although we do have
$\mathbb{Z}\left( \mathbf{N},D^{\prime } \right) =\mathbb{Z}^2$ for
that case; see also \cite{LaWa96} and \cite{GrMa92}. A digit set
$D^{\prime } $ such that $\mathbb{Z}\left( \mathbf{N},D^{\prime }
\right) =\mathbb{Z}^2$ is called \emph{primitive.} The $D$ from
Example \ref{JordanMatrix} is also primitive, but the Shark-Jawed
Parallelogram $\mathbf{T}\left( D\right) $ is a $\mathbb{Z}^2$-tile.

For the present examples with $\mathbf{N}=\left(
\begin{smallmatrix}
1 & 2 \\
-2 & 1
\end{smallmatrix}
\right) $ and $D_9$, $D_3$, $D_5$ (i.e., non-primitive cases), it can
be checked that the corresponding fractal $\mathbf{T}_3$ is a tile for
the lattice $\mathbb{Z}\left( \mathbf{N},D_3\right) =\mathbb{Z}\left[
\left(
\begin{smallmatrix}
2 \\
0
\end{smallmatrix}
\right) , \left(
\begin{smallmatrix}
0 \\
1
\end{smallmatrix}
\right) \right] $ which is strictly contained in $\mathbb{Z}^2$. So
$\mathbf{T}_3$ is still periodic, but not with $\mathbb{Z}^2$. We have
$\mathbb{Z}\left( \mathbf{N},D_3\right) =\mathbb{Z}\left[ \left(
\begin{smallmatrix}
2 \\
0
\end{smallmatrix}
\right) , \left(
\begin{smallmatrix}
0 \\
1
\end{smallmatrix}
\right) \right] $ (agreeing with Figure \ref{CloudThreeFigure}),
whereas $\mathbb{Z}\left( \mathbf{N},D_5\right) =\mathbb{Z}\left[
\left(
\begin{smallmatrix}
1 \\
0
\end{smallmatrix}
\right) , \left(
\begin{smallmatrix}
0 \\
2
\end{smallmatrix}
\right) \right] $ (see Figure \ref{CloudFiveFigure}). It then follows
that both $\mathbf{T}_3$ and $\mathbf{T}_5$ have Lebesgue measure
two. It is known from \cite{Ban91} that the boundary for all
$\mathbf{T}$'s has planar Lebesgue measure zero.

The present ``cloud'' examples are the first cases of an infinite
two-parameter family, analogous to the one-parameter family in
Subsection \ref{SubsecNew8.1}. For $\left( p,q\right) \in
\mathbb{Z}^2$, set
\begin{equation}
D\left( p,q\right) :=\left\{ 
\begin{pmatrix}
0 \\
0
\end{pmatrix}
, \pm
\begin{pmatrix}
1-p \\
2p
\end{pmatrix}
, \pm
\begin{pmatrix}
3q \\
1-q
\end{pmatrix}
\right\} \text{.}
\label{Dpq}
\end{equation}
Let the matrix $\mathbf{N}=\left( 
\begin{smallmatrix}
1 & 2 \\
-2 & 1
\end{smallmatrix}
\right) $ be fixed. We then get an example for each value of $\left(
p,q\right) \in \mathbb{Z}^2$, with $\left( 0,0\right) $ corresponding
to $D_1$, $\left( 1,0\right) $ to $D_3$, $\left( 0,1\right) $ to
$D_5$, and $\left( p,q\right) =\left( 1,1\right) $ to $D_9$. The
cycle-atom structure for these examples in $\mathbb{Z}^2$ is more
complicated than that for the Subsection-\ref{SubsecNew8.1} examples,
where $D\left( p\right) =\left\{ 0,p\right\} $, $p\in
\mathbb{N}_{odd}$, and
\begin{equation*}
B_\infty \left( p\right) =-\left\{ 0,1,\dots ,p\right\} \longleftrightarrow 
\mathbb{Z}\diagup p\mathbb{Z}\text{;}
\end{equation*}
cf. Remark \ref{RemNew8.2}. For the present examples, there is not a
simple analogue of this; but we can say something about the fixed
points
\begin{equation*}
B_\infty ^1\left( p,q\right) =\left\{ x\in B_\infty \left( p,q\right) 
\mid R_{p,q}\left( x\right) =x\right\}
\end{equation*}
where $R=R_{p,q}$ is the $\mathbb{Z}^2$-transformation from
(\ref{Trans}) associated with $\left( \mathbf{N},D\left( p,q\right)
\right) $:
\begin{enumerate}
\item \label{CloudProperty(i)}If $p$ is \emph{even,} then $B_\infty 
^1\left( p,q\right) =\left\{ \left(
\begin{smallmatrix}
0 \\
0
\end{smallmatrix}
\right) \right\} $, and so this accounts for the period-one structure
of examples $D_1$ and $D_5$.

\item \label{CloudProperty(ii)}If $p=2n+1$ $\left( n\in \mathbb{Z}\right) $,
i.e., $p$ \emph{odd,} then $B_\infty ^1\left( p,q\right) =\left\{ \left(
\begin{smallmatrix}
0 \\
0
\end{smallmatrix}
\right) , \pm \left(
\begin{smallmatrix}
2n+1 \\
n
\end{smallmatrix}
\right) \right\} $, and this accounts for the three fixed points for
each of the examples $D_3$ and $D_9$.
\end{enumerate}

To see this, let $E=\left\{ \left(
\begin{smallmatrix}
\alpha \\
\beta
\end{smallmatrix}
\right) \mid \alpha ,\beta \in \left\{ 0,\pm 1\right\} ,\alpha \cdot
\beta =0\right\} =D_1$. It follows from Section \ref{Monomial} that
$x\in B_\infty ^1\left( p,q\right) $ iff, for some $\varepsilon \in
E$, $x=\left( \openone -\mathbf{N}\right) ^{-1}\mathbf{M}\varepsilon
\in \mathbb{Z}^2$ where $\mathbf{M}=\left(
\begin{smallmatrix}
1-p & 3q \\
2p & 1-q
\end{smallmatrix}
\right) $. Indeed, writing $\varepsilon =\left( 
\begin{smallmatrix}
\alpha \\
\beta
\end{smallmatrix}
\right) $, we get 
\begin{equation*}
x=\left( \openone -\mathbf{N}\right) ^{-1}\mathbf{M}\varepsilon =
\begin{pmatrix}
p\alpha +\frac {\left( 1-q\right) \beta}2 \\
-\frac {\left( 1-p\right) \alpha +3q\beta }2
\end{pmatrix}
\text{.}
\end{equation*}
Referring back to the definition of $E\subset \mathbb{Z}^2$, the
conclusions in (\ref{CloudProperty(i)}) and (\ref{CloudProperty(ii)})
follow.

Of course, there is a similar argument for the set $B_\infty ^k\left(
p,q\right) $ of points of period $k$ more than one: indeed, $x\in
B_\infty ^k\left( p,q\right) $ iff there are points $\varepsilon
_1,\varepsilon _2,\dots ,\varepsilon _k\in E$ such that
\begin{equation}
x=\left( \openone -\mathbf{N}^k\right) ^{-1}\left( \mathbf{M}\varepsilon 
_1+\mathbf{N}\mathbf{M}\varepsilon _2+\dots 
+\mathbf{N}^{k-1}\mathbf{M}\varepsilon _k\right) \in \mathbb{Z}^2\text{.}
\label{Twist}
\end{equation}

So this is a two-dimensional matrix divisibility problem analogous to
the simpler one-dimensional one in Subsection \ref{SubsecNew8.1}, but
the analysis here is complicated considerably by the non-commutativity
of the two matrices $\mathbf{N}$ and $\mathbf{M}$. Indeed, if
\begin{equation*}
\mathbf{M}^{\prime } =\mathbf{N}\mathbf{M}\mathbf{N}^{-1}=\frac 15
\begin{pmatrix}
5+3p+2q & q-6p \\
4p-14q & 5-8p-7q
\end{pmatrix}
\end{equation*}
the entries of $\mathbf{M}^{\prime }$ are integral iff $p=q$. In that
case, $\mathbf{M}^{\prime } =\left(
\begin{smallmatrix}
1+p & -p \\
-2p & 1-3p
\end{smallmatrix}
\right) $ (${}\neq \mathbf{M}$ when $p\neq 0$).

For all the matrix examples $\left( \mathbf{N},D\right) $ in
$\mathbb{Z}^\nu $, the function $C:\mathbb{N}\rightarrow \left\{
0,1,2,\dots \right\} $ which counts the number of cycles $C\left(
k\right) $ of a given (minimal) period $k$, is an invariant for the
corresponding $R:\mathbb{Z}^\nu \rightarrow \mathbb{Z}^\nu $ from
(\ref{Trans}), and the corresponding Ruelle Zeta-function (see
\cite{Rue88}, \cite{Rue94}, \cite{LiMa}) computes out as
\begin{equation*}
\zeta _C\left( t\right) =\prod _k\frac 1{\left( 1-t^k\right) ^{C\left( k\right) }}\text{.}
\end{equation*}
\end{example}

When the formula
$\mathbf{M}\left( p,q\right) =\left(
\begin{smallmatrix}
1-p & 3q \\
2p & 1-q
\end{smallmatrix}
\right) $ is substituted in for the matrix $\mathbf{M}$ in
(\ref{Twist}), we arrive at a certain dichotomy for reflection
symmetry of cycles of the corresponding $D\left( p,q\right) $-examples
with transformation $R\left( p,q\right) $. We specialize to $p=1$ for
clarity, and note that, if $q$ is odd, then each cycle of (minimal)
period $k>1$ must be invariant under the reflection $x\mapsto -x$ in
$\mathbb{Z}^2$, while, if $q$ is even, this reflection interchanges
distinct cycles. Let $p=1$, and $q$ be odd: if the cycle is written in
the form $x^0,x^1,\dots ,x^{2m}$, $x^{2m}=x^0$, and $Rx^i=x^{i+1}$,
then $x^{i+m}=-x^i$, $i=0,\dots ,2m$. In other words, the period is
$k=2m$ and the second half from $x^m$ to $x^{2m-1}$ is the reflection
under $x\mapsto -x$ of the first half. Equivalently, these cycles $c$
are (closed) $R$-induced paths from $\mathbb{Z}\diagup k\mathbb{Z}$
into $\mathbb{Z}^2$ which, for $q$ odd, have the alternating property
relative to the index-$2$ subgroup of $\mathbb{Z}\diagup 2m\mathbb{Z}$
as indicated, i.e., $k=2m$, and $c\left( i+m\right) =-c\left( i\right) $.
\bigskip

\begin{figure}
\begin{minipage}{68bp}
\psfig{file=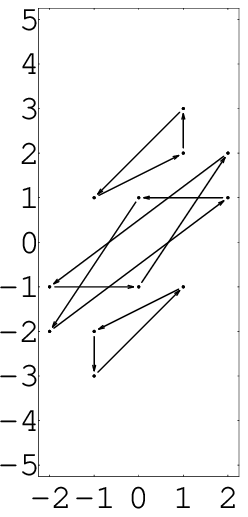}
\makebox[68bp]{(a) $p=1$, $q=2$}
\end{minipage}\hfill%
\begin{minipage}{94bp}
\psfig{file=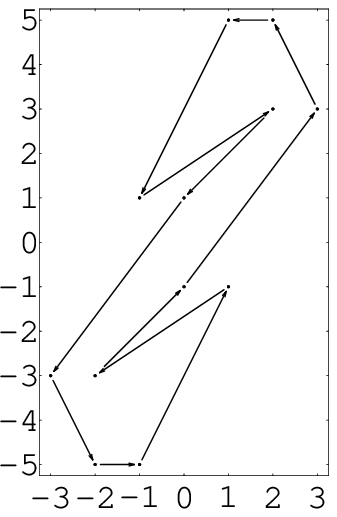}
\makebox[94bp]{(b) $p=1$, $q=3$}
\end{minipage}\hfill%
\begin{minipage}{172bp}
\psfig{file=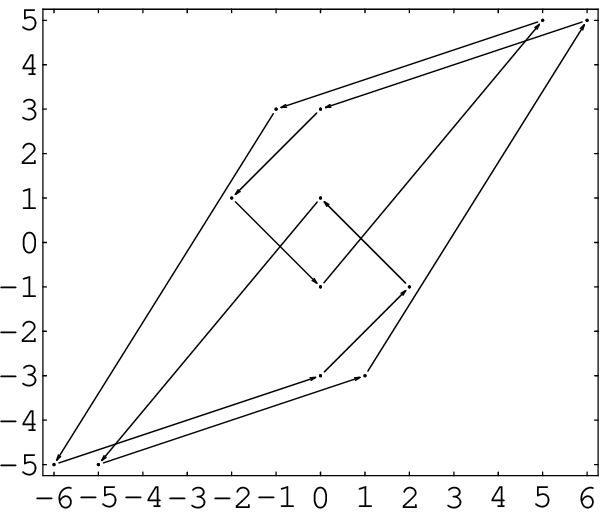}
\makebox[172bp]{(c) $p=1$, $q=5$}
\end{minipage}
\caption{Some non-singleton cycles in the $\left( p,q\right) $
family of Clouds}
\label{pqCycles}
\end{figure}

\noindent {\bfseries Cases:}
\begin{list}{}{\setlength{\leftmargin}{\leftmargini}
\setlength{\itemsep}{\customskipamount}
\setlength{\topsep}{\customskipamount}
\setlength{\labelsep}{2.7778pt}}
\item[\hss\llap{$p=$}]  $q=1$: $\mathbf{M}=\left(
\begin{smallmatrix}
0 & 3 \\
2 & 0
\end{smallmatrix}
\right) $, and we have a non-sin\-gle\-ton cycle with period $k=6$,
and $x^{i+3}=-x^i$, $i=0,\dots ,6$, $x^6=x^0$. (See Figure
\ref{CloudNineCycle}.)
 
\item[\hss\llap{$p=$}]  $1$, $q=2$: $\mathbf{M}=\left(
\begin{smallmatrix}
0 & 6 \\
2 & -1
\end{smallmatrix}
\right) $, and we have four non-sin\-gle\-ton cycles each with period
$3$. (See Figure \ref{pqCycles}(a).)

\item[\hss\llap{$p=$}]  $1$, $q=3$ or $q=5$: In each case, we have orbits of
(minimal) period $k=12$, and each satisfying $x^{i+6}=-x^i$, $i=0,\dots ,12$,
$x^{12}=x^0$. (See Figure \ref{pqCycles}(b) and (c).)
In each case, we may take $x^0=\left(
\begin{smallmatrix}
0 \\
1
\end{smallmatrix}
\right) $. If $q=5$, the $12$-cycle has winding number $3$ around $0=\left(
\begin{smallmatrix}
0 \\
0
\end{smallmatrix}
\right) $ in the plane. In fact, $\left( p,q\right) =\left( 1,3\right)
$ yields $x^{i+1}=Rx^i$, $i=0,\dots ,12$ as follows:
\begin{multline*}
\begin{pmatrix}
0 \\ 
1
\end{pmatrix}
\overset{\raise1pt\hbox{$\scriptstyle R$}}{\rightarrow }
\begin{pmatrix}
-3 \\ 
-3
\end{pmatrix}
\overset{\raise1pt\hbox{$\scriptstyle R$}}{\rightarrow }
\begin{pmatrix}
-2 \\ 
-5
\end{pmatrix}
\overset{\raise1pt\hbox{$\scriptstyle R$}}{\rightarrow }
\begin{pmatrix}
-1 \\ 
-5
\end{pmatrix}
\overset{\raise1pt\hbox{$\scriptstyle R$}}{\rightarrow }
\begin{pmatrix}
1 \\ 
-1
\end{pmatrix}
\overset{\raise1pt\hbox{$\scriptstyle R$}}{\rightarrow }
\begin{pmatrix}
-2 \\ 
-3
\end{pmatrix}
\\
\overset{\raise1pt\hbox{$\scriptstyle R$}}{\rightarrow }
\begin{pmatrix}
0 \\ 
-1
\end{pmatrix}
\overset{\raise1pt\hbox{$\scriptstyle R$}}{\rightarrow }
\begin{pmatrix}
3 \\ 
3
\end{pmatrix}
\overset{\raise1pt\hbox{$\scriptstyle R$}}{\rightarrow }
\begin{pmatrix}
2 \\ 
5
\end{pmatrix}
\overset{\raise1pt\hbox{$\scriptstyle R$}}{\rightarrow }
\begin{pmatrix}
1 \\ 
5
\end{pmatrix}
\overset{\raise1pt\hbox{$\scriptstyle R$}}{\rightarrow }
\begin{pmatrix}
-1 \\ 
1
\end{pmatrix}
\overset{\raise1pt\hbox{$\scriptstyle R$}}{\rightarrow }
\begin{pmatrix}
2 \\ 
3
\end{pmatrix}
\overset{\raise1pt\hbox{$\scriptstyle R$}}{\rightarrow }
\begin{pmatrix}
0 \\ 
1
\end{pmatrix}
\text{;}
\end{multline*}
while $\left( p,q\right) =\left( 1,5\right) $ yields:
\begin{multline*}
\begin{pmatrix}
0 \\ 
1
\end{pmatrix}
\overset{\raise1pt\hbox{$\scriptstyle R$}}{\rightarrow }
\begin{pmatrix}
-5 \\ 
-5
\end{pmatrix}
\overset{\raise1pt\hbox{$\scriptstyle R$}}{\rightarrow }
\begin{pmatrix}
1 \\ 
-3
\end{pmatrix}
\overset{\raise1pt\hbox{$\scriptstyle R$}}{\rightarrow }
\begin{pmatrix}
6 \\ 
5
\end{pmatrix}
\overset{\raise1pt\hbox{$\scriptstyle R$}}{\rightarrow }
\begin{pmatrix}
0 \\ 
3
\end{pmatrix}
\overset{\raise1pt\hbox{$\scriptstyle R$}}{\rightarrow }
\begin{pmatrix}
-2 \\ 
1
\end{pmatrix}
\\
\overset{\raise1pt\hbox{$\scriptstyle R$}}{\rightarrow }
\begin{pmatrix}
0 \\ 
-1
\end{pmatrix}
\overset{\raise1pt\hbox{$\scriptstyle R$}}{\rightarrow }
\begin{pmatrix}
5 \\ 
5
\end{pmatrix}
\overset{\raise1pt\hbox{$\scriptstyle R$}}{\rightarrow }
\begin{pmatrix}
-1 \\ 
3
\end{pmatrix}
\overset{\raise1pt\hbox{$\scriptstyle R$}}{\rightarrow }
\begin{pmatrix}
-6 \\ 
-5
\end{pmatrix}
\overset{\raise1pt\hbox{$\scriptstyle R$}}{\rightarrow }
\begin{pmatrix}
0 \\ 
-3
\end{pmatrix}
\overset{\raise1pt\hbox{$\scriptstyle R$}}{\rightarrow }
\begin{pmatrix}
2 \\ 
-1
\end{pmatrix}
\overset{\raise1pt\hbox{$\scriptstyle R$}}{\rightarrow }
\begin{pmatrix}
0 \\ 
1
\end{pmatrix}
\text{.}
\end{multline*}
\end{list}

By a theorem of Kronecker and the authors of \cite{GHJ}, all the
$\mathbb{Z}$-matrices $\mathbf{N}$ with norm satisfying $\left\|
\mathbf{N}\right\| <2$ can be described: the admissible values for the
norm must be $2\cos {\frac \pi n}$ where $n\in \mathbb{N}$, $n\geq
2$. Of the present examples, only Example \ref{Exa8ins1} satisfies
this: there $\left\| \mathbf{N}\right\| =\sqrt 2=2\cos \frac \pi 4$.

The authors of \cite{BaGe94} and \cite{Gel96} introduce a notion of
\emph{isomorphism} for reptiles $\mathbf{T}$, based on (\ref{Eq3ins1})
and some given affine system $\left( \mathbf{N},D\right) $ as
specified there; i.e., with the integral $\nu \times \nu $ matrix
$\mathbf{N}$ expansive, and $D$ a full set of residues for
$\mathbb{Z}^\nu \diagup \mathbf{N}\mathbb{Z}^\nu $. When a system is
given, we have the corresponding maps
\begin{equation*}
\tau _d^\mathbf{N}\left( x\right) :=\mathbf{N}^{-1}\left( x+d\right) \text{.}
\end{equation*}
Two given systems $\left( \mathbf{N},D\right) $ and $\left(
\mathbf{N}^{\prime },D^{\prime }\right) $ are said to be
\emph{isomorphic} if there is an invertible affine mapping $\alpha
:x\mapsto \mathbf{A}x+v$, i.e., $\mathbf{A}\in GL_\nu $ and $v\in
\mathbb{R}^\nu $, such that $\alpha \circ \tau _i^\mathbf{N}=\tau
_i^{\mathbf{N}^{\prime }}\circ \alpha $ where $\tau _i^\mathbf{N}$ and
$\tau _i^{\mathbf{N}^{\prime }}$ are defined from the systems for some
bijection $d_i\mapsto d_i^{\prime }$ between the respective digit sets
$D$ and $D^{\prime }$.

If $\alpha $ is such a mapping, defining the isomorphism, then it
follows from the uniqueness of the (\ref{Eq3ins1})-representation for
$\mathbf{T}$ that $\alpha \left( \mathbf{T}\right) =\mathbf{T}^{\prime
}$. But, more importantly, we also get that the respective parts $\tau
_{i_1}^\mathbf{N}\circ \dots \circ \tau _{i_k}^\mathbf{N}\left(
\mathbf{T}\right) $ in $\mathbf{T}$ are mapped by $\alpha $ onto the
corresponding parts $\tau _{i_1}^{\mathbf{N}^{\prime }}\circ \dots
\circ \tau _{i_k}^{\mathbf{N}^{\prime }}\left( \mathbf{T}^{\prime
}\right) $ in $\mathbf{T}^{\prime }$. Now, using the matrix
representation (\ref{Mpq}) for the infinite family of double-indexed
digit sets $D\left( p,q\right) $ in (\ref{Dpq}), it can be checked
that, when $\mathbf{N}=\left(
\begin{smallmatrix}
1 & 2 \\
-2 & 1
\end{smallmatrix}
\right) $ is fixed with $\left( p,q\right) $ varying over
$\mathbb{Z}^2$, we get a corresponding infinite family of
non-isomorphic reptile examples $\mathbf{T}\left( p,q\right) $.
\medskip

\subsection{Box spline matrices}\label{SubsecNew9.3}\label{BoxSpline}

\ \bigskip

Example \ref{Exa8ins2} can be put into a general class of examples of the
form 
\begin{equation}
\mathbf{N}=
\begin{pmatrix}
0 & 1 \\ 
N & 0
\end{pmatrix}
\text{,\quad }D=\left\{ 
\begin{pmatrix}
0 \\ 
s_i
\end{pmatrix}
\biggm|i=1,2,\dots ,N\right\} \text{,}  \label{EqNew8.13}
\end{equation}
where the $s_i$ are mutually incongruent modulo $N$. We will now describe a
method of reducing the analysis of the cycle-atom structure of this system
to that of the one-dimensional system defined by (\ref{Eq3.1}). We will
apply this analysis to Example \ref{Exa8ins2} at the end of this subsection,
and we find in that case that there are two cycles with one atom each, and
one cycle with two atoms. The fixed points for $R$ are $\left( 
\begin{smallmatrix}
0 \\
0
\end{smallmatrix}\right) $ and $\left( 
\begin{smallmatrix}
-1 \\
-1
\end{smallmatrix}\right) $ and the order $2$ orbit is $\left\{ \left( 
\begin{smallmatrix}
-1 \\
0
\end{smallmatrix}\right) ,\left( 
\begin{smallmatrix}
0 \\
-1
\end{smallmatrix}\right) \right\} $. We will also compute the cycle-atom
structure for some examples based on $\mathbf{N}=\left( 
\begin{smallmatrix}
0 & 1 \\
2 & 0
\end{smallmatrix}\right) $.

This analysis will more generally be applicable to $m\times m$ box spline
matrices 
\begin{equation}
\mathbf{N}=
\begin{pmatrix}
0 & 1 & 0 & 0 & \cdots & 0 & 0 \\ 
0 & 0 & 1 & 0 & \cdots & 0 & 0 \\ 
0 & 0 & 0 & 1 & \cdots & 0 & 0 \\ 
0 & 0 & 0 & 0 & \cdots & 0 & 0 \\ 
\vdots & \vdots & \vdots & \vdots & \ddots & \vdots & \vdots \\ 
0 & 0 & 0 & 0 & \cdots & 0 & 1 \\ 
N & 0 & 0 & 0 & \cdots & 0 & 0
\end{pmatrix}
\text{,\quad }d_i=
\begin{pmatrix}
0 \\ 
0 \\ 
0 \\ 
0 \\ 
\vdots \\ 
0 \\ 
s_i
\end{pmatrix}
\text{,}  \label{EqNew8.14}
\end{equation}
where $s_1,\dots ,s_N$ are mutually incongruent modulo $N$. Recall that $%
N=\left| \det \mathbf{N}\right| $. Note that $\mathbf{N}^m=N\openone_m$,
and the characteristic polynomial is $\lambda ^m-1$, so $\mathbf{N}$ has the
eigenvalues $N^{\frac 1m}e^{i\frac{2\pi k}m}$, $k=1,\dots ,m$, each with
multiplicity $1$. Thus Remark \ref{Rem8ins1} will apply.
However, the present analysis is based on a different method, which applies
also in the general setting of a branching function system from (\ref{Eq2.4}%
). Let us illustrate the method for $m=2$. Then, from a given branching
function system $\sigma _i$ with residual map $R$ defined as in Scholium \ref
{Sch2.4}, we define a new branching function system $\boldsymbol{\sigma }_i$
on $%
\mathbb{N}^2$ by 
\begin{equation}
\boldsymbol{\sigma }_i
\begin{pmatrix}
x \\ 
y
\end{pmatrix}
=
\begin{pmatrix}
0 & 1 \\ 
\sigma _i & 0
\end{pmatrix}
\begin{pmatrix}
x \\ 
y
\end{pmatrix}
=
\begin{pmatrix}
y \\ 
\sigma _i\left( x\right)
\end{pmatrix}
\text{,}  \label{EqNew8.15}
\end{equation}
and the corresponding residual map is then 
\begin{equation}
\mathbf{R}
\begin{pmatrix}
x \\ 
y
\end{pmatrix}
=
\begin{pmatrix}
0 & R \\ 
1 & 0
\end{pmatrix}
\begin{pmatrix}
x \\ 
y
\end{pmatrix}
=
\begin{pmatrix}
R\left( y\right) \\ 
x
\end{pmatrix}
\text{.}  \label{EqNew8.16}
\end{equation}
By Scholium \ref{Sch2.4}, the $\approx $-classes in $\mathbb{N}^2$ are now
given by the orbits of the semigroup $G_0$, and in this case these orbits
are given by 
\begin{equation}
\boldsymbol{\sigma }_{i_1}\cdots \boldsymbol{\sigma }_{i_{2k}}\mathbf{R}^{2k}
\begin{pmatrix}
x \\ 
y
\end{pmatrix}
=
\begin{pmatrix}
\sigma _{i_2}\sigma _{i_4}\cdots \sigma _{i_{2k}}R^k\left( x\right) \\ 
\sigma _{i_1}\sigma _{i_3}\cdots \sigma _{i_{2k-1}}R^k\left( y\right)
\end{pmatrix}
 \label{EqNew8.17}
\end{equation}
and 
\begin{equation}
\boldsymbol{\sigma }_{i_1}\cdots \boldsymbol{\sigma }_{i_{2k+1}}\mathbf{R}%
^{2k+1}
\begin{pmatrix}
x \\ 
y
\end{pmatrix}
=
\begin{pmatrix}
\sigma _{i_2}\sigma _{i_4}\cdots \sigma _{i_{2k}}R^k\left( x\right) \\ 
\sigma _{i_1}\sigma _{i_3}\cdots \sigma _{i_{2k+1}}R^{k+1}\left( y\right)
\end{pmatrix}
\text{.}  \label{EqNew8.18}
\end{equation}
From these formulae, it follows immediately that the $\approx $-equivalence
classes in $\mathbb{N}^2$ have the form 
\begin{equation*}
\begin{pmatrix}
a_1 \\ 
a_2
\end{pmatrix}
\text{,}
\end{equation*}
where $a_1$, $a_2$ are $\approx $-equivalence classes in $\mathbb{N}$ for
the old system. This argument works for a general $m$ as well. So we want to
construct a new branching function system $\boldsymbol{\sigma }_k$ of order $N$
from a given one $\sigma _k$. This is the following $m$'th root operation:
replace $\mathbb{N}$ by $\mathbb{N}^m$ and define 
\begin{equation}
\boldsymbol{\sigma }_k\left( x_1,\dots ,x_m\right) =\left( x_2,x_3,\dots
,x_m,\sigma _k\left( x_1\right) \right)  \label{EqNew8.19}
\end{equation}
for $x_1,\dots ,x_m\in \mathbb{N}$. Then the coding map for the new system
is 
\begin{equation}
\boldsymbol{\sigma }\left( x_1,\dots ,x_m\right) =\left( \sigma \left(
x_1\right) _1,\sigma \left( x_2\right) _1,\dots ,\sigma \left( x_m\right)
_1,\sigma \left( x_1\right) _2,\dots ,\sigma \left( x_m\right) _2,\sigma
\left( x_1\right) _3,\dots \right) \text{.}  \label{EqNew8.20}
\end{equation}
Thus the new coding map is injective if the old one is so. The new $R$-map
(of Scholium \ref{Sch2.4}) is given in terms of the old one by 
\begin{equation}
\mathbf{R}\left( x_1,\dots ,x_m\right) =\left( R\left( x_m\right)
,x_1,x_2,\dots ,x_{m-1}\right) \text{,}  \label{EqNew8.21}
\end{equation}
and hence 
\begin{multline}
\mathbf{R}^{km+l}\left( x_1,\dots ,x_m\right) \\ =
\Bigl( R^{k+1}\left(
x_{m+1-l}\right) ,R^{k+1}\left( x_{m+2-l}\right) ,\dots ,R^{k+1}\left(
x_m\right) , \\ R^k\left( x_1\right) ,R^k\left( x_2\right) ,\dots ,R^k\left(
x_{m-l}\right) \Bigr)  \label{EqNew8.22}
\end{multline}
for $k=0,1,\dots $ and $l=1,2,\dots ,m$.

Thus, by Corollary \ref{Cor2.4}, we have 
\begin{equation}
\left( x_1,\dots ,x_m\right) \approx \left( y_1,\dots ,y_m\right)
\label{EqNew8.23}
\end{equation}
if and only if 
\begin{equation}
x_i\approx y_i\text{\quad for }i=1,\dots ,m\text{.}  \label{EqNew8.24}
\end{equation}
It follows that the atoms in the new system are all $m$-tuples 
\begin{equation*}
\left( a_1,\dots ,a_m\right)
\end{equation*}
where $a_1,\dots ,a_m$ are atoms in the old system.

The new system of cycles is slightly more complicated to describe. By
Theorem \ref{Thm4.1}, it suffices to compute the action
$\boldsymbol{\tau }$ of $%
\mathbb{Z}$ on the new atoms.
By the formulae for $\mathbf{R}$ and $\boldsymbol{%
\sigma }_k$ this is given by 
\begin{equation}
\boldsymbol{\tau }\left( a_1,\dots ,a_m\right) =\left( \tau \left( a_m\right)
,a_1,a_2,\dots ,a_{m-1}\right)  \label{EqNew8.25}
\end{equation}
and hence 
\begin{equation}
\boldsymbol{\tau }^{-1}\left( a_1,\dots ,a_m\right) =\left( a_2,a_3,\dots
,a_m,\tau ^{-1}\left( a_1\right) \right) \text{.}  \label{EqNew8.26}
\end{equation}
The new cycles are the orbits under this action. One computes 
\begin{equation}
\boldsymbol{\tau }^{km+l}\left( a_1,\dots ,a_m\right) =\left( \tau ^{k+1}\left(
a_{m+1-l}\right) ,\dots ,\tau ^{k+1}\left( a_m\right) ,\tau ^k\left(
a_1\right) ,\dots ,\tau ^k\left( a_{m-l}\right) \right)  \label{EqNew8.27}
\end{equation}
for $k=0,1,\dots $ and $l=1,2,\dots ,m$. Thus the period $\limfunc{Per}%
\left( a_1,\dots ,a_m\right) $ is the smallest integer $km+l$ such that 
\begin{equation}
\begin{aligned} \tau ^{k+1}\left( a_{m+1-l}\right) &=a_1 \\
\text{\enspace} &\vdots \text{\enspace} \\ \tau ^{k+1}\left( a_m\right)
&=a_l \\ \tau ^k\left( a_1\right) &=a_{l+1} \\
\text{\enspace} &\vdots \text{\enspace} \\ \tau ^k\left( a_{m-l}\right)
&=a_m\text{,} \end{aligned}  \label{EqNew8.28}
\end{equation}
so the computation of the cycle-atom structure of the new system from that
of the old is in general complicated. One little observation from the case $%
l=m$ above is that $\limfunc{Per}\left( a_1,\dots ,a_m\right) $ divides $%
\func{lcm}\left( m,\limfunc{Per}\left( a_1\right) ,\dots ,\limfunc{Per}%
\left( a_m\right) \right) $. Thus $\limfunc{Per}\left( a_1,\dots ,a_m\right) 
$ is infinite if and only if $\limfunc{Per}\left( a_i\right) $ is infinite
form some $i$.

Let us finally use this to compute the cycle-atom structure in some specific
cases, so let us return to the box spline system defined by (\ref{EqNew8.14}%
). If now for example $s_i=i$ for $i=0,1,\dots ,N-1$, then the original
$\bmod{\,N}$
branching function system
has two cycles each containing one atom by the case of Proposition
\ref{ProNew8.2}.
Then the new system has $2^m$ atoms, and the number
of cycles is equal to the number of $m$-tuples of two elements, up to cyclic
permutations. For example, in Example \ref{Exa8ins2} we obtain two cycles
with one atom each, and one cycle with two atoms. The fixed points for $%
\mathbf{R}$ are $\left( \begin{smallmatrix}
0 \\
0
\end{smallmatrix}\right) $, $\left( 
\begin{smallmatrix}
-1 \\
-1
\end{smallmatrix}\right) $, and the order-$2$ orbit is $\left\{ \left( 
\begin{smallmatrix}
-1 \\
0
\end{smallmatrix}\right) ,\left( 
\begin{smallmatrix}
0 \\
-1
\end{smallmatrix}\right) \right\} $; see Figure \ref{UnitSquare}.

As a final example, look at 
\begin{equation}
\mathbf{N}=
\begin{pmatrix}
0 & 1 & 0 \\ 
0 & 0 & 1 \\ 
2 & 0 & 0
\end{pmatrix}
\text{,\quad }d_0=
\begin{pmatrix}
0 \\ 
0 \\ 
0
\end{pmatrix}
\text{,\quad }d_1=
\begin{pmatrix}
0 \\ 
0 \\ 
3
\end{pmatrix}
\text{.}  \label{EqNew8.29}
\end{equation}
By Proposition \ref{ProNew8.2}, the corresponding one-dimensional system $%
N=2 $, $s_0=0$, $s_1=3$ has two atoms, 
\begin{equation*}
a_1=\left\{ \dots ,-9,-6,-3\right\}
\end{equation*}
and 
\begin{equation*}
a_2=\left\{ 0,3,6,\dots \right\} \text{,}
\end{equation*}
that are fixed under $\tau $ and two others, 
\begin{equation*}
a_3=\left\{ \dots ,-5,-2,1,4,\dots \right\}
\end{equation*}
and 
\begin{equation*}
a_4=\left\{ \dots ,-4,-1,2,5,\dots \right\} \text{,}
\end{equation*}
that are intertwined under $\tau $. Thus the system above has the $4^3=64$
atoms 
\begin{equation*}
\begin{pmatrix}
a_{i_1} \\ 
a_{i_2} \\ 
a_{i_3}
\end{pmatrix}
\text{.}
\end{equation*}
Considering the more amenable case 
\begin{equation}
\mathbf{N}=
\begin{pmatrix}
0 & 1 \\ 
2 & 0
\end{pmatrix}
\text{,\quad }d_0=
\begin{pmatrix}
0 \\ 
0
\end{pmatrix}
\text{,\quad }d_1=
\begin{pmatrix}
0 \\ 
3
\end{pmatrix}
 \label{EqNew8.30}
\end{equation}
instead, we have $4^2=16$ atoms 
\begin{equation*}
\begin{pmatrix}
a_{i_1} \\ 
a_{i_2}
\end{pmatrix}
\end{equation*}
with $i_1,i_2\in \left\{ 1,2,3,4\right\} $. Let us compute the
$\boldsymbol{\tau 
}$ action on these using (\ref{EqNew8.25}), and the abbreviation $a_i=i$, so
that $\tau \left( 1\right) =1$, $\tau \left( 2\right) =2$, $\tau \left(
3\right) =4$, $\tau \left( 4\right) =3$. We have 
\begin{equation*}
\boldsymbol{\tau }
\begin{pmatrix}
i \\ 
j
\end{pmatrix}
=
\begin{pmatrix}
\tau \left( j\right) \\ 
i
\end{pmatrix}
\end{equation*}
so we get 
\begin{list}{}{\setlength{\leftmargin}{\leftmargini}
\setlength{\itemsep}{0\customskipamount}
\setlength{\topsep}{\customskipamount}}
\item[ ]  $2$ fixed points: $i=j\in \left\{ 1,2\right\} $,
 
\item[ ]  $1$ period-$2$ orbit: $i\neq j,\enspace i,j\in \left\{ 1,2\right\} $,
 
\item[ ]  $3$ period-$4$ orbits: $\left\{ i,j\right\}  \cap 
\left\{ 3,4\right\} \neq \emptyset $.
\end{list}
This kind of analysis for (\ref{EqNew8.29}) gives 
\begin{equation*}
\boldsymbol{\tau }
\begin{pmatrix}
i \\ 
j \\ 
k
\end{pmatrix}
=
\begin{pmatrix}
\tau \left( k\right) \\ 
i \\ 
j
\end{pmatrix}
\end{equation*}
and hence there are 
\begin{list}{}{\setlength{\leftmargin}{\leftmargini}
\setlength{\itemsep}{0\customskipamount}
\setlength{\topsep}{\customskipamount}}
\item[ ]  $2$ fixed points: $i=j=k\in \left\{ 1,2\right\} $,
 
\item[ ]  $1$ period-$2$ orbit: $i=k\neq j,\enspace i,j\in 
\left\{ 3,4\right\} $,
 
\item[ ]  $2$ period-$3$ orbits: $i,j,k$ not all equal, $i,j,k 
\in \left\{ 1,2\right\} $,
 
\item[ ]  $9$ period-$6$ orbits: All other combinations.
\end{list}
Thus this example has the same optimal property as the case $N=2$,
$p=2^{12}-1$ mentioned in Section \ref{SecNew8}, and the range of
the map $\sigma :\mathbb{Z}%
^2\rightarrow \bigcross_{k=1}^\infty \mathbb{Z}_2$ consists of \emph{all}
sequences with a periodic tail with period $6$ (i.e. the minimal period
divides $6$). The example (\ref{EqNew8.30}) also has the optimal property,
and there the range of the map $\sigma :\mathbb{Z}^2\rightarrow \bigcross%
_{k=1}^\infty \mathbb{Z}_2$ consists of all sequences with a periodic tail
of period $4$.

More generally, if $\mathbf{N}$ is the $m\times m$ matrix (\ref{EqNew8.14})
with $N=2$, and $d_i$ as given there with $d_0=0$, $d_1=2^n-1$, then it
follows from Proposition \ref{ProNew8.3} that the atoms of the corresponding
one-dimensional system can be indexed by $\left( i_1,\dots ,i_n\right) $, $%
i_k\in \left\{ 0,1\right\} $ with 
\begin{equation*}
\tau \left( i_1,\dots ,i_n\right) =\left( i_n,i_1,\dots ,i_{n-1}\right) 
\text{.}
\end{equation*}
Thus, an atom of the corresponding $\mathbf{N}$-system is a sequence 
\begin{equation*}
\left( \left( i_{1,1},\dots ,i_{1,n}\right) ,\left( i_{2,1},\dots
,i_{2,n}\right) ,\dots ,\left( i_{m,1},\dots ,i_{m,n}\right) \right)
\end{equation*}
in $\mathbb{Z}_2^{nm}$ with 
\begin{multline*}
\boldsymbol{\tau }\left( \left( i_{1,1},\dots ,i_{1,n}\right) ,\left(
i_{2,1},\dots ,i_{2,n}\right) ,\dots ,\left( i_{m,1},\dots ,i_{m,n}\right)
\right) = \\
\left( \left( i_{m,n},i_{m,1},\dots ,i_{m,n-1}\right) ,\left(
i_{1,1},\dots ,i_{1,n}\right) ,\dots ,\left( i_{m-1,1},\dots
,i_{m-1,n}\right) \right) \text{.}
\end{multline*}
It follows that $\boldsymbol{\tau }^{nm}=\limfunc{id}$, and all the $2^{nm}$
atoms have period under $\boldsymbol{\tau }$ dividing $nm$. It follows from
Corollary \ref{CorNew6.2} that the image of $\mathbb{Z}^m$ under $\sigma $
consists of all sequences in $\bigcross_{k=1}^\infty \mathbb{Z}_2$ having a
periodic tail with period dividing $nm$. Combining this with Proposition \ref
{ProNew8.3} again, we deduce

\begin{corollary}
\label{CorNew8.4}The branching function system defined on $\mathbb{Z}$
by \textup{(%
\ref{Eq3.1})} with $N=2$, $s_0=0$, $s_1=2^{nm}-1$ is isomorphic to the
branching function system defined on $\mathbb{Z}^m$ by
\textup{(\ref{Eq3.6bis})}
with 
\begin{equation*}
\mathbf{N}=
\begin{pmatrix}
0 & 1 & 0 & 0 & \cdots & 0 & 0 \\ 
0 & 0 & 1 & 0 & \cdots & 0 & 0 \\ 
0 & 0 & 0 & 1 & \cdots & 0 & 0 \\ 
0 & 0 & 0 & 0 & \cdots & 0 & 0 \\ 
\vdots & \vdots & \vdots & \vdots & \ddots & \vdots & \vdots \\ 
0 & 0 & 0 & 0 & \cdots & 0 & 1 \\ 
N & 0 & 0 & 0 & \cdots & 0 & 0
\end{pmatrix}
\text{,\quad }d_0=
\begin{pmatrix}
0 \\ 
0 \\ 
0 \\ 
0 \\ 
\vdots \\ 
0 \\ 
0
\end{pmatrix}
\text{,\quad }d_1=
\begin{pmatrix}
0 \\ 
0 \\ 
0 \\ 
0 \\ 
\vdots \\ 
0 \\ 
2^n-1
\end{pmatrix}
\text{.}
\end{equation*}
In both cases the image of $\mathbb{Z}$, respectively $\mathbb{Z}^m$, under $%
\sigma $ consists of all sequences in $\bigcross_{k=1}^\infty \mathbb{Z}_2$
having a periodic tail with period dividing $nm$. In particular, the two
corresponding representations of $\mathcal{O}_2$ are unitarily equivalent.
\end{corollary}

Of course the consequence pertaining to existence of sub-Cuntz states for
the corresponding representation of $\mathcal{O}_2$ is the same as in
Proposition \ref{ProNew8.3} (the $n$ there should be replaced by any number
dividing $nm$).

By Corollary \ref{CorNew8.4}, there exists a unique bijection $\phi :%
\mathbb{Z}\rightarrow \mathbb{Z}^m$ with the property 
\begin{equation*}
\phi \circ \sigma _k=\boldsymbol{\sigma }_k\circ \phi
\end{equation*}
for $k=0,1$. Let us indicate an algorithm for finding this bijection by
again considering the case $n=1$, $m=2$, where the image of both systems
under $\sigma $ in $\bigcross_{k=1}^\infty \mathbb{Z}_2$ consists of all
sequences with a periodic tail with period dividing $2$. Thus there are two
fixed points for $R$, $\mathbf{R}$: 
\begin{equation*}
\begin{alignedat}{c}{2}
\sigma _0\left( 0\right) &=0 & \qquad \boldsymbol{\sigma
}_0\begin{pmatrix}0 \\ 0\end{pmatrix} &=\begin{pmatrix}0 \\ 0\end{pmatrix}
\\ \sigma _1\left( -3\right) &=-3 & \qquad 
\boldsymbol{\sigma }_1\begin{pmatrix}-1
\\ -1\end{pmatrix} &=\begin{pmatrix}-1 \\ -1\end{pmatrix} \end{alignedat}
\end{equation*}
and one orbit of period-$2$ points 
\begin{equation*}
\begin{alignedat}{c}{2}
\sigma _0\left( -1\right) &=-2 &\qquad \boldsymbol{\sigma
}_0\begin{pmatrix}0 \\ -1\end{pmatrix} &=\begin{pmatrix}-1 \\ 0\end{pmatrix}
\\ \sigma _1\left( -2\right) &=-1 &\qquad 
\boldsymbol{\sigma }_1\begin{pmatrix}-1
\\ 0\end{pmatrix} &=\begin{pmatrix}0 \\ -1\end{pmatrix} \end{alignedat}
\end{equation*}
and hence 
\begin{equation*}
\phi \left( 0\right) =
\begin{pmatrix}
0 \\ 
0
\end{pmatrix}
\text{, }\phi \left( -3\right) =
\begin{pmatrix}
-1 \\ 
-1
\end{pmatrix}
\text{, }\phi \left( -1\right) =
\begin{pmatrix}
-1 \\ 
0
\end{pmatrix}
\text{, }\phi \left( -2\right) =
\begin{pmatrix}
0 \\ 
-1
\end{pmatrix}
\text{.}
\end{equation*}
One can now compute $\phi \left( n\right) $ for any $n\in \mathbb{Z}$ by
using the intertwining relation 
\begin{equation*}
\phi \sigma _{i_1}\cdots \sigma _{i_k}m=\boldsymbol{\sigma }_{i_1}\cdots 
\boldsymbol{%
\sigma }_{i_k}\phi \left( m\right)
\end{equation*}
for $m=0,-1,-3,$ since any $n$ has an expansion 
\begin{equation*}
n=\sigma _{i_1}\cdots \sigma _{i_k}m
\end{equation*}
for one of these $m$'s. For example, using 
\begin{equation*}
\sigma _{i_1}\sigma _{i_2}\cdots \sigma _{i_{2k}}0=3\left( i_1+2i_2+\dots
+2^{2k-1}i_{2k}\right)
\end{equation*}
one computes 
\begin{equation*}
\phi \left( 3\left( i_1+2i_2+\dots +2^{2k-1}i_{2k}\right) \right) =3
\begin{pmatrix}
i_2+2i_4+\dots +2^{k-1}i_{2k} \\ 
i_1+2i_3+\dots +2^{k-1}i_{2k-1}
\end{pmatrix}
\end{equation*}
for $i_j\in \left\{ 0,1\right\} $.

\section{\label{Sec9}The general $\bmod{\,\mathbf{N}}$ situation\label
{Appendix}}

Let $\mathbf{N}$ be a $\nu \times \nu $ matrix with integer entries, and $%
\left| \det \mathbf{N}\right| =N$. If $N\neq 0$ and $D=\left\{ d_1,\dots
,d_N\right\} $ is a set of $N$ points in $\mathbb{Z}^\nu $ which are
incongruent modulo $\mathbf{N}\mathbb{Z}^\nu $, we pointed out in the
Introduction and around (\ref{Eq3.7}) that in order that a point $x\in %
\mathbb{Z}^\nu $ should have unique expansions 
\begin{equation}
x=d_{i_1}+\mathbf{N}d_{i_2}+\dots +\mathbf{N}^{k-1}d_{i_k}%
\mod{\mathbf{N}^k\mathbb{Z}^\nu}  \label{Eq9.1}
\end{equation}
it is necessary and sufficient that 
\begin{equation}
\bigcap_{k=1}^\infty \mathbf{N}^k\mathbb{Z}^\nu =0\text{.}  \label{Eq9.2}
\end{equation}
In order to prove finiteness of the $\approx $-equivalence classes in
Corollary \ref{Cor3.9} we needed a stronger condition: 
\begin{equation}
\left| \lambda \right| >1  \label{Eq9.3}
\end{equation}
for all the (complex) eigenvalues of $\mathbf{N}$. The next result show that
this condition is strictly stronger than (\ref{Eq9.2}), so that, for
example, $\mathbf{N}=\left( 
\begin{smallmatrix}
3 & 1 \\
1 & 1
\end{smallmatrix}\right) $ satisfies (\ref{Eq9.2}) but not (\ref{Eq9.3}). To
formulate the result, recall from \cite{LeV96}, \cite{And94} that the ring $%
\mathbb{Z}\left[ \lambda \right] $ of all polynomials in $\lambda $ with
integer coefficients is a unique factorization domain, i.e., any polynomial
in $\mathbb{Z}\left[ \lambda \right] $ has a unique factorization (up to the
units $\pm 1$) in prime factors, up to the order of the factors. $\mathbb{Z}%
\left[ \lambda \right] $ is contained in the Euclidean domain $\mathbb{Q}%
\left[ \lambda \right] $, but is not itself Euclidean. However, the prime
factorization of a polynomial in $\mathbb{Z}\left[ \lambda \right] $ is also
a prime factorization in $\mathbb{Q}\left[ \lambda \right] $. A prime
polynomial in $\mathbb{Z}\left[ \lambda \right] $ will also be referred to
as an irreducible polynomial. The following proposition is also true for
integer matrices $\mathbf{N}$ with $\det \mathbf{N}=0$, so it is formulated
in that generality.

\begin{proposition}[David Handelman]
\label{Pro9.1}Let $\mathbf{N}$ be a $\nu \times \nu $ matrix with integer
entries. The following conditions are equivalent:

\begin{enumerate}
\item  \label{Pro9.1(i)}$\vphantom{\Bigl(}\displaystyle\bigcap\limits_{k=1}^%
\infty \mathbf{N}^k\mathbb{Z}^\nu =0$.

\item  \label{Pro9.1(ii)}If $\vphantom{\smash[b]{\Bigl(}}f$ is any prime
factor of the monic polynomial
$\det \left( \lambda \openone -\mathbf{N}\right) $ in 
$\mathbb{Z}\left[ \lambda \right] $, then $\left| f\left( 0\right) \right|
\neq 1$.
\end{enumerate}
\end{proposition}

\begin{remark}
The Proposition should be supplemented with the recent observation of
of Bandt and Gelbrich \cite{BaGe94} which states that, if the given
integral $\nu \times \nu $ matrix $\mathbf{N}$ in fact satisfies
(\ref{Eq9.3}), and if $N=\left| \det \mathbf{N}\right| $ is a prime
number, then the monic polynomial $\lambda \mapsto \det \left( \lambda
\openone -\mathbf{N}\right) $ is irreducible. Of course, in that case,
(ii) is satisfied, but for the obvious reason.

It was also noted by Gelbrich \cite{Gel96} that, if the number
$N=\left| \det \mathbf{N}\right| $ is fixed, but the matrix
$\mathbf{N}$ varies subject to (\ref{Eq9.3}), then the coefficients
$p_j=p_j\left( \mathbf{N}\right) $ in
\begin{equation*}
\det \left( \lambda \openone -\mathbf{N}\right) =\lambda ^\nu +p_{\nu 
-1}\lambda ^{\nu -1}+\dots +p_1\lambda \pm N\in \mathbb{Z}\left[ \lambda 
\right]
\end{equation*}
must satisfy $\left| p_j\right| \leq \binom {\nu }jN^{\nu -j}$. Hence,
for each fixed $N$ and $\nu $, there is only a finite number of
conjugacy classes of matrices $\mathbf{N}$ satisfying (\ref{Eq9.3}).
\end{remark}

\begin{proof}[Proof of Proposition \ref{Pro9.1} 
\textup{(}by David Handelman\/\textup{)}]Let $W=\mathbb{Z}^\nu $. We first
make the well known observation that if $A$ is a $\nu \times \nu $ integer
matrix and $Y$ is an $A$-invariant subgroup of $W$ such that $W\diagup Y$ is
torsion free, then there is a $\mathbb{Z}$-basis of $W$ with respect to which
the matrix of $A$ acting on columns is in block upper triangular form, with $%
A^{\prime \prime }$ acting as $A$ restricted to $Y$, 
\begin{equation*}
\begin{pmatrix}
A^{\prime } & X \\ 
0 & A^{{\prime \prime }}
\end{pmatrix}
\text{.}
\end{equation*}
(Simply complete a basis for $Y$ to a basis for $W$.)

Set $V=\bigcap_k\,\mathbf{N}^kW$. This is a sublattice of $W$, in particular
a free abelian group of rank at most $\nu $, and it is $\mathbf{N}$%
-invariant.

(\ref{Pro9.1(ii)})$\Rightarrow $(\ref{Pro9.1(i)}). Assume $V$ is not zero.
We will show that (\ref{Pro9.1(ii)}) fails under the additional assumption
that $\det \mathbf{N}\neq 0$, and then that it fails without this assumption.

Since $\det \mathbf{N}$ is not zero, $\mathbf{M}:V\to V$, defined as the
restriction of $\mathbf{N}$, is one-to-one, and it is clearly onto (for $z$
in $V$, for each $i$, there exists $w_i$ in $W$ such that $\mathbf{N}^iw_i=z$%
; as $\mathbf{N}$ is one-to-one, for all $i>1$, $\mathbf{N}^{i-1}w_i$ are
all equal to the same element, call it $w$; then $w$ belongs to $V$ and $%
\mathbf{N}w=z$). Hence $\left| \det \mathbf{M}\right| =1$. If $g$ is the
characteristic polynomial of $\mathbf{N}$ (in $\mathbb{Z}\left[ 
\lambda \right] $),
then the characteristic polynomial of $\mathbf{M}$, call it $h$, must
divide $g$. Since $\left| \det \mathbf{M}\right| =\left| h(0)\right| =1$, it
follows that for any irreducible $f$ dividing $h$, $\left| f(0)\right| =1$;
thus (\ref{Pro9.1(ii)}) fails.

Now drop the assumption that $\det \mathbf{N}$ is nonzero. If $\mathbf{N}$
is nilpotent, the result is trivial. Set $K$ to be the kernel of $\mathbf{N}%
^n$. As is well known, $K$ is the kernel of $\mathbf{N}^i$ for all $i\geq n-1
$, and of course, $W\diagup K$ is torsion free, and the endomorphism induced
on $W\diagup K$ by $\mathbf{N}$ is one-to-one. Applying the block upper
triangular block form from $K$, we obtain that $\mathbf{N}$ can be rewritten
in the form, 
\begin{equation*}
B:=
\begin{pmatrix}
A^{\prime } & X \\ 
0 & A^{{\prime \prime }}
\end{pmatrix}
\text{,}
\end{equation*}
where $A^{{\prime \prime }}$ is nilpotent; this forces
$(A^{{\prime \prime }})^n=0$. Then for all $k\geq n$, 
\begin{equation*}
B^k=
\begin{pmatrix}
(A^{\prime })^k & \sum_{i=0}^{n-1}(A^{\prime })^{k-i}X(A^{{\prime \prime }})^i
\\ 
0 & 0
\end{pmatrix}
\text{.}
\end{equation*}
Thus an element of $B^kW$ (for $k\geq n$) is of the form
$(A^{\prime })^{k-n}(w_k)$
for $w_k=(A^{\prime })^na_k+\sum_{j=1}^{n-1}{(A^{\prime })}%
^jX(A^{{\prime \prime }})^{n-j}b_k$ in the domain of $A^{\prime }$. Thus if $%
\bigcap B^kW$ is not zero then neither is its counterpart for $A^{\prime }$. By
the earlier result, (\ref{Pro9.1(ii)}) fails for $A^{\prime }$, and since
all the irreducible factors of the characteristic polynomial $A^{\prime }$
divide the characteristic polynomial of $\mathbf{N}$, we are done.

(\ref{Pro9.1(i)})$\Rightarrow $(\ref{Pro9.1(ii)}). Suppose there is an
irreducible factor $f$ of the characteristic polynomial $g$ of $\mathbf{N}$
such that $f(0)=\pm 1$. Set $Y$ to be the kernel of $f(\mathbf{N})$; then
restricted to $Y$, the minimal polynomial of $\mathbf{N}|_Y$ divides $f$ (so
must be $f$), and so $\mathbf{N}|_Y$ is invertible. It easily follows that $%
Y\subseteq V$, so $V$ is not zero.\TeXButton{End Proof}{\end{proof}}

(Since $\mathbb{Z}\left[ \lambda \right] $ is not a Euclidean field, it is
not possible to eliminate the blocks over the diagonal, but this does not
matter for the argument. It would be possible over $\mathbb{Q}\left[ \lambda
\right] $, but even then one could not of course write the matrix in normal
Jordan form as a genuine upper triangular matrix as over $\mathbb{C}\left[
\lambda \right] $.)

We next show that Corollary \ref{Cor3.9} does not extend to the situation
where some of the eigenvalues have modulus less than $1$.

\begin{proposition}
\label{Pro9.2}Let $\mathbf{N}$ be a $\nu \times \nu $ symmetric matrix with
integer entries. Assume that $\left| \det \left( \mathbf{N}\right) \right|
=N\geq 2$ and $\bigcap_k\mathbf{N}^k\mathbb{Z}^\nu =0$. Assume that $\mathbf{%
N}$ has at least one \textup{(}real\/\textup{)} eigenvalue of modulus
less than one. Let $%
D=\left\{ d_1,\dots ,d_N\right\} $ be a set of $N$ points in $\mathbb{Z}^\nu 
$ which are pairwise incongruent modulo $\mathbf{N}\mathbb{Z}^\nu $. Then
the associated permutative mul\-ti\-plic\-i\-ty-free representation of $%
\mathcal{O}_N$ has an infinite number of cycles with an infinite number of
atoms each, and at most finitely many cycles with a finite number of atoms
each \textup{(}the number of atoms of finite period being estimable
in terms of $%
\limfunc{const.}\cdot\; d^\nu $, where
$d=\max_i\left\| d_i\right\| $ and $\limfunc{const.}$ depends on $\mathbf{N}$
alone\/\textup{).}
\end{proposition}

\begin{remark}
\label{Hyperbolic}{\bfseries Hyperbolic transformations.}
It follows from the proof that the requirement that $\mathbf{N}$ is
symmetric can be relaxed somewhat. One sufficient condition for the
conclusion of the proposition is for example that $\mathbb{R}^\nu $ has a
linear decomposition $\mathbb{R}^\nu =V_1\oplus V_2$ into two nonzero $%
\mathbf{N}$-invariant subspaces $V_1$ and $V_2$ such that there is an
equivalent norm $\left\| \,\cdot \,\right\| $ on $\mathbb{R}^\nu $ and an $%
\varepsilon >0$ such that $\left\| \mathbf{N}x\right\| \leq 
\left( 1-\varepsilon
\right) \left\| x\right\| $ for $x\in V_1$ and $\left\| \mathbf{N}x\right\|
\geq \left( 1+\varepsilon \right) \left\| x\right\| $ for $x\in V_2$.
That is, the
(complex) eigenvalues of $\mathbf{N}|_{V_1}$ have all modulus less than one,
and those of $\mathbf{N}|_{V_2}$ have modulus larger than one. The proof is
the same as the present one: $R$ behaves like $\mathbf{N}^{-1}$ on a large
scale, so $V_1$ and $V_2$ are approximately ``unstable and stable
manifolds'' for $R$, except near $0$. (This
property is quite common from the theory of dynamical systems,
where it is called
hyperbolicity; see, e.g., \cite{PaTa95} and \cite{Rue94}.) .
\end{remark}

\TeXButton{Begin Proof}
{\begin{proof}[Proof of Proposition \ref{Pro9.2}]}The condition
$\bigcap_k\mathbf{N}^k%
\mathbb{Z}^\nu =0$ implies by Proposition \ref{Pro9.1} that $\pm 1$ cannot
be eigenvalues of $\mathbf{N}$. Since $\mathbf{N}$ is symmetric, it has a
spectral decomposition, so there is an orthonormal basis $\psi _1,\dots
,\psi _N$ for $\mathbb{R}^\nu $ and eigenvalues $\lambda _1,\dots ,\lambda
_\nu $ such that 
\begin{equation}
\mathbf{N}\psi _i=\lambda _i\psi _i\text{.}  \label{Eq9.4}
\end{equation}
Order the indices such that 
\begin{equation}
\left| \lambda _1\right| \leq \left| \lambda _2\right| \leq \dots \leq
\left| \lambda _\mu \right| <1<\left| \lambda _{\mu +1}\right| \leq \dots
\leq \left| \lambda _\nu \right|   \label{Eq9.5}
\end{equation}
and if $x\in \mathbb{Z}^\nu \subseteq \mathbb{R}^\nu $, define 
\begin{equation}
x_i=\left\langle x,\psi _i\right\rangle \text{.}  \label{Eq9.6}
\end{equation}
The map $R:\mathbb{Z}^\nu \rightarrow \mathbb{Z}^\nu $ is defined such that 
\begin{equation}
R\left( d_i+\mathbf{N}y\right) =y\text{.}  \label{Eq9.7}
\end{equation}
If $x=d_i+\mathbf{N}y$, then $y=\mathbf{N}^{-1}\left( x-d_i\right) $ so 
\begin{equation}
R\left( x\right) =\mathbf{N}^{-1}x+r\left( x\right)   \label{Eq9.8}
\end{equation}
where 
\begin{equation}
r\left( x\right) =-\mathbf{N}^{-1}\left( d_i\right)   \label{Eq9.9}
\end{equation}
is a remainder term depending on $x$. We have 
\begin{equation}
\left\| r\left( x\right) \right\| \leq \left\| \mathbf{N}^{-1}\right\|
d=\Delta \text{.}  \label{Eq9.10}
\end{equation}
If $x\in \mathbb{Z}^\nu $, we have by Corollary \ref{Cor2.4} to consider the
forward orbit $n\rightarrow R^n\left( x\right) $ of $x$. By (\ref{Eq9.8}),
we have 
\begin{equation}
R\left( x\right) _i=\lambda _i^{-1}x_i+r\left( x\right) _i  \label{Eq9.11}
\end{equation}
where 
\begin{equation}
-\Delta \leq r\left( x\right) _i\leq \Delta \text{.}  \label{Eq9.12}
\end{equation}
Consider first the case $i>\mu $. Then $\left| \lambda _i\right| >1$ and $%
\left| \lambda _i\right| ^{-1}<1$. Hence from (\ref{Eq9.11}) and (\ref
{Eq9.12}) we have 
\begin{equation*}
\limsup_{n\rightarrow \infty }\left| R^n\left( x\right) _i\right| \leq
\left| \lambda _i\right| ^{-1}\limsup_{n\rightarrow \infty }\left|
R^{n-1}\left( x\right) _i\right| +\Delta 
\end{equation*}
so 
\begin{equation}
\limsup_{n\rightarrow \infty }\left| R^n\left( x\right) _i\right| \leq
\Delta \left( 1-\left| \lambda _i\right| ^{-1}\right) ^{-1}\text{.}
\label{Eq9.13}
\end{equation}
This means that $R^n\left( x\right) $ is attracted to the ``hibachi'' \cite
{Str95} 
\begin{equation}
\left\{ y\in \mathbb{Z}^\nu \Bigm|\left| y_i\right| <\Delta \left( 1-\left|
\lambda _i\right| ^{-1}\right) ^{-1}+\varepsilon \text{,\quad }i=\mu +1,\dots 
,\nu
\right\}   \label{Eq9.14}
\end{equation}
as $n\rightarrow \infty $, i.e., $R^n\left( x\right) $ ultimately is
contained in the hibachi for each $\varepsilon >0$.

Next consider the case $i\leq \mu $. Then $\left| \lambda _i\right| <1$ and $%
\left| \lambda _i\right| ^{-1}>1$, and thus from (\ref{Eq9.11}) and (\ref
{Eq9.12}), 
\begin{equation}
\left| R^n\left( x\right) _i\right| \geq \left| \lambda _i\right|
^{-1}\left| R^{n-1}\left( x\right) _i\right| -\Delta \text{.}  \label{Eq9.15}
\end{equation}
From this we observe that if 
\begin{equation}
\left| R^{n-1}\left( x\right) _i\right| >\Delta \left( \left| \lambda
_i\right| ^{-1}-1\right)   \label{Eq9.16}
\end{equation}
for some $n$, then 
\begin{equation}
\left| R^n\left( x\right) _i\right| >\left| R^{n-1}\left( x\right) _i\right| 
\label{Eq9.17}
\end{equation}
and hence $\left| R^n\left( x\right) _i\right| $ will grow with $n$ for all
further $n$. Furthermore, from (\ref{Eq9.15}), this growth will
asymptotically be faster than $\left(
\vphantom{\left| \lambda \right| }\smash{\left| \lambda _i\right|
^{-1}-\varepsilon }\right) ^n$ for any $\varepsilon >0$. Thus $\limfunc{Per}%
\left( x\right) =\infty $, and $x^{\approx }$ will be contained in a cycle
with infinitely many points. Also, as $\left| R^n\left( x\right) _i\right|
-\left| R^{n-1}\left( x\right) _i\right| \rightarrow \infty $ as $%
n\rightarrow \infty $ we may find another $x$ giving rise to a disjoint $R$%
-orbit, and another cycle with infinitely many atoms. But distances between
these orbits will also grow exponentially, so one may find a third disjoint
orbit, and so on ad infinitum. This proves the first statement of
Proposition \ref{Pro9.2}. But the second statement follows from the fact
that if $x$ is part of a periodic orbit of $R$, then $x$ must be contained
in 
\begin{equation*}
\left\{ y\in \mathbb{Z}^\nu \Bigm|\left| y_i\right| \leq \Delta \left(
1-\left| \lambda _i\right| ^{-1}\right) ^{-1}\text{,\quad }i=\mu +1,\dots ,\nu
\right\} 
\end{equation*}
by (\ref{Eq9.14}), but $x$ also must be contained in 
\begin{equation*}
\left\{ y\in \mathbb{Z}^\nu \Bigm|\left| y_i\right| \leq \Delta \left(
\left| \lambda _i\right| ^{-1}-1\right) \text{,\quad }i=1,\dots ,\mu \right\} 
\end{equation*}
by (\ref{Eq9.16})--(\ref{Eq9.17}), otherwise $R^nx$ vanishes to infinity.
Thus $x$ must be contained in 
\begin{equation}
\left\{ y\in \mathbb{Z}^\nu \Biggm|
\begin{split}
& \left| y_i\right| \leq \Delta \left( 1-\left| \lambda _i\right|
^{-1}\right) ^{-1}\text{ for }i=\mu +1,\dots ,\nu  \\
& \text{ and }\left| y_i\right| \leq \Delta \left( \left| \lambda _i\right|
^{-1}-1\right) ^{-1}\text{ for }i=1,\dots ,\mu 
\end{split}
\right\}   \label{Eq9.18}
\end{equation}
But this box contains at most a finite number of points, so there are at
most a finite number of periodic orbits for $R$, and the second statement in
Proposition \ref{Pro9.2} is proved.\TeXButton{End Proof}{\end{proof}}

\begin{example}
\label{Insufficient}Let us give an example of a $2\times 2$ matrix $\mathbf{N%
}$ such that the conditions of Proposition \ref{Pro9.2} are fulfilled, but
not those of Corollary \ref{Cor3.9}. One such matrix is 
\begin{equation}
\mathbf{N}=
\begin{pmatrix}
3 & 1 \\ 
1 & 1
\end{pmatrix}
\text{,}  \label{Eq9.19}
\end{equation}
which has the eigenvalues and corresponding orthonormal set of eigenvectors 
\TeXButton{Equation 9.20}{} 
\begin{equation}
\begin{alignedat}{c}{2}
\lambda _1 &=2+\sqrt2 & \text{,\qquad} \psi _1
&=\frac1{\sqrt{2\left( 2+\sqrt2\right) }}\begin{pmatrix}\sqrt2+1 \\
1\end{pmatrix}\text{,} \\ \lambda _2 &=2-\sqrt2 & \text{,\qquad} \psi _2
&=\frac1{\sqrt{2\left( 2+\sqrt2\right) }}\begin{pmatrix}-1 \\
\sqrt2+1\end{pmatrix}\text{.} \end{alignedat}  \label{Eq9.20}
\end{equation}
In this case one computes that 
\begin{equation*}
\begin{pmatrix}
x_1 \\ 
y_1
\end{pmatrix}
\approx 
\begin{pmatrix}
x_2 \\ 
y_2
\end{pmatrix}
\mod{\mathbf{N}}
\end{equation*}
if and only if $x_1-y_1$ and $x_2-y_2$ have the same parity, so one may take 
\begin{equation}
D=\left\{ 
\begin{pmatrix}
0 \\ 
0
\end{pmatrix}
,
\begin{pmatrix}
0 \\ 
1
\end{pmatrix}
\right\} \text{.}  \label{Eq9.21}
\end{equation}
A little computation then shows 
\begin{equation}
R
\begin{pmatrix}
x \\ 
y
\end{pmatrix}
=\mathbf{N}^{-1}
\begin{pmatrix}
x \\ 
y
\end{pmatrix}
+\frac{p\left( x+y\right) }2
\begin{pmatrix}
1 \\ 
-3
\end{pmatrix}
\text{,}  \label{Eq9.22}
\end{equation}
where $p$ is the parity function.

We have 
\begin{equation}
\frac 12
\begin{pmatrix}
1 \\ 
-3
\end{pmatrix}
=\frac{1-\sqrt{2}}{2\sqrt{2+\sqrt{2}}}\psi _1-\frac{3+2\sqrt{2}}{%
2\sqrt{2+\sqrt{2}}}\psi _2=-a\psi _1-b\psi _2\text{.}  \label{Eq9.23}
\end{equation}
If now $x\in \mathbb{Z}^2$, define $a_0\left( x\right) $, $b_0\left(
x\right) $ by 
\begin{equation}
x=a_0\left( x\right) \psi _1+b_0\left( x\right) \psi _2  \label{Eq9.24}
\end{equation}
and $a_n\left( x\right) $, $b_n\left( x\right) $ by 
\begin{equation}
R^n\left( x\right) =a_n\left( x\right) \psi _1+b_n\left( x\right) \psi _2
\label{Eq9.25}
\end{equation}
and $p_n\left( x\right) $ as the parity of the sum of the two components of $%
R^n\left( x\right) $. Then 
\begin{equation*}
a_{n+1}\left( x\right) \psi _1+b_{n+1}\left( x\right) \psi _2=\mathbf{N}%
^{-1}\left( a_n\left( x\right) \psi _1+b_n\left( x\right) \psi _2\right)
+p_n\left( x\right) \left( -a\psi _1-b\psi _2\right) 
\end{equation*}
so 
\begin{equation}
\begin{aligned} a_{n+1}\left( x\right) &=\frac{a_n\left( x\right) }{\lambda
_1}-p_n\left( x\right) a\text{ and } \\ b_{n+1}\left( x\right) 
&=\frac{b_n\left(
x\right) }{\lambda _2}-p_n\left( x\right) b\text{.} \end{aligned}  
\label{Eq9.26}
\end{equation}
Let 
\begin{equation}
\begin{aligned} \overline{a}\left( x\right) &=\limsup_{n\rightarrow \infty
}a_n\left( x\right)\text{,} \\ \underline{a}\left( x\right)
&=\liminf_{n\rightarrow \infty }a_n\left( x\right)\text{.} \end{aligned}
\label{Eq9.12bis}
\end{equation}
We compute, from the recursion relations, 
\begin{equation}
\begin{aligned} \overline{a}\left( x\right) &\leq\frac{\overline{a}\left(
x\right) }{\lambda _1}\text{,} \\ \underline{a}\left( x\right)
&\geq\frac{\underline{a}\left( x\right) }{\lambda _1}-a\text{.} \end{aligned}
\label{Eq9.13bis}
\end{equation}
As $\lambda _1>1$, we thus obtain 
\begin{equation}
\begin{aligned} \overline{a}\left( x\right) & \leq 0\text{,} \\
\underline{a}\left( x\right) & \geq -a\left( 1-\frac 1{\lambda _1}\right)
^{-1}=-\frac 1{\left( 2+\sqrt{2}\right) ^{\frac 32}}\text{.} \end{aligned}
\label{Eq9.14bis}
\end{equation}
The first two of these inequalities mean that for all $x\in \mathbb{Z}^2$
and all $\varepsilon >0$, there is an $n_x\in \mathbb{N}$ such that if $n>n_x$,
then $R^n\left( x\right) $ is contained in the ``hibachi'' \cite{Str95} 
\begin{equation*}
\left\{ y\in \mathbb{Z}^2\biggm|-\frac 1{\left( 2+\sqrt{2}\right) ^{\frac
32}}-\varepsilon <\left\langle y,\psi _1\right\rangle <\varepsilon \right\} 
\text{,%
}
\end{equation*}
which is a narrow strip in the direction of $\psi _2$. This strip thus acts
as an attractor for the dynamical system defined by $R$. For vectors in this
hibachi which are far from the origin, $R$ acts approximately by
multiplication by $\frac 1{\lambda _2}=\frac 1{2-\sqrt{2}}$. More precisely,
from (\ref{Eq9.26}) one deduces 
\begin{equation*}
b_n\left( x\right) >b\left( \frac 1{\lambda _2}-1\right) 
^{-1}=\frac{3+2\sqrt{2}}{%
\sqrt{2\left( 2+\sqrt{2}\right) }}\Longrightarrow b_{n+1}\left( 
x\right) >b\left( \frac 1{\lambda
_2}-1\right) ^{-1}
\end{equation*}
and 
\begin{equation*}
b_n\left( x\right) <0\Longrightarrow b_{n+1}\left( x\right) <0\text{,}
\end{equation*}
so if for some $n$, $b_n\left( x\right) $ is outside the interval $\left[ 0,%
\frac{3+2\sqrt{2}}{%
\sqrt{2\left( 2+\sqrt{2}\right) }}\right] $, then $b_{n+1}\left( 
x\right) $ is also outside
this interval, and then $b_n\left( x\right) $ grows asymptotically like $%
\left( \lambda _2^{-1}-\varepsilon \right) ^n$.

It follows that there are infinitely many non-periodic tails of the
sequences $R^n\left( x\right) $ growing exponentially inside the hibachi,
and because of this exponential growth there must be infinitely many tails
up to translation. Thus there are infinitely many cycles with infinitely
many atoms. But in addition there might be a finite number of periodic tails
contained in the box 
\begin{figure}
\begin{center}
\ \psfig{file=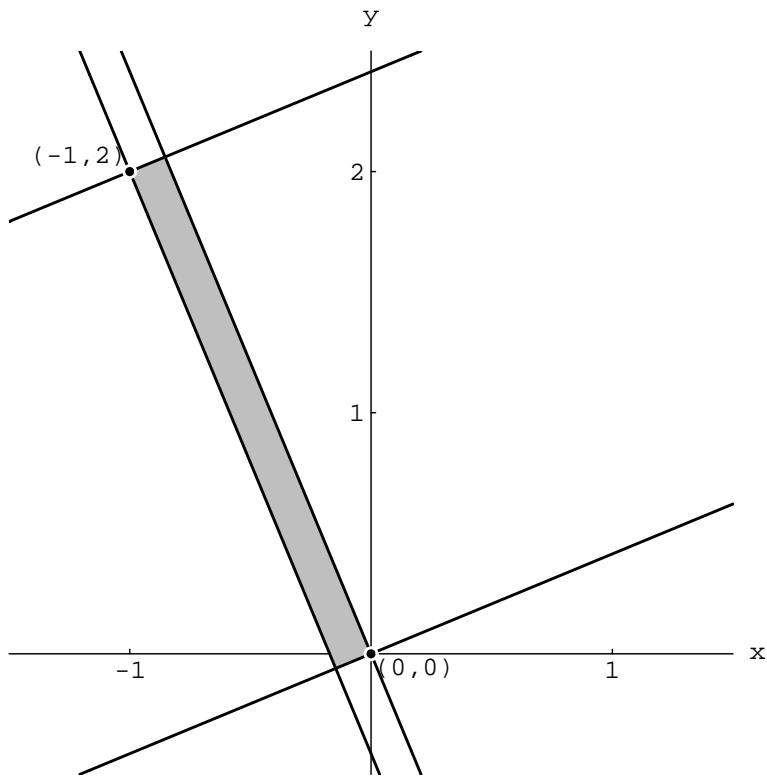}
\end{center}
\caption{Region containing periodic tails
(Example \protect \ref{Insufficient})}
\label{TailBox}
\end{figure}
\begin{align*}
&\left\{ y \in \mathbb{Z}^2\biggm|-\frac 1{\left( 2+\sqrt{2}\right) ^{\frac
32}}\leq \left\langle y,\psi _1\right\rangle \leq 0\text{ and }0\leq
\left\langle y,\psi _2\right\rangle \leq \frac{3+2\sqrt{2}}{%
\sqrt{2\left( 2+\sqrt{2}\right) }}\right\}  \\
&\qquad =\left\{ 
\begin{pmatrix}
x \\ 
y
\end{pmatrix}
\in \mathbb{Z}^2\Biggm| 
\begin{aligned}
\vphantom{\ }&1-\sqrt 2 \leq \left(\sqrt{2}+1\right)
x+y\leq 0 \\
\vphantom{\ }&\text{and }0 \leq -x+\left(\sqrt{2}+1\right)
y\leq 3+2\sqrt{2}
\end{aligned}
\right\} \\
&\qquad =\left\{ 
\begin{pmatrix}
0 \\ 
0
\end{pmatrix}
, 
\begin{pmatrix}
-1 \\ 
2
\end{pmatrix}
\right\} \text{.}
\end{align*}
See Figure \ref{TailBox}. But $R\left( \begin{smallmatrix}
0 \\ 
0
\end{smallmatrix}
\right) =\left( \begin{smallmatrix}
0 \\ 
0
\end{smallmatrix}
\right) $ and $R\left( \begin{smallmatrix}
-1 \\ 
2
\end{smallmatrix}
\right) =\left( \begin{smallmatrix}
-1 \\ 
2
\end{smallmatrix}
\right) $, thus this gives rise to two cycles with one atom
each. Conclusion: the representation associated with (\ref{Eq9.19})
and (\ref{Eq9.21}) has two cycles with one atom each, infinitely many
cycles with infinitely many atoms each, and no more cycles.
\end{example}

\section{\label{Sec11}Concluding remarks\label{Concluding}}

When finishing this paper, we became aware of the independent paper
\cite{DaPi96} with some overlapping results. Recall that the Toeplitz
algebra $\mathcal{T}_N$, for $N<\infty $, has $\mathcal{O}_N$ as a
quotient. We give a brief discussion of this below. For background
references on $\mathcal{T}_N$, $\mathcal{O}_N$, and the corresponding
representations, we use \cite{BEGJ}.

The paper \cite{DaPi96} describes representations of $\mathcal{O}_N$
which are called atomic representations (see next paragraph), and are
more general than our permutative ones: they include those of Remark
\ref{PermRepNotReg} and are explained as central extensions (see,
e.g., \cite{JoWe94}). As a result they get a general decomposition
theory for $\mathcal{O}_N$ but not for
$\limfunc{UHF}\nolimits_N$. Specifically, the results Theorem 6.4 and
Proposition 6.10 in \cite{DaPi96} show that, when the restricting
assumptions ``regular'' or ``mul\-ti\-plic\-i\-ty-free'' in our
Theorem \ref{Thm2.5} are removed, then the resulting orthogonal
representation decompositions get slightly more refined: there can be
continuous parts, or extra representations arising from central
extensions; these latter results are in turn closely related to
\cite{Jor96} and \cite{Rue88}.

Atomic representations of the Cuntz algebra $\mathcal{O}_N$ (or rather
the Toeplitz algebra $\mathcal{T}_N$) are representations such that
there exists an orthonormal basis $\left\{ e_n\right\} _{n=1}^\infty $
for the Hilbert space, maps $\sigma _k:\mathbb{N}\rightarrow
\mathbb{N}$ satisfying (\ref{Eq2.4a}) and (\ref{Eq2.4b}) and scalars
$\lambda _{i,k}\in \mathbb{T}$ such that
\begin{equation}
S_ie_k=\lambda _{i,k}e_{\sigma _i\left( k\right) }
\label{Eq11.1}
\end{equation}
(compare (\ref{Eq2.2})). In \cite{DaPi96}, a complete classification
up to unitary equivalence of irreducible representations of this kind
is given. Except for the left regular representation of the free
semigroup with $N$ generators, which corresponds to the unique
faithful irreducible representation of $\mathcal{T}_N$ (see, e.g.,
\cite{BEGJ}), all of these are representations of $\mathcal{O}_N$. In
\cite{DaPi96}, the representations are described in a different
language than ours, but translated to our setting, the classification
is as follows: for each $x\in \Omega =\bigcross _{k=1}^\infty
\mathbb{Z}_N$, and each $\lambda \in \mathbb{T}$, define a
representation on $\ell ^2\left( x^\sim \right) $ by $S_k\delta
_y=\lambda \delta _{ky}$ (see Section \ref{Universal}). This defines
an irreducible representation and, if $x^{\prime } \in \Omega $ and
$\lambda ^{\prime } \in \mathbb{T}$, the corresponding representations
are unitarily equivalent if and only if (\ref{Sec11(i)}) and
(\ref{Sec11(ii)}) below hold:
\begin{enumerate}
\item $x^{\prime } \sim x$ (i.e., $x^{\prime } $, $x$ have the same
tail up to translation).\label{Sec11(i)}

\item If $\limfunc{Per}%
\left( x\right) =m<+\infty $, then $\lambda ^m={\lambda ^{\prime } }^m$.
\label{Sec11(ii)}
\end{enumerate}
All atomic irreducible representations are of this form, and if $R$
has a finite minimal cyclic orbit $n_1,n_2,\dots ,n_m,n_{m+1}=n_1$,
and $j_1,j_2,\dots ,j_m,j_{m+1}=j_1$ is the corresponding coding
sequence then $\lambda $ can be taken to be any of the $m$ roots of
\begin{equation*}
\lambda ^m=\lambda _{i_1,n_1}\lambda _{i_2,n_2}\cdots \lambda _{i_m,n_m}
\text{.}
\end{equation*}
See \cite[Theorem 6.4]{DaPi96}.

In our setting of mul\-ti\-plic\-i\-ty-free function systems, we
obtained a discrete decomposition of any permutative representation
into mutually disjoint irreducible representations in Theorem
\ref{Thm2.5}. In the more general case of atomic representations,
there are some subrepresentations with a continuous decomposition
too. Without describing these in detail (see \cite[Proposition
6.10]{DaPi96}), a partial explanation of this phenomenon is as
follows: the decomposition is over the parameter $\lambda $ considered
above, but since the coding map is defined in terms of the maps
$\sigma _i$ alone, our condition of injectivity of the coding map
prevents two different subrepresentations with the same $\sigma _i$'s
(but different $\lambda $'s) from occurring. It should however be
noted that in the extension of Theorem \ref{Thm2.5} to general (not
necessarily mul\-ti\-plic\-i\-ty-free) branching function systems, one
must consider atomic representations with nontrivial $\lambda $'s in
the continuous part of the decomposition of a permutative
representation.

The present proof of Theorem \ref{Thm2.5} extends with obvious
modifications to the context of atomic representations with injective
coding maps, giving a discrete decomposition into the irreducible
representations described above, with only one $\lambda $ for each
$\sigma _i$-class.

Note that, if $x\in \Omega $ is asymptotically periodic with period
$m$, $I$ is a string of length $m$ from the tail of $x$, and $\lambda
\in \mathbb{T}$, then the irreducible representation described above
comes from the sub-Cuntz state determined by
\begin{equation*}
\xi =\lambda ^{-m}S_I\xi \text{;}
\end{equation*}
see Proposition \ref{Pro5.1}.

Although the overlap of \cite{DaPi96} and the present paper has been
emphasized above, the two papers are in other respects rather
different. The Gelfand pair aspect of the pair $\left(
\mathcal{O}_N,\limfunc{UHF}\nolimits_N\right) $ is only treated here,
as well as the connection to fractal properties of the dynamical
systems coming from the $\bmod{\,\mathbf{N}}$ situation.

As a final comment, we introduce operators that intertwine the two
pictures.  In \cite{DaPi96} the following class of representations of
$\mathcal{O}_N$ is considered: for $x=\left( j_1,j_2,\dots \right) \in
\Omega $ fixed, define $x_m$ as the word $\left( j_1,\dots ,j_m\right) $
and define $\mathcal{F}_nx^{-1}$ as the (finite) words in the free
\emph{group} generated by the $N$ symbols in $\mathbb{Z}_N$ of the
form $ux_m^{-1}$ for some $u\in \mathcal{F}_m$, where $\mathcal{F}_m$
is the free \emph{semigroup} generated by the $N$ symbols. If $\xi _v$,
$v\in \mathcal{F}_nx^{-1}$ is an orthonormal basis, the
representation is defined by
\begin{equation*}
S_i\xi _v=\xi _{iv}\text{.}
\end{equation*}
Now, note that the operation 
\begin{equation*}
ux_m^{-1}\rightarrow ux_m^{-1}x
\end{equation*}
defines an isometric intertwiner between this representation and our
representation defined in Section \ref{SecNew6} by the set of $y\in
\Omega $ having the same tail as $x$ up to translation. When $x$ does
not have a periodic tail, this intertwiner is a unitary operator,
establishing a unitary equivalence between the irreducible
representations. If, however, $x$ does have a periodic tail, the
intertwiner is merely an $\infty -1$ isometry, and hence the
corresponding \cite{DaPi96}-representation is an infinite multiple of
our irreducible representation, as is also proved in
\cite[Proposition 6.10(ii)]{DaPi96}.

It is interesting to note a connection to the construction of
unitaries $W$ on a Hilbert space tensor product $\mathcal{H}\otimes
\mathcal{H}$ subject to the pentagonal identity
\begin{equation}
W_{23}W_{12}=W_{12}W_{13}W_{23}\text{\quad on\quad }\mathcal{H}\otimes 
\mathcal{H}\otimes \mathcal{H}
\label{Pentagonal}
\end{equation}
in the theory of quantum groups (recall the subscripts refer to the
usual leg-notation for tensors, that is, $W_{12}=W\otimes \openone
_\mathcal{H}$, $W_{23}=\openone _\mathcal{H}\otimes W$, and $W_{13}$
acts on the tensor product of the first and last tensor factor of
$\mathcal{H}\otimes \mathcal{H}\otimes \mathcal{H}$). Recently Saad
Baaj and Georges Skandalis have used a method that is similar to ours
in the starting point in their construction of such unitaries $W$ on
$\mathcal{H}\otimes \mathcal{H}$ which permute the elements in a
tensor product basis. Very recently they have also considered the
modification to the effect that the action of $W$ on those basis
vectors gets modified with a cocycle factor \emph{\`{a} la}
(\ref{Eq11.1}). We are indebted to S.L. Woronowicz for this
information: see \cite{Wor96} and \cite{Cun93} for more on the
pentagonal identity (\ref{Pentagonal}).

\subsection*{Acknowledgements}
Both the statement and the proof of Proposition \ref{Pro9.1} are entirely
due to David Handelman, and likewise the statement and proof of Lemma \ref
{Lem3.4} are due to Helge Tverberg. We are indebted to these two colleagues
both for permission to publish these results here, and for other useful
remarks. This paper was started when P.E.T.J. visited Oslo November 1995
with support from NFR, and continued in Oberwolfach January 1996 where O.B.
was supported by NFR. We would like to thank Brian Treadway for an excellent
typesetting job, helpful suggestions, and for producing Mathematica-generated
graphics to perfection. P.E.T.J. also had helpful
discussions with C. Bandt,
P. Muhly, R. Curto, R. Strichartz, S. Pedersen, and Yang Wang.
P.E.T.J. was also supported by the NSF and a University of Iowa
fellowship.

\subsection*{Note on the graphics}
The figures of the various reptiles $\mathbf{T}$ in this paper were executed
by Brian Treadway using Mathematica and the formulas (\ref{Eq3ins1}) and
(\ref{Eq3ins2}) in the form
\begin{equation*}
\mathbf{T}=\bigcup _{d\in D}\mathbf{N}^{-1}\left( d+\mathbf{T}\right) \text{.}
\end{equation*}
There is of course a lower limit to the size of the objects drawn, so
the series in (\ref{Eq3ins2}) is truncated in the figures. The
question of what graphic object to put at the position of each point
in the set given by (\ref{Eq3ins2}) with a finite upper limit on $i$
is slightly subtle. The best choice for the reptiles illustrated in
this paper seems to be a parallelogram (a square in the case of the
Fractal Red Cross, Figure \ref{FractalRedCross}) that is of the right
shape and just large enough to touch the neighboring parallelograms in
the picture without overlapping. Compared with the use of a circular
dot, this choice enhances the appearance that $\mathbf{T}$ is
``compact with non-empty interior'', avoids ``aliasing'' on the scale
of the laser-printer resolution, and makes the graphic easier to
rescale without changing its appearance. In this paper, Cloud Three
(Figure \ref{CloudThreeFigure}) is drawn with circular dots, but all
the other fractal figures are drawn with parallelograms.

\bibliographystyle{amsalpha}
\bibliography{jorgen}

\end{document}